\documentclass[10pt, a4paper]{TheseAlinou}

\usepackage{array}
\usepackage{amsthm}
\usepackage{amssymb}
\usepackage[all,dvips]{xy}
\usepackage{amsmath}
\usepackage{amsfonts}
\usepackage{latexsym}
\usepackage[francais]{babel}
\usepackage{ucs}
\usepackage{epsfig}
\usepackage{psfrag}
\usepackage{graphicx}
\usepackage[utf8x]{inputenc}
\usepackage[T1]{fontenc}
\usepackage{url}
\usepackage{listings}
\usepackage{fancyhdr}

\def\fract#1/#2{\hbox{\leavevmode
\kern.1em \raise .5ex \hbox{\the\scriptfont0 $#1$}\kern-.1em }/
\hbox{\kern-.15em \lower .25ex \hbox{\the\scriptfont0 $#2$}}
}

\newtheorem{thm}{Th\'eor\`eme}[section]
\newtheorem{prop}[thm]{Proposition}
\newtheorem{lem}[thm]{Lemme}
\newtheorem{cor}[thm]{Corollaire}
\newtheorem{defn}[thm]{D\'efinition}

\newtheorem{fait}[thm]{Fait}
\newtheorem{nota}[thm]{Notation}
\newtheorem{ques}{Question}
\newtheorem{rem}[thm]{Remarque}

\theoremstyle{remark}
\newtheorem{exmp}[thm]{Exemple}

\newcommand{\OP}{\mathcal{O}_{S,P}}
\newcommand{\OC}{\mathcal{O}_{S,C}}
\newcommand{\mP}{\mathfrak{m}_{S,P}}
\newcommand{\mC}{\mathfrak{m}_{S,C}}

\newcommand{\hOP}{\widehat{\mathcal{O}}_{S,P}}
\newcommand{\hOC}{\widehat{\mathcal{O}}_{S,C}}

\newcommand{\p}{\partial}
\newcommand{\w}{\wedge}
\newcommand{\F}{\mathbf{F}}
\newcommand{\res}{\textrm{res}}
\newcommand{\supp}{\textrm{Supp}}

\def\pick{{\textrm{Pic}_{k}}}
\def\divk{{\textrm{Div}_{k}}} 

\def\div{{\textrm{Div}_{\mathbf{F}_q}}} 
\def\val{{\textrm{val}}}
\def\N{{\mathbf N}}

\def\Z{{\mathbf Z}}
\def\P{{\mathbf P}}
\def\C{{\mathbf C}}
\def\L{{\mathcal L}}

\def\ev{{\textrm{ev}}}

\title{Résidus de $2$-formes différentielles sur les surfaces algébriques et applications aux codes correcteurs d'erreurs}
\author{Alain \textsc{Couvreur}}
\date{Le \today}

\begin{document}
\renewcommand{\proofname}{\textsc{Preuve}}
\renewcommand{\labelenumi}{(\theenumi)}
\renewcommand{\labelenumii}{(\theenumii)}

\thispagestyle{empty}

\begin{center}
\Huge \textbf{Th\`ese de \\  l'Universit\'e  de  Toulouse}
\end{center}

\vspace{.75cm}

\begin{center}
  {\Large \bf  R\'esidus de $2$-formes diff\'erentielles sur les surfaces alg\'ebriques
    et
    application aux codes correcteurs d'erreurs}
\end{center}

\begin{center}
  { par Alain Couvreur}
 \bigbreak
 Soutenue le lundi 8 d\'ecembre \`a 16h\\
 dans l'amphith\'e\^atre Schwartz   
\end{center}

\vspace{.5cm}

\noindent \hrulefill $\quad${\sffamily\Large\textbf{Jury}}  $\ $ \hrulefill

\vspace{.2cm}

\begin{center}
\begin{tabular}{lllll}

 {\bfseries Emmanuel Hallouin} &  $\ \ $ & Universit\'e de Toulouse II & $\ \ $ & Examinateur\\
 {\bfseries Gilles Lachaud } &  & Universit\'e d'Aix-Marseille II  & & Examinateur\\
 {\bfseries Marc Perret} & & Universit\'e de Toulouse II  & & Directeur\\
 {\bfseries Marc Reversat} & & Universit\'e de Toulouse III  & & Directeur\\
 {\bfseries Felipe Voloch} & & University of Texas  & & Rapporteur\\
 {\bfseries Gilles Zemor} & & Universit\'e de Bordeaux I  & & Examinateur\\
\end{tabular}
\end{center}

\vspace{.2cm}

\noindent \hrulefill

\vspace{.5cm}

\vfill

\begin{center}
  \includegraphics[width=6cm, height=4.5cm]{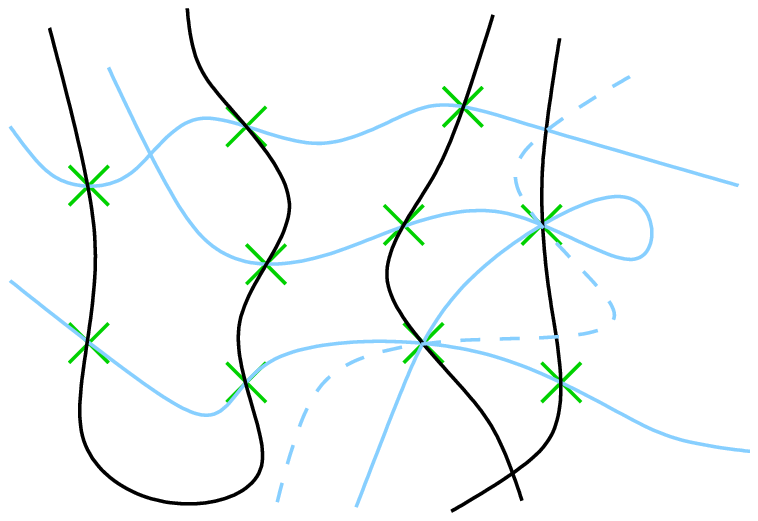}
\end{center}

\vfill

\begin{center}
  Institut de Math\'ematiques de Toulouse, UMR 5219, UFR MIG\\
  Laboratoire \'Emile Picard,\\
  Universit\'e Paul Sabatier 31062 TOULOUSE C\'edex 9
\end{center}


\newpage
\thispagestyle{empty}
\null

\newpage

\thispagestyle{empty}
\null

\vfill

\begin{flushright}
\textit{Ce qui nous rassure du sommeil,
c'est qu'on en sort,\\
et qu'on en sort inchangé,
puisqu'une interdiction bizarre \\ nous empêche de rapporter avec nous
l'exact résidu de nos songes.}

\medbreak

Marguerite Yourcenar\\
\textit{Mémoires d'Hadrien}

\end{flushright}

\vfill

\null
\newpage
\thispagestyle{empty}
\null
\newpage

\chapter*{Remerciements}

Quelques mois avant de commencer ma thèse, une jeune chercheuse qui me vantait les mérites de la recherche avait formulé cette mise en garde: ``une thèse, c'est difficile, on est très seul.''
Je suis convaincu qu'une écrasante majorité de doctorant(e)s - dont je fais partie - a éprouvé ce sentiment au moins une fois.
Pourtant, maintenant que la rédaction de ma thèse touche à sa fin, je réalise à quel point rien de tout ce que j'ai fait n'aurait été possible si j'avais réellement été seul.

Très souvent, les premières lignes d'un livre sont celles que l'auteur a écrit en dernier. 
Cette thèse n'échappe pas à la règle.
De la même manière, mes premiers remerciements vont vers les acteurs du dénouement.

\paragraph{Ceux sans qui cette thèse ne pourrait se terminer.}
Je remercie les six personnes qui ont accepté de composer mon jury. Tout d'abord, Marc Reversat sans qui mon intégration dans l'ex-laboratoire Émile Picard n'aurait pas été possible.
Ma gratitude va ensuite à Gilles Lachaud, directeur de mon directeur, pour son accueil chaleureux à l'IML au début de cette année.
Un grand merci également à Gilles Zemor pour tout ce qu'il m'a appris lors de ma visite à Télécom Paris ainsi que pour les pertinentes corrections qu'il a suggéré pour ce manuscrit.
Je remercie ensuite Felipe Voloch pour son accueil et son encadrement durant mon séjour à l'Université du Texas, je le remercie également d'avoir accepté de rapporter ma thèse.
Je remercie bien sûr mon directeur Marc Perret, mais j'y reviendrai un peu plus loin...
Enfin, je remercie Emmanuel Hallouin pour ses nombreux conseils tout au long de ma thèse. Je garderai un excellent souvenir des discussions de politique ou d'algèbre locale que nous avons eues ensemble.

Je remercie par ailleurs Dino Lorenzini pour toutes ses remarques et pertinentes suggestions sur mon manuscrit. Je lui suis également reconnaissant d'avoir accepté de rapporter ma thèse. Pour la même raison, je remercie sincèrement le troisième rapporteur de cette thèse Michael Tsfasman.

\paragraph{L'équipe de la dernière ligne droite.}
Merci à Matthieu et Aurélie qui ont sans aucune hésitation accepté que leur appartement soit envahi et leur cuisine dévastée, l'espace d'un week-end. Merci également à Tiphaine qui, durant la préparation du pot, n'a jamais rechigné à nettoyer tout couvert dégoulinant de chocolat fondu.

\paragraph{Ceux qui murmurent à l'oreille des processeurs.} Face aux nombreux soucis dus à un ordinateur récalcitrant, j'ai toujours bénéficié d'une assistance mail ou téléphonique d'une efficacité remarquable.
À ce titre, je remercie Maxime, grâce à qui ma version d'Emacs est si performante qu'elle serait presque capable de commander une pizza ou de cirer un parquet.
Merci à Gâchette pour son objectivité dans le débat Ubuntu vs Madrake et enfin merci à Benjamin qui nous a démontré par l'exemple que son ordinateur cherchait à devenir le maître du monde.

\paragraph{Ceux qui font avancer à grands pas.}
Durant ces trois dernières années, un certain nombre de conversations orales ou écrites avec d'autres chercheurs m'ont indiscutablement aidé à avancer. Je remercie tout d'abord Julien Duval, Steven L. Kleiman et Gerhard Frey pour leurs précieuses explications.
Merci ensuite à Antoine Ducros pour son aptitude à dégainer un contre exemple plus vite que son ombre et à Joseph Tapia pour son excellente synthèse orale du \textit{Residues and Duality} d'Hartshorne.
Enfin, je tiens à remercier sincèrement Tom H\o holdt pour sa gentillesse, son dévouement, son hospitalité ainsi ses indiscutables qualités de pédagogues qui ont fait de ma visite à la DTU (Université Technique du Danemark) un séjour aussi fructueux qu'agréable.

\paragraph{Ceux qui font la même chose... ou presque.} Un grand merci à tous les doctorants (et docteurs) Toulousains que j'ai côtoyés durant ces trois années.
Je remercie tout particulièrement Anne pour les ``quatre heures'' du mois d'août, Cécile pour sa maîtrise des blagues Carambar, Tony pour ses imitations du Schtroumpf grognon et Landry pour sa capacité à animer un débat en évitant systématiquement le consensus\footnote{Je sais, tu n'es pas d'accord avec ce que je viens de dire.}.
Merci également à Fred Protin pour ses mails aux jeux de mots décapants, à Fred Pitoun pour son dynamisme décontracté qui tonifiait nos jeudis, à Julien Roques qui partageait mon avis sur Bono et les performances de Randy Marsh.
Je remercie enfin mes cobureaux, Charef, Ghada et Matthieu, leur présence me fut toujours des plus agréables.

\paragraph{Les uns et les hôtes.} Lors de mes missions hors de Toulouse, j'ai fréquemment apprécié la qualité de mon accueil. Je remercie à ce titre Delphine Boucher de Rennes, Sylvain Duquesne et Louise Nyssen de Montpellier
 et Pierre-Louis Cayrel de Limoges.
Je remercie également Frederic Edoukou, Adnen Sboui, leur ex-directeur de thèse François Rodier et Christophe Ritzenthaler de m'avoir acueilli durant une semaine à l'IML durant ma seconde année de thèse.

Je tiens également à remercier les membres de l'IML et du département de mathématiques de Luminy qui m'ont si bien accueilli et intégré dans leur équipe cette année. Remerciements particuliers à tous les membres de l'équipe ATI. 

\paragraph{Parce qu'il n'y a pas que la recherche dans la thèse...} Je remercie tous les enseignants et chercheurs du département de mathématiques de l'Université du Mirail. Mes trois années d'enseignement dans cette université furent un plaisir, tant pour le contact des étudiants que pour celui des collègues. 
Je tiens tout particulièrement à remercier Julien Labetaa pour qui j'ai donné des TD's et avec qui la collaboration fut des plus agréables.

\paragraph{... et qu'il n'y a pas que le travail dans la vie.} Parmi les nombreux souvenirs que je garderai de ces trois années, il y aura les nombreuses fins d'après midi ensoleillées en terrasse.
Merci à ceux avec qui j'ai partagé ces si agréables moments.
Merci à Émilie, Gustavo, Romain, Seb, Solenn, Soazig, Tanguy et Perrine. 
Un grand merci également à Cécile qui m'a expliqué la différence entre un bus et un TUB, à Erwan parce qu'il comprend le 229$^{\textrm{e}}$ degré et à Xavier pour son sens de l'orientation en situation critique.

\paragraph{Ceux qui simplifient la vie.} Merci à Véronique Fabris, Agnès Requis et Jocelyne Picard pour leur disponibilité et leur patience en toute circonstance.

\paragraph{Celle que j'ai croisé.} Merci à Lara pour tous les conseils qu'elle a pu me donner lors des nombreuses conversations que l'on a eues ensemble. C'est toujours un plaisir pour moi de trouver son nom sur la liste des participants d'une conférence à laquelle je me rends.

\paragraph{Ceux qui m'ont accueilli.} Merci à tous les membres de l'ex-Grimm de m'avoir intégré parmi eux. Remerciements particuliers à Thierry Henocq et Christian Maire.

\paragraph{Ceux sans qui cette thèse n'aurait pas eu lieu.} Je remercie Jean-Marc Couveignes qui fut mon premier contact à l'ex-Grimm et qui consacra une énergie particulière à la régularisation de ma complexe situation administrative.
Merci également à Arnaud Debussche et Michel Pierre de m'avoir si bien conseillé en fin de Master.

\paragraph{Au chef.} Merci à Marc Perret de m'avoir si bien encadré durant ces trois années. J'ai particulièrement apprécié son investissement, sa patience, sa capacité a expliquer en des termes simples les faits et objets mathématiques les plus complexes. Connaissant son extrême modestie, je préfère ne pas en dire plus de peur de le mettre mal à l'aise, mais si tout était à refaire, je lui conseillerais de ne rien changer.

\paragraph{Celles qui étaient là dès le début.} Merci à mes s{\oe}urs Nadine et Sylvie.

\paragraph{Ceux sans qui je ne serais pas là.} À mes parents, merci pour tout.

\paragraph{Celle qui est tout pour moi.} Merci Gwenola. 



\newpage
\thispagestyle{empty}
\null
\newpage

\tableofcontents

\newpage
\thispagestyle{empty}
\null
\newpage
\chapter*{Introduction}

Cette thèse est composée de deux parties.
La première porte sur les notions de résidus de $2$-formes différentielles rationnelles sur une surface algébrique.
La seconde partie utilise les résultats de la première en vue d'applications aux codes correcteurs d'erreurs. 
Ce travail de recherche est parti d'une constatation simple. En théorie des codes géométriques construits à partir de courbes algébriques, on distingue deux types de constructions.
La construction \textit{fonctionnelle}  qui, comme son nom l'indique, utilise des fonctions et la construction \textit{différentielle} qui utilise des formes différentielles.
Cependant, tous les travaux de recherche abordant l'étude des codes géométriques construits à partir de variétés de dimension supérieure ou égale à $2$ font systématiquement appel à une construction de type fonctionnelle.
De cette observation est née une question: \textit{peut-on généraliser la construction différentielle en dimension supérieure ou égale à $2$?}
Seule la généralisation aux surfaces sera abordée, nous justifierons ce choix un peu plus loin dans cette introduction.

\section*{Historique des codes géométriques}

La première construction de codes correcteurs d'erreurs par des méthodes issues de la géométrie algébrique a été présentée par Goppa dans \cite{goppa}.
Peu après, dans \cite{TVZ}, Tsfasman, Vl\u{a}du\c{t} et Zink, utilisaient cette approche géométrique pour construire des familles de codes  dont les performances asymptotiques dépassaient celles de toutes les familles de codes connues jusque là.
Ces résultats ont été la principale motivation du développement de la théorie des codes géométriques.

\subsubsection*{Codes sur les courbes algébriques} Dès la fin des années 80, la théorie des codes géométriques était devenue un thème de recherche  extrêmement dynamique.
Plusieurs centaines d'articles ont été publié sur l'étude de ces codes, que ce soit sur la recherche de bons codes, de bonnes familles de codes ou encode d'algorithmes de décodage.
Il serait donc difficile de fournir une bibliographie complète sur le sujet. Signalons tout de même les quelques publications présentant un point de vue général sur la théorie.  
Le premier article de synthèse sur la question est dû à Lachaud \cite{lachsembou}, il y est présenté toutes propriétés théoriques connues sur les codes géométriques. Pour des références plus détaillées, on peut consulter le livre de Goppa \cite{livregoppa} ou celui de Tsfasman et Vl\u{a}du\c{t} \cite{TV}. Enfin, pour une synthèse sur les algorithmes de décodage de codes géométriques on pourra se référer à l'article de synthèse de H{\o}holdt et Pellikaan \cite{hohpel} pour les travaux connus avant 1995 et au chapitre de \cite{bouquindiego} écrit par Beelen et H{\o}holdt pour les travaux plus récents.

\newpage
\thispagestyle{myheadings}
\markboth{}{} 
\subsubsection*{Codes sur les variétés en dimension supérieure} Si le sujet des codes géométriques sur les courbes a été étudié de façon très approfondie, la recherche sur les codes construits à partir de variétés de dimension supérieure ou égale à $2$ est restée nettement plus marginale. Historiquement, le premier à avoir donné une construction de codes correcteurs d'erreurs à partir de variétés de dimension quelconque est Manin dans \cite{manin}. 
Par la suite, un certain nombre d'articles est paru sur la question. La liste de références qui suit n'est pas exhaustive.\label{refbib}

On dénombre au moins trois publications fournissant des résultats généraux sur les codes géométriques construits à partir de variétés algébriques de dimension supérieure ou égale à $2$.
Dans \cite{lachaudbis} et \cite{lachaud}, Lachaud fournit une minoration de la distance minimale des codes construits sur une variété projective lisse quelconque.
Dans \cite{soH}, S\o ren Have Hansen étudie les paramètres des codes construits à partir de variétés algébriques lisses quelconques et propose des exemples issus des variétés de Deligne-Lustzig.
Enfin, dans \cite{bouganis}, Bouganis étudie les codes construits sur des surfaces algébriques lisses quelconques  puis étudie le comportement asymptotique de certaines familles de tels codes.

Pour le reste, la plupart des autres travaux publiés portent sur l'estimation des paramètres de codes construits à partir de variétés appartenant à une classe particulière.
Les codes sur les surfaces toriques ont été étudiés par Hansen dans \cite{hansen}. Ses résultats ont ensuite été généralisés en dimension quelconque par Ruano dans \cite{diego}.
Les codes construits sur des Grassmanniennes ont d'abord été étudiés par Nogin dans \cite{nogin}, puis par Ghorpade et Lachaud dans \cite{gl}.
Les codes construits à partir de variétés Hermitiennes ont été abordés pour la première fois par Chakravarti dans \cite{chak2}, ensuite par Hirschfeld, Tsfasman et Vl\u{a}du\c{t} dans \cite{HTV}, puis par S{\o}rensen dans sa thèse \cite{sorensen} et enfin par Edoukou dans \cite{fred}. Notons que les variétés Hermitiennes et Grassmaniennes peuvent être vues comme des variétés drapeaux. Ce point de vue unifié est discuté par Rodier dans \cite{rodier}.
 La distance minimale des codes sur les variétés quadriques de dimension quelconque est étudiée par Aubry dans \cite{aubry}.
Le cas des surfaces quadriques est approché de façon pus détaillée par Edoukou dans \cite{fred2}. 
Enfin, Zarzar a traité le cas des surfaces dont le rang du Groupe de Néron-Sévéri arithmétique est petit dans \cite{zarzar}. Il propose ensuite dans un travail commun avec F. Voloch \cite{agctvoloch}, une approche de décodage utilisant un algorithme de décodage itératif proposé par Luby et Mitzenmacher \cite{luby}.

Enfin, signalons qu'une excellente synthèse sur les travaux connus sur les codes construits sur des variétés de dimension supérieure est présentée dans une prépublication de Little (voir \cite{little}).

\medbreak

À présent, rappelons que, comme indiqué au début de ce chapitre introductif, en théorie des codes sur les courbes on distingue deux méthodes de construction de codes respectivement appelées construction fonctionnelle et différentielle.
Cependant, en dimension supérieure, on ne dispose que de la construction fournie par Manin dans \cite{manin}. Cette dernière est une généralisation naturelle de la construction fonctionnelle sur les courbes.
Tous les travaux cités ci-dessus s'appuient sur cette construction et aucune généralisation de la construction différentielle n'a été proposée jusque là.
Notons d'ailleurs que Little signale dans l'introduction de son article de synthèse \cite{little} une obstruction majeure à une telle généralisation.

\medbreak
\textit{``In a sense, the first major difference between higher dimensional varieties and curves is that points on $X$ of dimension $\geq 2$ are subvarieties of codimension $\geq 2$, not divisors. This means that many familiar tools used for Goppa codes (e.g. Riemann-Roch theorems, the theory of differentials and residues etc.) do not apply exactly in the same way. ''}

\medbreak

En quelques mots, l'objectif de cette thèse est, après avoir mis en place le matériel théorique nécessaire, de fournir une construction différentielle de codes sur les surfaces, puis de l'appliquer à l'étude des codes géométriques.

\section*{Pourquoi des codes différentiels sur les surfaces?}

Outre la volonté de généralisation en vue d'une harmonisation des théories entre le cas des courbes et celui des variétés de dimension supérieure, plusieurs arguments motivent cette question.

\paragraph*{Un intérêt historique.}
Les codes géométriques ont été introduits pour la première fois par V.D Goppa en 1981 \cite{goppa}.
Dans cet article, la construction présentée était différentielle.
Aussi, même si les codes fonctionnels sont plus populaires chez les spécialistes des codes géométriques, la construction historique est de type différentielle.

\paragraph*{L'intérêt d'une nouvelle construction géométrique.}
Le second argument réside dans l'intérêt de disposer d'une construction géométrique de codes. Pour comprendre en quoi une telle construction est avantageuse, commençons par réfléchir aux différentes façons de décrire un code. La manière la plus simple est de s'en donner une base, c'est-à-dire une matrice génératrice.
Cependant, une telle description n'est pas du tout adaptée à la résolution de problèmes tels que la minoration de la distance minimale ou la recherche d'un algorithme de décodage efficace.
Par conséquent, on cherche en général à résoudre ces problèmes pour des classes de codes admettant une \textit{réalisation} par des objets appartenant à une autre branche des mathématiques, comme l'arithmétique ou la géométrie. 
C'est par exemple le cas des codes de Reed-Solomon qui font appel à des polynômes en une variable, des codes de Reed-Müller qui se construisent à partir de polynômes à plusieurs variables ou encore des codes de résidus quadratiques dont la construction et l'étude font appel à de l'arithmétique des corps finis.
Par ce biais, les problèmes de minoration de la distance minimale et de recherche d'algorithmes de décodage peuvent être traduits sous forme de problèmes d'algèbre ou de géométrie.
On se ramène donc à un contexte comportant une \textit{structure} (arithmétique ou géométrique par exemple) et dans lequel on dispose de davantage d'outils mathématiques pour résoudre un problème donné.
En conclusion, \textbf{il est toujours intéressant de disposer d'une réalisation géométrique d'un code pour l'étudier}. À ce titre, la construction de codes correcteurs à partir de formes différentielles sur des surfaces est une voie que l'on se doit d'explorer.

\paragraph*{Des codes en relation avec les codes fonctionnels.}
En théorie des codes géométriques construits à partir de courbes, on dispose de relations entre codes fonctionnels et codes différentiels.
\begin{enumerate}
\item[$(\mathbf{R1})$]\label{r1} Un code différentiel sur une courbe est toujours l'orthogonal d'un code fonctionnel construit à partir de la même courbe et associé aux mêmes diviseurs.
\item[$(\mathbf{R2})$]\label{r2} Tout code différentiel sur une courbe se réalise comme un code fonctionnel construit à partir de la même courbe mais associé à des diviseurs différents.
\end{enumerate}

La relation ($\mathbf{R1}$) est une conséquence de la formule des résidus et du théorème de Riemann-Roch.
Cette propriété d'orthogonalité est de plus un ingrédient utilisé dans de nombreux algorithmes de décodage (voir \cite{hohpel}).
D'une façon générale, disposer d'une réalisation géométrique de l'orthogonal ou d'un sous-code de l'orthogonal d'un code correcteur peut être fort utile pour le décodage.
La relation ($\mathbf{R2}$) est une conséquence du théorème d'approximation faible (\cite{sti} I.3.1).
Elle implique que les codes fonctionnels et les codes différentiels construits à partir de courbes algébriques, bien qu'obtenus par des constructions différentes, appartiennent à la même classe.
On peut donc restreindre l'étude générale de ces codes à celle de codes provenant d'une seule des deux constructions.
Le plus souvent, c'est la construction fonctionnelle qui est adoptée.
Ce choix vient sans doute de ce que, pour beaucoup de mathématiciens, la notion d'évaluation d'une fonction en un point est plus intuitive et manipulable que celle d'évaluation du résidu d'une forme différentielle.

Ainsi, après s'être interrogé sur la possibilité d'étendre aux surfaces la construction différentielle de codes, il est naturel de réfléchir aux perspectives d'extension aux surfaces des propriétés ($\mathbf{R1}$) et ($\mathbf{R2}$). De tels résultats contribueraient en effet à approfondir nos connaissances des codes géométriques construits à partir de surfaces.
Nous détaillerons les résultats obtenus dans ce sens en page \pageref{resultats1}.

\paragraph*{D'intéressants développements théoriques.}
Nous allons voir que la construction et l'étude des codes différentiels construits sur des surfaces a nécessité de nombreux résultats théoriques concernant les formes différentielles sur les surfaces. 
Les résultats énoncés dans le premier chapitre ne sont pas réellement nouveaux.
En géométrie algébrique, la notion de résidu en dimension supérieure à $2$ a été abordée par Grothendieck et Hartshorne dans \cite{harRD} ainsi que par Lipman dans \cite{lip}.
Cependant, à la différence de ces références, la notion de résidu présentée dans le chapitre \ref{chapres} provient d'une construction explicite ne faisant appel à aucun raisonnement de type fonctoriel.
La volonté de construire des codes différentiels sur des surfaces algébriques a donc permis l'élaboration d'une introduction au résidus sur des surfaces par une approche plus explicite et constructive que celles qui existaient jusque-là\footnote{Peu de temps après l'envoi de la seconde version de ce manuscrit, Oleg Osipov du Steklov Mathematical Institute, m'a contacté après avoir consulté une prépublication de mes résultats sur ArXiv (voir \cite{oim}).
Il m'a alors signalé qu'une approche similaire avait été donnée par Par\v{s}in dans \cite{parshin}. J'ignorais l'existence de cet article peu connu et rarement cité lorsque j'ai travaillé sur ces questions.}\label{parshin}.

\medbreak

Avant de passer à une présentation plus détaillée des différentes parties de la thèse. Signalons que le contenu des chapitres \ref{chapres} et \ref{chapdiff} en version \textit{condensée} a donné lieu à la rédaction d'un article \cite{oim}.

\section*{Présentation de la première partie}

Si la notion de résidu est bien connue dans le cas des $1$-formes différentielles sur une courbe algébrique et qu'une unique définition de cet objet fait l'unanimité dans la littérature, en dimension supérieure la situation est nettement moins claire.
Par exemple, en géométrie algébrique complexe, la définition énoncée dans l'ouvrage \cite{gh} de Griffiths et Harris  diffère de celle du livre  \cite{bhpv} de Bath, Peters, Hulek et Van de Ven.
Pour le premier, un résidu est un élément du corps de base (le corps des complexes) obtenu à partir de la donnée d'une $n$-forme méromorphe $\omega$ définie sur une variété complexe $X$ de dimension $n$, d'un point $P$ de $X$ et d'une famille ordonnée de $n$ diviseurs de cette variété s'intersectant en $P$.
Pour le second, étant donnée une variété complexe $X$ et une sous-variété $Y$ de codimension un dans $X$, le résidu d'une $r$-forme méromorphe sur $X$ le long de $Y$ est la donnée d'une $(r-1)$-forme sur $Y$.
Notons dès maintenant que ces ouvrages se placent dans le contexte des variétés complexes, contexte dans lequel on peut calculer les résidus avec l'aide de la formule de Cauchy.
En d'autres termes, les résidus peuvent être obtenus en intégrant une forme différentielle sur une sous-variété réelle.
Ce point de vue utilise le fait qu'une variété complexe de dimension $n$ peut être vue comme une variété réelle de dimension $2n$. 
Un tel point de vue ne peut évidemment pas s'étendre à un autre cadre comme par exemple celui des variétés sur un corps fini.

Dans un contexte plus général, on trouve dans \cite{harRD} un objet appelé \textit{résidu de Grothendieck} qui ressemble à l'objet défini par Griffiths et Harris en ce sens qu'il associe à une forme différentielle de degré maximal un élément du corps (ou de l'anneau) de base. 
Cet objet est cependant plus fortement relié à un système de coordonnées locales et sa construction nécessite un important arsenal d'objets et de raisonnements fonctoriels.

\medbreak

Dans la première partie, qui est composée du seul chapitre \ref{chapres}, on introduira les notions de $1$-résidu qui correspondront à la définition de \cite{bhpv} et de $2$-résidu qui correspondront à la définition de \cite{gh}.
Nous étudierons également les relations qui lient ces objets.
Pour ce faire, nous étudierons les développements de fonctions et de $2$-formes différentielles en séries de Laurent de deux variables. 
Le $2$-résidu sera l'objet qui suscitera le plus notre attention. Il permet d'extraire un élément du corps de base à partir de la donnée d'une $2$-forme rationnelle $\omega$ sur une surface, d'une courbe $C$ plongée dans cette surface et d'un point rationnel $P$ de $C$. On le notera
$$
\res^2_{C,P}(\omega).
$$

\subsection*{Présentation des résultats de la première partie}

Les travaux effectués dans la première partie (donc le premier chapitre) aboutissent à deux types de résultats.

\subsubsection*{Invariance des $2$-résidus}
Le premier résultat majeur est le théorème \ref{inv2res} qui assure que l'application $\res^2_{C,P}$ est bien définie. En d'autres termes, le $2$-résidu en un point $P$ le long d'une courbe $C$ d'une $2$-forme rationnelle $\omega$ ne dépend pas d'un choix de coordonnées locales.

\subsubsection*{Formules de sommation}
L'objectif principal de ce travail est d'obtenir des formules du type: ``\textit{la somme des résidus de $\omega$ est nulle}'', en vue de relations d'orthogonalité entre codes dans la seconde partie. 
Dans la section \ref{sommation} du chapitre \ref{chapres}, on fournira trois formules de sommation\footnote{Comme signalé dans la note au bas de la page \pageref{parshin}, une partie des résultats présentés dans la première partie de cette thèse avaient en fait déjà été démontrées dans \cite{parshin} par des méthodes similaires. C'est par exemple le cas des deux premières formules de sommation de résidus, à savoir les théorèmes \ref{FR1} et \ref{FR2}}.

\medbreak

\noindent \textbf{Théorème \ref{FR1}} (Première formule des résidus)\textbf{.}
\textit{Soit $S$ une surface projective irréductible lisse définie sur un corps algébriquement clos. Soient $C$ une courbe projective irréductible plongée dans $S$ et $\omega$ une $2$-forme rationnelle sur $S$. On a
$$
\sum_{P\in C} \res^2_{C,P}(\omega)=0.
$$}   

\medbreak

\noindent \textbf{Théorème \ref{FR2}} (Deuxième formule des résidus)\textbf{.}
\textit{Soit $S$ une surface quasi-projective irréductible lisse définie sur un corps algébriquement clos. 
Soient $P$ un point de $S$ et $\mathcal{C_{S,P}}$ l'ensemble des germes de courbes irréductibles tracées sur $S$ et contenant $P$. Pour toute $2$-forme $\omega$ rationnelle sur $S$, on a
$$
\sum_{C\in \mathcal{C}_{S,P}} \res^2_{C,P}(\omega)=0.
$$
}

\medbreak

La troisième formule de sommation nécessite la définition de $2$-résidu en un point le long d'un diviseur. Nous renvoyons le lecteur à la définition \ref{resdiv} page \pageref{resdiv}.

\medbreak

\noindent \textbf{Théorème \ref{FR3}} (Troisième formule des résidus, \cite{lip} chap. 12)\textbf{.}
\textit{Soit $S$ une surface projective irréductible lisse définie sur un corps algébriquement clos. 
Soient $D_a$ et $D_b$ deux diviseurs sur $S$ dont l'intersection des supports est un ensemble fini $Z$. Soit $\Omega^2(-D_a-D_b)$ le faisceau de $2$-formes rationnelles vérifiant localement
$$
(\omega) \geq -D_a-D_b.
$$  
Alors, pour toute section globale $\omega$ du faisceau $\Omega^2 (-D_a-D_b)$, on a 
$$
\sum_{P\in S} \res^2_{D_a,P}(\omega)=\sum_{P\in Z} \res^2_{D_a,P}(\omega)=0.
$$}

\medbreak

La troisième formule des résidus est le résultat que nous utiliserons dans le chapitre \ref{chapdiff} pour obtenir un résultat d'orthogonalité entre codes.
Elle se démontre à l'aide des deux autres formules de sommation  énoncées (les théorèmes \ref{FR1} et \ref{FR2}).
Cette troisième formule des résidus n'est pas nouvelle, on en trouve un énoncé similaire dans le chapitre 12 de \cite{lip} qui est valable en toute dimension et pas seulement sur les surfaces.
Nous insistons une fois de plus sur le fait que la démonstration donnée dans cette thèse a l'intérêt de faire appel à des constructions plus explicites et plus constructives que celles utilisées dans les démonstrations connues de ce résultat. 
Achevons notre argumentation à ce sujet par une citation justement extraite de \cite{lip}, afin de légitimer (``\textit{more or less}'') le choix que nous avons fait de présenter et démontrer ces résultats de manière nouvelle et plus accessible.

\medbreak

\textit{``Statements 0.3A and 0.3B, are consequences (more or less) of (\cite{harRD} page 383 corollary 3.4). However, one of our main purposes in this paper is to provide a proof of 0.3 for which loc. cit. is not a prerequisite. The other main purpuse is to describe the connection between local and global duality, via residues (c.f. \cite{harRD} page 386 prop 3.5).''}

\medbreak

Avant de passer à la présentation de la seconde partie, finissons par une remarque.
Il peut sembler naturel de se demander pourquoi les résultats énoncés dans cette thèse ne portent principalement que sur les surfaces et non sur les variétés de dimension supérieure.
Différentes raisons ont motivé ce choix.
La première est que, même s'il est fort probable que les constructions et les résultats présentés dans la première partie admettent une généralisation en dimension supérieure à $2$, tout travail dans cette direction aurait entraîné d'importantes lourdeurs dans les notations.
Nous avons donc choisi de nous restreindre au cas déjà non trivial des surfaces, sachant que, pour ce type de problème de géométrie algébrique, le passage de la dimension $1$ à $2$ est l'étape difficile à franchir.
Enfin, l'objectif étant de travailler sur les codes correcteurs, il semblait déjà fort intéressant de ne considérer que le cas des surfaces, ce dernier n'ayant été que rarement exploré.
Il nous a donc semblé inutile de partir vers de telles généralités alors que le monde des surfaces algébriques offrait déjà de si nombreuses perspectives.

\section*{Présentation de la seconde partie}

La seconde partie contient les chapitres \ref{chapdiff} à \ref{chapldpc}. Elle concerne les codes géométriques et plus précisément les codes différentiels construits sur des surfaces algébriques.

\subsection*{Présentation des résultats de la seconde partie}

Les codes différentiels sur une surface algébrique sont définis dans le chapitre \ref{chapdiff} (définition \ref{defdiff} page \pageref{defdiff}).
On se donne dans tout ce chapitre une surface projective lisse géométriquement intègre $S$ sur $\F_q$, un diviseur $G$ sur $S$ et une famille de points rationnels $P_1, \ldots, P_n$ de $S $ qui évitent le support de $G$.
On note $\Delta$ le $0$-cycle
$$
\Delta:=P_1+ \cdots +P_n.
$$

Dans tout ce qui suit et jusqu'à la fin de cette introduction, les codes fonctionnels seront notés ``$C_L$'' et les codes différentiels ``$C_{\Omega}$''.
Les définitions respectives de ces codes sont données en sections \ref{codesfonc} et \ref{codesdiff}.

\subsubsection*{Codes différentiels sur les surfaces}\label{resultats1}

La construction de ces codes nécessite l'introduction d'une paire de diviseurs $(D_a,D_b)$.
Pour obtenir une relation d'orthogonalité on définit la notion de paires de diviseurs $\Delta$-convenables (voir définition \ref{deltac} page \pageref{deltac}).
Il s'agit de paires de diviseurs qui sont en un certain sens \textit{reliées} au $0$-cycle $\Delta$.
Le premier résultat majeur de ce chapitre est une relation d'othogonalité qui est plus faible que la propriété ($\mathbf{R1}$) dans le cas des courbes puisqu'il ne s'agit plus que d'une inclusion au lieu d'une égalité.

\medbreak

\noindent \textbf{Théorème \ref{orthocode}} (Théorème d'orthogonalité)\textbf{.}
\textit{Soient $(D_a,D_b)$ une paire $\Delta$-convenable de diviseurs et $D:=D_a+D_b$. On a alors,
$$
C_{\Omega,S}(\Delta, D_a, D_b, G) \subseteq C_{L,S}(\Delta, G)^{\bot}.
$$}

\medbreak

L'inclusion réciproque est en général fausse, comme le montre le contre-exemple donné en section \ref{P1P1}. Plus précisément, on présente l'exemple d'une surface (le produit de deux droites projectives) sur laquelle l'orthogonal d'un code fonctionnel ne se réalise sous la forme d'un code différentiel pour aucun choix de paire de diviseurs $\Delta$-convenable $(D_a,D_b)$.

\medbreak

Nous étudions ensuite la possibilité d'étendre aux codes sur les surfaces la propriété ($\mathbf{R2}$).
À la différence de ($\mathbf{R1}$), cette seconde relation s'étend parfaitement au cas des surfaces.

\medbreak

\noindent \textbf{Théorème \ref{diff=fonc}.}
\textit{Soient $(D_a,D_b)$ une paire $\Delta$-convenable de diviseurs et $D:=D_a+D_b$, alors il existe un diviseur canonique $K$ tel que
$$
C_{\Omega}(D_a,D_b,G)=C_L(\Delta, K-G+D).
$$   }

\medbreak

\noindent \textbf{Théorème \ref{fonc=diff}.}
\textit{ Étant donné un diviseur $G$ sur $S$, il existe un diviseur canonique $K$ et une paire $\Delta$-convenable $(D_a,D_b)$ telle que
$$
C_L(\Delta, G)=C_{\Omega}(D_a,D_b,K-G+D).
$$}

\medbreak

Le chapitre \ref{chapdiff} se termine par une discussion en section \ref{heuris} autour des raisons du défaut d'inclusion réciproque dans le théorème d'othogonalité \ref{orthocode}. Cette discussion est consécutive à la présentation d'un contre-exemple à l'inclusion réciproque du théorème \ref{orthocode} donnée en section \ref{P1P1}.
Par ailleurs, ce contre-exemple permet de conclure le second chapitre sur une importante constatation.
Il assure en effet que les codes fonctionnels construits sur une surface algébrique et leurs orthogonaux appartiennent en général à une classe différente.
C'est un phénomène qui différencie fondamentalement le cas des courbes de celui des surfaces.
Notons que cette asymétrie entre les codes fonctionnels et leurs orthogonaux avait déjà été signalée par Voloch et Zarzar dans \cite{agctvoloch}.

\medbreak
\textit{``It is interesting to note that Goppa codes coming from curves are seldom LDPC since their duals are also Goppa codes coming from curves and, as such, have a large minimal distance, whereas the dual of an LDPC has a small minimal distance by definition.''}

\medbreak

Pour le reste, cette observation ouvre un intéressant axe de recherche, celui de l'étude de l'orthogonal d'un code fonctionnel sur une surface. C'est ce qui donnera lieu au chapitre \ref{chaporth}, nous y reviendrons plus loin.

\subsubsection*{Théorème de réalisation}

Dans le chapitre \ref{chapreal} on montre comment, sous certaines conditions sur la surface $S$ et le diviseur $G$, on peut réaliser l'orthogonal d'un code fonctionnel non pas comme un code différentiel mais comme une somme de codes différentiels.
L'énoncé du théorème fait appel à la notion de sous-$\Delta$-convenance définie en section \ref{secsous} (définition \ref{sousDconv} page \pageref{sousDconv}).

\medbreak

\noindent \textbf{Théorème \ref{thmreal}} (Théorème de réalisation)\textbf{.}
\textit{Soient $S$ une surface lisse géométriquement intègre et intersection complète dans un espace projectif $\P^r_{\F_q}$ et $G$ un diviseur sur $S$ linéairement équivalent à une section de $S$ par une hypersurface de $\P^r$.
On se donne également un $0$-cycle $\Delta$ qui est la somme de $n$ points rationnels de $S$ évitant le support de $G$.
Soit $c$ un mot du code $C_{L,S} (\Delta, G)^{\bot}$. Alors, il existe une paire de diviseurs $(D_a, D_b)$ et une $2$-forme $\omega$ appartenant à l'espace des sections globales $\Gamma(S,\Omega^2(G-D_a-D_b))$, tels que
$$
c=\res^2_{D_a, \Delta}(\omega).
$$
}

\medbreak

\noindent \textbf{Remarque.}
\textit{  Le théorème de réalisation dit en fait un peu plus que ça, il fournit également des informations sur les structures géométriques et les classes d'équivalence linéaires des diviseurs $D_a$ et $D_b$ (voir page \pageref{thmreal}).}

\medbreak

\noindent \textbf{Corollaire \ref{correal}.}
\textit{Sous les hypothèses du théorème de réalisation,
il existe une famille finie $(D_a^{(1)}, D_b^{(1)}), \ldots ,(D_a^{(s)}, D_b^{(s)})$ de paires de diviseurs sous-$\Delta$-convenables telles que
$$
C_{L,S}(\Delta, G )^{\bot} = \sum_{i=1}^s C_{\Omega,S} (\Delta, D_a^{(i)}, D_b^{(i)}, G).
$$}

\medbreak

La démonstration du théorème de réalisation utilise un théorème ``à la Bertini'' sur les corps finis démontré par Poonen en 2004 dans \cite{poon}.
On termine le chapitre en montrant qu'un argument ``à la Bertini'' de type différent pourrait permettre d'obtenir d'intéressantes informations sur la distance minimale d'un code fonctionnel sur une surface. 
Ce problème reste ouvert, nous en discuterons de nouveau page \pageref{openprob}.

\subsubsection*{Étude de l'orthogonal d'un code fonctionnel}

Le chapitre \ref{chaporth} explore la voie ouverte par le chapitre \ref{chapdiff}, à savoir l'étude de cette nouvelle classe de codes que sont les orthogonaux de codes fonctionnels sur une surface.
La première section de ce chapitre se place en fait dans un contexte plus général, celui des variétés de dimension quelconque.
Son objectif est de minorer la distance minimale de l'orthogonal d'un code fonctionnel sur une telle variété à l'aide de méthodes d'algèbre linéaire.
On obtient un résultat de minoration.

\medbreak

\noindent \textbf{Théorème \ref{minor}.}
\textit{ On suppose $N$ supérieur ou égal à $2$.
  Soit $m$ un entier tel que $G\sim mL_X$, alors
  \begin{enumerate}
  \item  la distance minimale $d^{\bot}$ du code $C_{L,X} (\Delta, G)^{\bot}$ vérifie $$d^{\bot} \geq m+2$$ et il y a égalité si et seulement si le support de $\Delta$ contient $m+2$ points alignés;
  \item sinon, si le support de $\Delta$ ne contient pas $m+2$ points alignés, alors $$d^{\bot}\geq 2m+2$$ et il y a égalité si et seulement si le support de $\Delta$ contient $2m+2$ points sur une même conique plane.  
  \end{enumerate}}

\medbreak

On conclut cette première section en donnant quelques applications de ce résultat. On montre par exemple que si $X$ est une courbe plane, alors pour certaines valeurs de $m$, la borne fournie par le théorème \ref{minor} (\ref{minor1}) est meilleure que la distance construite\footnote{Le terme de ``distance construite'' a été choisi par l'auteur comme traduction de \textit{designed minimal distance}.} de Goppa (\cite{sti} def II.2.4).

\medbreak

La deuxième section du chapitre \ref{chaporth} présente une méthode de minoration de la distance minimale de l'orthogonal d'un code fonctionnel sur une surface, sous réserve de disposer d'un résultat ``à la Bertini'' que l'on énonce.
Cette partie ne fournit donc pas de résultat à proprement parler mais motive un problème ouvert que l'on énoncera à la fin de cette introduction (voir question  \ref{qpoon}G page \pageref{qbertintro}).

\subsubsection*{Codes LDPC et décodage itératif}

Le chapitre \ref{chapldpc} porte sur l'étude de certains codes fonctionnels construits sur des surfaces. Cette question a déjà été abordée par Voloch et Zarzar dans \cite{agctvoloch}.

Le chapitre commence par une série de prérequis concernant les codes LDPC (\textit{Low Density Parity Check}, ce sont les codes admettant une matrice de parité \textit{creuse}). On y rappelle les notions de graphe de Tanner et présente un algorithme de décodage itératif.

Dans un second temps on étudie la possibilité de construire une matrice de parité creuse pour certains codes fonctionnels sur des surfaces et on applique à ces codes l'algorithme de décodage itératif présenté en première partie de chapitre.
Ce chapitre présente un volet plus expérimental de ce travail de thèse, en décrivant des calculs effectués avec le logiciel \textsc{Magma}.

\section*{Problèmes ouverts}\label{openprob}

Dans ce qui précède, nous avons signalé à plusieurs reprises l'existence de problèmes ouverts posés par ce travail de thèse. Nous concluons cette introduction en énonçant les plus importants.

\subsection*{Sur l'orthogonal d'un code fonctionnel}

Dans le chapitre \ref{chapreal}, on montre que sous certaines hypothèses sur la surface $S$ et le diviseur $G$, l'orthogonal du code fonctionnel se réalise comme somme de codes différentiels.
On remarque ensuite par l'étude d'un exemple (page \pageref{P1bu}) que les conditions que doivent vérifier $S$ et $G$ dans l'énoncé du théorème de réalisation sont suffisantes mais pas nécessaires.

\medbreak

\noindent \textbf{Question \ref{hypless}.}
\textit{Le résultat du théorème de réalisation (théorème \ref{thmreal}) reste-t-il vrai si l'on élimine les hypothèses sur $S$ et $G$ dans l'énoncé?}

\medbreak

Une autre question naturelle se pose concernant le théorème de réalisation, ou plutôt le corollaire \ref{correal}.

\medbreak

\noindent \textbf{Question \ref{bornsom}.}
\textit{Sous les conditions du corollaire \ref{correal}, peut-on estimer le nombre minimal de codes différentiels dont la somme est égale à l'orthogonal d'un code fonctionnel en fonction d'invariants géométriques de la surface?}

\medbreak

\subsection*{Sur les théorèmes ``à la Bertini''}

Une question majeure est posée à la fin du chapitre \ref{chapreal} et une variante de cette dernière est posée à la fin du chapitre \ref{chaporth}.
Une réponse à ce problème pourrait fournir des minorations de la distance minimale de codes fonctionnels construits sur des surfaces et d'orthogonaux de tels codes.

\medbreak

\noindent \textbf{Question \ref{qpoon}} (Arithmétique)\textbf{.}
\textit{Soient $X$ une variété projective lisse géométriquement intègre sur un corps fini $\F_q$ et $P_1, \ldots, P_n$, une famille de points fermés de $X$.
Peut-on évaluer explicitement ou majorer de façon précise le plus petit entier $d$ tel qu'il existe au moins une hypersurface définie sur $\F_q$ de degré inférieur ou égal à $d$ qui interpole tous les $P_i$ et dont l'intersection schématique avec $X$ soit une sous-variété lisse géométriquement intègre de codimension $1$?}

\medbreak

\noindent \textbf{Question \ref{qpoon}}\label{qbertintro} (Géométrique)\textbf{.}
\textit{  Soit $X$ une variété projective irréductible lisse définie sur $\overline{\F}_q$ et $P_1, \ldots, P_n$ une famille de points de $X$. Peut-on évaluer explicitement ou majorer de façon précise le plus petit entier $d$ tel qu'il existe au moins une hypersurface $H$ de degré inférieur ou égal à $d$ contenant tous les $P_i$ et telle que $H\cap X$ soit une sous-variété lisse de codimension un de $X$?}

\medbreak

Une présentation plus complète des questions et problèmes ouverts posés par cette thèse sera faite dans la conclusion page \pageref{chapconclu}.


\newpage
\thispagestyle{empty}
\null

\part{R\'esidus}

\chapter{R\'esidus de 2-formes sur une surface}\label{chapres}

\begin{flushright}
\begin{tabular}{p{8cm}}
{\small \textbf{Résidu.} n.m (lat. residuum). Matière qui subsiste après une opération physique ou chimique, un traitement industriel etc... Syn. Débris, déchet, rebut, reste.}

\end{tabular}
\end{flushright}

Ce chapitre est relativement différent de ceux qui vont suivre.
Il est en effet le seul dont le contenu ne soit pas directement relié à la théorie des codes correcteurs d'erreurs.
L'objectif est de fournir le matériel théorique nécessaire à la construction et l'étude de codes différentiels construits sur des surfaces algébriques.

La notion centrale de ce premier chapitre est celle de résidu. 

\section{Notations}
 
Soit $X$ une variété algébrique définie sur un corps $k$, on note $k(X)$ le corps des fonctions rationnelles sur $X$. De même, on note $\Omega^i_{k(X)/k}$ le $k(X)$-espace vectoriel des $i$-formes différentielles rationnelles sur $X$.
Soit $Y$ une sous-variété irréductible de $X$, on dira qu'une fonction (resp. une forme différentielle) rationnelle sur $X$ est \textit{régulière au voisinage de} $Y$, si et seulement si elle est régulière sur un ouvert dont l'intersection avec $Y$ est non vide\footnote{Dans le langage des schémas, cela revient à dire que la fonction (resp. la forme différentielle) est régulière au voisinage du point générique de $Y$.}.
L'anneau local des fonctions régulières au voisinage de $Y$ et son idéal maximal sont respectivement notés $\mathcal{O}_{X,Y}$ et $\mathfrak{m}_{X,Y}$.
On rappelle que le corps résiduel de cet anneau est le corps $k(Y)$ des fonctions rationnelles sur $Y$.
Si $u$ est un élément de $\mathcal{O}_{X,Y}$, on note $u_{|Y}$ sa restriction à $Y$.
Si par ailleurs il n'y a pas d'ambiguïté concernant la sous-variété $Y$ le long de laquelle on restreint notre fonction cette restriction pourra être notée $\bar{u}$.
Enfin, le complété $\mathfrak{m}_{X,Y}$-adique de cet anneau est noté $\widehat{\mathcal{O}}_{X,Y}$

\section{Cadre}\label{cadre}

Dans ce chapitre, sauf mention contraire, $k$ désigne un corps quelconque (donc de caractéristique quelconque) et $S$ une surface algébrique quasi-projective \textbf{lisse} géométriquement intègre\footnote{C'est-à-dire que sur tout ouvert affine $U$ de $S$, l'anneau de coordonnées de $U\times_k \bar{k}$ est intègre. En d'autres termes, la surface $S$ est absolument réduite et absolument irréductible.}  définie sur $k$. De plus, sauf mention contraire, $C$ désigne une courbe irréductible absolument réduite définie sur $k$ et plongée dans $S$ et $P$ un point rationnel lisse de $C$.
Notons que, comme $S$ est supposée lisse, $C$ est non contenue dans le lieu singulier de cette surface.
Par conséquent, l'anneau $\mathcal{O}_{S,C}$ est de valuation discrète.
De plus, la valuation $\mathfrak{m}_{S,C}$-adique de cet anneau s'étend en une valuation discrète $\textrm{val}_C$ sur $k(S)$.

\subsection*{Sur la notion de variété}

Dans toute cette thèse, nous parlerons de \textit{variétés}, or il s'avère que ce terme n'est pas réellement standard. Il est donc nécessaire de commencer par fixer une définition de cette notion.

\begin{defn}
  Une variété $X$ sur un corps $k$ est un schéma noethérien de type fini sur $k$.
\end{defn}

Pour les définitions de schéma noethérien et de type fini voir \cite{H} II.3.

\section{Résidus en codimension 1  et 2}\label{secres}

Il est signalé dans l'introduction, qu'en dimension supérieure à $1$, différents objets portent le nom de \textit{résidu} dans la littérature.
Nous allons introduire ces objets et étudier les relations qui les relient.
La définition de résidu la plus simple à introduire est celle de résidu en codimension $1$.
Rappelons que l'on se place sous les hypothèses énoncées en section \ref{cadre}.

\begin{prop}\label{1-res_classic}
Soit $v$ une uniformisante\footnote{C'est à dire une fonction de valuation $1$ le long de $C$.} de l'anneau $\mathcal{O}_{S,C}$.
Soit $\omega$ une $2$-forme rationnelle de valuation supérieure ou égale à $-1$ le long de $C$. Alors, il existe $\eta_1 \in \Omega^1_{k(S)/k}$ et $\eta_2 \in \Omega^2_{k(S)/k}$, toutes deux régulières au voisinage de $C$ et telles que
\begin{equation}\label{decomp}
\omega=\eta_1\wedge \frac{\displaystyle dv}{\displaystyle v}+\eta_2 .
\end{equation}

\noindent De plus, la forme différentielle  ${\eta_1}_{|C}\in \Omega_{k(C)/k}^1$ est unique et ne dépend ni du choix de l'uniformisante $v$ ni du choix de la décomposition (\ref{decomp}).
\end{prop}

\begin{defn}\label{def-1-res_classic}
On appelle cette $1$-forme sur $C$ le $1$-résidu de $\omega$ le long de $C$ et on la note
$$\res^1_{C}(\omega):={\eta_1}_{|C}.$$
\end{defn}

Un analogue de la proposition \ref{1-res_classic} est énoncé et démontré dans \cite{bhpv} au début de la section II.4. Notons que ladite référence se place dans un cadre sensiblement différent, à savoir celui des formes holomorphes sur les variétés complexes.
Toutefois, la preuve d'invariance ne fait en aucun cas appel à des propriétés spécifiques des variétés complexes. Elle s'étend de fait aisément au cadre dans lequel nous travaillons.
Nous donnerons en section \ref{resformel} une preuve de cette proposition-définition dans un contexte plus général (voir lemme \ref{valCV}).

\begin{defn}\label{2-res_classic}
Sous les hypothèses de la proposition \ref{1-res_classic}, soit $P$ un point $k$-rationnel lisse de $C$. Le $2$-résidu de $\omega$ en $P$ le long de $C$ est le résidu en $P$ du $1$-résidu de $\omega$ le long de $C$. On le note
$$
\res^2_{C,P}(\omega):= \res_P (\res^1_C (\omega)).
$$  
\end{defn}

\begin{rem}\label{corps}
Étant donné que la $2$-forme $\omega$ est $k$-rationnelle sur $S$, que la courbe $C$ est définie sur $k$ et que $P$ est un point $k$-rationnel de $C$, ce $2$-résidu est un élément de $k$.  
\end{rem}

Notons que, comme le corps de base n'est pas supposé algébriquement clos, il peut sembler logique de se placer dans un cadre plus général, à savoir que $P$ est un point fermé lisse de $C$.
Cependant, la motivation de ce chapitre est d'aboutir à des formules de sommation de $2$-résidus, dont l'une (le théorème \ref{FR3}) peut être vue comme une version en dimension $2$ de la formule des résidus bien connue en dimension $1$. 
Pour parvenir à ces formules, nous avons trouvé plus confortable d'adopter une approche géométrique. Ainsi, dans la section \ref{sommation} qui concerne ces formules de sommation, le corps de base est supposé algébriquement clos.

D'un autre côté, nous aurons tout de même besoin dans les chapitres suivants d'un résultat de type arithmétique, à savoir la remarque \ref{corps}.
En effet, l'objectif étant de construire des codes par évaluation de résidus, il faut s'assurer que les mots de code construits sont bien à coefficients dans un corps fixé.

Ainsi, le compromis adopté est le suivant. Étant donné que tout point géométrique de $S$ est un point rationnel de cette surface après une certaine extension des scalaires, on travaillera toujours avec des points rationnels. Dans un second temps lorsqu'il s'agira d'énoncer des résultats de sommation, on se placera dans $S\times_k \bar{k}$ de façon à pouvoir considérer sans distinction tous les points géométriques de $S$.

Avant de passer à la section suivante, donnons quelques exemples et remarques pour commencer à développer une certaine intuition des résidus.

\begin{rem}
Dans les deux définitions précédentes on a supposé que $\omega$ n'avait pas de pôle multiple le long de $C$.
Dans ce qui va suivre, nous verrons que les $2$-résidus sont bien définis même si l'on retire cette hypothèse.
Cependant, cette condition sur la valuation de $\omega$ le long de $C$ est indispensable pour la bonne définition des $1$-résidus le long de $C$. C'est ce que montre l'exemple \ref{exval}.
\end{rem}

\begin{exmp}\label{exval}
Supposons que $S$ est le plan affine complexe $\mathbf{A}^2_{\C}$ muni d'un système de coordonnées affines $(x,y)$. Soient $C$ la droite d'équation $y=0$ et $P$ l'origine du plan affine.
Considérons la $2$-forme
$$
\omega:= x dx\w \frac{dy}{y^2}.
$$ 
Une généralisation naturelle de la notion de $1$-résidu serait d'extraire de $\omega$, la restriction à $C$ du terme en $dy/y$. Dans l'expression ci-dessus on obtiendrait un $1$-résidu nul. 
Effectuons maintenant le changement de variables, $x:=u+y$. L'expression de $\omega$ devient
$$
\omega=  (u+y) du\w \frac{dy}{y^2} =udu \w \frac{dy}{y^2} + du \w \frac{dy}{y} 
$$
et on obtiendrait dans ce cas un $1$-résidu égal à $d\bar{u}$.
\end{exmp}

\begin{rem}
Il faut insister dès à présent sur le fait que l'on ne peut pas parler de résidu d'une $2$-forme en un point mais de résidu d'une $2$-forme, \textbf{le long d'une courbe} $C$ en un point $P$. Cela peut sembler étrange, mais le calcul présenté dans l'exemple \ref{courbe_indisp} permet de se convaincre du fait que cette spécification est incontournable.   
\end{rem}

\begin{exmp}\label{courbe_indisp}
On reprend $S=\mathbf{A}^2_{\C}$ et les mêmes $C$ et $P$ que dans l'exemple \ref{exval}. Soit
$$\omega:=\frac{dx}{x}\w \frac{dy}{y}.$$
Ici nous sommes dans un cas sympathique, la $2$-forme $\omega$ n'a que des pôles simples au voisinage de $P$. 
On a
$\res^1_C(\omega)={\frac{d\bar{x}}{\bar{x}}}$ et donc
$$\res^2_{C,P}(\omega)=1.$$
À présent, posons $C':=\{x=0\}$. L'anticommutativité du produit extérieur entraîne que $\res^1_{C'}(\omega)
=-{\frac{d\bar{y}}{\bar{y}}}$. De fait, 
$$
\res^2_{C',P}(\omega)=-1.
$$
Enfin, si on appelle $C''$ la droite d'équation $\{x=y\}$, en posant $v=y-x$, on obtient,
$$
\omega=\frac{\displaystyle dx}{\displaystyle x}\w \frac{\displaystyle dv}{\displaystyle v+x}.
$$
On développe alors en série de Laurent en la variable $v$,
$$
\omega=\frac{dx}{x} \w \frac{dv}{x\left(1+\frac{v}{x} \right)}=
\left(1-\frac{v}{x}+\frac{v^2}{x^2}-\cdots \right) \frac{dx}{x^2}\w dv.
$$
Par conséquent, il n'y a pas de terme en $\frac{dv}{v}$, donc $\res^1_{C''}(\omega)=0$ et
$$
\res^2_{C'',P}(\omega)=0.
$$  
\end{exmp}

Ce dernier exemple motive les constructions introduites dans la section suivante.
En effet, le calcul effectué correspond à un développement du coefficient de cette $2$-forme en une série de Laurent appartenant à $k((x))((v))$.
Par ailleurs, les séries de Laurent étant l'objet utilisé en théorie des courbes algébriques pour calculer des résidus il semble naturel d'en introduire une généralisation en dimension $2$.

\section{Complétions et séries de Laurent en deux variables}\label{secseclolo}

\subsection{Problématique}
En un point $k$-rationnel lisse $Q$ d'une courbe algébrique $X$, il est aisé de décrire le complété $\mathfrak{m}_{X,Q}$-adique de $k(X)$. Il s'identifie au corps des séries de Laurent $k((u))$, où $u$ est un paramètre local en $Q$.
Ici, le fait que $k(X)$ contienne le corps résiduel de $\widehat{k(X)}$, à savoir $k$, permet d'obtenir un unique plongement $k(X)\hookrightarrow k((T))$ envoyant $u$ sur $T$.
On dispose en particulier d'une méthode explicite pour décomposer une fonction en séries de Laurent en la variable $u$, et calculer le résidu d'une $1$-forme en $Q$.

Dans le cas d'un corps de fonctions de dimension $2$, la situation se complique lourdement. Si $Y$ est une surface irréductible sur $k$,
les anneaux de valuation discrète de $k(Y)$ sont ses sous-anneaux de la forme $\mathcal{O}_{Y,C}$, où $C$ est une courbe irréductible absolument réduite contenue dans le complémentaire du lieu singulier de $Y$ (ou d'une surface birationnelle à $Y$).
Soit $C$ une telle courbe et $v$ une uniformisante de $\OC$. Le corps résiduel de cet anneau local est le corps $k(C)$ des fonctions $k$-rationnelles sur $C$.
De fait, l'anneau $\OC$ est de valuation discrète, de même caractéristique que son corps résiduel (ils contiennent tous deux $k$) et contient un corps.
D'après le théorème de structure de Cohen (voir \cite{eis} théorème 7.7 ou \cite{cohen} théorème 9 pour une référence historique), l'anneau $\hOC$ est isomorphe à $k(C)[[v]]$ et le complété $\mathfrak{m}_{Y,C}$-adique de $k(Y)$ est isomorphe à $k(C)((v))$.
Cette description peut sembler commode, elle a toutefois un défaut qui la rend difficile à exploiter: en général $k(Y)$ ne contient pas $k(C)$. L'exemple suivant illustre ce phénomène.

\begin{exmp} Soient $S=\mathbf{P}^2_k$ et $C\subset S$ une courbe elliptique. Alors, il existe $x\in k(S)$ telle que $k(C)$ est une extension quadratique de $k(x)$ et $k(S)$ une extension transcendante pure de $k(x)$.
D'après le théorème de Lur\"oth, $k(S)$ ne peut pas contenir $k(C)$.   
\end{exmp}

Une autre approche consiste à considérer l'anneau local $\OP$ et à le compléter $\mathfrak{m}_{S,P}$-adiquement.
Si l'on se donne un système de coordonnées locales $(u,v)$ en $P$, l'anneau $\hOP$ est isomorphe à $k[[u,v]]$ et la décomposition en série de Taylor d'un élément de $\OP$ est explicitement calculable (voir \cite{sch1} II.2.2).
Malheureusement, si le corps des fractions de $k[[t]]$ est isomorphe à $k((t))$, on ne dispose pas d'une description aussi agréable du corps des fractions de $k[[u,v]]$.
Ces constatations motivent le travail qui va être effectué dans cette section. Il s'agit de plonger les complétés $\mathfrak{m}_P$ et $\mathfrak{m}_C$-adiques de $k(S)$ dans un corps ``plus gros''.
Deux approches vont être proposées. Moralement, la première utilise la structure de $\OP$ et la seconde celle de $\OC$.

\subsection{Développements en séries de Laurent, première approche}\label{seclolo}
Rappelons que $C$ est supposée être une courbe irréductible sur $k$ plongée dans $S$ et $P$ un point rationnel lisse de $C$.
Dans la section \ref{secres}, nous avons vu que les $2$-résidus d'une forme différentielle dépendaient d'une courbe et d'un point de celle-ci. De fait nous allons introduire un type de système de coordonnées locales \textit{relié} à $P$ et $C$.

\newpage

\begin{defn}[$(P,C)$-paires fortes]\label{PCfort}
On dit qu'une paire $(u,v)$ d'éléments de $\OP$ est une $(P,C)$-paire forte, si elle vérifie les deux conditions suivantes.
\begin{enumerate}
\item Le couple $(u,v)$ est un système de coordonnées locales en $P$.
\item La fonction $v$ est une équation locale de $C$ au voisinage de $P$.
\end{enumerate}
\end{defn}

\begin{lem}\label{lolo}
Soit $(u,v)$ une $(P,C)$-paire forte, alors il existe un morphisme $\phi: k(S) \hookrightarrow k((u))((v))$ qui envoie $u$ et $v$ sur eux-mêmes et tel que l'image de $\OP$ est contenue dans $k[[u,v]]$ et celle de $\OC$ dans $k((u))[[v]]$.  
\end{lem}

\begin{rem}
  La proposition \ref{lolo2} de la section \ref{seclolo2} entraînera qu'un tel morphisme est unique. 
\end{rem}

\begin{proof}
Comme $k(S)$ est le corps des fractions de $\OC$, il suffit de montrer l'existence d'un morphisme $\phi_0: \OC \hookrightarrow k((u))[[v]]$, qui envoie $u$ et $v$ sur eux-mêmes et injecte $\OP$ dans $k[[u,v]]$. Le lemme s'en déduira en appliquant la propriété universelle des corps de fractions.

Commençons par montrer que $\OC$ est isomorphe à ${\OP}_{(v)}$. 
Soit $U$ un voisinage affine de $P$ tel que $v$ soit une fonction régulière sur $U$ dont le lieu d'annulation sur cet ouvert soit exactement $C \cap U
$. Un tel ouvert existe étant donné que $v$ est une équation locale de $C$ au voisinage de $P$.
Notons $k[U]$ l'anneau des fonctions régulières sur $U$.

Les anneaux $\OP$ et $\OC$ s'identifient respectivement aux localisés $k[U]_{\mathfrak{m}_P}$ et $k[U]_{\mathfrak{m}_C}$ où $\mathfrak{m}_P$  
et $\mathfrak{m}_C$ correspondent respectivement à $P$ et $C$.
De plus, l'idéal $\mathfrak{m}_C$ est principal et engendré par $v$.
De fait, comme $\mathfrak{m}_C \subset \mathfrak{m}_P$, on a 
$$
{\OP}_{(v)} \cong {\left( k[U]_{\mathfrak{m}_P} \right)}_{\mathfrak{m}_C}
\cong k[U]_{\mathfrak{m}_C} \cong \OC.
$$

\noindent Ensuite, la complétion $\mP$-adique de $\OP$ fournit un morphisme injectif $\OP \hookrightarrow k[[u,v]]$, qui à une fonction régulière au voisinage de $P$ associe sa série de Taylor en les variables $u$ et $v$.
On considère alors le diagramme
\begin{equation}\label{diag1}
\xymatrix{\relax \OP  \ar@{^{(}->}[d] \ar[r]_{\textrm{loc}} &
\OC \ar[r]_{\textrm{comp}} \ar@{^{(}->}[d]_{\exists !}
& \hOC \ar[d]_{\exists !} \\
 k[[u,v]] \ar[r]^{\textrm{loc}} 
 & k[[u,v]]_{(v)} \ar[r]^{\textrm{comp}}
 & \widehat{k[[u,v]]}_{(v)}.
}
\end{equation}
Les deux premières flèches horizontales du carré de gauche sont des localisations. Celles du carré de droite sont des complétions $(v)$-adiques.

Pour finir il ne nous reste qu'à montrer que $\widehat{k[[u,v]]}_{(v)}$ est isomorphe à $k((u))[[v]]$.
Pour ce faire, on commence par montrer que le corps résiduel de l'anneau $\widehat{k[[u,v]]}_{(v)}$ est $k((u))$. En effet,
$$
\widehat{k[[u,v]]}_{(v)}/(v)\cong \textrm{Frac}\left( k[[u,v]]/(v)\right)\cong\textrm{Frac}\left( k[[u]]\right). 
$$

\noindent On invoque ensuite le théorème de structure de Cohen. L'anneau $\widehat{k[[u,v]]}_{(v)}$ est complet, de même caractéristique que son corps résiduel et contient un corps. Il est donc isomorphe à l'anneau $k((u))[[v]]$.
\end{proof}

\begin{rem}
Noter que les variables $u$ et $v$ ne jouent pas un rôle symétrique, par exemple la série 
$$
f:=\sum_{n=0}^{\infty} \frac{v^n}{u^n}
$$ 
est un élément de $k((u))((v))$ mais pas de $k((v))((u))$.
Cette asymétrie n'a rien de choquant étant donné que, dans la définition de $(P,C)$-paire forte, les fonctions $u$ et $v$ elles-mêmes jouent des rôles asymétriques.
\end{rem}

\subsection{Développements en séries de Laurent, seconde approche}\label{seclolo2}

Dans ce paragraphe, nous allons introduire une autre approche du développement en série de Laurent.
Pour cette nouvelle approche, nous nous placerons dans le contexte des $(P,C)$-paires faibles (voir définition \ref{PCfaible}), moins restrictif que celui des $(P,C)$-paires fortes.

La principale motivation de cette seconde construction est que si l'on prend une fonction rationnelle $f$ sur $S$, elle admet un développement en série de Laurent que l'on peut mettre sous la forme
$$
f=\sum_{n\geq l} f_i(u)v^i,
$$
où l'entier $l$ désigne la valuation $\mathfrak{m}_{S,C}$-adique de $f$. Les coefficients $f_i$ sont des éléments de $k((u))$. La série $f_{l}(\bar{u})$ est le développement $\bar{u}$-adique au voisinage de $P$ de la restriction à $C$ de la fonction $v^{-l}f$. Il s'agit donc du développement en série de Laurent en la variable $\bar{u}$ d'une fonction rationnelle sur $C$. Une question se pose: \textit{en est-il de même pour les autres coefficients $f_i$?}

Soit $(u,v)$ une $(P,C)$-paire forte. Comme nous l'avons signalé dans l'introduction de cette section, d'après le théorème de structure de Cohen, l'anneau $\hOC$ est isomorphe à $k(C)[[v]]$. Malheureusement cet isomorphisme n'est en aucun cas unique.
En effet, d'après \cite{cohen} théorème 10(c), si $k$ est de caractéristique positive, il y a une infinité de sous-corps de $\hOC$ qui sont envoyés sur le corps résiduel $k(C)$ via le morphisme de réduction modulo $\mC$.
Plus précisément, ce défaut d'unicité d'un représentant du corps résiduel est lié au fait que ce dernier n'est pas parfait.
D'une certaine manière, le choix de $u$ permet de contourner les éventuels problèmes d'inséparabilité.
De ce fait, pour utiliser le théorème de Cohen, nous allons choisir un représentant du corps $k(C)$ qui sera en un certain sens \textit{relié} à la fonction $u$.

\begin{prop}[Le corps $\mathcal{K}_u$]\label{Ku}
Soit $u\in \OC$ une fonction dont la restriction $\bar{u}$ à $C$ est un élément séparant\footnote{Voir \cite{sti} III.9 pour une définition d'\textit{élément séparant}.} de $k(C)$ au-dessus de $k$. Alors, il existe un unique sous-corps $\mathcal{K}_u$ de $\hOC$ contenant $k(u)$ et isomorphe à $k(C)$ via le morphisme de réduction modulo $\mC$.
De plus, ce corps et une extension monogène de $k(u)$ engendrée par un élément $y$ de $\hOC$.
\end{prop}

\begin{proof}
\textbf{Existence.}  
Par hypothèse, l'extension de corps $k(C)/k(\bar{u})$ est une extension finie séparable. D'après le théorème de l'élément primitif, il existe une  fonction $\bar{y}$ rationnelle sur $C$ qui engendre $k(C)$ sur $k(\bar{u})$.
D'après le lemme de Hensel, $\bar{y}$ se relève en un unique élément $y$ de $\hOC$ dont le polynôme minimal sur $k(u)$ est celui de $\bar{y}$ sur $k(\bar{u})$.
Soit $\mathcal{K}_u$, le sous-anneau de $\hOC$ engendré par $k(u)$ et $y$, c'est-à-dire
$$
\mathcal{K}_u:=k(u)[y].
$$
On obtient ainsi une copie de $k(C)$ qui contient $k(u)$ et s'envoie isomorphiquement sur $k(C)$ via la réduction modulo $\mC$.

\medbreak

\noindent \textbf{Unicité.} 
Soit $\mathcal{K}'$ un corps distinct de $\mathcal{K}_u$ et vérifiant les mêmes propriétés.
Il existe donc un élément de l'un de ces corps qui n'appartient pas à l'autre.
Supposons par exemple qu'il existe $z\in \mathcal{K}'$ tel que $z\notin \mathcal{K}_u$. La classe de $z$ modulo $\mathfrak{m}_{S,C}$ est une fonction $\bar{z}\in k(C)$. Cette dernière admet un unique relevé $z'$ dans $\mathcal{K}_u$. De fait, soit $R\in k(\bar{u})[T]$ le polynôme minimal unitaire de $\bar{z}$ au dessus de $k(\bar{u})$. Alors les éléments $z$ et $z'$ de $\hOC$ sont tous deux solution du problème suivant,
$$
\left\{\begin{array}{cccc}
Z & \equiv & \bar{z}  & \textrm{mod}\ \mathfrak{m}_{S,C} \\
R(u,Z) & = & 0 .& 
\end{array}\right.
$$
Ce problème admet une solution unique d'après le lemme de Hensel (\cite{eis} théorème 7.3) ce qui contredit l'hypothèse que $z$ n'appartient pas à $\mathcal{K}_u$.
\end{proof}


\begin{cor}\label{Ku((v))}
Soit $u$ une fonction rationnelle sur $S$ régulière au voisinage de $C$ dont la restriction $\bar{u}$ à $C$ est un élément séparant de $k(C)/k$.
Alors, toute fonction rationnelle $f$ sur $S$ admet un unique développement dans $\mathcal{K}_u ((v))$. 
\end{cor}

\begin{rem} 
En réalité, le résultat énoncé dans le corollaire \ref{Ku((v))} est valable pour tout élément du complété $\mathfrak{m}_{S,C}$-adique du corps $k(S)$. 
\end{rem}

Notons que, pour décrire ce corps $\mathcal{K}_u$ nous avons eu besoin de conditions plus faibles sur $u$ que celles qui sont exigées dans la définition de $(P,C)$-paire forte. C'est ce qui motive la définition suivante.

\begin{defn}[$(P,C)$-paires faibles]\label{PCfaible}
Une $(P,C)$-paire faible est une paire $(u,v)$ d'éléments de $\OC$ vérifiant les conditions suivantes.
\begin{enumerate}
\item La restriction de $u$ à $C$ est une uniformisante de $\mathcal{O}_{C,P}$.
\item La fonction $v$ est une uniformisante de $\OC$.
\end{enumerate}
\end{defn}

\begin{rem}
  Dans \cite{parshin}, le contexte décrit page 699 revient exactement à se donner une $(P,C)$-paire faible.
\end{rem}

\noindent Il va de soi qu'une $(P,C)$-paire forte est faible, mais la réciproque est fausse. En effet, en ce qui concerne $u$, le fait que sa restriction à $C$ soit régulière au voisinage de $P$ ne signifie pas que $u$ l'est. Quant à $v$, la condition: \textit{être une équation locale de $C$ au voisinage de $P$} est plus forte que celle d'\textit{être une uniformisante de $\OC$}. L'exemple qui suit permet de s'en convaincre.

\begin{exmp}
Supposons que $S$ soit le plan affine complexe muni de coordonnées affines $x$ et $y$. Soient $C$ la droite d'équation $y=0$ et $P$ l'origine du plan affine.
Posons
$$
u:=\frac{(x+y)(x-y)}{x} \quad \textrm{et} \quad v:=xy .
$$
Alors, le couple $(u,v)$ est une $(P,C)$-paire faible qui n'est pas forte. En effet, la fonction $u$ n'est pas régulière en $P$ et la fonction $v$ est dans $\mathfrak{m}_{S,P}^2$, elle n'est donc pas une équation locale de $C$ au voisinage de $P$.
\end{exmp}

\noindent Nous pouvons maintenant présenter le second procédé de décomposition en séries de Laurent.

\begin{prop}\label{lolo2}
Soit $(u,v)$ une $(P,C)$-paire faible, il existe un unique morphisme $\varphi: k(S) \hookrightarrow k((u))((v))$ qui envoie $\OC$ sur $k((u))[[v]]$ et envoie $u,v$ sur eux-mêmes. 
\end{prop}

\begin{proof}
\textbf{Existence.}
Tout comme dans la preuve du lemme \ref{lolo}, il suffit de prouver l'existence d'un morphisme $\varphi_0: \OC \hookrightarrow k((u))[[v]]$ envoyant $u$ et $v$ sur eux-mêmes, puis d'appliquer la propriété universelle des corps de fractions.
La courbe $C$ est supposée absolument réduite. Donc, d'après \cite{mumford} proposition II.4.4 (i), l'extension $k(C)/k$ est séparable, donc admet une base de transcendance séparante.
Par ailleurs, la fonction $\bar{u}$ est une uniformisante de $\mathcal{O}_{C,P}\subset k(C)$, donc sa différentielle $d\bar{u}\in \Omega^1_{k(C)/k}$ est non nulle et d'après \cite{bou} V.16.7 théorème 5, c'est un élément séparant de $k(C)/k$.

D'après le corollaire \ref{Ku((v))}, on dispose d'une injection $\OC \hookrightarrow \mathcal{K}_u [[v]]$ et $\mathcal{K}_u$ est isomorphe à $k(C)$ via le morphisme de réduction modulo $\mC$. De plus, comme $\bar{u}$ est une uniformisante de $\OP$, le complété $\mathfrak{m}_{C,P}$-adique de $k(C)$ est isomorphe à $k((\bar{u}))$. On dispose donc d'une injection $\mathcal{K}_u \hookrightarrow k((u))$ qui s'étend coefficient par coefficient en un morphisme $\mathcal{K}_u [[v]] \hookrightarrow k((u))[[v]]$. On en déduit l'existence de l'application $\varphi_0: \OC \hookrightarrow k((u))[[v]]$ recherchée.

\medbreak

\noindent \textbf{Unicité.}
Soit $\varphi_0': \OC \rightarrow k((u))[[v]]$, un autre morphisme d'anneaux envoyant $u$ et $v$ sur eux-mêmes. Nous allons montrer que le diagramme suivant est commutatif.
$$
\xymatrix{
\relax
\OC \ar[r] \ar[rd]_{\varphi_0'} \ar@/^2pc/[rrr]^{\varphi_0}&
\hOC \ar[r]^-{\sim} &
\mathcal{K}_u [[v]] \ar[r]^-{r} & k((u))[[v]] \ar@{.>}[lld]^{\textrm{id}}\\ 
 & k((u))[[v]]
}
$$
Comme $\varphi_0'$ envoie $v$ sur lui-même, on en déduit que c'est un morphisme local non ramifié. La propriété universelle du complété, implique l'existence et l'unicité d'un morphisme $\hat{\varphi}_0'$ qui fait commuter le diagramme suivant.
$$
\xymatrix{
\relax
\OC \ar[r] \ar[rd]_{\varphi_0'} \ar@/^2pc/[rrr]^{\varphi_0}&
\hOC \ar[r]^-{\sim}  \ar[d]^{\hat{\varphi}_0'}&
\mathcal{K}_u [[v]] \ar[ld]^{r'} \ar[r]^-{r} & k((u))[[v]] \\ 
 & k((u))[[v]]
}
$$
Le morphisme $r'$ est la composée du morphisme inverse de $\xymatrix{ \relax\hOC \ar[r]^-{\sim}& \mathcal{K}_u((v))}$ et de $\hat{\varphi}_0$.
Il reste à montrer que $r=r'$. Un morphisme local de $\mathcal{K}_u [[v]]$ dans $k((u))[[v]]$ est entièrement déterminé par les images de $u$, $v$ et $y$. Il suffit donc de montrer que $r(y)=r'(y)$.
Remarquons dès à présent que, d'après la construction de $\varphi_0$ et donc de $r$, on a $r(y)=\psi(u)$, où $\psi(\bar{u})$ est le développement en série de Laurent en $P$ de $\bar{y}\in k(C)$.
Soit $F\in k(\bar{u})[T]$, le polynôme minimal unitaire de $\bar{y}$ sur $k(\bar{u})$. L'élément $y$ de $\hOC$ vérifie $F(u,y)=0$ et comme $r$ et $r'$ sont des morphismes d'anneau, on en déduit
$$
F(u,r(y))=0 \quad  \textrm{et} \quad F(u,r'(y))=0
\quad \textrm{dans}\ k((u))[[v]].
$$
De plus, par passage au quotient modulo $v$, on a
$$
r(\bar{y})\equiv r'(\bar{y}) \equiv \psi(\bar{u})\ \mod (v).
$$
Ainsi $r(y)$ et $r'(y)$ sont tous deux solution du problème suivant.
$$
\left\{
\begin{array}{cccc}
F(u,Z) & = & 0 & \\
Z & \equiv & \psi(u) & \mod \ (v).
\end{array}
\right.
$$
D'après le lemme de Hensel, ce problème admet une unique solution qui est $\psi(u)$. Ce dernier étant égal à $r(y)$, cela conclut la preuve.
\end{proof}

\begin{rem}
  La principale différence entre les résultats de ce chapitre et ceux de la première partie de \cite{parshin} est que ce dernier suppose que le corps de base est parfait, alors que nous ne nous sommes donnés aucune restriction sur $k$ dans ce chapitre. On trouve dans cet article la démonstration d'un énoncé analogue à celui de la proposition \ref{lolo2}. Cette dernière se trouve de fait simplifiée grâce à cette hypothèse supplémentaire sur $k$.
\end{rem}

Ainsi, nous avons montré que si $(u,v)$ est une $(P,C)$-paire forte, les deux approches fournissent les mêmes développements en série de Laurent. Par ailleurs nous avons obtenu une réponse à la question posée à la fin de la section \ref{seclolo}. Cela donne lieu au corollaire suivant.

\begin{cor}\label{coeffts}
Soit $(u,v)$, une $(P,C)$-paire faible.
Alors, toute fonction $f\in k(S)$ admet un unique développement en séries de Laurent
$$
f=\sum_{j\geq l}f_j(u)v^j\ \in k((u))((v)).
$$  
De plus, pour tout $j\geq l$, la série de Laurent $f_j(\bar{u})$ est une fonction rationnelle sur $C$.
\end{cor}



\subsection{Changement de variables}

Les séries de Laurent ont été introduites de façon à montrer que l'on peut définir le $2$-résidu d'une $2$-forme $\omega \in \Omega^2_{k(S)/k}$ en $P$ le long de $C$, sans aucune condition sur la valuation de $\omega$ le long de $C$.
Nous allons donc donner une définition générale des $2$-résidus en utilisant les séries de Laurent.
Ensuite, il faudra prouver que cet objet ne dépend pas du choix d'une $(P,C)$-paire.
C'est la raison pour laquelle nous devons introduire les changements de $(P,C)$-paires.

\begin{lem}\label{lemcv}
  Soient $(u,v)$ et $(x,y)$ deux $(P,C)$-paires faibles. Les fonctions $u$ et $v$ se décomposent en séries de Laurent en les variables $x$ et $y$ et leurs développements sont de la forme suivante
\begin{equation}\label{cv}\tag{CV}
\left\{
\begin{array}{rclcccl}
u & = & f(x,y) &\ \textrm{avec} & 
\ f(x,0) & \in &  xk[[x]]\smallsetminus x^2k[[x]]\\
v & = & g(x,y) &\ \textrm{avec} & 
\ g(x,y) & \in & yk((x))[[y]]\smallsetminus y^2k((x))[[y]].     
\end{array}
\right.
\end{equation}
De plus, si $(u,v)$ et $(x,y)$ sont des $(P,C)$-paires fortes, alors $f$ et $g$ sont des séries de Taylor, c'est-à-dire des éléments de $k[[x,y]]$.
\end{lem}

\begin{proof}
  Les fonctions $u$ et $v$ sont des éléments de $\OC$. D'après la proposition \ref{lolo2}, leurs développements respectifs en séries de Laurent $f(x,y)$ et $g(x,y)$ sont dans $k((x))[[x]]$.
De plus, si $(u,v)$ et $(x,y)$ sont des $(P,C)$-paires fortes, alors ces fonctions sont des éléments de $\OP$.
Or, d'après le lemme \ref{lolo}, les fonctions $\bar{u}$ et $\bar{x}\in k(C)$ sont toutes deux des uniformisantes de $\mathcal{O}_{C,P}$, donc $\bar{u}=f(\bar{x},0)$ est une série de Taylor de valuation $(\bar{x})$-adique $1$.
Les fonctions $v$ et $y$ sont de valuation $\mC$-adique $1$ le long de $C$, donc leur quotient $v/y$ est un inversible de $\OC$. Par conséquent, $G(x,y):=v/y$ est un élément de $k((x))[[y]]$ est de valuation $(y)$-adique nulle. Ainsi, comme $g=yG$, on en déduit que $g$ est de valuation $(y)$-adique $1$.
\end{proof}

Avant de passer à la suite, faisons un courte remarque sur ce changement de variables. Soient $(u,v)$ et $(x,y)$ deux $(P,C)$-paires fortes, on dispose donc d'un changement de variables de la forme (\ref{cv}),
$$
\left\{
\begin{array}{rcl}
u & = & f(x,y) \\    
v & = & g(x,y).
\end{array}
\right.
$$ 

\noindent De plus, les séries $f$ et $g$ sont des séries de Taylor et vérifient
\begin{equation}\label{condf}
f=\sum_{i,j \geq 0} f_{i,j}x^i y^j \ \ \textrm{avec} \ \ f_{0,0}=0 \ \ 
\textrm{et}\ \ f_{1,0}\neq 0
\end{equation}
et
\begin{equation}\label{condg}
g=\sum_{i,j \geq 0} g_{i,j}x^i y^j \ \ \textrm{avec} \ \ \forall k \in \N,\ g_{k,0}=0 \ \ 
\textrm{et}\ \ g_{0,1}\neq 0.
\end{equation}

\noindent De toutes la assertions ci-dessus seule ``$g_{0,1}\neq 0$'' n'est pas complètement évidente. Supposons que $g_{0,1}=0$, alors, comme $g_{k,0}$ est nul pour tout entier naturel $k$, on en déduit que $g(x,y)$ est dans l'idéal $((x,y))^2$, ce qui contredit le fait que la paire $(u,v)$ est une $(P,C)$-paire forte.

Regardons à présent la matrice jacobienne de ce changement de variables.
$$
\textrm{Jac}\left( \frac{f,g}{x,y}\right)=
\left(
\begin{array}{cc}
\frac{\p f}{\p x}(0,0) & \frac{\p f}{\p y}(0,0) \\
\frac{\p g}{\p x}(0,0) & \frac{\p g}{\p y}(0,0)
\end{array}
\right)
=
\left(
\begin{array}{cc}
f_{1,0} & f_{0,1}\\
g_{1,0} & g_{0,1}
\end{array}
\right)
=
\left(
\begin{array}{cc}
f_{1,0} & f_{0,1}\\
0 & g_{0,1}
\end{array}
\right).
$$
D'après (\ref{condf}) et (\ref{condg}), le produit $f_{1,0}g_{0,1}$ est non nul, donc que cette matrice est inversible, ce qui est normal puisque $(u,v)$ est un système de coordonnées locales. 
On voit ainsi que les changement de $(P,C)$-paires fortes sont des changements de variables dont la jacobienne en $P$ est triangulaire supérieure et inversible.
On peut donner une interprétation géométrique à ce fait.
Une matrice triangulaire supérieure est la matrice d'un endomorphisme qui préserve un drapeau.
Le changement de variables (\ref{cv}) préserve le \textit{drapeau géométrique} $(P,C)$.
\subsection{Objets rationnels et formels}\label{ratformel}

Dans la section \ref{secdefgen}, nous manipulerons fréquemment des séries de Laurent.
Cependant, l'objectif de ce travail n'est pas d'obtenir des résultats sur les séries formelles mais sur des objets géométriques, en l'occurrence les $2$-formes rationnelles sur $S$.
Aussi, les séries de Laurent ne sont qu'un outil pour arriver à nos fins.
Elles nous permettront de traduire certains problèmes géométriques sous forme de problèmes purement combinatoires. 
Dans ce qui suit, outre les séries de Laurent nous allons manipuler de formes différentielles formelles, c'est à dire des éléments des espaces de différentielles relatives $\Omega^i_{k((u))((v))/k}$. On renvoie le lecteur au chapitre IX de \cite{matsumura} pour une définition de ces espaces. Le lemme qui suit nous permet d'obtenir une description agréable des ces modules de différentielles relatives.

\begin{lem}\label{reldif}
Soit $(u,v)$ une $(P,C)$-paire faible,
on a les isomorphismes
$$
\begin{array}{ccclc}
& \Omega^i_{k((u))((v))/k} & \cong & \Omega^i_{k(S)/k} \otimes_{k(S)} k((u))((v)), & \textrm{pour} \ i\in \{1,2\}\\
 & & & & \\
\textrm{et} & \Omega^1_{k((u))/k} & \cong & \Omega^1_{k(C)/k} \otimes_{k(C)} k((u)).
\end{array}
$$  
\end{lem}

\begin{proof}
  Voir annexe \ref{ann_dif_lolo}.
\end{proof}


\begin{lem}\label{defcv}
Soient $(u,v)$ deux éléments de $k((x))((y))$ liés aux variables $(x,y)$ par un changement de variables de la forme\footnote{Voir lemme \ref{lemcv}.} (\ref{cv}).
Alors, ce changement de variables induit un isomorphisme de corps locaux $k((u))((v)) \rightarrow k((x))((y))$.
C'est-à-dire qu'il envoie un série de Laurent de valuation $(v)$-adique $m\in \Z$ sur une série de valuation $(y)$-adique $m$. De même, il induit un isomorphisme $\Omega^2_{k((u))((v))/k}\rightarrow \Omega^2_{k((x))((y))/k}$ qui préserve les valuations.
\end{lem}

\begin{proof}
Voir annexe \ref{ann_topo}.
\end{proof}

\section{Définition générale des résidus}\label{secdefgen}

En utilisant les notions introduites dans la section \ref{secseclolo}, nous allons pouvoir donner une définition plus générale de résidus. 

\subsection{Invariance des $2$-résidus}\label{resformel}

Dans ce qui suit nous allons travailler exclusivement avec des objets formels. Ensuite, en section \ref{secgeom}, on appliquera les résultats obtenus dans le cadre formel aux différentielles rationnelles.
Noter que, le but étant d'obtenir des informations sur les $2$-formes rationnelles, nous aurions pu énoncer un résultat géométrique.
Cependant, la preuve du théorème \ref{inv2res}, qui est le point clé de cette section, consiste uniquement en des manipulations sur les coefficients de séries formelles.
Surtout, nous aurons absolument besoin de la version formelle de ce résultat pour démontrer la proposition \ref{1resprop2} (voir section \ref{secgeom}). C'est pourquoi nous avons choisi de l'énoncer dans ce contexte.

\begin{nota}\label{indices}
Dans tout ce qui suit, lorsque nous aurons affaire à une série de Laurent $f\in k((u))((v))$ ou $k((x))((y))$, nous adopterons le système d'indices suivant. L'indice ``$i$'' sera lié à la première variable ($u$ ou $x$) et l'indice ``$j$'' à la seconde ($v$ ou $y$). De fait, $f$ s'écrit,
$$
f=\sum_{j\geq l}f_j(u)v^j, \ \ \textrm{avec}\ \ f_j(u)=\sum_{i\geq l_j}f_{i,j}u^i\ \in k((u)).
$$
\end{nota}

\newpage

\begin{defn}\label{defres}
Soit $\omega=h(u,v)du \w dv$ avec $h=\sum_j h_j(u)v^i \in k((u))((v))$, une $2$-forme formelle, on définit les objets suivants.
\begin{enumerate}
\item Le $(u,v)$-$1$-résidu de $\omega$ est défini par
$$(u,v)\res^1(\omega):= h_{-1}(u)du\ \in \Omega^1_{k((u))/k}.$$
\item Le $(u,v)$-$2$-résidu de $\omega$ en $P$ le long de $C$ est défini par
$$(u,v)\res^2(\omega):=h_{-1,-1} \in k.$$
\end{enumerate}
\end{defn}

\noindent Le théorème qui suit est la clé de la définition des $2$-résidus.

\begin{thm}\label{inv2res}
Soit $(x,y)$ une paire d'éléments du corps $k((u))((v))$ liée aux fonctions $(u,v)$ par un changement de variables de la forme (\ref{cv}). Alors, pour toute $2$-forme formelle $\omega=h(u,v)du\w dv \in \Omega^2_{k((u))((v))/k}$, on a
$$
(u,v)\res^2(\omega)=(x,y)\res^2(\omega).
$$
\end{thm}

La preuve de ce théorème nécessite les lemmes \ref{lem_invCV1} et \ref{valCV} qui seront énoncés plus loin.
Tout d'abord, considérons de nouveau le changement de variables (\ref{cv}).
\begin{equation}\tag{CV}
\left\{
\begin{array}{rclcl}
u & = & f(x,y) &\ \textrm{avec} & 
\ f(x,0)\in xk[[x]]\smallsetminus x^2k[[x]]\\
v & = & g(x,y) &\ \textrm{avec} & 
\ g\in yk((x))[[y]]\smallsetminus y^2k((x))[[y]].     
\end{array}
\right.
\end{equation}

\noindent Ce changement de variables peut être appliqué en deux étapes. On peut dans un premier temps passer de $(u,v)$ à $(u,y)$ puis de $(u,y)$ à $(v,y)$. C'est-à-dire,
$$
\begin{array}{cccc}
\textrm{d'abord} &
\textrm{(CV1)}
\left\{\begin{array}{rcl}
u & = & u \\   
v & = & \gamma(u,y)
\end{array}\right.
,& \textrm{ensuite} &
\textrm{(CV2)}
\left\{
\begin{array}{rcl}
u & = & f(x,y)\\
y & = & y
\end{array}
\right.
\end{array},
$$
où $\gamma(x,y)$ de valuation $(y)$-adique $1$.
Nous allons montrer successivement que les $2$-résidus sont invariants sous l'action de (CV1), puis de (CV2).

\begin{lem}[Invariance des $1$-résidus sous l'action de (CV1)]\label{lem_invCV1}
Soit $\omega=h(u,v)du\w dv$ une $2$-forme formelle.
Pour tout $y$ lié à $(u,v)$ par un changement de variables (CV1): v=g(u,y), on a 
$$
(u,v)\res^1(\omega)=(u,y)\res^1(\omega).
$$  
\end{lem}

\begin{proof}
En appliquant (CV1), on obtient
$$
\omega=h(u,g(u,y))\frac{\p g}{\p y} du\w dy.
$$  
On peut voir le corps $k((u))((v))$ comme un corps de séries de Laurent à une variable au-dessus de $k((u))$. De ce point de vue, $y$ est une autre uniformisante de ce corps.
D'après \cite{sti} proposition IV.2.9, le coefficient en $v^{-1}$ de $h(u,v)$ est égal au coefficient en $y^{-1}$ de $h(u,g(u,y))\p g/\p y$.
\end{proof}

\begin{rem}\label{invCV1}
  Dans tout le chapitre IV de \cite{sti} , le corps de base est supposé parfait.
Or si le corps $k$ est de caractéristique positive, le corps $k((u))$ n'est pas parfait. Cependant la démonstration de la proposition IV.2.9 de cet ouvrage est purement formelle et ne nécessite en aucun cas un corps de base parfait.
\end{rem}

\noindent Nous avons donc vu que le changement de variables (CV1) n'avait pas d'influence sur les $1$-résidus.
Il n'en aura donc à fortiori pas sur les $2$-résidus.
En ce qui concerne le changement de variables (CV2), ce dernier peut avoir une influence sur les $1$-résidus, c'est ce que montrait l'exemple \ref{exval}.
Nous allons cependant montrer qu'il n'a pas d'influence sur les $2$-résidus.
Pour ce faire, nous aurons besoin du lemme qui suit. 

\begin{lem}\label{valCV}
Pour toute $2$-forme formelle $\omega=h(u,v)du\w dv\in \Omega^2_{k((u))((v))/k}$, de valuation $(v)$-adique supérieure ou égale à $-1$, on a 
$$
(u,v)\res^1(\omega)=(x,y)\res^1(\omega).
$$      
\end{lem}

\begin{rem}
La proposition \ref{1-res_classic} est une conséquence immédiate de ce lemme.
\end{rem}

\begin{proof}
  D'après le lemme \ref{invCV1}, on a $(u,v)\res^1(\omega)=(u,y)\res^1(\omega)$. 
De fait, nous n'avons qu'à étudier le comportement de $\omega$ sous l'action de (CV2). Soit donc $u=f(x,y)$. Isolons dans $\omega$ le terme en $y^{-1}$:
$$
\omega= \frac{h_{-1}(u)}{y} du\w dy + \left( \sum_{j\geq 0}h_j(u)y^j  \right)du \w dy = \omega_{-1}+\omega_+ .
$$
La forme $\omega_+$ a une valuation $(y)$-adique positive, donc d'après le lemme \ref{defcv}, le changement de variables (CV2) n'a pas d'influence sur cette valuation.
De fait, le $(x,y)$-$1$-résidu de $\omega$ est celui de $\omega_{-1}$ et
$$
\omega_{-1}=\frac{h_{-1}(f(x,y))}{y}\frac{\p f}{\p x}dx\w dy.
$$
D'après le lemme \ref{defcv}, la valuation $(y)$-adique de $h_{-1}(f(x,y))$ est égale à la valuation $(v)$-adique de $h_{-1}(u)$, à savoir $0$. Ainsi,
$$
(x,y)\res^1(\omega)=h_{-1}(f_0(\bar{x}))f'_0(\bar{x})d\bar{x}
=h_{-1}(f_0(\bar{x}))d(f_0(\bar{x})),
$$
où $f_0(x):=f(x,0)$. Cette
 $1$-forme formelle est égale à $(u,y)\res^1_{C,P}(\omega)=h_{-1}(\bar{u})d\bar{u}$. Il suffit pour s'en convaincre d'appliquer à $h_{-1}(\bar{u})d\bar{u}$ le changement de variables
$$\bar{u}=f(\bar{x},0).$$
\end{proof}

\noindent Dans la preuve du théorème \ref{inv2res} nous aurons besoin du lemme suivant dont la preuve est donnée en annexe.

\begin{lem}\label{jacob}
  Soient $A,B$ deux séries de Laurent appartenant à $k((u))((v))$. Pour toute paire de séries $(x,y)$ liée à $(u,v)$ par un changement de variables de la forme (\ref{cv}), on a
$$
(x,y)\res^2\left(dA\w dB \right)=0.
$$
\end{lem}

\begin{proof}
Voir annexe \ref{demojacob}.
\end{proof}

\noindent Nous disposons à présent de tous les outils nécessaires à la démonstration du théorème \ref{inv2res}. Dans un premier temps, nous allons le démontrer  dans le cas où la caractéristique du corps de base $k$ est nulle.

\begin{proof}[\textsc{Preuve du théorème \ref{inv2res} si $k$ est de caractéristique nulle}]\label{preuveinv}
Dans l'in\-té\-gra\-li\-té de cette preuve,
les $(x,y)$-$1$- et $2$-résidus sont toujours en $P$ le long de $C$. Aussi,
pour alléger la rédaction, nous omettrons de signaler les ``en $P$ le long de $C$''.
 
Commençons par décomposer $\omega=hdu \w dv$ en isolant les termes de valuation $(y)$-adique inférieure ou égale à $-2$.
$$
\omega=\sum_{j=-l}^{-2} h_j(u)y^j du\w dy + \sum_{j\geq -1}h_j(u)y^j du\w dy=\omega_{-}+\omega_{\textrm{inv}} .
$$
Noter que l'extraction d'un $2$-résidu est une application $k$-linéaire.
Aussi, d'après les lemmes \ref{invCV1} et \ref{valCV}, il suffit d'étudier le comportement des $2$-résidus de $\omega_{-}$ sous l'action de (CV2).
De plus, toujours pour des raisons de linéarité, on peut scinder le problème et restreindre notre étude aux $2$-formes de la forme:
$$
\omega=\phi(u)du\w \frac{dy}{y^n}\quad \textrm{avec}\quad \phi \in k((u))\quad \textrm{et}\quad n\geq 2.
$$
Le $(u,y)$-$1$-résidu de la $2$-forme formelle ci-dessus est nul, il est en donc de même pour son $(u,y)$-$2$-résidu. Avant d'appliquer (CV2), nous allons continuer à travailler $\omega$ au corps.
Isolons le terme en $u^{-1}$  de la série de Laurent $\phi \in k((u))$.
$$\phi(u)=\widetilde{\phi}(u)+\frac{\phi_{-1}}{u},\  \textrm{où}\ \widetilde{\phi}_i=
\left\{
\begin{array}{lll}
\phi_i & \textrm{si} & i\neq -1 \\
0 & \textrm{si} & i= -1.
\end{array}
\right.$$ 

\noindent Comme $k$ est supposé de caractéristique nulle, la série $\widetilde{\phi}(u)$ a une primitive formelle $\widetilde{\Phi}(u)$. De même, posons $s:=\frac{1}{(1-n)y^{n-1}}$. C'est une primitive formelle de $1/y^n$. On a alors
$$
\omega= d\widetilde{\Phi} \w ds+ \phi_{-1} \frac{du}{u}\w ds =\omega_{r}+\phi_{-1}\omega_{-1}.
$$
D'après le lemme \ref{jacob}, la forme $\omega_r$ a un $2$-résidu nul et indépendant du choix de $(u,y)$. Il reste à étudier la $2$-forme
$$
\omega_{-1}=\frac{du}{u}\w \frac{dy}{y^n}.
$$ 
En lui appliquant (CV2), on obtient
$$
\omega_{-1}=\frac{df(x,y)}{f(x,y)}\w \frac{dy}{y^n}.
$$ 

\noindent Rappelons que $f$ est de la forme: $\sum_{j\geq 0} f_j(x)y^j$ avec:
$$f_0(x)=f_{1,0}x+f_{2,0}x^2+\cdots \quad \textrm{et} \quad f_{1,0}\neq 0.$$ 
On peut donc factoriser $f_0$ en
$$
f_0(x)=f_{1,0}x\left(1+\frac{\displaystyle f_{2,0}}{\displaystyle f_{1,0}}x+\cdots \right).
$$
Posons
$$
\begin{array}{lrcll}
 & r(x) & := & \frac{\displaystyle f_{2,0}}{\displaystyle f_{1,0}}x+\frac{\displaystyle f_{3,0}}{\displaystyle f_{1,0}}x^2 +\cdots & \in k[[x]]\\
\textrm{et} & \mu(x,y) & := & \frac{\displaystyle f_1(x)}{\displaystyle f_0(x)}y+\frac{\displaystyle f_2(x)}{\displaystyle f_0(x)}y^2 +\cdots & \in k((x))[[y]].
\end{array}
$$
La série $f$ se factorise donc en
\begin{equation}\label{factf}
f(x,y)=f_{1,0}x(1+r(x))(1+\mu(x,y)). 
\end{equation}

\noindent Par ailleurs, pour toute série $S$ appartenant à $xk[[x]]$ (resp. à $yk((x))[[y]]$), on définit le logarithme formel de $1+S$ par
$$
log(1+S):=\sum_{k=0}^{+\infty} (-1)^{k+1}\frac{S^k}{k}.
$$

\noindent Cette définition a un sens puisque $k$ est supposé de caractéristique nulle. De plus, cette série converge pour la topologie $(x)$-adique (resp. $(y)$-adique). Enfin, on a
$$d\log (1+S)=\frac{d(1+S)}{1+S}.$$

\noindent En utilisant la factorisation \ref{factf}, on obtient:
$$
\omega_{-1}  =  \underbrace{\frac{ dx}{ x} \w \frac{ dy}{ y^n}}_{\mu_1} +\underbrace{d\log(1+r)\w ds}_{\mu_2} + \underbrace{d\log(1+\mu)\w ds}_{\mu_3}.   
$$
D'après le lemme \ref{jacob}, les $(x,y)$-$2$-résidus des formes $\mu_2$ et $\mu_3$ sont nuls. 
La forme $\mu_1$ a un $(x,y)$-$1$-résidu nul (elle n'a pas de terme en $dy/y$), son $(x,y)$-$2$-résidu est également nul. On en conclut que 
$$
(x,y)\res^2_{C,P}(\omega_{-1})=0.
$$
\end{proof}
\label{finpreuve}

\begin{proof}[\textsc{Esquisse de preuve du théorème \ref{inv2res} en caractéristique positive}]
La preu\-ve est semblable à celle de l'invariance des résidus de $1$-formes sur des courbes (c.f. \cite{sti} IV.2.9 ou \cite{ser} II.7 proposition 5). On montre que le $(x,y)$-$2$-résidu de $\omega$ est une expression polynomiale en certains coefficients de $f$. Ce polynôme ne dépend ni de $f$ ni du corps de base. D'après le théorème de prolongement des identités algébriques (\cite{bou} IV.3 proposition 9) et le travail effectué en caractéristique nulle, on conclut que ce polynôme est nul.
Une preuve détaillée est donnée en annexe \ref{demochiante}.
\end{proof}

Ainsi, à partir de maintenant, lorsque nous parlerons de $2$-résidus en un point le long d'un courbe, nous n'aurons plus à préciser la $(P,C)$-paire.

\subsection{Le cadre géométrique}\label{secgeom}

À présent nous allons introduire les notions de $1$- et $2$-résidus pour des $2$-formes rationnelles. Tout comme en section \ref{secres}, les $2$-résidus seront associés à un point $P$ et une courbe $C$ contenant $P$ et les $1$-résidus seront associés à une courbe $C$. Ces derniers, seront également associés à un autre paramètre si $\omega$ a un pôle multiple le long de $C$.  
Commençons par définir les $2$-résidus qui sont \textit{plus intrinsèques}.

\begin{defn}
Soient $(u,v)$ une $(P,C)$-paire faible et $\omega$ une $2$-forme rationnelle sur $S$. Il existe une fonction $h\in k(S)$ telle que $\omega=hdu\w dv$ et cette fonction admet un unique développement en série de Laurent $H(u,v)$. On appelle $2$-résidu de $\omega$ en $P$ le long de $C$ le coefficient $H_{-1, -1}$ de $x^{-1}y^{-1}$ de $H$. On le note
$$
\res^2_{C,P}(\omega):=H_{-1,-1}.
$$
\end{defn}

\noindent En d'autres termes et pour faire le lien avec ce qui précède, le $2$-résidu en $P$ le long de $C$ de $\omega$ est le $(u,v)$-$2$-résidu de $\omega$ vue comme une forme formelle.
Le travail effectué en section \ref{resformel} nous assure que l'objet est bien défini et ne dépend pas du choix de la paire $(u,v)$.

Passons maintenant à la définition de résidus en codimension $1$. Notre objectif est d'associer à une $2$-forme différentielle rationnelle $\omega$ sur $S$ une $1$-forme rationnelle $\mu$ sur $C$.

\begin{prop}\label{1resgenprop}
Soit $u$ une fonction rationnelle sur $S$ régulière au voisinage de $C$ dont la restriction $\bar{u}$ à $C$ est un élément séparant de $k(C)/k$.
Soient $\omega$ une $2$-forme rationnelle sur $S$ et $v$ une uniformisante de l'anneau $\OC$. On rappelle que, d'après la proposition \ref{Ku} et son corollaire \ref{Ku((v))}, il existe une unique série de Laurent
$$f=\sum_{j\geq -l}f_jv^j \in \mathcal{K}_u((v))\quad \textrm{telle que} \quad \omega = fdu\w dv.$$
De plus, la $1$-forme $\bar{f}_{-1} d\bar{u}$ est rationnelle sur $C$ et ne dépend pas du choix de $v$.
\end{prop}

\begin{defn}\label{1resgendef}
Sous les hypothèses de la proposition \ref{1resgenprop}, on appelle $(u)$-$1$-résidu de $\omega$ le long de $C$ la $1$-forme rationnelle sur $C$ définie par
$$
(u)\res^1_{C}(\omega):=\bar{f}_{-1}d\bar{u}.
$$ 
\end{defn}

\begin{proof}[\textsc{Preuve de la proposition \ref{1resgenprop}}]
Rappelons que $f_{-1}$ est un élément de $\mathcal{K}_u$. Ainsi, sa classe $\bar{f}_{-1}$ modulo $\mathfrak{m}_{S,C}$ est un élément de $k(C)$, la $1$-forme $\bar{f}_{-1}d\bar{u}$ est donc rationnelle. 
L'indépendance de cette $1$-forme par rapport au choix de $v$ se démontre de la même façon que le lemme  \ref{lem_invCV1}. Si $w$ est une autre uniformisante de $\OC$, d'après \cite{sti} IV.2.9, le coefficient en $v^{-1}$  de la série de Laurent $f \in \mathcal{K}_u((v))$ est égal à celui en $w^{-1}$ de $f\frac{\p v}{\p w}$.
\end{proof}

\begin{rem}
L'exemple \ref{exval} confirme la nécessité de faire intervenir la fonction $u$ dans la définition. On ne peut espérer obtenir un objet qui ne dépende que de $\omega$ et de la courbe $C$.
\end{rem}

Maintenant que nous disposons d'une définition générale de $1$-résidu, il serait souhaitable que cette dernière soit compatible avec la définition \ref{def-1-res_classic}. De plus, étant donné un point rationnel $P$ de $C$, il serait intéressant de savoir quelle relation lie le $(u)$-$1$-résidu de $\omega$ le long de $C$ et son $2$-résidu en $P$ le long de $C$. C'est le but du lemme \ref{1resprop1} et de la proposition \ref{1resprop2}.

\begin{lem}\label{1resprop1}
Sous les conditions de la proposition \ref{1resgenprop}, soit $\omega$ une $2$-forme rationnelle sur $S$ de valuation supérieure ou égale à $-1$ le long de $C$. Alors le $(u)$-$1$-résidu de $\omega$ le long de $C$ coïncide avec le résidu de $\omega$ le long de $C$ de la définition \ref{1-res_classic}.
\end{lem}

\begin{proof}
Il existe une unique série de Laurent $f$ appartenant à $\mathcal{K}_u((v))$ de valuation $-1$ telle que 
$$
\omega = fdu \w dv =\sum_{j\geq -1} f_j v^j du \w dv.
$$
Soit $\varphi$ une fonction rationnelle sur $S$ régulière au voisinage de $C$ et dont la restriction à $C$ est égale à $\bar{f}_{-1}$, on pose
$$
\eta_1:=\varphi du \quad \textrm{et} \quad \eta_2:=\omega - \varphi du \w \frac{dv}{v}.
$$
La $2$-forme $\eta_2$ est régulière le long de $C$ et $\omega$ se décompose en
$$
\omega = \eta_1 \w \frac{dv}{v}+\eta_2.
$$
D'après la définition \ref{def-1-res_classic}, la $1$-forme ${\eta_1}_{|C}$ sur $C$ est le $1$-résidu de $\omega$ le long de $C$. Or, ${\eta_1}_{|C}$ est égale à $\bar{f}_{-1}d\bar{u}$.
\end{proof}

\begin{prop}\label{1resprop2}
Sous les conditions de la proposition \ref{1resgenprop}. Soient $\omega$ une $2$-forme rationnelle sur $S$ et $P$ un point rationnel lisse de $C$, alors
$$
\res_P ((u)\res^1_C(\omega))=\res^2_{C,P}(\omega).
$$
\end{prop}

\begin{rem}
Si $\omega$ est de valuation supérieure ou égale à $-1$ le long de $C$, la proposition \ref{1resprop2} entraîne que le $2$-résidu de $\omega$ en $P$ le long de $C$ défini dans cette section coïncide avec celui de la section \ref{secres}.
\end{rem}

\begin{rem}
  La condition ``$P$ est un point lisse de $C$'' pourra être supprimée dès que l'on saura définir des $2$-résidus le long de $C$ en des points singuliers de cette courbe (voir section \ref{secsing2}).
\end{rem}

\begin{proof}[\textsc{Preuve de la proposition \ref{1resprop2}}]
Soient $P$ un point rationnel de $C$ et $v$ une uniformisante de $\OC$.
Commençons par noter que, si $\bar{u}$ est une uniformisante de $\mathcal{O}_{C,P}$, alors le résultat est évident d'après la définition du $2$-résidu en $P$ le long de $C$. En effet, dans cette situation, pour toute uniformisante $v$ de $\OC$, le couple $(u,v)$ est une $(P,C)$-paire faible et
$$
\res^2_{C,P}(\omega)=\res_P \left( (u)\res^1_C (\omega) \right)
.$$
Si maintenant $\bar{u}$ n'est pas une uniformisante de $\mathcal{O}_{C,P}$, alors quatre situations peuvent survenir. Dans ce qui suit, $t$ désigne la fonction $\frac{1}{u}$.
\begin{enumerate}
\item\label{oo} La fonction $\bar{u}$ est régulière en $P$ et $\bar{u}':=\bar{u}-\bar{u}(P)$ est une uniformisante de $\mathcal{O}_{C,P}$.
\item\label{on} La fonction $\bar{u}$ est régulière  en $P$ et $\bar{u}':=\bar{u}-\bar{u}(P)$ n'est pas une uniformisante de $\mathcal{O}_{C,P}$.
\item\label{no} La fonction $\bar{u}$ n'est pas régulière en $P$ et $\bar{t}':=\bar{t}-\bar{t}(P)$ est une uniformisante de $\mathcal{O}_{C,P}$.
\item\label{nn} La fonction $\bar{u}$ n'est pas régulière en $P$ et $\bar{t}':=\bar{t}-\bar{t}(P)$ n'est pas une uniformisante de $\mathcal{O}_{C,P}$.
\end{enumerate}

\noindent Remarquons que l'on peut donner une interprétation géométrique simple des deux premières situations si $u$ est régulière\footnote{Noter que le fait que $\bar{u}$ soit régulière en $P$ ne signifie en rien que $u$ l'est. Par exemple, sur $\mathbf{A}^2$, la fonction $u:= (x+y)/(x-y)$ n'est pas régulière à l'origine. Par contre, sa restriction à la courbe $C:=\{y=0\}$ est la fonction constante et égale à $1$ qui est régulière à l'origine.} en $P$. La première situation signifie que la ligne de niveau $u=u(P)$ intersecte $C$ transversalement en $P$.
Dans la seconde situation, cette ligne de niveau est singulière en $P$ ou tangente à $C$ en $P$.

Nous allons à présent traiter successivement ces quatre situations. On rappelle qu'il existe une unique fonction rationnelle $f$ sur $S$ telle que $\omega=fdu \w dv$ et que $f$ se décompose de façon unique en série de Laurent
$$
f=\sum_{j\geq -l}f_j v^j.
$$ 
On pose également $\mu:=(u)\res^1_C(\omega)$.
\medbreak

\noindent \textbf{Situation \ref{oo}.}
Le sous-corps $\mathcal{K}_u$ de $\hOC$ défini dans la proposition \ref{Ku} est l'unique sous-corps de $\hOC$ qui contienne $k(u)$ et soit isomorphe à $k(C)$ via le morphisme de réduction modulo $\mC$. La fonction $u_0:=u-\bar{u}(P)$ engendre le même sous-corps de $\OC$. En d'autres termes, $k(u)=k(u_0)$. De ce fait, les sous-corps $\mathcal{K}_u$ et $\mathcal{K}_{u_0}$ de $\hOC$ sont égaux.
De plus, $du=du_0$. Par conséquent,
$$\omega=fdu \w dv =fdu_0 \w dv$$
et le $(u_0)$-$1$-résidu de $\omega$ le long de $C$ est $\bar{f}_{-1}d\bar{u}'$ qui est égal à $\mu$. On déduit que 
$$
\res^2_{C,P}(\omega)=\res_P ((u_0)\res^1_C(\omega))=\res_P(\mu).
$$

\medbreak

\noindent \textbf{Situation \ref{on}.}
Soit $x$ une fonction rationnelle sur $S$ telle que $(x,v)$ soit une $(P,C)$-paire faible. La fonction $\bar{x}$ est donc une uniformisante de $\mathcal{O}_{C,P}$. De fait, il existe une série formelle $\phi \in k[[T]]$ telle que $\bar{u}'=\phi(\bar{x})$.
Soit $\sigma$ le relevé de Hensel de $\bar{x}$ dans $\hOC$, alors c'est un élément de $\mathcal{K}_u$ et dans ce corps, on a la relation
$$
u_0=\phi(\sigma).
$$
On en déduit la nouvelle expression de $\omega$
\begin{equation}\label{formexpr}
\omega=\sum_{j\geq -l} f_j v^j \phi'(\sigma)d\sigma\w dv,
\end{equation}

\noindent où $\phi'$ désigne la dérivée formelle de $\phi$.
Notons que $\sigma$ n'est à priori pas une fonction. Le second membre de l'expression (\ref{formexpr}) est à priori une $2$-forme formelle appartenant à $\Omega^2_{k((u))((v))/k}$ (voir section \ref{ratformel}).
À présent, rappelons que $\sigma$ est un élément de $\hOC$ qui est congru à $x$ modulo $\mC$. Par conséquent, $\sigma$ se décompose dans $k((x))[[v]]$ de la façon suivante
$$
\sigma=x+\sigma_1(x)v+\sigma_2(x)v^2+\cdots
$$
La paire $(\sigma, v)$ est liée à la paire $(x,v)$ par un changement de variables de la forme (\ref{cv}). D'après le théorème \ref{inv2res}, l'expression (\ref{formexpr}) fournit le même $2$-résidu que la décomposition de $\omega$ dans $k((x))((v))$. On en conclut que
$$
\res^2_{C,P}(\omega)=\res_P (\bar{f}_{-1}\phi'(\bar{\sigma})d\bar{\sigma})
=\res_P (\bar{f}_{-1}\phi'(\bar{x})d\bar{x}).
$$
Or, $\phi'(\bar{x})d\bar{x}$ est égal à $d\phi(\bar{x})=d\bar{u}'$, donc $\bar{f}_{-1}\phi'(\bar{x})d\bar{x}$ est égal à $\mu$.

\medbreak

\noindent \textbf{Situation \ref{no}.} Tout comme dans la situation \ref{oo}, on remarque que, comme $u$ est égal à $\frac{1}{t}$, on a
$$
k(u)=k(t), \quad \textrm{ce qui implique} \quad \mathcal{K}_u=\mathcal{K}_t.
$$

\noindent De fait, le développement de $\omega$ à coefficient dans $\mathcal{K}_t((v))$ s'écrit
$$
\omega = \sum_{j\geq -l}f_j v^j \left(-\frac{dt}{t^2} \right) \w dv.
$$

\noindent Par conséquent, on a
$$
(t)\res^1_C (\omega)=-\bar{f}_{-1} \frac{d\bar{t}}{\bar{t}^2} = \bar{f}_{-1}d\bar{u}=(u)\res^1_C (\omega).
$$

\noindent Comme $\bar{t}$ est par hypothèse une uniformisante de $\mathcal{O}_{C,P}$, le couple $(t,v)$ est une $(P,C)$-paire faible et donc
$$
\res^2_{C,P}(\omega)=\res_P ((t)\res^1_C(\omega))=\res_P ((u)\res^1_C(\omega)).
$$ 

\medbreak

\noindent \textbf{Situation \ref{nn}.} D'après la situation \ref{no}, le $(t)$-$1$-résidu de $\omega$ le long de $C$ est égal à son $(u)$-$1$-résidu. On reprend alors le travail effectué dans la situation \ref{on} en remplaçant $u$ par $t$ et on en déduit le résultat.
\end{proof}

En conclusion, les objets définis dans cette section sont bien une généralisation de ceux introduits en section \ref{secres}.

\section{Propriétés des résidus}\label{secprop}

Sur une courbe algébrique irréductible lisse $X$, une $1$-forme régulière en un point $P$ a un résidu nul en ce point. On dispose donc d'une condition nécessaire pour que le résidu d'une $1$-forme en un point soit non nul. Le lemme qui suit fournit un énoncé analogue pour les $2$-résidus. On rappelle que l'on se place toujours dans le cadre donné en section \ref{cadre}.

\begin{lem}\label{oukisont}
  Soit $\omega$ une $2$-forme rationnelle sur $S$ admettant $C$ comme pôle (éventuellement multiple) et $P$ un point rationnel lisse de $C$. Si $\omega$ n'a pas d'autre pôle que $C$ au voisinage de $P$, alors
$$
\res^2_{C,P}(\omega)=0.
$$
\end{lem}

\begin{proof}
Soient $(u,v)$ une $(P,C)$-paire forte et $n$ un entier tel que la valuation de $\omega$ le long de $C$ soit égale à $-n$. Par hypothèse, l'entier $n$ est positif. Par ailleurs, il existe une fonction rationnelle $h$, régulière au voisinage de $C$, telle que 
$$
\omega=hdu\w \frac{dv}{v^n}.
$$
De plus, comme $C$ est le seul pôle de $\omega$ au voisinage de $P$, on en déduit que $h$ est régulière au voisinage de $P$, elle se développe donc en série de Taylor
$$
h=\sum_{j\geq 0} h_j(u)v^j \quad \textrm{où}\ h_j\in k[[u]]\ \textrm{pour tout entier}\ j.
$$
Par conséquent,
$$
(u)\res^1_C(\omega)=h_{n-1}(\bar{u})d\bar{u}.
$$
Cette $1$-forme sur $C$ est régulière au voisinage de $P$, son résidu en $P$ est donc nul.
\end{proof}

\begin{rem}
C'est pour démontrer ce type d'énoncé que la notion de $(P,C)$-paire forte est très utile.
\end{rem}

En d'autres termes, une $2$-forme sur $S$ admet un $2$-résidu non nul en un point $P$ le long d'une courbe $C$, seulement si plusieurs $\omega$ a des pôles autres que $C$ au voisinage de $P$.

\subsection{Influence d'un éclatement sur les résidus}

Soient $P$ un point rationnel lisse de $C$ et $(u,v)$ une $(P,C)$-paire faible. On note $\pi:\ \widetilde{S}\rightarrow S$ l'éclatement de $S$ en $P$. On note $E$, le diviseur exceptionnel de $\widetilde{S}$.
La transformation stricte par $\pi$ d'une courbe $X$ passant par $P$ sera notée $\widetilde{X}$. On rappelle que la transformation stricte d'une courbe est l'adhérence de Zariski dans $\widetilde{S}$ de la courbe
$
\pi^{-1}(X)\smallsetminus E
$.

\begin{lem}\label{1resBU}
Soit $\omega$ une $2$-forme rationnelle sur $S$, on a
$$
(\pi^* u)\res^1_{\widetilde{C}}(\pi^*\omega)=\pi^* \left((u)\res^1_C (\omega)\right).
$$ 
\end{lem}

\begin{proof}
L'application $\pi$ induit un isomorphisme entre un ouvert de $C$ et un ouvert de sa transformée stricte.
Les $1$-formes $(\pi^* u)\res^1_{\widetilde{C}}(\pi^*\omega)$ et $(u)\res^1_C (\omega)$ sont tirées en arrière l'une de l'autre par cet isomorphisme.
\end{proof}

\begin{cor}\label{2resBU}
Soit $\omega$ une $2$-forme rationnelle sur $U$ et $Q$ le point d'intersection\footnote{La courbe $C$ est supposée lisse en $P$. Par conséquent, elle intersecte $E$ en un unique point.} de $E$ avec $\widetilde{C}$. On a 
$$
\res^2_{\widetilde{C},Q}(\pi^* \omega) = \res^2_{C,P}(\omega).
$$   
\end{cor}

\subsection{Le cas des points singuliers d'une courbe}\label{secsing2}

Les deux énoncés qui précèdent permettent de généraliser la définition $2$-résidu d'une $2$-forme en $P$ le long de $C$ au cas où $P$ est un point singulier de $C$. Dans ce qui suit, $P$ désigne un point rationnel éventuellement singulier de $C$.

\begin{prop}
  Soit $\pi:  \widetilde{S} \rightarrow S$ un morphisme birationnel provenant d'une suite finie d'éclatements de $S$ induisant une résolution de la singularité de $C$ en $P$. Soit $\widetilde{C}$ la transformée stricte de $C$ par $\pi$, alors, la somme
$$
\sum_{Q \rightarrow P} \res^2_{\widetilde{C},Q} (\pi^* \omega)
$$
ne dépend pas de $\pi$. La notation ``$Q \rightarrow P$'' signifie que $Q$ est un point de $\widetilde{C}$ envoyé sur $P$ par $\pi$.
\end{prop}

\begin{proof}
  Soient $\pi_1: \widetilde{S}_1 \rightarrow S$ et $\pi_2: \widetilde{S}_2 \rightarrow S$ deux tels morphismes et notons $\widetilde{C}_1$ et $\widetilde{C}_2$ les transformées strictes respectives de $C$ par ces applications.
Comme les deux applications induisent un résolution de la singularité de $C$ en $P$, le point $P$ a le même nombre d'antécédents par ${\pi_1}_{|\widetilde{C}_1}$ et ${\pi_2}_{|\widetilde{C}_2}$. Notons respectivement $P_{1,1}, \ldots, P_{1,n}$ et $P_{2,1}, \ldots, P_{2,n}$ ces antécédents. Il existe alors deux ouverts $U_1 \subseteq \widetilde{C}_1$ et $U_2 \subseteq \widetilde{C}_2$ contenant respectivement $P_{1,1}, \ldots, P_{1,n}$ et $P_{2,1}, \ldots, P_{2,n}$ et un isomorphisme $\varphi:U_1 \rightarrow U_2$ tel que ${\pi_1}_{|U_1}= {\pi_2}_{|U_2} \circ \varphi$.
De plus, quitte à ré ordonner les indices, $\varphi$ envoie $P_{1,i}$ sur $P_{2,i}$ pour tout $i$ appartenant à $\{1, \ldots, n\}$.

On se donne alors une fonction $u\in \OC$ dont la restriction à $C$ est un élément séparent de $k(C)/k$. D'après le corollaire \ref{1resBU}, les $1$-formes sont tirées en arrière l'une de l'autre par $\varphi$. Par conséquent, on a
$$
\forall i \in \{1, \ldots, n\}, \quad \res^2_{\widetilde{C}_1,P_{1,i}}(\pi_1^* \omega) = \res^2_{\widetilde{C}_2,P_{2,i}}(\pi_2^* \omega).
$$
 
\end{proof}

\begin{defn}\label{singdef}
Soit $P$ un point rationnel singulier de $C$ et $\pi: \widetilde{S} \rightarrow S$ un morphisme birationnel défini par une séquence finie d'éclatements induisant une résolution de la singularité de $C$ en $P$. Soit $\omega$ une $2$-forme rationnelle sur $S$, on définit le $2$-résidu de $\omega$ en $P$ le long de $C$ par
$$
\res^2_{C,P}(\omega):=\sum_{Q\rightarrow P} \res^2_{\widetilde{C},Q}(\pi^* \omega).
$$  
\end{defn}

\begin{rem}
Les points $Q$ au-dessus de $P$ peuvent être définis sur une extension finie de $k$. Cependant, on peut facilement se convaincre du fait que la somme qui définit $\res^2_{C,P}(\omega)$ est invariante sous l'action du groupe de Galois absolu de $k$. Le $2$-résidu reste donc un élément de $k$.  
\end{rem}

\begin{rem}\label{sparad}
  Muni de cette définition on montre aisément que l'énoncé de la proposition \ref{1resprop2} reste valable dans le cas où le point $P$ est un point singulier de la courbe $C$. Pour ce faire, il suffit de réaliser le même raisonnement que dans la preuve de cette proposition mais en l'appliquant à la courbe $\widetilde{C}$ de la définition ci-dessus.
\end{rem}

Une autre façon de définir et de calculer les $2$-résidus d'une $2$-forme le long d'une courbe $C$ en un point singulier $P$ de cette dernière est d'étendre à cette situation la définition de $(P,C)$-paire faible. Ce point de vue nous sera utile dans le chapitre \ref{chapreal}.

\begin{defn}\label{PCfsing}
Soient donc $C$ une courbe irréductible absolument réduite plongée dans $S$ et $P$ un point rationnel éventuellement singulier de $S$. On appelle $\pi:\widetilde{C}\rightarrow C$ la normalisation de $C$.
Une $(P,C)$-paire faible est un couple de fonctions $(u,v)$ appartenant à l'anneau local $\OC$ et vérifiant les conditions suivantes.
\begin{enumerate}
\item[$(i)$] Pour tout point fermé $Q$ de $\widetilde{C}$ au-dessus de $P$,
la fonction $\pi^*\bar{u} \in k(\widetilde{C})$ est une uniformisante de $\mathcal{O}_{\widetilde{C},Q}$.
\item[$(ii)$] La fonction $v$ est une uniformisante de $\OC$.
\end{enumerate}
\end{defn}
 
\begin{rem}
Notons que cette définition est une généralisation naturelle de la définition \ref{PCfaible} de $(P,C)$-paires faibles.
\end{rem}

\begin{rem}
Une $(P,C)$-paire faible existe toujours. En effet, d'après le théorème d'approximation faible (\cite{sti} théorème I.3.1), il existe une fonction $z$ dans $k(\widetilde{C})$ qui est une uniformisante de $\mathcal{O}_{\widetilde{C},Q}$ pour tout point fermé $Q$ au-dessus de $P$. On descend $z$ en une fonction $\bar{u}$ sur $C$ que l'on relève ensuite en une fonction $u$ dans $\OC$.  
\end{rem}

\begin{exmp}
Prenons $S=\mathbf{A}^2_{\C}$ et $C$ la cubique \textit{cuspidale} d'équation affine $y^2=x^3$ et $P$ l'origine. On pose alors
$$
u:=\frac{y}{x}\quad \textrm{et} \quad v:=y^2-x^3.
$$

\noindent On appelle $L_x$ et $L_y$ les droites d'équations respectives $x=0$ et $y=0$.
Éclatons $\mathbf{A^2}$ en $P$. On considère dons la surface plongée dans $\mathbf{A}^2 \times \P^1$ définie par l'équation $xu=yv$ où $(u:v)$ est un système de coordonnées homogènes de $\P^1$. On pose $z:=\frac{u}{v}$ et dans la carte affine ${v \neq 0}$ on obtient des équations pour $\widetilde{C}$,
$$
\widetilde{C}=\{y=xz\}\cap \{z^2=x\}.
$$

\begin{center}
\includegraphics[width=4.6cm, height=8cm]{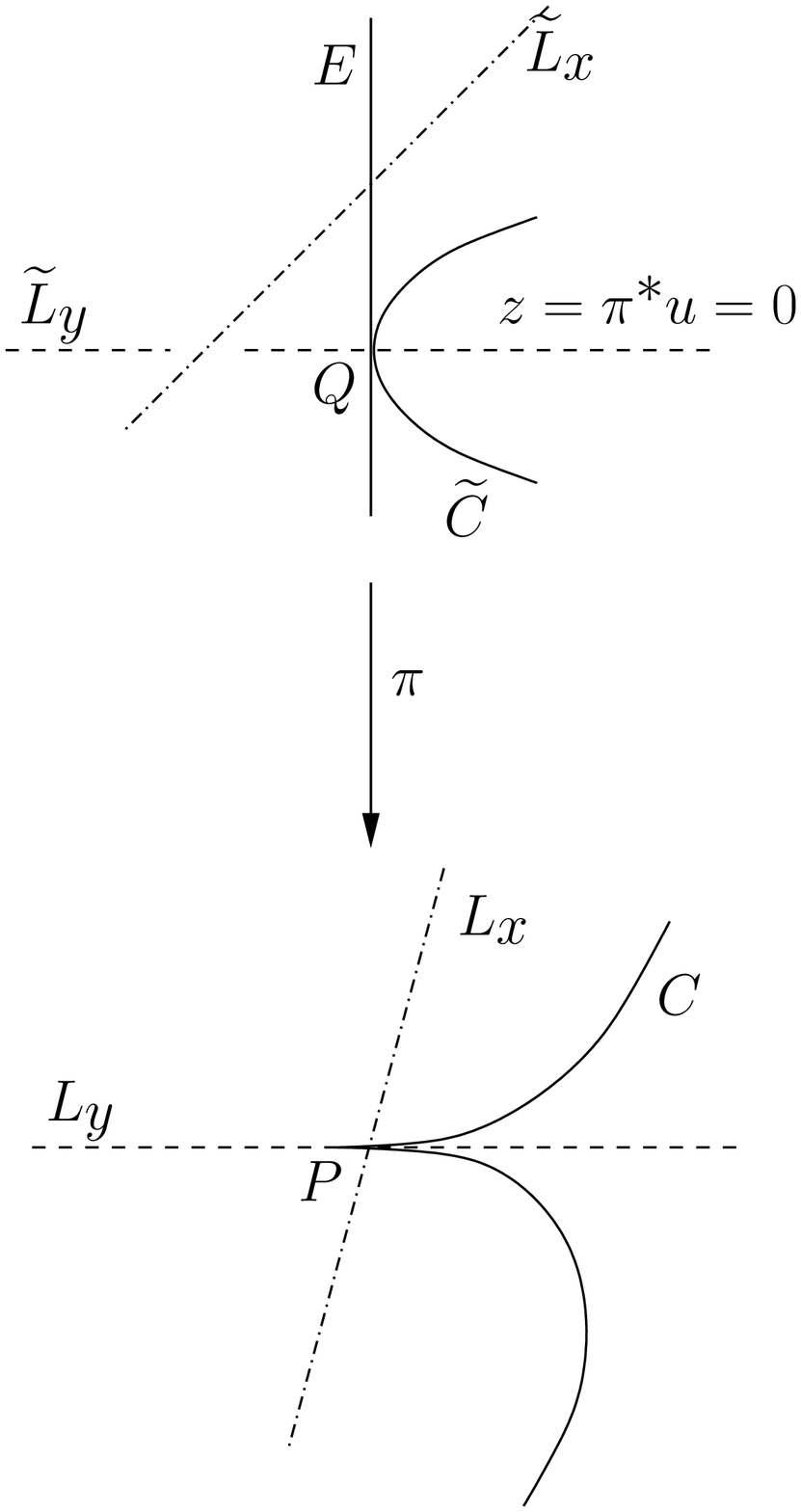}
\end{center}

On appelle $Q$ le point de $\widetilde{C}$ au-dessus de $P$. La fonction $\pi^*u$ est égale à $z$ et son lieu d'annulation au voisinage de $Q$ est $\widetilde{L}_y$ qui intersecte $\widetilde{C}$ en $Q$ avec multiplicité $1$. De fait, $\pi^* \bar{u}=\bar{z}$ est une uniformisante de $\mathcal{O}_{\widetilde{C},Q}$. En conclusion, la couple $(u,v)$ est bien une $(P,C)$-paire faible au sens de la définition \ref{PCfsing}.
\end{exmp}

Terminons cette section sur un commentaire élémentaire.
Soient $C$ une courbe irréductible absolument réduite plongée dans $S$ et $P$ un point rationnel singulier de $C$.
Soient alors $(u,v)$ une $(P,C)$-paire faible et $\omega$ une $2$-forme sur $S$.
Comme $\bar{u}$ est un élément séparant de $k(C)/k$, on peut définir un $(u)$-$1$-résidu de $\omega$ le long de $C$ qui s'identifie à une $1$-forme sur la normalisée $\widetilde{C}$ de $C$. La somme des résidus de cette $1$-forme sur $\widetilde{C}$ en tous les points au-dessus de $P$ est évidemment égal au $2$-résidu de $\omega$ en $P$ le long de $C$ de la définition \ref{singdef}.

\section{Formules de sommation}\label{sommation}

Les résultats énoncés dans cette section, en particulier le théorème \ref{FR3}, sont ceux que nous utiliserons dans les chapitres suivants qui portent sur les codes correcteurs. On rappelle que notre objectif est de construire des codes à partir de $2$-formes sur une surface et d'obtenir des relations d'orthogonalité entre ces codes et les codes fonctionnels.
Dans le cas des courbes, une partie de la démonstration de cette relation d'othogonalité consiste à utiliser la formule des résidus. Aussi, il semble intéressant de pouvoir disposer de formules de sommation de $2$-résidus d'une $2$-forme sur une surface. Nous allons fournir trois relations de sommations. La troisième est celle qui nous sera la plus utile dans ce qui suit.

\medbreak

\noindent \textbf{Attention!} Tout comme dans le cas des courbes, les formules de sommation qui suivent font intervenir tous des points géométriques de la surface. Aussi pour plus de confort, le corps de base $k$ sera supposé \textbf{algébriquement clos} dans cette section.

\begin{thm}[Première formule des résidus]\label{FR1}
Soit $S$ une surface quasi-projective lisse géométriquement intègre définie sur $k$.
Soient $C$ une courbe projective irréductible plongée dans $S$ et $\omega$ une $2$-forme rationnelle sur $S$. On a
$$
\sum_{P\in C} \res^2_{C,P}(\omega)=0.
$$   
\end{thm}

\begin{rem}\label{sensthm}
 D'après le lemme \ref{oukisont}, la somme ci-dessus a un support fini, l'énoncé a donc un sens.
\end{rem}

\begin{proof}
  Commençons par supposer que $C$ est lisse. Soit $u$ une fonction rationnelle sur $S$, régulière au voisinage de $C$ et dont la restriction $\bar{u}$ à $C$ est un élément séparant de $k(C)/k$. Soit $\mu$ le $(u)$-$1$-résidu de $\omega$ le long de $C$. D'après la proposition \ref{1resprop2} et la remarque \ref{sparad}, on a 
$$
\sum_{P\in C} \res^2_{C,P}(\omega)= \sum_{P\in C} \res_P(\mu)
$$
et cette dernière somme est nulle d'après la formule des résidus sur une courbe.
Si maintenant $C$ est singulière, on considère un morphisme birationnel $\pi: \widetilde{S} \rightarrow S$ obtenu par une séquence finie d'éclatements et tel que la transformée stricte $\widetilde{C}$ de $C$ soit lisse. D'après le lemme \ref{1resBU} et son corollaire \ref{2resBU}, on a
$$
\sum_{P\in C} \res^2_{C,P}(\omega)=
\sum_{Q\in \widetilde{C}} \res^2_{\widetilde{C},Q}(\omega)
=\sum_{Q\in \widetilde{C}} \res_Q (\pi^* \mu).
$$ 
On conclut de nouveau en appliquant la formule des résidus à la courbe $\widetilde{C}$ et la $1$-forme $\pi^* \mu$.
\end{proof}

\begin{rem}
  Noter que si la valuation de $\omega$ le long de $C$ est supérieure ou égale à $-1$, alors le résultat est évident. En effet, il suffit d'appliquer la formule des résidus au $1$-résidu de $\omega$ le long de $C$. La partie non évidente de la preuve ci-dessus est l'étude du cas où $\omega$ a un pôle multiple le long de $C$. Le travail sur les $(u)$-$1$-résidus effectué dans les sections \ref{secdefgen} et \ref{secprop} avait pour principal objectif de fournir les outils nécessaires à la preuve de ce résultat.
\end{rem}

\begin{thm}[Deuxième formule des résidus]\label{FR2}
Soit $S$ une surface quasi-projective lisse géométriquement intègre définie sur $k$.
Soient $P$ un point de $S$ et $\mathcal{C_{S,P}}$ l'ensemble des germes courbes irréductibles tracées sur $S$ et contenant $P$. Pour toute $2$-forme $\omega$ rationnelle sur $S$, on a
$$
\sum_{C\in \mathcal{C}_{S,P}} \res^2_{C,P}(\omega)=0.
$$
\end{thm}

\begin{rem}
La somme ci-dessus a également un support fini (c.f. remarque \ref{sensthm}).
\end{rem}

\begin{proof}
Soient $\omega$ une $2$-forme rationnelle sur $S$ et $C_1,\ldots , C_n$ les composantes irréductibles du lieu des pôles de $\omega$ au voisinage de $P$.

\medbreak

\noindent \textbf{Étape 1.} Dans un premier temps, nous allons supposer que les pôles $C_1,\ldots, C_n$ de $\omega$ sont lisses en $P$ et se croisent deux à deux transversalement en ce point.

\begin{itemize}

\item[\textbullet]
Si $n=1$, d'après le lemme \ref{oukisont}, le $2$-résidu de $\omega$ en $P$ le long de $C_1$ est nul. Le résultat est donc immédiat.

\item[\textbullet]
Si $n=2$, soient $u_1$ et $u_2$ des équations locales respectives des courbes $C_1$ et $C_2$ au voisinage de $P$.
Par hypothèse, $C_1$ et $C_2$ se croisent transversalement en $P$, donc $(u_1,u_2)$ est un système de coordonnées locales en ce point.
De fait, $(u_1,u_2)$ est une $(P,C_2)$-paire forte et $(u_2,u_1)$ une $(P,C_1)$-paire forte. Soient $-n_1$ et $-n_2$ les valuations respectives de $\omega$ le long de $C_1$ et $C_2$.
Il existe une fonction $h$ régulière au voisinage de $P$ telle que
$$
\omega = h \frac{du_1}{u_1^{n_1}}\w \frac{du_2}{u_2^{n_2}}.
$$
On développe $h$ en série de Taylor 
$$
h=\sum_{i,j\neq 0} h_{i,j}u_1^i u_2^j.
$$
Le $2$-résidu de $\omega$ en $P$ le long de $C$ est égal à $h_{n_1-1,n_2-1}$ et, comme le produit extérieur est anticommutatif, le $2$-résidu de $\omega$ le long de $C_1$ est égal à $-h_{n_1-1,n_2-1}$. Leur somme est donc nulle.

\item[\textbullet]
Si $n\geq 2$, soit $\pi: \widetilde{S} \rightarrow S$ l'éclatement de $S$ en $P$. Le diviseur exceptionnel est noté $E$, la transformée stricte d'une courbe $C_i$ est notée $\widetilde{C}_i$.
On rappelle que, par hypothèse les courbes $C_i$ sont lisses en $P$ et s'y croisent deux à deux transversalement. Donc, pour tout $i$, la courbe $\widetilde{C}_i$ intersecte $E$ en un unique point que l'on appelle $Q_i$ et les points $Q_1, \ldots, Q_n$ sont deux à deux distincts. La figure \ref{bam} résume cette situation.

\begin{figure}[!t]
\begin{center}
\includegraphics[width=7cm, height=11cm]{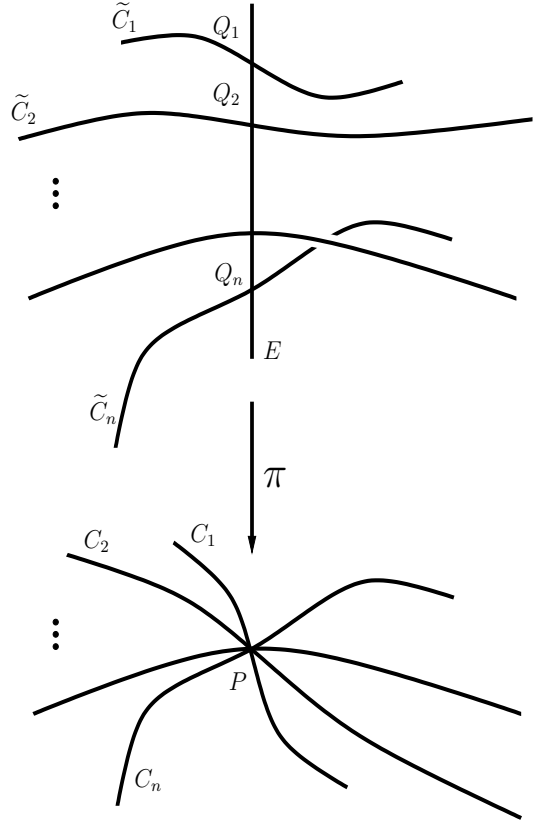}
\end{center}
\caption{Le cas $n\geq 2$ dans l'étape 1 de la preuve.}\label{bam}
\end{figure}

La courbe $E$ est projective et les $\widetilde{C}_i$ sont les seuls pôles de $\pi^* \omega$ qui intersectent $E$.
Soit $i\in \{1, \ldots, n\}$, d'après le cas $n=2$, on a 
\begin{equation}\label{excep}
\res^2_{\widetilde{C}_i, Q_i} (\pi^* \omega)
=-\res^2_{E,Q_i} (\pi^* \omega).
\end{equation}
Ainsi, en appliquant le lemme \ref{1resBU} et la relation (\ref{excep}) ci-dessus, on en déduit que
$$
\sum_{i=1}^n \res^2_{C_i,P}(\omega) = \sum_{i=1}^n \res^2_{\widetilde{C}_i,Q_i}
(\pi^* \omega) = -\sum_{i=1}^n \res^2_{E,Q_i}(\pi^* \omega).
$$
Soit $Q$ un point de $E$ autre que $Q_1, \ldots, Q_n$, la courbe $E$ est le seul pôle la $2$-forme $\pi^* \omega$ au voisinage de ce point. Par conséquent, d'après le lemme \ref{oukisont} le $2$-résidu de $\pi^* \omega$ en $Q$ le long de $E$ est nul.
En reprenant (\ref{excep}), on en déduit que
$$
\sum_{i=1}^n \res^2_{C_i,P}(\omega) = -\sum_{i=1}^n \res^2_{E,Q_i}(\pi^* \omega) = -\sum_{Q\in E} \res^2_{E,Q}(\pi^* \omega)
$$
et cette somme est nulle d'après le théorème \ref{FR1}.
\end{itemize}
 
\medbreak

\noindent \textbf{Étape 2.} Dans le cas général, les courbes $C_1,\ldots, C_n$ peuvent être singulières en $P$ et la multiplicité d'intersection de deux d'entre elles peut être supérieure ou égale à $2$. On réalise alors une désingularisation à \textbf{croisements normaux} de la courbe $C_1\cup \ldots \cup C_n$.
C'est-à-dire qu'à partir d'une séquence finie d'éclatements on obtient un morphisme birationnel $\pi: \widetilde{S} \rightarrow S$ tel que la surface $\widetilde{S}$ vérifie les propriétés suivantes.
\begin{enumerate}
\item[$(i)$] Les courbes $\widetilde{C}_1,\ldots, \widetilde{C}_n$ sont lisses en tout point $Q$ tel que $\pi(Q)=P$.
\item[$(ii)$] Par un point $Q$ tel que $\pi(Q)=P$ passe au plus une courbe $\widetilde{C}_i$.
\item[$(iii)$] L'intersection d'une courbe $\widetilde{C}_i$ avec la courbe $\pi^{-1}(\{P\})$ est de multiplicité un.
\end{enumerate}

\noindent On appelle arbre de résolution l'image réciproque par $\pi$ du point $P$. Il s'agit d'une réunion de courbes projectives de genre nul. Les relations d'incidence entre ces diviseurs se représentent sous la forme d'un arbre que l'on notera $\mathcal{A}$.
$$
\xymatrix{
E_{s,1} \ar@{.}[rd] & \cdots \ar@{.}[rd] & \cdots \ar@{.}[d] &
\cdots \ar@{.}[ld]& \cdots\ar@{.}[rd] & \cdots \ar@{.}[d] & E_{s',n_s'} \ar@{.}[ld] \\
 & E_{2,1} \ar@{-}[rrd] & E_{2,2} \ar@{-}[rd] & \cdots \ar@{-}[d] & \cdots \ar@{-}[ld] & E_{2,n_2} \ar@{-}[lld]  & \\
 & & & E_1 & & & 
}
$$


\noindent Noter que les feuilles de cet arbre ne sont pas forcément toutes au même étage, même si le diagramme ci-dessus laisse supposer le contraire. C'est la raison pour laquelle les indices des deux feuilles (extrémités supérieures) représentées sont différents ($s$ et $s'$).

À présent, nous allons appliquer les résultats de l'étape précédente aux sommets l'arbre $\mathcal{A}$, en partant de ses feuilles et en remontant à sa racine.
Dans ce qui suit, nous illustrerons notre travail de la façon suivante. Si $\widetilde{C}$ est un pôle de $\pi^* \omega$ et $Q$ un point de $\widetilde{C}$, alors le $2$-résidu $r$ de $\pi^* \omega$ en $Q$ le long de $\widetilde{C}$ apparaîtra dans un dessin sous la forme suivante.

\begin{center}
\includegraphics[width=5cm, height=1.5cm]{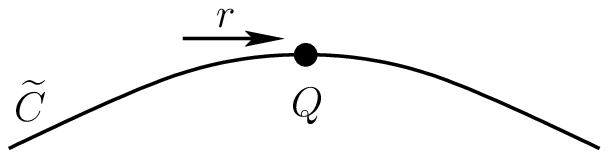}
\end{center}

Pour tout diviseur $E$ correspondant à un sommet autre que la racine de l'arbre
$\mathcal{A}$, on appelle $T$ le point d'intersection de $E$ avec son ascendant et $\mathcal{P}_E^j$ l'ensemble des points de $\widetilde{C}_j$ qui intersectent $E$
ou un de ses ascendants dans l'arbre $\mathcal{A}$. On note
$$
\sigma_E:= \res^2_{E,T_E} (\pi^* \omega).
$$
Commençons par montrer que pour tout diviseur $E$ correspondant à un sommet de l'arbre autre que sa racine, on a
\begin{equation}\label{inc}\tag{R}
\sigma_E= \sum_{j=1}^n \ \sum_{Q\in \mathcal{P}_E^j} \res^2_{\widetilde{C}_j, Q}(\pi^* \omega).
\end{equation}
Nous allons démontrer cette relation par récurrence sur les étages de l'arbre.

\medbreak

\noindent \textbf{Étape 2.a.} Soit $E$ un diviseur correspondant à une feuille de l'arbre. Il admet un unique ascendant dans l'arbre que l'on note $E'$ et qui intersecte $E$ en un point $T$. Par ailleurs, quitte à réordonner les indices des courbes, ce diviseur $E$ intersecte les courbes $\widetilde{C}_1, \ldots, \widetilde{C}_k$ en les points $T_{1,1}, \ldots, T_{1,l_1},\ldots , T_{k,1}, \ldots, T_{k,l_k}$. On note enfin
$$
r_{i,j}:=\res^2_{\widetilde{C}_{i}, T_{i,j}}(\pi^* \omega).
$$

\noindent La situation peut être représentée par la figure suivante.

\begin{center}
\includegraphics{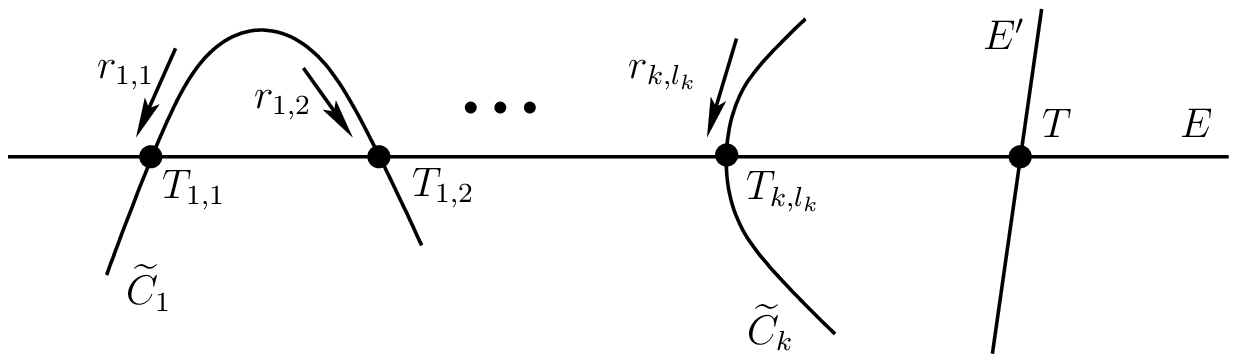}
\end{center}

\noindent D'après le travail effectué dans l'étape $1$, on sait que pour tout couple $(i,j)$
$$
\res^2_{\widetilde{C}_i,T_{i,j}} (\pi^* \omega) = r_{i,j} = -\res^2_{E,T_{i,j}} (\pi^* \omega).
$$

\noindent De plus, d'après le lemme \ref{oukisont}, la $2$-forme $\pi^* \omega$ a des résidus non nuls seulement en les points $T_{i,j}$ et $T$. Donc,
d'après la première formule des résidus (théorème \ref{FR1}), on obtient
$$
\sigma_E=\res^2_{E,T}(\pi^* \omega)=-\sum_{i,j} \res^2_{E,T_{i,j}} (\pi^* \omega)=\sum_{i,j} r_{i,j}.
$$
C'est-à-dire la relation (\ref{inc}) pour une feuille de l'arbre $\mathcal{A}$.
Le schéma suivant résume le travail qui vient d'être effectué. 

\begin{center} 
\includegraphics{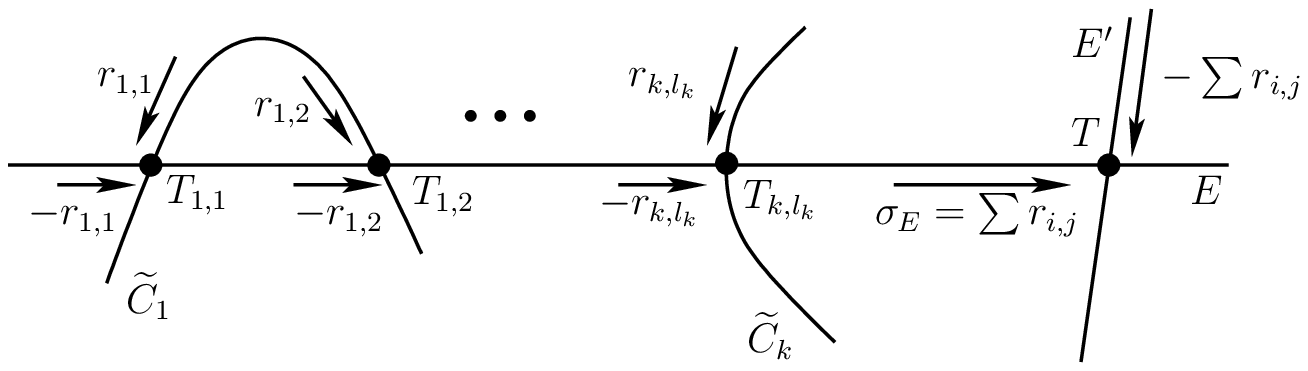}
\end{center}

\noindent Notons que, comme cela apparaît sur le dessin ci-dessus, en appliquant de nouveau le travail effectué dans l'étape $1$, on obtient
$$
\res^2_{E',T}(\pi^* \omega) = -\res^2_{E,T} (\pi^* \omega) =-\sum_{i,j} \res^2_{\widetilde{C}_i,T_{i,j}} (\pi^* \omega).
$$

\noindent \textbf{Étape 2.b.} 
Soit $E$ un diviseur correspondant à un sommet intermédiaire de l'arbre, c'est-à-dire un sommet qui n'est ni une feuille ni la racine. Par récurrence, supposons que la relation (\ref{inc}) est vérifiée par tous les descendants (directs ou indirects) de $E$ dans $\mathcal{A}$.
Notons $E'$ l'ascendant direct de $E$ dans $\mathcal{A}$ et $D_1,\ldots, D_r$ ses descendants directs. Soient également $\widetilde{C}_1, \ldots ,\widetilde{C}_q$ les courbes\footnote{Quitte à ré indicer les courbes.} qui intersectent $E$.
On désigne par $T$, le point d'intersection de $E$ avec $E'$. Les points d'intersection de $E$ avec $D_1,\ldots, D_r$ sont notés $U_1,\ldots, U_r$ et les points d'intersection de $E$ avec $\widetilde{C}_1, \ldots ,\widetilde{C}_q$ sont notés $T_{1,1}, \ldots ,T_{1,l_1}, \ldots ,T_{q,1},\ldots ,T_{Q,l_q}$. On reprend la notation
$$
r_{i,j}:=\res^2_{\widetilde{C}_j,T_{i,j}}(\pi^* \omega).
$$
D'après l'hypothèse de récurrence on a 
$$
\forall i\in \{1, \ldots, r\},\quad \res^2_{D_i,U_i} (\pi^* \omega)=\sigma_{D_i}.
$$
Donc d'après le cas $n=2$ de l'étape 1, on a 
$$
\forall i\in \{1, \ldots, r\},\quad \res^2_{E,U_i} (\pi^* \omega)=-\sigma_{D_i}
$$
et
$$
\forall j\in \{1, \ldots, q\}, \forall i\in \{1, \ldots, l_q\},\quad
\res^2_{E,T_{i,j}}(\pi^*\omega)=-r_{i,j}.
$$

\noindent Ainsi, d'après le théorème \ref{FR1} appliqué à $E$ et $\pi^* \omega$, on a 
$$
\sigma_E=\res^2_{E,T}(\pi^* \omega)=-\sum_{k=1}^r \res^2_{E,U_k}(\pi^* \omega)-\sum_{i,j}
\res^2_{E,T_{i,j}}(\pi^* \omega)=\sum_{k=1}^r \sigma_{D_k} +\sum_{i,j}r_{i,j}$$
et cette dernière somme n'est autre que $\sigma_E$.
La relation (\ref{inc}) est donc vérifiée par $E$. Le dessin suivant résume le travail effectué.
\begin{center}
\includegraphics{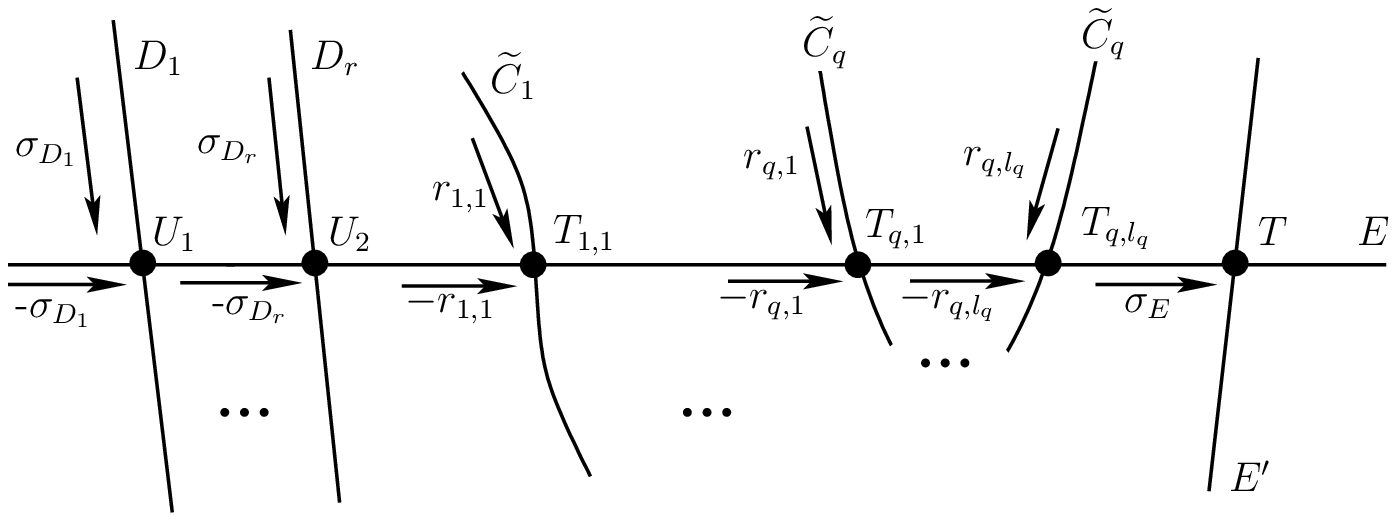}
\end{center}

\medbreak

\noindent \textbf{Étape 2.c.} Considérons maintenant $E_1$ la racine de l'arbre $\mathcal{A}$. Ce dernier n'a pas d'ascendant. Reprenons les notations de l'étape précédente en ce qui concerne ses descendants directs et les courbes $\widetilde{C}_j$ qu'il intersecte.
Par un raisonnement analogue à celui qui a été effectué dans l'étape précédente, en appliquant le théorème \ref{FR1} à $E_1$ et $\pi^* \omega$ on obtient
\begin{equation}\label{sommecool}
-\sum_{k=1}^r \sigma_{D_k} -\sum_{i,j}r_{i,j}=0.
\end{equation}

\noindent Or cette somme n'est autre que la somme
$$
\sum_{i=1}^n \sum_{Q\in \mathcal{P}_{E_1}^i} \res^2_{\widetilde{C}_i, Q}(\pi^* \omega)=0,
$$
où l'on rappelle que  $\mathcal{P}_{E_1}^i$ est l'ensemble des points de $\widetilde{C}_i$ qui intersecte $E_1$ ou l'un de ses antécédents. C'est donc l'ensemble des points des courbes $\widetilde{C}_i$ qui sont envoyés sur $P$ par $\pi$. On conclut en rappelant que, d'après la définition \ref{singdef}, on a
$$
\sum_{Q\in \mathcal{P}_{E_1}^i} \res^2_{\widetilde{C}_i, Q}(\pi^* \omega)=\res^2_{C_i,P}(\omega).
$$
En combinant cette relation avec l'équation (\ref{sommecool}) on obtient le résultat attendu, à savoir
$$
\sum_{i=1}^n \res^2_{C_i,P}(\omega)=0.
$$
\end{proof}

\begin{rem}
Noter que, dans la démonstration du théorème \ref{FR2} qui précède, on n'a appliqué le théorème \ref{FR1} qu'à des diviseurs appartenant à l'image réciproque de $P$. Or, il a été signalé que la preuve du théorème \ref{FR1} est évidente si
la valuation de $\omega$ le long de $C$ est supérieure ou égale à $-1$. Aussi, il est important de signaler que, dans la preuve qui précède, la valuation de $\pi^* \omega$ le long d'un diviseur exceptionnel $E$ provenant d'un éclatement peut être inférieure à $-2$. Le lemme suivant fournit une formule explicite de la valuation de $\pi^* \omega$ le long d'un tel diviseur exceptionnel.
\end{rem}

\begin{lem}\label{valex}
Soit $\omega$ une $2$-forme rationnelle sur $S$ et $P$ un point de $S$.
On note $\mathcal{C}_P$, l'ensemble de toutes les courbes irréductibles contenues dans  $S$ et contenant $P$.
Soit $\pi: \widetilde{S} \rightarrow S$ l'éclatement de $S$ en $P$, alors la valuation de $\pi^* \omega$ le long du diviseur exceptionnel $E$ est donnée par la formule:
$$
\textrm{val}_E(\pi^* \omega)=1+\sum_{C \in \mathcal{C}_P} m_P(C) \textrm{val}_{C} (\omega),
$$ 
où $m_P(C)$ est la multiplicité de $C$ en $P$.
\end{lem}

\begin{proof}
D'après \cite{H} proposition V.3.3, on a l'égalité de diviseurs
$$
(\pi^* \omega)=\pi^* (\omega)+E.
$$
Ensuite, d'après \cite{H} proposition V.3.6, la valuation du diviseur $\pi^* (\omega)$ le long de $E$ est
$$
\textrm{val}_E (\pi^* (\omega))=\sum_{C \in \mathcal{C}_P} m_P(C) \textrm{val}_{C} (\omega).
$$
\end{proof}

Remarquons également que la théorème \ref{FR2} permet d'étendre le lemme \ref{oukisont} au cas des courbes singulières voire même réductibles.
C'est ce qui fera l'objet du corollaire \ref{oukisont2}.
Pour énoncer ce dernier, nous avons besoin de la convention et la définition ci-dessous.

\medbreak

\noindent \textbf{Convention.}
Soit $C$ une courbe quelconque plongée dans $S$ et $P$ un point de $S$ n'appartenant pas à $C$, on dit alors que 
$$
\res^2_{C,P}(\omega)=0.
$$

\begin{defn}[$2$-résidu en un point d'une courbe réductible]
Soit $C$ une courbe réductible plongée dans $S$ et $C_1, \ldots, C_d$ ses composantes irréductibles. Soit $P$ un point de $C$ et $\omega$ une $2$-forme rationnelle sur $S$. On définit le $2$-résidu de $\omega$ en $P$ le long de $C$ par
$$
\res^2_{C,P}(\omega):= \sum_{i=1}^d \res^2_{C_i,P}(\omega).
$$
\end{defn}

\begin{cor}\label{oukisont2}
Soit $C$ une courbe quelconque plongée dans $S$ et $P$ un point de $C$. Soit une $2$-forme rationnelle $\omega$ dont le lieu des pôles au voisinage de $P$ est contenu dans $C$, alors
$$
\res^2_{C,P}(\omega)=0.
$$ 
\end{cor}

\begin{proof}
On applique la deuxième formule des résidus (théorème \ref{FR2}) à $\omega$, ses pôles au voisinage de $P$ faisant partie des composantes irréductibles de $C$.
\end{proof}

Pour finir ce chapitre, nous allons énoncer la troisième formule des résidus, qui est celle que nous appliquerons aux codes correcteurs dans le chapitre suivant. 
Noter que, dans le chapitre suivant nous manipulerons fréquemment des diviseurs. C'est ce qui motive la définition suivante.
\medbreak

\begin{defn}[$2$-Résidu en un point le long d'un diviseur]\label{resdiv}
Soit $D$ un diviseur sur $S$. Pour toute $2$-forme rationnelle $\omega$ et tout point $P$ de $S$, on appelle $2$-résidu de $\omega$ en $P$ le long de $D$ et on notera $\res^2_{D,P}(\omega)$ le $2$-résidu de $\omega$ en $P$ le long du \textbf{support} de $D$.  
\end{defn}

\noindent \textbf{Attention!} Noter qu'il  ne s'agit pas exactement d'une extension de la définition par linéarité. D'une certaine façon, la définition \ref{resdiv} ci-dessus autorise un abus de langage pour éviter d'avoir à parler de \textit{$2$-résidu en un point le long du support du diviseur $D$}.
En particulier, il faut faire attention au fait que, selon cette définition, le résidu d'une $2$-forme $\omega$ en un point $P$ le long d'un diviseur $D$ est par exemple égal à celui de $\omega$ en $P$ le long du diviseur $2D$.

\begin{thm}[Troisième formule des résidus, \cite{lip} chap. 12]\label{FR3}
Soit $S$ une surface projective lisse géométriquement intègre. Soient $D_a$ et $D_b$ deux diviseurs sur $S$ dont l'intersection des supports est un ensemble fini $Z$. Soit $\Omega^2(-D_a-D_b)$ le faisceau de $2$-formes vérifiant localement
$$
(\omega) \geq -D_a-D_b.
$$  
Alors, pour toute section globale $\omega$ du faisceau $\Omega^2 (-D_a-D_b)$, on a 
$$
\sum_{P\in S} \res^2_{D_a,P}(\omega)=\sum_{P\in Z} \res^2_{D_a,P}(\omega)=0.
$$
\end{thm}

\begin{proof}
La $2$-forme $\omega$ n'a pas de pôles hors du support de $D_a+D_b$.
Donc, d'après le corollaire \ref{oukisont2}, les $2$-résidus de $\omega$ le long de $D_a$ sont nuls en tout point $P$ n'appartenant pas à $Z$, ce qui nous donne la première égalité.
La seconde égalité vient de la première formule des résidus (théorème \ref{FR1}). En effet, soient $D_{a,1},\ldots, D_{a,m_a}$ les composantes irréductibles du support de $D_a$. Le théorème \ref{FR1} entraîne
$$
\forall i\in \{1, \ldots, m_a\},\quad \sum_{P\in D_{a,i}}\res^2_{D_{a,i},P}(\omega)=0. 
$$
De fait, en sommant ces $m_a$ relations, on obtient
$$
\sum_{P\in \supp D_a} \res^2_{D_a,P} (\omega)=0.
$$
Et il revient au même de sommer sur tous les points de $S$ puisque les $2$-résidus de $\omega$ le long de $D_a$ en un point hors du support de $D_a$ sont nuls par convention.
\end{proof}

\begin{rem}
Sous les conditions du théorème \ref{FR3}, soit $\omega$ une section globale du faisceau $\Omega^2 (-D_a-D_b)$.  
D'après la deuxième formule des résidus (thm \ref{FR2}) on a 
$$
\forall P \in S, \quad \res^2_{D_a,P}(\omega)=-\res^2_{D_b,P}(\omega).
$$  
Par conséquent l'énoncé du théorème \ref{FR3} est symétrique, c'est-à-dire qu'il reste vrai si l'on échange $D_a$ et $D_b$.
\end{rem}


\newpage
\thispagestyle{empty}
\null

\part{Codes géométriques}

\chapter{Codes diff\'erentiels sur une surface}\label{chapdiff}

\begin{flushright}
\begin{tabular}{p{8cm}}
\begin{flushright}
{\small \textit{``Je vous jure d'être décent et de ne pas dire un seul gros mot ni rien qui blesse les convenances.''}}
\medbreak
{\small Musset}\\
{\small \textit{Il ne faut jurer de rien}} 
\end{flushright}
\end{tabular}
\end{flushright}

Dans ce chapitre, nous allons appliquer les résultats du chapitre \ref{chapres}.
Le but est de construire des codes correcteurs d'erreurs en évaluant les $2$-résidus de $2$-formes différentielles en des points rationnels d'une surface algébrique.

\section{Langage et Notations}

Soit $X$ une variété géométriquement intègre de dimension $n$ définie sur un corps $k$. On note $\mathcal{O}_X$ son faisceau structural. L'ensemble des diviseurs de Weil sur $X$ sera noté $\divk (X)$. Étant donné un diviseur $D$ sur $S$ on note respectivement $D^+$ et $D^-$ ses parties effectives et non effectives. Les diviseurs $D^+$ et $D^-$ sont tous deux effectifs et $D$ s'écrit
$$
D=D^+-D^- .
$$
 L'équivalence linéaire sera notée ``$\sim$''. Si $X$ est lisse, le groupe $\divk (X)/\sim$ s'identifie au groupe de Picard de $X$ que l'on notera $\pick (X)$.
Si $X$ est une surface lisse et que 
les supports de deux diviseurs $D$ et $D'$ n'ont pas de composante irréductible commune, leur multiplicité d'intersection en un point $P$ est notée $m_P(D,D')$.
Étant donné un diviseur $G$ sur $X$, on note $\L (G)$ (resp. $\Omega^n (G)$) le 
faisceau inversible des fonctions rationnelles (resp. des $n$-formes rationnelles) sur $X$ qui vérifient localement
$$
(f)\geq -G \quad (\textrm{resp.}\ (\omega)\geq G).
$$
L'ensemble des sections globales d'un faisceau $\mathcal{F}$ sera noté $\Gamma(X, \mathcal{F})$ et $\mathcal{F}_P$ désignera sa fibre en un point $P$.
En ce qui concerne les faisceaux $\L (G)$ (qui seront fréquemment utilisés), on utilisera la notation standard $L(G)$ pour $\Gamma (X, \L (G))$.
Noter que, dans la littérature, le symbole $\Omega^n(G)$ peut désigner un faisceau inversible ou l'espace des sections globales du faisceau en question.
Insistons donc sur le fait que, dans ce qui suit $\Omega^n (G)$ désignera toujours un faisceau de $n$-formes.

Pour finir, soit $\bar{k}$ la clôture algébrique de $k$, on note $\overline{X}$ la variété
$$
\overline{X}:=X\times_{k} \bar{k}
$$
et, étant donné un faisceau $\mathcal{F}$ sur $X$, on note $\overline{\mathcal{F}}$ le tiré en arrière de $\mathcal{F}$ sur $\overline{X}$.

\paragraph{Important.} Dans tout ce qui suit, sauf mention contraire, si $D_1$ et $D_2$ sont deux diviseurs sur une surface lisse $S$ dont les supports n'ont pas de composante irréductible commune, alors ``$D_a \cap D_b$'' signifiera \textbf{intersection au sens de la théorie des schémas}.
Il s'agit donc d'une intersection tenant compte des multiplicités et non d'une intersection ensembliste.

\section{Rappels sur les codes construits à partir de courbes}\label{codescourbes}

Dans cette section, $X$ désigne une courbe algébrique projective lisse au-dessus de $\F_q$. On se donne également un diviseur $\F_q$-rationnel $G$ sur $X$ et une famille de points $P_1, \ldots, P_n$ rationnels sur $X$ et qui évitent le support de $G$. On note $D$ le diviseur
$$
D:=P_1+\cdots +P_n.
$$

\subsection{Codes fonctionnels et différentiels}\label{courbes}

La donnée des diviseurs $G$ et $D$ permet de construire deux codes différents. Ces codes sont respectivement appelés \textit{codes fonctionnels} et \textit{codes différentiels}. Les premiers sont construits par évaluation de fonctions en des points rationnels de $X$ et les seconds par évaluation de résidus de formes différentielles en ces mêmes points.

\paragraph{Le code fonctionnel.}
Soit $\ev_D$ l'application,
$$
\textrm{ev}_{D}:\left\{
\begin{array}{ccc}
L(G) & \rightarrow & \F_q^n \\
f & \mapsto & (f(P_1),\ldots ,f(P_n)).
\end{array}
\right. 
$$
L'image de cette application est appelée code fonctionnel associé aux diviseurs $D$ et $G$ et notée $C_{L,X}(D,G)$.

\medbreak

\paragraph{Le code différentiel.}
Soit $\textrm{res}_D$ l'application,
$$
\textrm{res}_{D}:\left\{
\begin{array}{ccc}
\Gamma(X,\Omega^1(G-D)) & \rightarrow & \F_q^n \\
\omega & \mapsto &
(\textrm{res}_{P_1}(\omega),\ldots ,\textrm{res}_{P_n}(\omega)).
\end{array}
\right.
$$
L'image de cette application est appelée code différentiel associé aux diviseurs $D$ et $G$ et notée $C_{\Omega,X}(D,G)$.

\begin{exmp}
Supposons que $X$ soit la droite projective $\P^1_{\F_q}$ et que le diviseur $G$ soit de la forme $kP$ où $P$ est un point rationnel de $\P^1_{\F_q}$.
Alors, en prenant pour $D$ une somme de points rationnels autres que $P$, le code fonctionnel est un code de Reed-Solomon et le code différentiel un code de Goppa classique (cf \cite{hohpel} exemples 3.3 et 3.4).
\end{exmp}

\subsection{Paramètres de ces codes}

L'un des intérêts de ces codes est que l'on dispose d'outils simples provenant de la théorie des courbes algébriques pour en évaluer les paramètres.
\begin{itemize}
\item[\textbullet] Le théorème de Riemann-Roch permet de minorer, voire d'évaluer exactement la dimension de ces codes.
\item[\textbullet] Le fait que le degré d'un diviseur principal soit nul permet d'obtenir de façon élémentaire une borne inférieure pour la distance minimale d'un code fonctionnel. Un raisonnement analogue sur le degré d'un diviseur canonique fournit une méthode de minoration de la distance minimale d'un code différentiel.
Nous renvoyons le lecteur aux références \cite{step}, \cite{sti} et\cite{TV} pour plus de détails sur ces propriétés (la liste n'est bien sûr pas exhaustive).

\end{itemize}

\noindent Concernant la distance minimale $d_{\textrm{min}}$, 
on obtient pour les codes fonctionnels la minoration
$$
n-k+1-g \leq d_{\textrm{min}},
$$
où $k$ désigne la dimension du code et $g$ le genre de la genre de la courbe. Ainsi on sait que ces codes sont toujours \emph{à $g$ de la borne de singleton}. C'est l'une des raisons qui a motivé les nombreux travaux effectués autour de l'étude des codes géométriques durant les 25 dernières années.

\subsection{Relation d'orthogonalité et décodage}

Un autre intérêt de ces codes est que l'on dispose d'une relation d'orthogonalité entre le code fonctionnel et le code différentiel, à savoir
$$
C_{\Omega,X}(D,G)=C_{L,X}(D,G)^{\bot}.
$$
La preuve de cette relation est une conséquence de la formule des résidus pour l'inclusion $\subseteq$ et du théorème de Riemann-Roch qui fournit une égalité de dimension entre ces codes entraînant l'inclusion réciproque. Nous renvoyons le lecteur à \cite{TV} théorème 3.1.44 ou \cite{sti} théorème II.2.8 pour plus de détails sur ce résultat.

Noter que, cette relation d'orthogonalité est l'outil de base de la majorité des algorithmes de décodage (voir \cite{hohpel} ou \cite{handhp} ).

\subsection{Deux constructions distinctes mais une seule classe de codes}

Un dernier résultat bien connu concernant ces codes est que tout code différentiel est un code fonctionnel associé à d'autres diviseurs et réciproquement. Plus précisément, étant donnés deux diviseurs $D$ et $G$ comme précédemment, il existe un diviseur canonique $K$ tel que
$$
C_{\Omega,X}(D,G)=C_{L,X}(D, K-G+D).
$$
Le diviseur $K$ est celui d'une forme différentielle dont les résidus en les points du support de $D$ sont tous égaux à $1$.
L'existence d'un tel diviseur est une conséquence du théorème d'approximation faible dans les corps de fonctions (\cite{sti} théorème I.3.1).

Aussi, si l'on souhaite étudier ces codes, on peut sans perte de généralité se restreindre à l'étude de codes issus d'une seule des deux constructions.
Généralement, on se focalise sur les codes fonctionnels dont la construction semble plus accessible, les fonctions étant un objet plus intuitif que les formes différentielles.

\begin{rem}
Il est toutefois intéressant de noter que la première construction de codes géométriques fut donnée par V.D. Goppa en 1981 dans l'article \cite{goppa}. Dans cet article, les codes introduits sont des codes différentiels. Sans doute parce qu'il sont une généralisation des codes de Goppa classiques. 
\end{rem}

\section{Codes géométriques construits à partir de surfaces algébriques}

Si de telles constructions sont possibles sur les courbes, il est naturel de s'interroger sur les perspectives d'extension de ces constructions à des variétés de dimension supérieure. 

\subsection{Cadre}\label{diffcadre}

Dans ce qui suit et jusqu'à la fin de ce chapitre, $S$ désignera une surface algébrique projective lisse définie au-dessus d'un corps fini $\F_q$.
Sauf mention contraire, lorsque l'on parlera de courbe pongée dans $S$, il s'agira de courbe définie sur $\F_q$.
On se donne également un diviseur $\F_q$-rationnel $G$ sur $S$ et une famille de points rationnels $P_1,\ldots ,P_n$ de $S$ qui évitent le support de $G$.
On appelle $\Delta$ le $0$-cycle sur $S$ défini par
$$
\Delta:= P_1+ \cdots + P_n.
$$
Notons que $\Delta$ joue plus ou moins le rôle du diviseur $D$ de la section précédente. Nous avons cependant choisi de le noter avec une lettre grecque car ce n'est plus un diviseur mais un $0$-cycle.
D'une façon générale, dans tout ce qui suit, les lettres latines majuscules désigneront des diviseurs et les lettres grecques majuscules des $0$-cycles.
Enfin, signalons dès à présent que la différence de dimension entre $\Delta$ et $G$ sera à l'origine de la plupart des difficultés que posent la construction et l'étude des codes différentiels sur des surfaces.

\subsection{Codes fonctionnels}\label{codesfonc}
Comme nous l'avons dit précédemment, la construction fonctionnelle présentée dans le cas des courbes se généralise à des variétés de dimension quelconque (voir \cite{manin} I.3.1). Pour ce faire, on définit l'application 
$$
\ev_{\Delta}:\left\{
\begin{array}{ccc}
L(G) & \rightarrow & \F_q^n \\
f & \mapsto & (f(P_1), \ldots, f(P_n)).
\end{array}
\right.
$$
L'image de cette application est appelée code fonctionnel sur $S$ associé à $\Delta$ et $G$ et est notée $C_{L,S}(\Delta,G)$.
L'étude des paramètres de ce type de code est nettement plus ardue que dans le cas des courbes.

\paragraph{Sur la longueur du code.}
Si l'on veut étudier l'asymptotique des codes construits sur des surfaces algébriques, on doit disposer de moyens d'évaluer le nombre de points rationnels d'une surface sur un corps fini. 
La borne de Weil-Deligne (voir \cite{deligne}) valable pour des variétés de dimension quelconque permet de majorer ce nombre de points. Lachaud et Tsfasman ont donné des estimations plus précises de ce nombre de points via des formules explicites dans \cite{LT}.

\paragraph{Sur la dimension.}
  En dimension supérieure ou égale à $2$, le théorème de Riemann-Roch se complique. Aussi, l'évaluation de la dimension de ce type de code est en général plus ardue. Cependant, dans les exemples étudiés dans la littérature, le diviseur $G$ est presque toujours un diviseur très ample obtenu par intersection de $S$ avec une hypersurface.
Dans ce cas, l'évaluation de la dimension du code se ramène au calcul élémentaire de la dimension d'espaces de polynômes en plusieurs variables et de degré total borné.

\paragraph{Sur la distance minimale.}
Alors que l'on disposait facilement de la distance minimale construite de Goppa (\textit{designed minimal distance}) dans le cas des courbes (voir \cite{sti} def II.2.4), dans le cas des variétés de dimension supérieure, 
la  minoration de la distance minimale d'un code fonctionnel devient un problème infiniment plus complexe.
Elle revient à majorer le nombre maximal de points rationnels du lieu d'annulation d'un élément de $L(G)$.
Les références citées dans l'introduction (page \pageref{refbib}) portent principalement sur la résolution de ce problème lorsque $S$ appartient à une classe spécifique de surfaces.

\begin{rem}
Notons que Lachaud propose dans \cite{lachaud2} une construction sensiblement différente du code fonctionnel. Cette approche a été reprise par un certain nombre des auteurs précédemment cités tels que Aubry, Edoukou et S{\o}rensen.
L'annexe \ref{annexelachaud}
est consacrée à cette autre construction et au moyen de la relier à celle présentée ci-dessus.  
\end{rem}

Passons maintenant à une première approche de la construction de codes différentiels sur une surface.

\subsection{Codes différentiels}\label{codesdiff}
Pour réaliser une construction analogue à celle qui a été présentée en section \ref{courbes}, il faut plus que la donnée de $G$ et $\Delta$.
En effet, on souhaiterait évaluer les $2$-résidus de $2$-formes ayant des pôles prescrits.
Il faut donc introduire un nouveau diviseur que l'on notera $D$ et qui, d'une certaine manière, jouera le rôle\footnote{On a signalé en section \ref{diffcadre} que $\Delta$ jouait le rôle du diviseur $D$ dans le cas des courbes. En réalité, dans le cas des codes différentiels sur des surfaces deux objets de dimension différente endossent le rôle joué par le diviseur $D$ dans le cas des courbes. Il y a d'un côté le $0$-cycle $\Delta$ et d'un autre la paire de diviseurs $(D_a,D_b)$.} du diviseur du même nom dans la construction de codes différentiels sur une courbe.
Il faut également que l'on évalue les $2$-résidus le long de certains pôles de $\omega$ mais 
pas tous. En effet, d'après le théorème \ref{FR2} si l'on évalue le $2$-résidu de $\omega$ en un point le long de tous ses pôles, on obtient zéro. Il faut donc décomposer le diviseur $D$ en deux parties distinctes. C'est ce qui motive la définition suivante.

\begin{defn}\label{defdiff}
Soient $D_a$ et $D_b$ deux diviseurs sur $S$ dont les supports n'ont pas de composante irréductible commune et soit $D$ la somme de ces deux diviseurs. On définit l'application
$$
\res^2_{D_a,\Delta}:
\left\{
\begin{array}{ccc}
\Gamma (S, \Omega^2 (G-D)) & \rightarrow & \F_q^n \\
\omega & \mapsto & (\res^2_{D_a,P_1}(\omega), \ldots ,
\res^2_{D_a,P_n}(\omega)).
\end{array}
\right.
$$
L'image de cette application est appelée code différentiel associé à $\Delta$, $D_a$, $D_b$ et $G$. On le note $C_{\Omega, S}(\Delta, D_a,D_b,G)$.
\end{defn}

\begin{rem}\label{sym}
On peut également construire une application $\res^2_{D_b, \Delta}$ en échangeant $D_a$ et $D_b$ dans l'énoncé de la définition \ref{defdiff}. Le théorème \ref{FR2} entraîne 
$$
\res^2_{D_a, \Delta}=-\res^2_{D_b, \Delta}.  
$$
De ce fait, les applications sont différentes mais ont même image. La construction du code ne dépend donc pas de l'ordre des éléments dans le couple $(D_a,D_b)$. 
\end{rem}

\begin{rem}
L'application $\res^2_{D_a, \Delta}$ peut en fait être définie sur $\Omega^2_{\F_q(S)/\F_q}$ tout entier. Aussi, on s'autorisera à l'appliquer à des $2$-formes quelconques de $\Omega^2_{\F_q(S)/\F_q}$ voire de
$\Omega^2_{\overline{\F}_q(\overline{S})/\overline{\F}_q}$. Cet abus de notation sera par exemple utilisé dans la définition \ref{deltac}.   
\end{rem}

La définition \ref{defdiff} n'est pas complètement satisfaisante, car il n'y a pas de lien entre le couple $(D_a,D_b)$ et $\Delta$. De fait, il se peut par exemple que les supports de $D_a$ et $D_b$ ne se croisent en aucun point du support de $\Delta$, ce qui, d'après le corollaire \ref{oukisont2}, donnerait un code nul.
Nous allons donc introduire une nouvelle notion permettant de relier un $0$-cycle sur $S$ à une paire de diviseurs.
Notons que cette définition (définition \ref{deltac}) pourra sembler inutilement compliquée au premier abord. Cependant, les commentaires en section \ref{discute} justifieront à posteriori la pertinence de ce choix.

\subsection{Paires de diviseurs $\Delta$-convenables}

Commençons par se donner un cahier des charges. On souhaite disposer des propriétés suivantes.
\begin{enumerate}
\item[$(a)$] On aimerait que les codes différentiels construits à partir de $\Delta$, $G$ et du couple $(D_a, D_b)$ n'aient pas une coordonnée systématiquement nulle.
On souhaiterait donc que pour tout point $P$ appartenant au support de $\Delta$, il existe une section de $\Omega^2 (G-D)$ qui n'annule pas l'application $\res^2_{D_a,P}$.
\item[$(b)$] Notre but est également d'obtenir une relation d'orthogonalité entre les codes $C_{L,S}(\Delta,G)$ et $C_{\Omega,S}(\Delta, D_a,D_b,G)$. Pour ce faire, nous allons utiliser la troisième formule des résidus (théorème \ref{FR3}) et adopter une démarche proche de celle qui est utilisée dans le cas des courbes.
\end{enumerate}

\noindent La définition suivante répond à ce cahier des charges.

\begin{defn}\label{deltac}
Soient $D_a$ et $D_b$ deux diviseurs sur $S$ dont les supports n'ont pas de composante irréductible commune et soit $D$ le diviseur $D:=D_a+D_b$.
La paire $(D_a,D_b)$ est dite $\Delta$-convenable si elle vérifie les conditions suivantes.
\begin{enumerate}
\item[$(i)$]\label{deltac1} Pour tout point $P$ de $\overline{S}$, l'application 
$
\res^2_{D_a,P}:\overline{\Omega^2(-D)}_P \rightarrow \overline{\F}_q
$
est $\mathcal{O}_{\overline{S},P}$-linéaire. On rappelle que $\overline{\Omega^2(-D)}_P$ désigne la fibre en $P$ du tiré en arrière sur $\overline{S}$ du faisceau $\Omega^2(-D)$
\item[$(ii)$]\label{deltac2} L'application $\res^2_{D_a,P}$ définie ci-dessus est surjective pour tout point $P$ appartenant au support de $\Delta$ et nulle pour tout autre point de $\overline{S}$.
\end{enumerate}
\end{defn}

\noindent \textbf{Attention.}
Même si les propriétés requises dans la définition \ref{deltac} sont d'ordre géométriques, c'est-à-dire qu'elles concernent $\overline{S}$, les diviseurs $D_a$ et $D_b$ sont \textbf{rationnels}, c'est-à-dire définis sur $\F_q$.

\begin{rem}
La structure de $\mathcal{O}_{\overline{S},P}$-module de $\overline{\F}_q$ est induite par l'application d'évaluation
$
f \rightarrow f(P)
$.  
Aussi, la condition $(i)$ signifie que pour toute fonction $f$ régulière au voisinage de $P$ et tout germe de $2$-forme $\omega$ appartenant à $\overline{\Omega^2(-D)}_P$, on a
$$
\res^2_{D_a,P}(f\omega)=f(P)\res^2_{D_a,P}(\omega).
$$
Notons également que si $(D_a,D_b)$ vérifie $(i)$, alors l'application $\res^2_{D_a,P}$ s'annule sur $\mathfrak{m}_{\overline{S},P}\overline{\Omega^2(-D)}_P$.
\end{rem}

\begin{rem}\label{DasymDb}
Par un raisonnement analogue à celui qui est utilisé dans la remarque \ref{sym}, on montre aisément que si $(D_a,D_b)$ est $\Delta$-convenable, alors l'application $\res^2_{D_b,\Delta}$ vérifie les mêmes propriétés de $\mathcal{O}_{\overline{S}}$-linéarité que $\res^2_{D_a, \Delta}$.
Par conséquent la notion de $\Delta$-convenance est symétrique.
$$
(D_a,D_b)\ \textrm{est}\ \Delta \textrm{-convenable}\ \Longleftrightarrow
(D_b,D_a)\ \textrm{est}\ \Delta \textrm{-convenable}.
$$
\end{rem}

Nous allons maintenant donner une critère de $\Delta$-convenance faisant intervenir des propriétés d'intersection entre les composantes des diviseurs $D_a$ et $D_b$.

\begin{prop}[Critère de $\Delta$-convenance]\label{crit}
Soit $(D_a, D_b)$ une paire de diviseurs dont les supports n'ont pas de composante irréductible commune et soit $D$ la somme de ces deux diviseurs. Si $D_a$ et $D_b$ vérifient les conditions suivantes, alors la paire $(D_a,D_b)$ est $\Delta$-convenable.
\begin{enumerate}
\item\label{crit1} Pour tout $P$ appartenant au support de $\Delta$, il existe une courbe irréductible $C$ définie sur $\F_q$, lisse en $P$ telle que, sur un voisinage $U$ de $P$ on ait
$$
{D_a}^+_{| U}=C \cap U \quad \textrm{ou} \quad {D_b}^+_{| U}=C \cap U \quad \textrm{et}\quad
m_P(C,D-C)=1.
$$
\item\label{crit2} Pour tout point géométrique $P$ de $\overline{S}$ n'appartenant pas au support de $\Delta$, alors l'un des diviseurs $D_*=D_a$ ou $D_*=D_b$ vérifie les conditions suivantes. Pour toute composante $\overline{\F}_q$-irréductible $\overline{C}$ de $D_*^+$ contenant $P$, on a: 
  \begin{enumerate}
  \item\label{crit21} la courbe $\overline{C}$ est lisse en $P$;
  \item\label{crit22} la courbe $\overline{C}$ apparaît dans la décomposition de $D_*$ en combinaison $\Z$-linéaire de composantes $\overline{\F}_q$-irréductibles avec le coefficient $1$;
  \item\label{crit23} $m_P(\overline{C}, D-\overline{C})\leq 0$.
  \end{enumerate}
\end{enumerate}
\end{prop}

\begin{rem}
  Ce critère, quoique technique présente un avantage majeur, il permet de construire des paires de diviseurs $\Delta$-convenables. La preuve de la proposition \ref{Deltacex} fournit un algorithme de construction d'une paire $\Delta$-convenable étant donné un $0$-cycle $\Delta$ (voir aussi remarque \ref{remcrit}).
\end{rem}

Avant de fournir une démonstration de cette proposition, nous allons faire quelques remarques. Nous donnerons ensuite quelques illustrations pour tenter de se développer une intuition des conditions exigées par le critère. 

\begin{rem}\label{Det}
Dans la condition (\ref{crit2}) de la proposition \ref{crit}, le fait que $D_*$ soit $D_a$ ou $D_b$ dépend du point $P$. En d'autres termes, en chaque point $P$ de $\overline{S}$ évitant le support de $\Delta$, les conditions (\ref{crit21}), (\ref{crit22}) et (\ref{crit23}) doivent être vérifiées soit par $D_a$, soit par $D_b$. 
Par ailleurs, 
le diviseur $D_*$ peut-être nul au voisinage de $P$ (c'est d'ailleurs ce qui arrive en presque tout point de $\overline{S}$). Dans ce cas, les conditions (\ref{crit21}), (\ref{crit22}) et (\ref{crit23})  sont trivialement vérifiées. De fait, si au voisinage d'un point $P$ de $\overline{S}$, l'un des diviseurs $D_a$ ou $D_b$ est nul, c'est celui que l'on choisit pour jouer le rôle de $D_*$. Cela permet de se ramener à un nombre fini de vérifications.
\end{rem}

Dans tout ce qui suivra nous utiliserons le code de couleurs suivant. 

$$
\begin{array}{c>{\centering\arraybackslash}m{1cm}>{\centering\arraybackslash}m{1.5cm}c>{\centering\arraybackslash}m{1cm}}
D_a^+ & \includegraphics{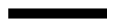} &  & D_a^- & \includegraphics{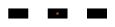}\\
D_b^+ & \includegraphics{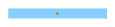} &  & D_b^- & \includegraphics{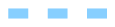}
\end{array}
$$

\noindent Nous allons illustrer les conditions du critère.
Pour ce faire, nous allons représenter des situations dans lesquelles ces situations sont vérifiées et d'autres dans lesquelles elles ne le sont pas.
Ces conditions sont locales.
Nous allons donc présenter deux séries de figures. 
La première série correspond au voisinage d'un point du support de $\Delta$ et la seconde au voisinage d'un point géométrique de $S$ non contenu dans le support de $\Delta$.

\medbreak

\noindent \textbf{En un point $P$ du support de $\Delta$.} Dans le tableau qui suit, les figures de la colonne de gauche représentent des situations où la condition (\ref{crit1}) du critère est vérifiée. 
Dans ce tableau ainsi que dans celui qui suit, on suppose que les courbes représentées apparaissent dans l'expression de $D_a$ (resp. de $D_b$) avec coefficient $1$.

$$
\begin{array}{>{\centering\arraybackslash}m{5cm}>{\centering\arraybackslash}m{1.5cm}>{\centering\arraybackslash}m{5cm}}
\textrm{\textbf{Vérifiées}} & \quad & \textrm{\textbf{Non vérifiées}} \\
 & & \\
\includegraphics{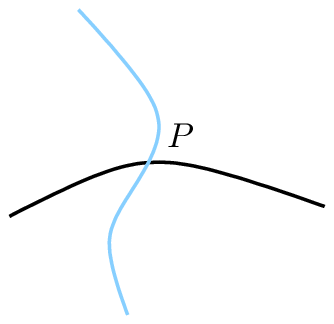} & & \includegraphics{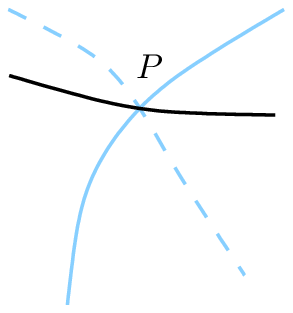} \\
\end{array}
$$
\medbreak
$$
\begin{array}{>{\centering\arraybackslash}m{5cm}>{\centering\arraybackslash}m{1.5cm}>{\centering\arraybackslash}m{5cm}}
\includegraphics{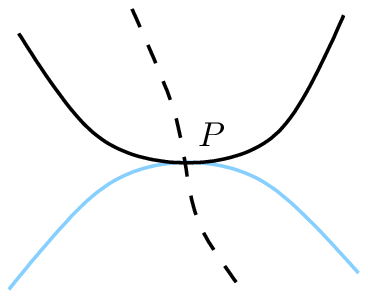} & & \includegraphics{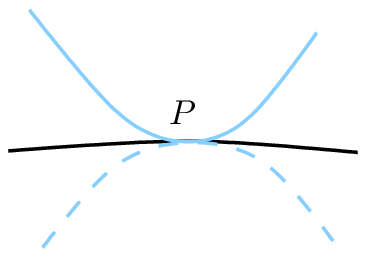} \\
\includegraphics{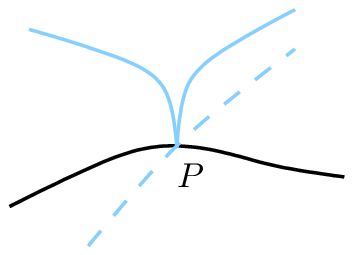} & & \includegraphics{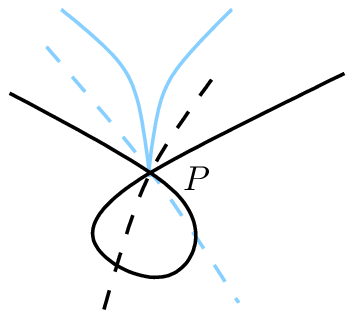}
\end{array}
 $$

\medbreak

\noindent \textbf{En un point $P$ n'appartenant pas au support de $\Delta$.}
Les figures de la colonne de gauche représentent des situations ou les conditions (\ref{crit21}), (\ref{crit22}) et (\ref{crit23}) sont vérifiées. Dans la colonne de droite elle ne le sont pas.

$$
\begin{array}{>{\centering\arraybackslash}m{5cm}>{\centering\arraybackslash}m{1.5cm}>{\centering\arraybackslash}m{5cm}}
\textrm{\textbf{Vérifiées}} & \quad & \textrm{\textbf{Non vérifiées}} \\
 & & \\
\includegraphics{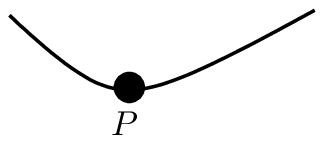} & & \includegraphics{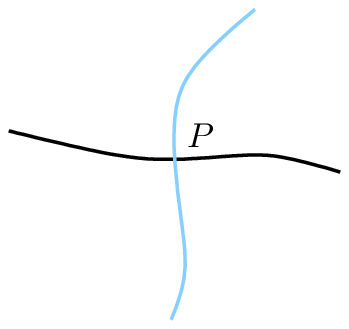} \\
\includegraphics{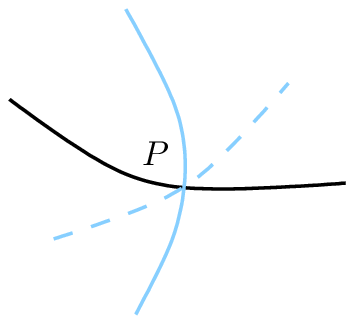} & & \includegraphics{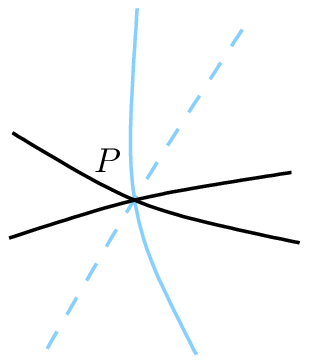} \\
\end{array}$$

$$
\begin{array}{>{\centering\arraybackslash}m{5cm}>{\centering\arraybackslash}m{1.5cm}>{\centering\arraybackslash}m{5cm}}
\includegraphics{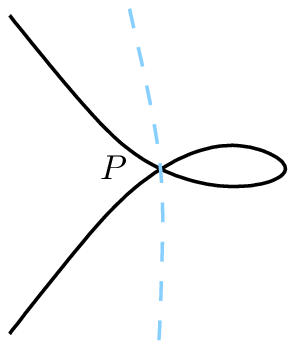} & & \includegraphics{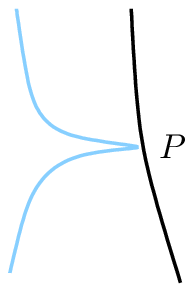} \\
\includegraphics{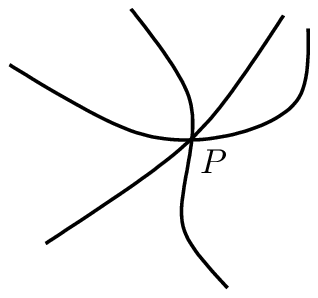} & & \includegraphics{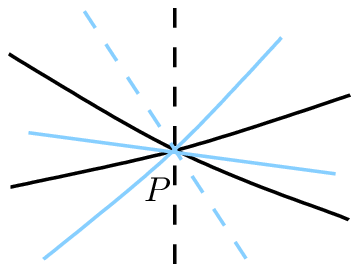} 
\end{array}
$$

\begin{rem}
  Le dernier exemple de la colonne de gauche peut surprendre. Pour le comprendre, on pourra se référer à la remarque \ref{Det}.
\end{rem}

\noindent La démonstration de la proposition \ref{crit} nécessite le lemme et le corollaire qui suivent.

\begin{lem}\label{valres}
Soient $\overline{C}$ une courbe $\overline{\F}_q$-irréductible plongée dans $\overline{S}$ et $P$ un point lisse de $\overline{C}$. Soit $\omega$ une $2$-forme sur $\overline{S}$ admettant un pôle simple le long de $\overline{C}$. Alors, le $1$-résidu de $\omega$ le long de $\overline{C}$ vérifie
$$
\textrm{val}_P (\res^1_{\overline{C}}(\omega))=m_P(\overline{C}, (\omega) + \overline{C}).
$$ 
\end{lem}

\begin{rem}
On déduit de ce lemme que le diviseur $(\res^1_{\overline{C}}(\omega))$ sur $\overline{C}$ est égal à $i^* ((\omega)+\overline{C})$, où $i$ désigne l'injection naturelle $i: \overline{C} \hookrightarrow \overline{S}$. Cette relation est utilisée par Serre dans \cite{ser} IV.8 lemme 2, pour démontrer la formule d'adjonction. 
\end{rem}

\begin{proof}[\textsc{Preuve du lemme \ref{valres}}]
Soient $\varphi, \psi$ et $v$ des équations locales respectives des diviseurs $\left((\omega)+\overline{C}\right)^+,\ \left((\omega)+\overline{C}\right)^-$ et $\overline{C}$ au voisinage de $P$. Soit $u$ un élément de $\overline{\F}_q(\overline{S})$ tel que $(u,v)$ soit une $(P,\overline{C})$-paire forte. Il existe une fonction $h$ sur $\overline{S}$ régulière et inversible au voisinage de $P$ telle que 
$$
\omega=h\frac{\varphi}{\psi} du \w \frac{dv}{v}.
$$
Le $1$-résidu de $\omega$ le long de $\overline{C}$ est $\bar{h}\bar{\varphi} \bar{\psi}^{-1} d\bar{u}$ et $\bar{h}$ est une fonction sur $\overline{C}$ régulière et inversible au voisinage de $P$. Par conséquent la valuation en $P$ de $\bar{h}d\bar{u}$ est nulle et
$$
\textrm{val}_P(\res^1_{\overline{C}}(\omega))=\textrm{val}_P(\bar{\varphi})-\textrm{val}_P (\bar{\psi}).
$$
De plus,
\begin{equation}\label{multip}
m_P(\overline{C}, (\omega)+\overline{C})  =  m_P(\overline{C}, ((\omega)+\overline{C})^+)  -  m_P(\overline{C}, ((\omega)+\overline{C})^-).
\end{equation}
On utilise ensuite la définition de la multiplicité d'intersection,
\begin{equation}\label{suce}
m_P(\overline{C}, ((\omega)+\overline{C})^+) = \dim_{\overline{\F}_q} \mathcal{O}_{\overline{S},P}/(\varphi,v) = \dim_{\overline{\F}_q} \mathcal{O}_{\overline{C},P}/(\bar{\varphi}) = \textrm{val}_P (\bar{\varphi}).
\end{equation}

\noindent De même, 
\begin{equation}\label{boule}
  m_P(\overline{C}, ((\omega)+\overline{C})^-) = \textrm{val}_P (\bar{\psi}).
\end{equation}

\noindent En injectant les résultats de (\ref{suce}) et (\ref{boule}) dans l'expression (\ref{multip}), on obtient le résultat recherché.
\end{proof}

\begin{cor}\label{OPlin}
Soient $\overline{C}$ une courbe $\overline{\F}_q$-irréductible plongée dans $\overline{S}$ et $P$ un point lisse de $C$. Soit $\omega$ une $2$-forme sur $\overline{S}$ telle que
$$
\textrm{val}_{\overline{C}}(\omega)\geq -1 \quad \textrm{et} \quad m_P(\overline{C},(\omega)+\overline{C})\geq -1.
$$ 
Alors, pour toute fonction rationnelle $f$ sur $\overline{S}$ régulière au voisinage de $P$, on a
$$
\res^2_{\overline{C},P}(f\omega)=f(P)\res^2_{\overline{C},P}(\omega).
$$
\end{cor}

\begin{proof}
  Soient $(u,v)$ une $(P,\overline{C})$-paire forte et $f$ une fonction rationnelle sur $\overline{S}$ régulière au voisinage de $P$. Il existe une fonction rationnelle $\psi$ sur $\overline{S}$ régulière au voisinage de $\overline{C}$ telle que 
$$
\omega = \psi du \w \frac{dv}{v}.
$$

\noindent Posons
$$\mu:= \res^1_{\overline{C}}(\omega)=\bar{\psi}d\bar{u}.$$

\noindent Comme $\omega$ est de valuation supérieure à $-1$ et $f$ de valuation positive le long de $\overline{C}$, on a
$$
\res^1_{\overline{C}}(f\omega) = \bar{f}\mu.
$$
D'après le lemme \ref{valres}, la valuation de $\mu$ en $P$ est supérieure ou égale à $-1$, donc
$$
\res_P(\bar{f}\mu)=\bar{f}(P)\res_P(\mu)\quad \Longrightarrow \quad \res^2_{\overline{C},P}(f\omega)=f(P)\res^2_{\overline{C},P}(\omega).
$$
\end{proof}

\begin{proof}[\textsc{Preuve de la proposition \ref{crit}}]
Soit $(D_a,D_b)$ une paire de diviseurs vérifiant le critère, c'est-à-dire les conditions (\ref{crit1}), (\ref{crit21}), (\ref{crit22}) et (\ref{crit23}) de la proposition \ref{crit}. Montrons qu'elle vérifie alors les conditions $(i)$ et $(ii)$ de la définition de $\Delta$-convenance. 

\medbreak

\noindent \textbf{Condition $(i)$.} Soient $P$ un point appartenant au support de
$\Delta$ et $\omega$ un germe de $2$-forme appartenant à $\overline{\Omega^2(-D)}_P$. D'après la condition (\ref{crit1}), il existe une courbe irréductible $C$ qui est égale à $D_a^+$ ou $D_b^+$ au voisinage de $P$. D'après la remarque \ref{DasymDb}, on peut supposer sans perte de généralité que $C$ est égale à $D_a^+$ au voisinage de $P$.
De fait,
$$
\res^2_{D_a,P}(\omega)=\res^2_{C,P}(\omega).
$$
De plus, la multiplicité d'intersection en $P$ de $C$ et $D_b$ est inférieure à $1$ donc
$$
m_P (C, (\omega)+C) \geq -1.
$$
Ainsi, comme la $2$-forme $\omega$ est de valuation supérieure ou égale à $-1$ le long de $C$, d'après le corollaire \ref{OPlin}, l'application $\res^2_{C,P}$ (donc $\res^2_{D_a,P}$) restreinte à $\overline{\Omega^2(-D)}_P$ est $\mathcal{O}_{\overline{S},P}$-linéaire.

\medbreak

\noindent \textbf{Condition $(ii)$.}
Soit $P$ un point de $\overline{S}$ hors du support de $\Delta$. Encore d'après la remarque \ref{DasymDb}, on peut supposer que le $D_a$ est le diviseur $D_*$ de la condition (\ref{crit2}) de la proposition \ref{crit}. 
Par conséquent, toute composante $\overline{\F}_q$-irréductible $\overline{C}$ du support de $D_a^+$ contenant $P$ est lisse en ce point, apparaît dans $D_a$ avec le coefficient $1$ et vérifie
\begin{equation}\label{interneg}
m_P(\overline{C},D-\overline{C}) \leq 0.
\end{equation}
Soient $\overline{C}$ une telle composante et $\omega$ un germe de $2$-forme appartenant à $\overline{\Omega^2(-D)}_P$. D'après la condition (\ref{crit22}), $\omega$ est de valuation supérieure ou égale à $-1$ le long de $\overline{C}$.
D'après le lemme \ref{valres}, l'inégalité (\ref{interneg}) entraîne que le $1$-résidu $\res^1_{\overline{C}}(\omega)$ de $\omega$ le long de $\overline{C}$ est de valuation positive en $P$. De fait, le $2$-résidu de $\omega$ en $P$ le long de $\overline{C}$ est nul.
Cette assertion est valable pour toute composante $\overline{\F}_q$-irréductible de $D_a^+$ au voisinage de $P$, d'après la définition \ref{resdiv}, on en déduit que 
$$
\res^2_{D_a,P}(\omega)=0.
$$
\end{proof}

\subsection{Exemples de diviseurs $\Delta$-convenables.}\label{Dconvex}

\subsubsection{Le plan projectif}\label{P2exdelta}
Si $S$ est le plan projectif $\P^2_{\F_q}$ et $X,Y,Z$ des coordonnées homogènes sur $S$. Soient $U$ la carte affine $U:=\{Z\neq 0\}$ et $\Delta$ la somme formelle de tous rationnels de $U$.
Pour tout élément $\alpha$ appartenant à $\F_q$, on définit les droites
$$
L_{a,\alpha}:=\{X=\alpha\} \quad \textrm{et}\quad L_{b,\alpha}:=\{Y=\alpha\},
$$
puis les diviseurs
$$
D_a:=\sum_{\alpha \in \F_q} L_{a,\alpha} \quad \textrm{et}\quad D_b:=
\sum_{\alpha \in \F_q} L_{b,\alpha}.
$$

\noindent La figure suivante représente cette paire de diviseurs dans le cas où le corps de base est $\F_3$. La droite en pointillés fins est la droite ``à l'infini'' et les $\bullet$ représentent les éléments du support de $\Delta$.

\begin{center}
    \includegraphics[width=7cm, height=7cm]{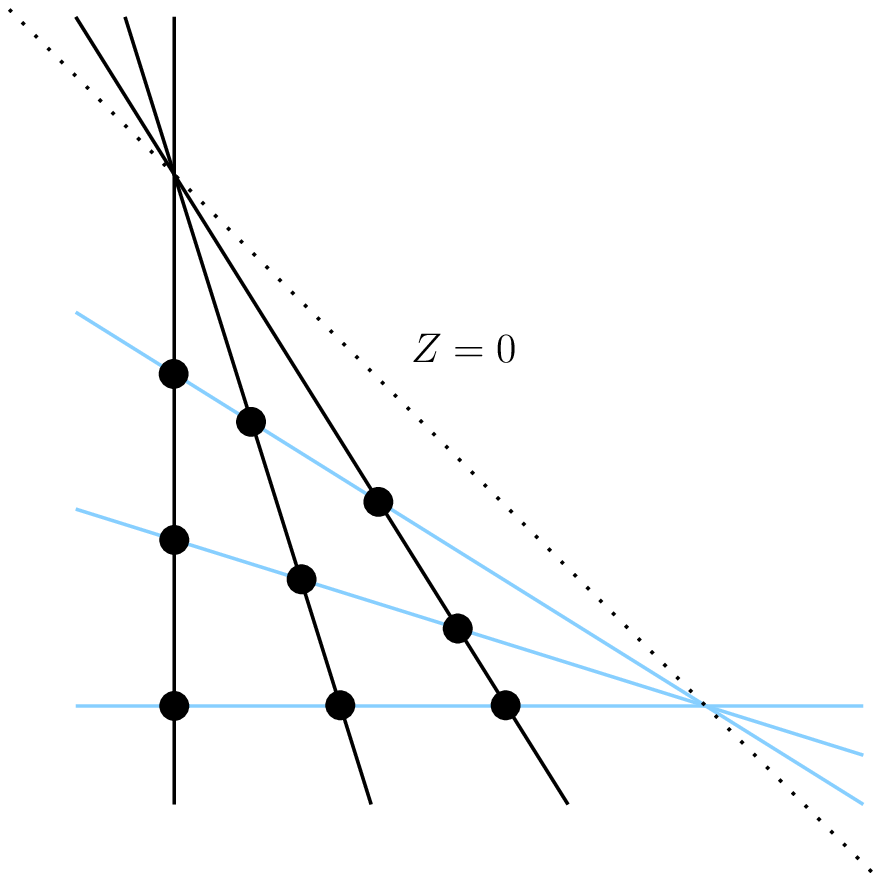}
\end{center}

La paire $(D_a,D_b)$ vérifie le critère de la proposition \ref{crit}. En effet, soit $P$ un point du support de $\Delta$ de coordonnées homogènes $(\alpha :\beta :1)$, alors $L_{a,\alpha}$ (resp. $L_{b,\beta}$) est la seule composante de $D_a$ (resp. $D_b$) passant par $P$, ces deux composantes sont lisses en $P$ et s'intersectent avec multiplicité $1$. En un point géométrique $P$ de $\overline{S}$ n'appartenant pas au support de $\Delta$, au moins l'un des diviseurs $D_a$ ou $D_b$ ne possède pas $P$ dans son support, on lui fait donc jouer le rôle de $D_*$ (voir remarque \ref{Det}). 

\subsubsection{Le produit de deux droites projectives}\label{P1prodexdelta}
Supposons que $S$ est la variété $\P^1_{\F_q} \times \P^1_{\F_q}$. Soit $((U,V),((X,Y))$ un système de coordonnées bihomogènes sur $S$. Soient $W$ la carte affine $W:=\{V \neq 0\}\cap\{Y\neq 0\}$ et $\Delta$ la somme des points rationnels de $W$.
Pour tout $\alpha$ appartenant à $\F_q$, on introduit les droites
$$
L_{a,\alpha}:=\{U=\alpha\} \quad \textrm{et}\quad L_{b,\alpha}:=\{X=\alpha\},
$$
puis les diviseurs
$$
D_a:=\sum_{\alpha \in \F_q} L_{a,\alpha} \quad \textrm{et}\quad D_b:=
\sum_{\alpha \in \F_q} L_{b,\alpha}.
$$

\noindent La figure suivante représente cette paire de diviseurs dans le cas où le corps de base est $\F_3$. Les droites en pointillés fins correspondent aux deux droites ``à l'infini'' et les $\bullet$ représentent les éléments du support de $\Delta$.

\begin{center}
\includegraphics[width=5.5cm, height=5.5cm]{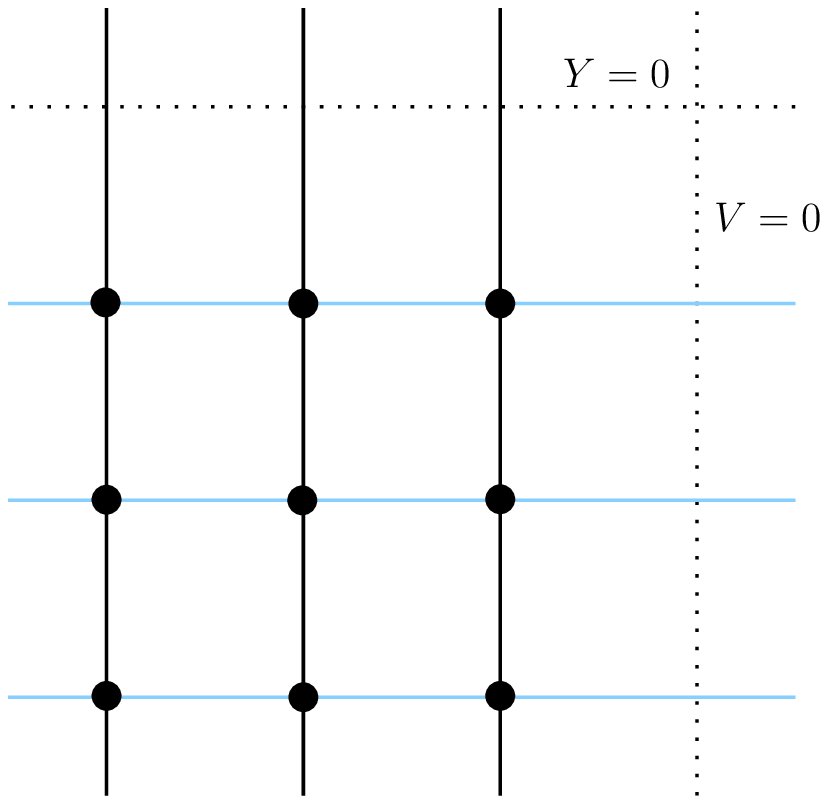}
\end{center}

\noindent On montre aisément que le couple de diviseurs ainsi construit satisfait le critère de la proposition \ref{crit}, c'est donc un couple de diviseurs $\Delta$-convenable.

\begin{rem}
La construction ci-dessus s'étend aisément au cas où $S$ est un produit de deux courbes admettant toutes deux des points rationnels.   
\end{rem}

\subsubsection{Quadriques lisses de $\P^3$}\label{qq}
Sur un corps algébriquement clos, une quadrique lisse s'obtient à partir de $\P^2$ en éclatant deux points $P$ et $Q$ puis en contractant la transformée stricte l'unique droite contenant ces deux points. De fait, une quadrique lisse de $\P^3$ est toujours géométriquement isomorphe à un produit de deux droites projectives. 
Si le corps de base est un corps fini $\F_q$, on distingue deux classes d'isomorphisme de quadriques lisses dans $\P^3$. Les quadriques hyperboliques sont $\F_q$-isomorphes à $\P^1\times \P^1$ et correspondent au cas où les points $P$ et $Q$ sont rationnels. Les quadriques elliptiques sont $\F_{q^2}$-isomorphes à $\P^1\times \P^1$ et correspondent au cas où les points $P$ et $Q$ sont définis sur $\F_{q^2}$ et conjugués sous l'action de $\textrm{Gal}(\F_{q^2}/\F_q)$.

De ce fait, on peut remarquer une relation entre les couples de diviseurs $\Delta$-convenables des deux exemples précédents. Partons de l'exemple où $S$ est le plan projectif.
Appelons $P$ et $Q$ les points de concours respectifs des composantes de $D_a$ et $D_b$ et $D$ la droite qui relie ces deux points.
Alors, le processus d'éclatements et contractions décrit ci-dessus permet d'obtenir le couple de diviseurs $\Delta$-convenables de l'exemple où $S$ est $\P^1\times \P^1$ à partir de celui où $S$ est $\P^2$, par le procédé suivant

$$\xymatrix{\relax
\includegraphics[width=4cm, height=4cm]{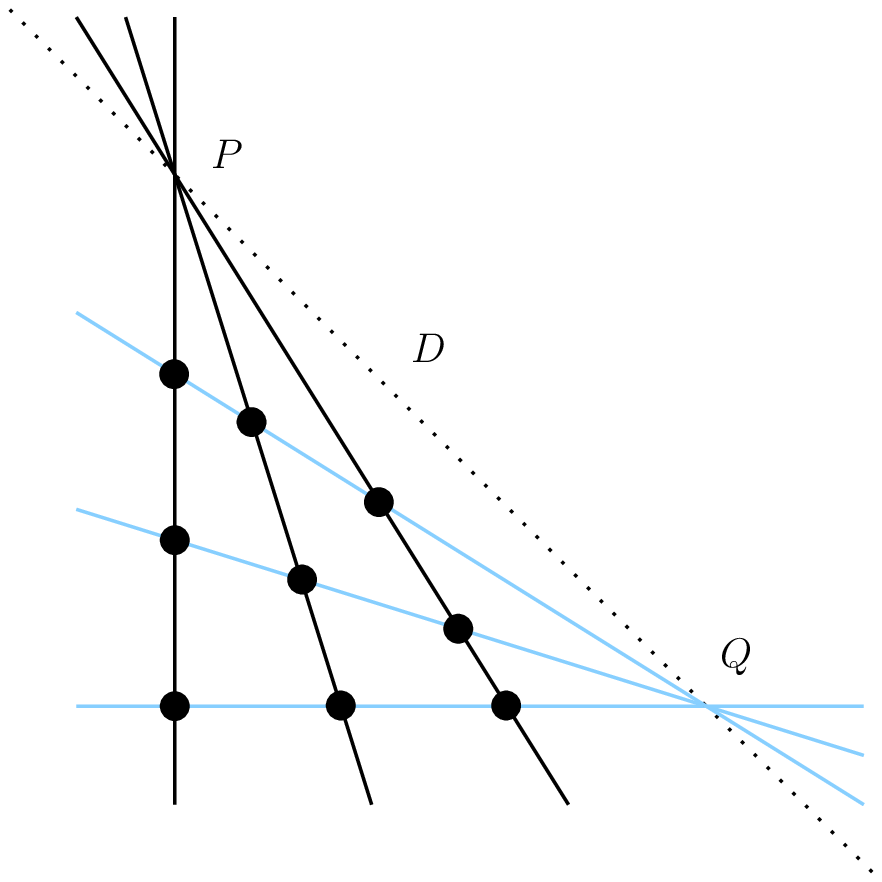}
\ar[r] &
\includegraphics[width=4cm, height=4cm]{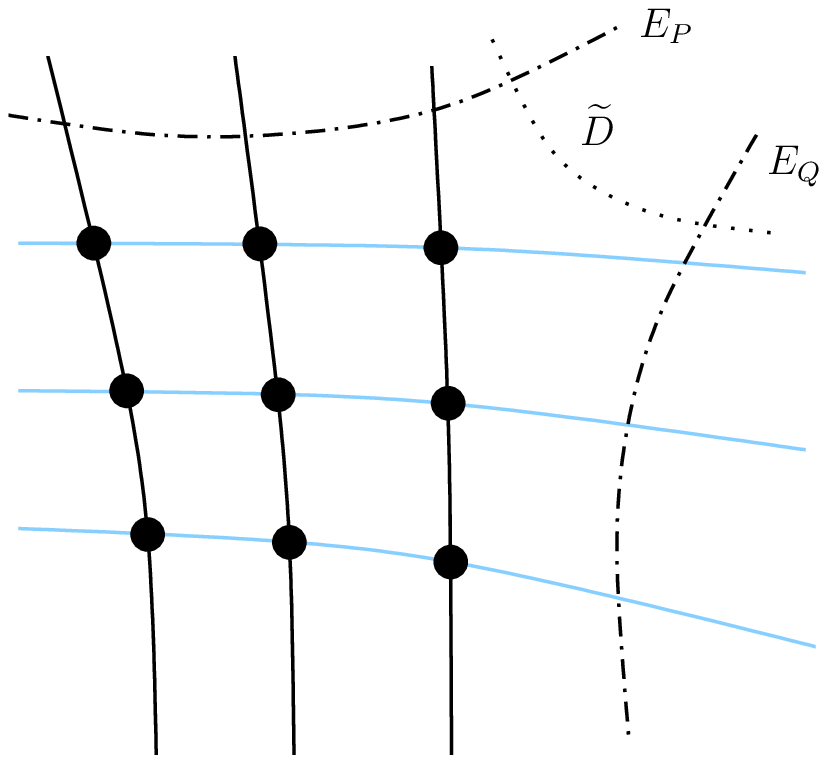}
\ar[r] &
\includegraphics[width=4cm, height=4cm]{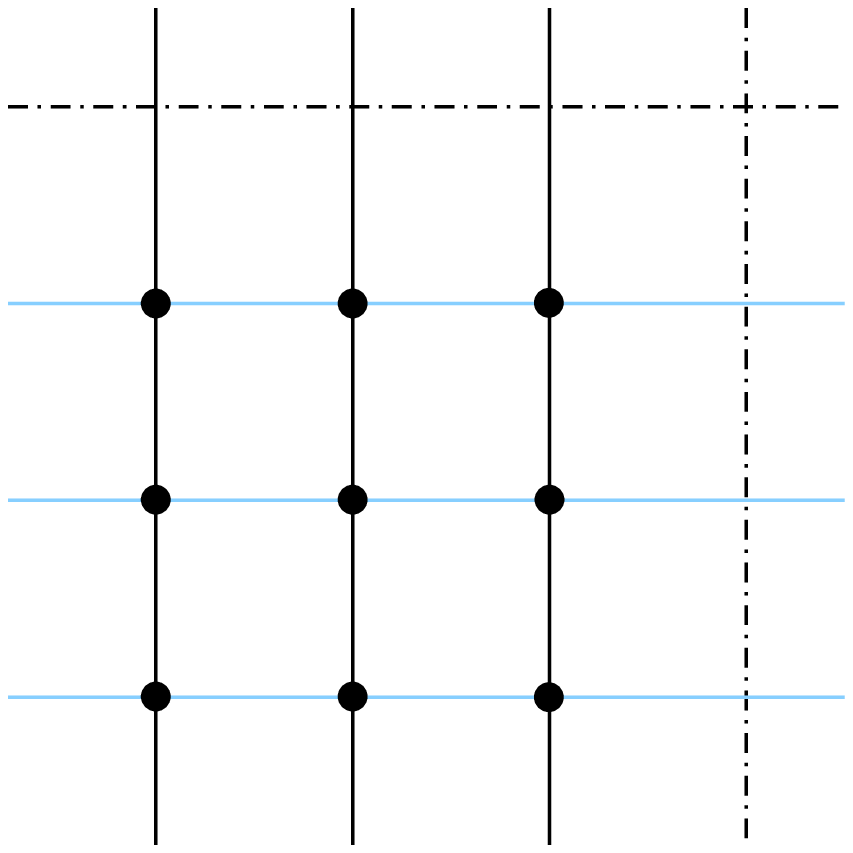}
}$$

\noindent Les courbes $E_P$ et $E_Q$ de la figure centrale sont les diviseurs exceptionnels correspondant respectivement à $P$ et $Q$. Dans la dernière figure, les courbes en pointillés sont les images de $E_P$ et $E_Q$ après contraction de $\widetilde{D}$.
\newpage 

\paragraph{Construction d'un diviseur $\Delta$-convenable sur une quadrique elliptique de $\P^3$.}

Dans cet exemple, on suppose que le corps de base $\F_q$ est de \textbf{caractéristique différente de $\mathbf{2}$}.
On considère une quadrique elliptique $Q$ plongée dans $\P^3$.
Soit $P_{\infty}$ un point rationnel de $Q$ et $\Delta$, la somme de tous les points rationnels de $Q$ sauf $P_{\infty}$.

D'après les commentaires sur les quadriques donnés page \pageref{qq}, la surface $Q$ se construit à partir de $\P^2$ en éclatant un point fermé de degré $2$ puis en contractant la transformée stricte de l'unique droite rationnelle contenant ce point.

Nous allons construire notre paire de diviseurs $\Delta$-convenable à partir de $\P^2$.
On considère le plan projectif muni d'un système de coordonnées homogènes $(X,Y,Z)$. Soit $D$, la droite d'équation $Z=0$ et $P$ un point fermé de $D$ de degré $2$. On notera $p$ et $p^{\sigma}$ les points correspondants après extension des scalaires de degré $2$, l'exposant $\sigma$ représente le conjugué sous l'action du Fr\"obenius.
On pose $\Delta_1$, la somme des points rationnels de $\P^2 \smallsetminus D$.
Après éclatement de $P$ et contraction de la transformée stricte $\widetilde{D}$ de $D$, on obtient une quadrique elliptique, l'image de $\Delta_1$ par cette opération est $\Delta$.

Soient $L$ la droite d'équation $X=0$ et $s$ la symétrie d'axe $L$ définie par
$$
s:(x:y:z) \mapsto (x:-y:z)
.$$

\noindent On rappelle que la caractéristique du corps $\F_q$ est supposée différente de deux dans cet exemple. De ce fait, l'application ci-dessus n'est pas l'identité.

\medbreak

\noindent \textbf{Construction de $\mathbf{D_a}$.} Soit $\alpha \in \F_{q^2}\smallsetminus \F_q$ dont la trace $\alpha+\alpha^q$ sur $\F_q$ est nulle. Soit $v\in \F_{q^2}^2$, le vecteur $v:=(1,\alpha)$. Le conjugué $v^{\sigma}$ de $v$ est égal au symétrique de $v$ par $s$.
On considère l'ensemble de coniques rationnelles contenant $P$ et admettant $\{v, v^{\sigma}\}$ comme vecteur tangent en ce point.
Le système linéaire correspondant est de dimension $1$, il y a donc $q+1$ coniques rationnelles satisfaisant ces contraintes. L'une d'entre elles est la droite double $2D$, les autres sont notées $C_1^a, \ldots, C_q^a$.

\begin{center}
\includegraphics[width=5cm, height=8.75cm]{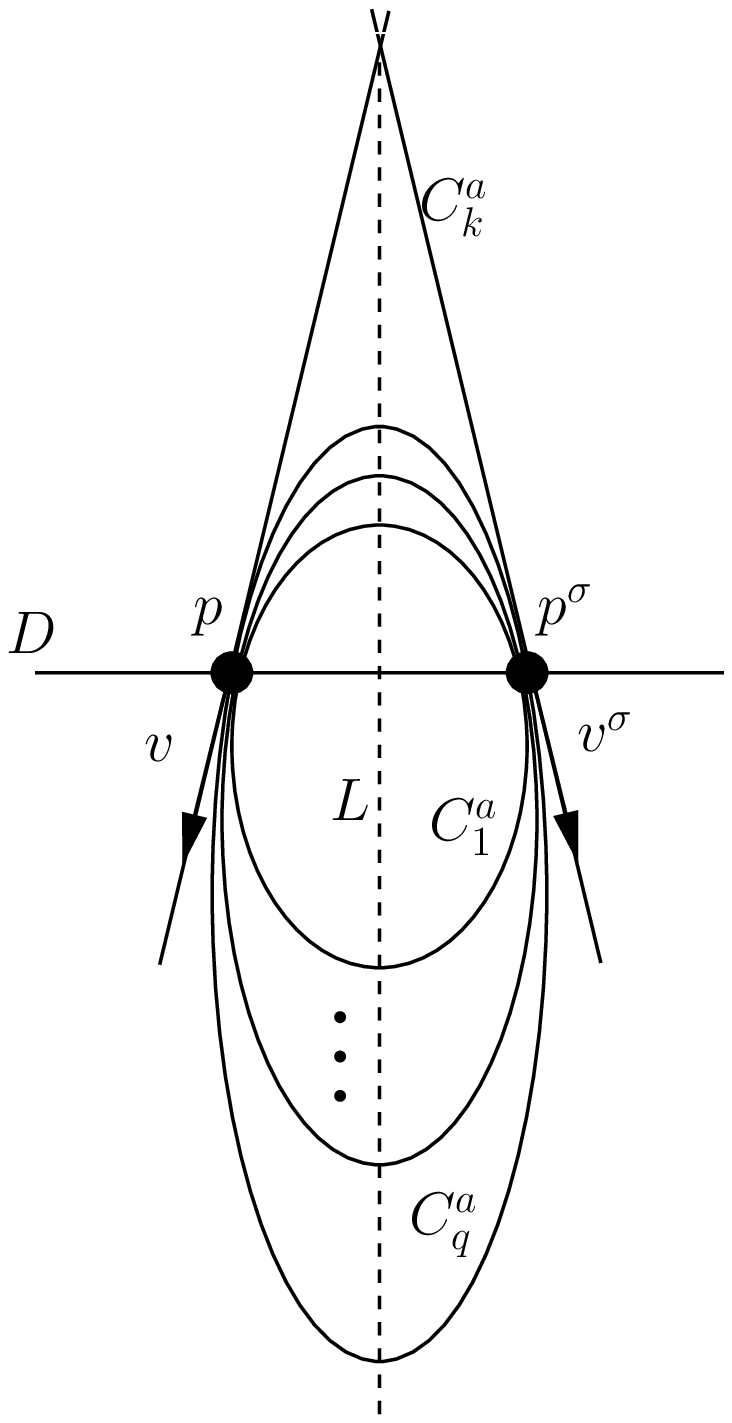}
\end{center}

\noindent La figure est de plus invariante par $s$.

\begin{rem}
  Comme cela apparaît, l'une des coniques notée $C_k^a$ dans la figure est dégénérée, elle est réunion de deux droites quadratiques conjuguées.
\end{rem}

On pose
$$
D_a:=\sum_{i=1}^q C_{i}^a.
$$

\medbreak

\noindent \textbf{Construction de $\mathbf{D_b^+}$.}
On se donne un autre vecteur $w\in \F_{q^2}^2\smallsetminus \F_q^2$ dont le conjugué coïncide avec son image par $s$ et on construit une seconde famille de coniques $C_1^b, \ldots, C_q^b$ comme dans l'étape précédente. On pose
$$
D_b^+:=\sum_{i=1}^q C_{i}^b.
$$
Ce diviseur est également invariant sous l'action de la symétrie $s$.

 Les diviseurs $D_a$ et $D_b^+$ sont tous deux linéairement équivalents à $2qL$ où $L$ est la classe d'équivalence linéaire d'une droite quelconque de $\P^2$. Le produit d'intersection de ces diviseurs est donc $4q^2$. Il va donc falloir éliminer des points en ajoutant à $D_b$ une partie négative. 
Prenons le temps de décrire le $0$-cycle d'intersection $D_a \cap D_b^+$.
\begin{enumerate}
\item Tous les points du support de $\Delta$ y apparaissent. Les points doubles correspondent à une tangence entre un élément du support de $D_a$ et un du support de $D_b$.
L'invariance de ces deux diviseurs sous l'action de $s$ entraîne que les points doubles sont sur l'axe de symétrie $L$.
On a donc dans ce $0$-cycle les $q^2$ points du support de $\Delta_1$ dont $q$ points de $L$ qui apparaissent avec coefficient $2$.
\item La multiplicité d'intersection de $D_a$ et $D_b^+$ en $p$ (resp. $p^{\sigma}$) est $q^2$, donc le point $P$ apparaît $q^2$ fois dans ce $0$-cycle.
\item Il reste donc $q^2-q$ points géométriques à identifier dans ce $0$-cycle. Ils s'agit en fait de $(q^2-q)/2$ points de degré $2$ provenant de l'intersection d'un élément de $\supp (D_a)$ et d'un élément de $\supp (D_b)$.
\end{enumerate}

\medbreak

\noindent \textbf{Construction de $\mathbf{D_b^-}$.} 
Pour construire $D_b^-$, nous allons avoir besoin de donner des équations explicites pour $D_a$ et $D_b^+$.
Soit $a\in \F_q^{\times} \smallsetminus {\F_q^{\times}}^2$ et $\alpha\in \F_{q^2}$ une racine carrée de $a$.
On peut supposer que le point $p$ est de coordonnées $(1:\alpha:0)$. On peut alors se convaincre du fait que
les deux équations suivantes fournissent de bons candidats pour $D_a$ et $D_b^+$.
$$
H_a =  \prod_{t\in \F_q} (x^2+y^2+tz^2)
$$

\noindent et 
$$
H_b =  \prod_{t\in \F_q} ((x-z)^2+y^2+tz^2)=  \prod_{t'\in \F_q} (x^2+y^2-2xz+t'z^2).
$$

Trouver un point d'intersection dans le plan affine de $D_a$ et $D_b^+$ revient à trouver un point
d'intersection entre une conique du support de $D_a$ et une conique du support 
de $D_b$. Ce qui revient à résoudre le système:
$$
\left\{
  \begin{array}{rcl}
    x^2+y^2+t & = & 0 \\
    x^2+y^2-2x+u & = & 0
  \end{array}
\right.
\Longleftrightarrow
\left\{
  \begin{array}{rcl}
    x^2+y^2+t & = & 0 \\
    x & = & \frac{t-u}{2}
  \end{array}
\right.
\Longleftrightarrow
\left\{
  \begin{array}{rcl}
     x & = & \frac{t-u}{2}\\
     y^2 & = & t - (\frac{t-u}{2})^2.
  \end{array}
\right.
$$

On vérifie ensuite que l'application $(t,u)\rightarrow t - (t-u)^2/4$ est surjective de $\F_q\times \F_q$ dans $\F_q$.
En d'autres termes les points d'intersections dans le plan affine de deux telles coniques sont soit des points $\F_q$-rationnels du plan affine soit des points de la forme $(s,\tau)$ où $s$ est un élément de $\F_q$ et $\tau$ un élément de $\F_{q^2}\smallsetminus \F_q$ dont le carré est dans $\F_q$. 

Rappelons que les points doubles dans l'intersection de $D_a$ et $D_b^+$ sont tous sur l'axe de symétrie, à savoir la droite d'équation $y$. Aussi, le diviseur effectif d'équation
$$
H_c:=y \prod_{s\in \F_q^{\times} \smallsetminus {\F_{q}^{\times}}^2} (y^2-sz^2)
$$
fournit un bon candidat pour le diviseur $D_b^-$.

\noindent On pose alors
$
D_b:=D_b^+-D_b^-.
$

\medbreak

On définit enfin sur $Q$ les diviseurs $\widetilde{D}_a$ et $\widetilde{D}_b$ construits comme étant les transformées strictes des diviseurs du même nom par l'opération d'éclatement/contraction.
Le $0$-cycle d'intersection de ces diviseurs est exactement $\Delta$, car l'opération d'éclatement a ``séparé'' les diviseurs $D_a$ et $D_b$ en $P$.
On vérifie alors que cette paire de diviseurs vérifie le critère de la proposition \ref{crit}, elle est donc $\Delta$-convenable.

\begin{rem}
  Le groupe de Picard d'une quadrique elliptique est libre de rang $1$ et engendré par la classe d'une section plane $L_Q$.
On Peut aisément montrer que les diviseurs ainsi construits vérifient
$$
\widetilde{D}_a \sim \widetilde{D}_b^+ \sim \widetilde{D}_b^- \sim qL_Q.  
$$
\end{rem}

\subsubsection{Autres exemples}

D'une façon générale le calcul de paires de diviseurs $\Delta$-convenables est ardue.
Toutefois, le lemme \ref{Deltacex} et la remarque \ref{remcrit} stipulent que des paires de diviseurs $\Delta$-convenables
sont explicitement calculables via des méthodes d'interpolation implémentables sur ordinateur.
Un programme appelé \textsc{DeltaConv} permettant de calculer des paires de diviseurs $\Delta-convenables$ à l'aide du logiciel \textsc{magma} est proposé en annexe \ref{prgmDconv}.
Avec l'aide de ce programme nous avons calculé des paires de diviseurs $\Delta$-convenables pour quelques exemples moins triviaux que ceux qui précèdent.

\paragraph{Sur une surface Hermitienne.} On considère la surface Hermitienne sur $\F_4$ plongée dans $\P^3$ d'équation
$
X^3+Y^3+Z^3+T^3=0$.
On la munit du $0$-cycle égal à la somme de ses points rationnels dans la carte affine $\{Z \ne 0 \}$. 
Le programme nous retourne les résultats suivants.

\medbreak

\begin{lstlisting}[frame=single]
> 
> F<w>:=FiniteField(4);
> P3<x,y,z,t>:=ProjectiveSpace(F,3);
> Herm:=Scheme(P3,x^3+y^3+z^3+t^3);
> 
> 
> A:=DeltaConv(Herm,Points,1000,1000,10);
D_a a ete trouve de maniere deterministe.
D_b^+ a ete trouve de maniere deterministe. 
D_b^- a ete trouve de maniere deterministe.


Le diviseur d'equation 
x*y^3 + x*z^3 + y^4 + y*t^3 + z^4 + z*t^3
 convient pour D_a. 

Le diviseur d'equation 
w^2*x^4 + w^2*x*t^3 + w*y^4 + w*y*t^3 + z^4 + z*t^3
 convient pour D_b^+. 

Le diviseur d'equation 
x^2 + w^2*y^2 + w*z^2
 convient pour D_b^-. 
\end{lstlisting}

\paragraph{Sur une surface quartique.} On considère la surface d'équation
$X^4+Y^4+Z^4+T^4=0$ définie sur $\F_3$. On prend comme $0$-cycle $\Delta$, la somme des points rationnels de la carte affine $\{Z=0\}$ de cette surface. On obtient le résultat suivant.

 \medbreak

\begin{lstlisting}[frame=single]
> S:=Scheme(P3,x^4+y^4+z^4+t^4);
> A:=DeltaConv(Herm,Points,1000,1000,10);

Le diviseur d'equation 
x^2 + y^2 + z^2 + t^2
 convient pour D_a. 

Le diviseur d'equation 
z^3 + 2*z*t^2
 convient pour D_b^+. 

Le diviseur d'equation 
x^3*z + 2*x^3*t + 2*x^2*y^2 + x^2*y*z + x^2*y*t +
     x^2*z^2 + x^2*z*t + x^2*t^2 + x*y^2*z + 2*x*y^2*t +
     x*y*z^2 + 2*x*z^2*t + 2*x*z*t^2 + 2*x*t^3 + y^4 +
     y^3*z + y^3*t + y^2*z^2 + y^2*z*t + y*z^2*t +
     2*y*z*t^2 + y*t^3 + z^3*t + 2*z^2*t^2 + z*t^3 +
     2*t^4
 convient pour D_b^-. 
\end{lstlisting}

\subsection{Discussion sur la $\Delta$-convenance et le critère}\label{discute}

Au vu des exemples du plan projectif et du produit de deux droites projectives, on est tenté d'envisager une définition nettement plus simple.
Un couple $(D_a, D_b)$ serait $\Delta$-convenable si et seulement si les diviseurs $D_a$ et $D_b$ étaient effectifs et le $0$-cycle d'intersection $D_a \cap D_b$ vérifiait
$$
D_a \cap D_b =\Delta.
$$

On montre aisément qu'une telle définition implique en fait la $\Delta$-convenance (il implique le critère de la proposition \ref{crit}).
Cependant, elle n'est pas vraiment intéressante en ce sens où, étant donné un $0$-cycle $\Delta$ sur $S$, une paire de diviseurs vérifiant de telles conditions n'existe pas en général.
Par exemple, si $S$ est le plan projectif et $\Delta$ la somme de trois points rationnels non alignés on ne peut construire une paire $(D_a, D_b)$ vérifiant ces conditions.
En effet, supposons qu'une telle paire existe et soit $L$ une droite de $S=\P^2$, alors la classe de $L$ engendre la groupe de Picard de $S$ et il existe deux entiers positifs $n_a$ et $n_b$ tels que $D_a \sim n_a L$ et $D_b \sim n_b L$.
De plus, le produit d'intersection de $D_a$ et $D_b$ est $3$ par hypothèse.
Donc, comme l'auto-intersection de $L$ est égale à $1$, on en déduit que $n_a$ ou $n_b$ est égal à $1$, ce qui contredit le fait que les éléments du support de $\Delta$ sont non alignés.

\section{Relations entre codes fonctionnels et différentiels sur une surface}

Nous étudions dans cette section l'extension aux surfaces de relations connues en théorie des codes construits à partir de courbes. On rappelle que l'on se place toujours dans le cadre décrit en section \ref{diffcadre}.

\subsection{Relation d'orthogonalité}

Le théorème qui suit est celui qui a motivé l'introduction de la notion de $\Delta$-convenance.

\begin{thm}[Théorème d'orthogonalité]\label{orthocode}
Soient $(D_a,D_b)$ une paire $\Delta$-convenable de diviseurs et $D:=D_a+D_b$. On a alors,
$$
C_{\Omega,S}(\Delta, D_a, D_b, G) \subseteq C_{L,S}(\Delta, G)^{\bot}.
$$
\end{thm}

\begin{proof}
  Soient $c$ un mot de $C_{L,S}(\Delta, G)$ et $c'$ un mot de $C_{\Omega,S}(\Delta, D_a, D_b, G)$. Il existe respectivement une fonction $f$ appartenant à $L(G)$ et une $2$-forme $\omega$ appartenant à $\Gamma (S, \Omega^2(G-D))$ telles que
$$
c=\ev_{\Delta}(f)\quad \textrm{et} \quad c'= \res^2_{D_a, \Delta}(\omega).
$$
On note ``$\langle \ \ \!, \ \rangle$'' la forme bilinéaire canonique sur $\F_q^n$. On a donc
$$
\langle c,c' \rangle=\sum_{i=1}^n f(P_i) \res^2_{D_a,P_i}(\omega).
$$
Comme le support diviseur $G$ est supposé éviter celui du $0$-cycle $\Delta$, on en déduit que $f$ est régulière au voisinage de tout point $P_i$ appartenant au support de $\Delta$. De plus, la $2$-forme $f\omega$ appartient à $\Gamma(S,\Omega^2 (-D))$. Donc, d'après la définition de $\Delta$-convenance, on a
$$
\forall P\in \overline{S},\ \res^2_{D_a,P}(f\omega)=\left\{
  \begin{array}{ccc}
    0 & \textrm{si} & P \notin \supp (\Delta) \\
f(P)\res^2_{D_a, P}(\omega) & \textrm{si} & P\in \supp (\Delta).
  \end{array}
\right.
$$
Par conséquent, 
$$
\langle c,c' \rangle=\sum_{i=1}^n \res^2_{D_a,P_i}(f\omega) = \sum_{\P \in \overline{S}}
\res^2_{D_a,P}(f\omega)
$$
et cette dernière somme est nulle d'après la troisième formule des résidus (théorème \ref{FR3}).
\end{proof}

En section \ref{P1P1} l'étude d'un exemple simple nous montrera que l'inclusion réciproque est en général fausse. Elle peut d'ailleurs être fausse pour tout choix de paire $\Delta$-convenable de diviseurs. Nous discuterons de ce défaut d'inclusion réciproque en section \ref{heuris}.

Avant de passer à la suite faisons une courte remarque sur les notations adoptées.
\medbreak

\noindent \textbf{Allègement des notations.}
Les notations de codes fonctionnels et différentiels sur $S$ sont respectivement $C_{L,S}(\Delta,G)$ et $C_{\Omega,S}(\Delta, D_a,D_b,G)$.
Dans ce qui suit, s'il n'y a pas d'ambiguïté sur la variété sur laquelle on travaille (en l'occurrence dans ce chapitre on travaille systématiquement sur $S$), on s'autorisera à ne pas la signaler en indice. On parlera alors de $C_L(\Delta, G)$ et $C_{\Omega}(\Delta, D_a, D_b, G)$.

\subsection{Un code différentiel est fonctionnel}

Nous avons vu en section \ref{codescourbes} qu'un code différentiel sur une courbe vérifiait deux propriétés intéressantes. La première est qu'il est l'orthogonal d'un code fonctionnel.
La seconde est que ce code s'identifie à un code fonctionnel associé à d'autres diviseurs.
Dans ce qui précède, nous avons cherché un analogue de la première propriété.
Nous allons maintenant en chercher un pour la seconde
et voir qu'à la différence de la première, cette seconde propriété s'étend parfaitement aux codes construits sur des surfaces.

Dans le cas des codes sur les courbes, ce résultat est une conséquence du théorème d'approximation faible aussi appelé théorème d'indépendance des valuations (voir \cite{sti} I.3.1).
De fait, nous allons avoir besoin d'un résultat analogue en dimension $2$. 
La proposition qui suit est sensiblement différente d'un énoncé de théorème d'indépendance de valuations, mais elle va nous fournir exactement le résultat nécessaire pour la suite.
 
\begin{prop}\label{approx}
Soit $C$ une courbe irréductible plongée dans $S$.
Soient $P_1,\ldots, P_r$ une famille de points fermés de $S\smallsetminus C$ et $Q_1, \ldots, Q_S$ une famille de points fermés de $C$.
Alors, il existe une uniformisante $v\in \mathcal{O}_{S,C}$ telle que le support du diviseur principal $(v)$ évite les points $P_1, \ldots, P_r$ et celui du diviseur $(v)-C$  évite les points $Q_1, \ldots, Q_s$.
\end{prop}

\begin{rem}
Une autre façon de formuler le résultat consiste à dire qu'il existe une fonction $v$ qui est une équation locale de $C$ au voisinage des points $Q_1, \ldots, Q_s$ et qui n'a ni zéro ni pôle en les points $P_1, \ldots, P_r$. 
\end{rem}

\begin{proof}
Soit $v_0$ une uniformisante de $\OC$. Alors, le diviseur de $v_0$ est de la forme
$$
(v_0)=C+D,
$$
où $D$ est un diviseur dont le support ne contient pas $C$.
D'après le \emph{moving lemma} (\cite{sch1} III.1.3 théorème 1), il existe un diviseur $D'$ linéairement équivalent à $D$ dont le support évite les points $P_1, \ldots,\ P_r,\ Q_1, \ldots,\ Q_r$. Ainsi, il existe une fonction $f$ rationnelle sur $S$ telle que 
$$
D'=D+(f).
$$ 
La fonction $v:=fv_0$ est solution du problème.
\end{proof}

\begin{rem}\label{remsch1}
  Dans le premier volume du livre \cite{sch1} de Shafarevich, le corps de base est supposé algébriquement clos. Cependant, l'étude de la preuve du ``\textit{moving lemma}'' permet de constater que cette hypothèse n'est pas utile pour prouver ce résultat.
Une preuve directe de la proposition \ref{approx} est donnée en annexe \ref{annexeapprox}. La logique de cette preuve est sensiblement la même que celle du moving lemma.
\end{rem}

\begin{cor}\label{omega0}
  Soit $(D_a,D_b)$ une paire $\Delta$-convenable de diviseurs et $D:=D_a+D_b$. Alors, il existe une $2$-forme $\omega_0$ rationnelle sur $S$ qui vérifie les propriétés suivantes.
\begin{enumerate}
\item\label{omega1} Il existe un ouvert $U$ contenant le support de $\Delta$ et tel que 
$
({\omega_0}_{|U})=-D_{|U}.
$
\item\label{omega2} Pour tout point $P$ appartenant au support de $\Delta$, on a 
$
\res^2_{D_a,P}(\omega_0)=1.
$
\item\label{omega3} Pour tout point $P$ appartenant au support de $\Delta$ et pour toute fonction $f$ régulière au voisinage de $P$, on a $\res^2_{D_a,P}(f\omega_0)=f(P)\res^2_{D_a,P}(\omega_0)$.
\end{enumerate}
\end{cor}

\begin{proof}
Soient $X_1, \ldots,\ X_r$ et $Y_1, \ldots,\ Y_r$ les composantes irréductibles respectives des supports de $D_a$ et $D_b$. Il existe des entiers $m_1, \ldots , m_r$ et $n_1, \ldots, n_s$ tels que
$$
D_a=m_1 X_1+\cdots +m_r X_r \quad \textrm{et} \quad D_b= n_1 Y_1+\cdots +n_s Y_s.
$$
D'après la proposition \ref{approx}, il existe un voisinage $U$ de $\supp (\Delta)$ et des fonctions $u_1, \ldots, u_r$ et $v_1, \ldots, v_s$ régulières sur $U$ telles que pour tout $i$ (resp. tout $j$), la fonction ${u_i}_{|U}$ est une équation de $X_i \cap U$ (resp. ${v_j}_{|U}$ est une équation de $Y_j \cap U$).

Soit $\mu$ une $2$-forme rationnelle sur $S$ n'ayant ni zéro ni pôle au voisinage du support de $\Delta$. Une telle $2$-forme existe, car d'après le \textit{moving lemma} (voir \cite{sch1} III.1.3 thm1 et remarque \ref{remsch1}), il existe un diviseur canonique dont le support évite celui de $\Delta$. 
Quitte à remplacer $U$ par un voisinage plus petit de $\supp (\Delta)$ on peut supposer que la $2$-forme $\mu$ restreinte à $U$ n'a ni zéro ni pôle.
On pose alors
$$
\omega:= \frac{\mu}{uv}.
$$
On a donc 
$$
(\omega_{|U})=-D_{|U}
$$
et d'après la définition de $\Delta$-convenance, pour tout $P$ appartenant au support de $\Delta$, on a
$$
\res^2_{D_a,P}(\omega)=a_P \neq 0.
$$
Par interpolation, on peut construire une fonction $g$ régulière au voisinage du support de $\Delta$ et telle que pour tout $P$ dans ce support,
$$
g(P)=a_P^{-1}.
$$
Quitte a réduire encore la taille de $U$, on peut supposer que $g$ n'a ni zéro ni pôle sur $U$. On pose alors
$$
\omega_0:=g\omega.
$$
Comme $g$ n'a ni zéro ni pôle sur $U$, les $2$-formes $\omega$ et $\omega_0$ restreintes à $U$ ont même diviseur, c'est-à-dire
$$
({\omega_0}_{|U})=-D_{|U}.
$$ 
La $\Delta$-convenance du couple $(D_a,D_b)$ permet de conclure que $\omega_0$ vérifie les propriétés requises.
\end{proof}

\begin{thm}\label{diff=fonc}
Soient $(D_a,D_b)$ une paire $\Delta$-convenable de diviseurs et $D:=D_a+D_b$, alors il existe un diviseur canonique $K$ tel que
$$
C_{\Omega}(\Delta,D_a,D_b,G)=C_L(\Delta, K-G+D).
$$   
\end{thm}

\begin{proof}
Soit $\omega_0$ une $2$-forme rationnelle sur $S$ vérifiant les propriétés du corollaire \ref{omega0} et soit $K$ son diviseur. D'après la propriété \ref{omega1} du corollaire \ref{omega0}, le diviseur $K$ est de la forme
$$
K=-D+R,
$$
où le support de $R$ évite celui de $\Delta$.
Soit $\omega$ une $2$-forme appartenant à $\Gamma (S, \Omega^2 (G-D))$, il existe une unique fonction $f$ dans $L(K-G+D)$ telle que 
$$
\omega=f\omega_0.
$$
Notons que
$$
K-G+D=-G+R.
$$
Aussi, les éléments de $L(K-G+D)$ sont des fonctions régulières au voisinage de $\supp (\Delta)$.
Soit $P$ un point du support de $\Delta$, d'après les propriétés \ref{omega2} et \ref{omega3} du corollaire \ref{omega0}, on a
$$
\res^2_{D_a,P}(\omega)=\res^2_{D_a,P}(f\omega_0)=f(P)\underbrace{\res^2_{D_a,P}(\omega_0)}_{=1}.
$$ 

\noindent On en déduit la relation
$$
\res^2_{D_a,\Delta}(\omega)=\ev_{\Delta}(f).
$$
\end{proof}

Nous avons montré que tout code différentiel est en fait un code fonctionnel associé à d'autres diviseurs.
Notons à ce stade que, dans le cas des courbes, la réciproque est élémentaire, à savoir: \textit{tout code fonctionnel est différentiel}.
Dans le cas des surfaces, cette réciproque est moins évidente.
En effet, étant donné un code fonctionnel $C_L(\Delta, G)$, si l'on veut prouver que ce code se réalise sous la forme d'un code différentiel, il faut d'abord disposer d'une paire $\Delta$-convenable de diviseurs. 

\subsection{Réciproque, un code fonctionnel est différentiel}

\begin{lem}[Existence d'une paire $\Delta$-convenable pour tout $\Delta$]\label{Deltacex}
Soient $S$ une surface algébrique projective lisse géométriquement intègre
et $Q_1, \ldots, Q_m$ une famille de points rationnels de $S$. Posons $\Delta:=Q_1+\cdots +Q_m$. Alors, il existe une infinité de paires $\Delta$-convenables qui vérifient le critère de la proposition \ref{crit}.
\end{lem}

\begin{rem}\label{remcrit}
 La démonstration qui suit est constructive. De fait, elle accentue l'intérêt de la proposition \ref{crit}. En effet, même si le critère qui y est énoncé est extrêmement technique, les paires qui le vérifient peuvent être calculées explicitement.
\end{rem}

\begin{proof}

\noindent \textbf{Étape 1: Construction de $\mathbf{D_a}$.} 
On choisit une courbe réduite $C$, éventuellement réductible qui contienne tout le support de $\Delta$ et qui soit régulière en chaque point de ce dernier. 
Assurons nous de l'existence d'une telle courbe. Soit $U$ un voisinage affine du support de $\Delta$, on note $\mathfrak{m}_{P_1}, \ldots, \mathfrak{m}_{P_n}$ les idéaux maximaux de $\F_q [U]$ correspondant aux points $P_1, \ldots, P_n$.
Il s'agit de choisir un élément de l'idéal produit $\mathfrak{m}_{P_1} \cdots \mathfrak{m}_{P_n}$ qui n'appartienne à aucun des idéaux $\mathfrak{m}_{P_i}^2$. Les idéaux $\mathfrak{m}_{P_i}^2$ étant deux à deux étrangers, d'après le théorème chinois, on a
l'isomorphisme
$$
\xymatrix{
\relax
\F_q [U] / \mathfrak{m}_{P_1}^2 \cdots \mathfrak{m}_{P_n}^2
\ar[r]^{\sim}
&
\prod_{i=1}^n \F_q [U] / \mathfrak{m}_{P_i}^2. }
$$ 
On choisit un élément $(\bar{a}_1, \ldots, \bar{a}_n)$ dans $\prod_i \mathfrak{m}_{P_i}/ \mathfrak{m}_{P_i}^2$ dont aucune des coordonnées $\bar{a}_i$ n'est nulle. On relève cet élément en une fonction $a$ de $\F_q[U]$ par le biais de l'isomorphisme ci-dessus. La fonction obtenue appartient bien à l'idéal produit $\mathfrak{m}_{P_1} \cdots \mathfrak{m}_{P_n}$ et n'est dans aucun des $\mathfrak{m}_{P_i}^2$.
La fermeture projective du lieu d'annulation de $a$ est une courbe $C$ vérifiant les conditions exigées. De plus, si l'on remplace $a$ par $a+m$ où $m$ est un élément de $\mathfrak{m}_{P_1}^2 \cdots \mathfrak{m}_{P_n}^2$ n'appartenant pas à l'idéal engendré par $a$, on obtient une courbe $C'$ distincte de $C$ et vérifiant les mêmes conditions.
On en déduit l'existence d'une infinité de courbes interpolant le support de $\Delta$ et lisse en chaque point de ce dernier.
Soit donc $C$ une telle courbe, on pose alors
$$D_a:=C_1+\cdots +C_k,$$
où les $C_i$ sont les composantes irréductibles de $C$.

\medbreak

\noindent \textbf{Étape 2: Construction de $\mathbf{D_b}$.}
On choisit, toujours par interpolation, un diviseur effectif $D'$ interpolant tous les points de $\supp (\Delta)$ et n'ayant pas de composante irréductible commune avec $D_a$.
Soit $\Theta$, le $0$-cycle obtenu par l'intersection au sens de la théorie des schémas des diviseurs $D_a$ et $D'$.
On a donc
$$
\Theta= \Delta + \Delta'
$$
où $\Delta'$ est un $0$-cycle effectif.
On choisit alors un diviseur $D''$ tel que
$$
D_a \cap D'' = \Delta' +\Delta''
$$
où le support de $\Delta''$ est disjoint de celui de $\Delta$. Le diviseur $D''$ se construit également par interpolation.
On pose enfin
$$
D_b:=D'-D''. 
$$
On montre aisément que la paire ainsi construite vérifie le critère de la proposition \ref{crit}. Elle est donc $\Delta$-convenable.
Comme il existe une infinité de façons de construire $D_a$ (et $D_b$) on en déduit qu'il existe une infinité de paires $\Delta$-convenables.
\end{proof}

\begin{thm}\label{fonc=diff}
  Étant donné un diviseur $G$ sur $S$, il existe un diviseur canonique $K$ et une paire $\Delta$-convenable $(D_a,D_b)$ telle que
$$
C_L(\Delta, G)=C_{\Omega}(\Delta,D_a,D_b,K-G+D).
$$
\end{thm}

\begin{proof}
Le lemme \ref{Deltacex} assure l'existence d'une paire $\Delta$-convenable $(D_a, D_b)$.
À partir de cette paire, on construit une $2$-forme $\omega_0$ en utilisant corollaire \ref{omega0}. On pose
$$
K:=(\omega_0).
$$  
D'après le théorème \ref{diff=fonc} on a 
$$
\begin{array}{rcl}
C_{\Omega}(\Delta,D_a,D_b,K-G+D) & = & C_L(\Delta, K-(K-G+D)+D)\\
 & = & C_L(\Delta,G).
\end{array}$$
Ce qui conclut la démonstration.
\end{proof}

\section{Défaut d'inclusion réciproque pour le théorème d'orthogonalité}\label{defaut}

Nous allons présenter deux exemples de surfaces, qui sont en l'occurrence les exemples les plus simples que l'on connaisse, à savoir $\P^2$ et $\P^1 \times \P^1$.
Nous allons voir que l'on dispose d'une inclusion réciproque systématique pour le premier exemple (le plan projectif) en choisissant une paire de diviseurs $\Delta$-convenable extrêmement simple.
Ensuite, nous observerons que dans le second exemple (le produit de deux droites projectives), l'inclusion réciproque pour le théorème d'orthogonalité n'a jamais lieu, et ce quel que soit la paire $\Delta$-convenable choisie.

\subsection{Codes sur le plan projectif}\label{exmpP2}

On reprend les notations de la section \ref{P2exdelta}.
On rappelle que $X,Y$ et $Z$ désignent des coordonnées homogènes sur $\P^2$, que l'ouvert $U$ est le complémentaire de la droite d'équation $Z=0$ et que le $0$-cycle $\Delta$ est la somme de tous les points rationnels de l'ouvert $U$.
Pour construire des codes fonctionnels on doit également introduire un diviseur $G$.
Soient $m$ un entier positif et $L_{\infty}$ la droite d'équation $Z=0$, on pose
$$
G_m:= mL_{\infty}.
$$
Notons que, comme la classe d'équivalence linéaire de la droite $L_{\infty}$ engendre le groupe de Picard de $\P^2$, on peut sans perte de généralité considérer que le diviseur $G$ est de la forme $G_m$.

\begin{rem}
  On a supposé que l'entier naturel $m$ était positif, on aurait pu omettre cette hypothèse. Cependant, si $m<0$, alors l'espace $L(G_m)$ est nul et le code fonctionnel le sera également. Nous avons donc choisi d'éviter cette situation totalement inintéressante. 
\end{rem}


\subsubsection{Codes fonctionnels}

Commençons par rappeler que les codes fonctionnels $C_L (\Delta, G_m)$ ne sont autre que des codes de Reed-Müller affines. 
Posons
$$
x:=\frac{X}{Z} \quad \textrm{et} \quad y:= \frac{Y}{Z}.
$$
L'espace vectoriel $L(G_m)$ s'identifie à l'espace $\F_q [x,y]_{\leq m}$ des polynômes en $x$ et $y$ de degré total inférieur ou égal à $m$. 
On rappelle que, par convention, si l'entier $m$ est strictement négatif, alors l'espace $\F_q [x,y]_{\leq m}$ est nul.
Ainsi, le code $C_L (\Delta, G_m)$ s'obtient par évaluation en tous les points du plan affine $U$ des éléments de l'espace vectoriel $\F_q [x,y]_{\leq m}$, c'est donc le code de Reed-Müller $RM_q (2,m)$.
Pour plus d'informations sur les codes de Reed-Müller, voire \cite{slmc1} chap 13 pour les codes binaires et \cite{delsarte} pour le cas général.

\subsubsection{Codes différentiels}

On reprend la paire $\Delta$-convenable $(D_a, D_b)$ de la section \ref{P2exdelta}. C'est à dire que $D_a$ (resp. $D_b$) est la somme de toutes les droites d'équations $x=\alpha$ (resp. $y=\alpha$) avec $\alpha \in\F_q$. 
Pour procéder à l'étude des codes différentiels  de la forme $C_{\Omega}(\Delta,D_a,D_b,G_m)$, nous allons utiliser le théorème \ref{diff=fonc} et chercher à quels codes fonctionnels ils s'identifient.
Pour ce faire, nous allons introduire explicitement une $2$-forme rationnelle $\omega_0$ vérifiant les propriétés du corollaire \ref{omega0}. Soit donc
$$
\omega_0:= \frac{dx}{\prod_{\alpha \in \F_q} (x-\alpha)} \w
\frac{dy}{\prod_{\beta \in \F_q} (y- \beta)}.
$$

\paragraph{Calcul du diviseur de $\omega_0$.}
D'après \cite{sch1} III.6.4, on sait qu'un diviseur canonique sur $\P^2$ est linéairement équivalent à $-3L_{\infty}$. De plus,
$$
({\omega_0}_{|U})= -D_{|U}.
$$

\noindent De fait, le diviseur canonique $(\omega_0)$ est de la forme
$$
(\omega_0)=kL_{\infty} -D
$$ 
où $k$ est un entier à déterminer. On sait également que, par construction, le diviseur $D$ est linéairement équivalent à $2q L_{\infty}$.
On en déduit la relation
$$
-3L_{\infty} \sim (k-2q)L_{\infty}, \quad \textrm{donc} \quad k=2q-3.
$$

\noindent En conclusion,
\begin{equation}\label{miche}
(\omega_0)=(2q-3)L_{\infty}-D=G_{2q-3}-D.
\end{equation}

\paragraph{Propriétés vérifiées par $\omega_0$.}
Maintenant que l'on connaît le diviseur $(\omega_0)$ et que l'on sait qu'il coïncide avec $-D$ sur le voisinage $U$ du support de $\Delta$, il reste à vérifier que $\omega_0$ vérifie les deux autres propriétés du corollaire \ref{omega0}. 
Soient $P$ un point appartenant au support de $\Delta$ et $x_P, y_P$ ses coordonnées affines dans $U$.
On appelle $C$, la droite d'équation $y=y_P$ et on calcule le $1$-résidu de $\omega_0$ le long de cette droite. On obtient
$$
\res^1_{C}(\omega_0)= \frac{1}{\prod_{\beta \neq y_P} (y_P -\beta)}\  \frac{d\bar{x}}{\prod_{\alpha \in \F_q} (\bar{x}-\alpha)}   .
$$  

\noindent On remarque que le produit $\prod_{\beta \neq y_P} (y_P -\beta)$ est égal au produit de tous les éléments de $\F_q^{\times}$, il est donc égal à $-1$. 

\noindent Calculons à présent le $2$-résidu en $P$ le long de $C$ de $\omega_0$, c'est à dire le résidu en $P$ de la $1$-forme sur $C$ ci-dessus. On obtient
$$
\res^2_{C,P}(\omega_0) = -\frac{ 1}{ \prod_{\alpha \neq x_P} (x_P-\alpha)} =1.
$$

\noindent La propriété \ref{omega3} est une conséquence immédiate de la $\Delta$-convenance de $(D_a, D_b)$.

\paragraph{Identification à des codes fonctionnels et orthogonalité.}
La $2$-forme $\omega_0$ vérifie les propriétés du corollaire \ref{omega0}. On en déduit que pour tout entier $m$, le code $C_{\Omega}(\Delta,D_a,D_b, G_m)$ s'identifie au code $C_L (\Delta, (\omega_0)-G_m+D)$. D'après le calcul du diviseur $(\omega_0)$ en (\ref{miche}), on conclut que pour tout entier $m$, on a 
$$
C_{\Omega}(\Delta,D_a,D_b,G_m)=C_L (\Delta, G_{2q-3-m}).
$$
D'après le théorème \ref{orthocode}, on a l'inclusion
$$
C_{\Omega}(\Delta,D_a,D_b,G_m) \subseteq C_L (\Delta, G_m)^{\bot}.
$$
On peut évaluer les dimensions de ces codes et montrer que l'inclusion réciproque est vérifiée. Ce résultat n'a  absolument rien de nouveau. Il est en effet connu que l'orthogonal d'un code de Reed-Müller est encore un code de Reed-Müller (voir \cite{slmc} ch 13 et \cite{delsarte} 3.2). 

En conclusion, l'orthogonalité parfaite entre code fonctionnel et code différentiel est obtenue dans ce cas élémentaire et très particulier. Il ne s'agit malheureusement pas d'un fait général.
L'exemple suivant, qui est pourtant presque aussi élémentaire, montre qu'en général on doit se contenter d'une inclusion stricte. 

\subsection{Codes sur un produit de deux droites projectives}\label{P1P1}

On reprend les notations de la section \ref{P1prodexdelta}. On rappelle que $((U,V),(X,Y))$ est un système de coordonnées bihomogènes sur $\P^1 \times \P^1$. On note $E$ et $F$ les droites d'équations respectives $V=0$ et $Y=0$.
On rappelle également que l'ouvert $U$ est le complémentaire de $E \cup F$.
Enfin, pour tout couple d'entiers $(m,n)$ on définit le diviseur $G_{m,n}$ par
$$
G_{m,n}:=mE+nF.
$$

\noindent Tout comme dans l'exemple précédent, on sait que l'on peut sans perte de généralité supposer que le diviseur $G$ intervenant dans  la construction du code fonctionnel est de la forme $G_{m,n}$. En effet, les classes d'équivalence linéaires de $E$ et $F$ engendrent le groupe de Picard de $\P^1 \times \P^1$.

\subsubsection{Codes fonctionnels}
Nous allons montrer tout d'abord que les codes fonctionnels de la forme $C_L (\Delta, G_{m,n})$ sont des produits tensoriels de codes de Reed-Solomon. Posons
$$
u:=\frac{U}{V} \quad \textrm{et} \quad x:= \frac{X}{Y}.
$$
L'espace vectoriel $L(G_{m,n})$ s'identifie au sous-espace de $\F_q [u,x]$ des polynômes de degré en $u$ inférieur ou égal à $m$ et de degré en $x$ inférieur ou égal à $n$. En d'autres termes on a l'identification
$$
L(G_{m,n}) \cong \F_q[u]_{\leq m} \otimes_{\F_q} \F_q[x]_{\leq n}.
$$
On note $RS_q (n)$ le code de Reed-Solomon de longueur $q$ obtenu par évaluation en tous les éléments de $\F_q$ des polynômes de $\F_q [t]_{\leq n}$. Le code fonctionnel sur $\P^1 \times \P^1$ est donc de la forme
$$
C_L (\Delta, G_{m,n})= RS_q (m) \otimes_{\F_q} RS_q (n).
$$

\paragraph{L'orthogonal ne peut être différentiel.}
D'après le théorème \ref{diff=fonc}, il suffit de montrer que l'orthogonal du code fonctionnel $C_L (\Delta, G_{m,n})$ n'est pas un code fonctionnel sur $\P^1 \times \P^1$. Un tel résultat entraînerait, qu'il n'existe aucun couple $\Delta$-convenable $(D_a, D_b)$ tel que le code $C_L (\Delta, G_{m,n})^{\bot}$ soit égal à $C_{\Omega}(\Delta,D_a,D_b, G_{m,n})$. On a vu dans le paragraphe précédent que le code $C_L (\Delta, G_{m,n})$ était égal au produit tensoriel des codes $RS_q (m)$ et $RS_q (n)$. Ces codes de Reed-Solomon sont non triviaux, si et seulement si
$$
0 \leq m \leq q-2 \quad \textrm{et} \quad 0 \leq n \leq q-2.
$$
Si les entiers $m$ et $n$ vérifient les encadrements ci-dessus, alors l'orthogonal du code $C_L (\Delta, G_{m,n})$ ne peut être fonctionnel. En effet, si $C_L (\Delta, G_{m,n})^{\bot}$ était un code fonctionnel $C_L(\Delta, G)$, alors $G$ serait linéairement équivalent à un certain $G_{a,b}$. De fait, le code fonctionnel $C_L (\Delta, G)$ serait isométrique\footnote{Au sens de la métrique de Hamming.} à $C_L (\Delta, G_{a,b})$ et cette isométrie serait représentée par une matrice diagonale dans la base canonique de $\F_q ^{q^2}$.
En regardant $\F_q^{q^2}$ comme le produit tensoriel de deux copies de $\F_q^q$, le code $C_L (\Delta, G_{a,b})$ est un produit tensoriel de deux codes et cette propriété est invariante sous l'action d'une isométrie diagonale. Ainsi, l'orthogonal $C_L (\Delta, G)$ de $C_L (\Delta, G_{m,n})$ serait un produit tensoriel de deux codes. Ce qui est impossible d'après le lemme \ref{tens} énoncé en annexe \ref{linalg}.

En conclusion, pour tout couple $\Delta$-convenable $(D_a, D_b)$ et tout couple d'entiers $(m,n)$ tous deux compris entre $0$ et $q-2$, on a
$$
C_{\Omega} (\Delta, D_a, D_b, G) \varsubsetneq C_L (\Delta, G)^{\bot}.
$$

\begin{rem}
Par un raisonnement identique, on peut montrer que ce défaut d'inclusion réciproque a lieu pour toute surface $S$ qui est un produit de deux courbes.
\end{rem}

\paragraph{Une réalisation de l'orthogonal.} 

D'après le lemme \ref{orthotens}, l'orthogonal du code $C_L (\Delta, G_{m,n})$ est une somme de deux produits tensoriels. À savoir
\begin{equation}\label{prout}
C_L (\Delta, G_{m,n})^{\bot}=RS_q (m)^{\bot}\otimes \F_q^q + \F_q^q \otimes RS_q(n)^{\bot}.
\end{equation}
L'orthogonal d'un code de Reed-Solomon étant encore un code de Reed-Solomon, les deux termes de la somme ci-dessus ($RS_q (m)^{\bot}\otimes \F_q^q$ et $\F_q^q \otimes RS_q(n)^{\bot}$) sont des produits tensoriels de codes de Reed-Solomon. Ce sont donc des codes fonctionnels sur $\P^1 \times \P^1$. Nous allons tenter de les réaliser sous forme de codes différentiels.

Pour tout $\alpha \in \F_q$, on appelle $L_{d,\alpha}$ la droite d'équation $u-x-\alpha$. Les droites $(L_{d, \alpha})_{\alpha \in \F_q}$ forment un famille de droites \textit{diagonales parallèles} dans l'ouvert $U$ elles sont concourantes en le point $Q$ d'intersection des droites à l'infini $E$ et $F$. Elles sont également deux à deux tangentes en ce point. On définit le diviseur $D_d$ par
$$
D_d:= \sum_{\alpha \in \F_q} L_{d, \alpha}. 
$$
La figure suivante est une tentative de représentation du support de $D_d$ dans le cas où le corps de base est $\F_3$. Si les droites ne ressemblent plus à des droites, nous avons par contre cherché à représenter les points rationnels de $U$ que ces \textit{droites} interpolent. 
 
\begin{center}
\includegraphics[height=8cm, width=8cm]{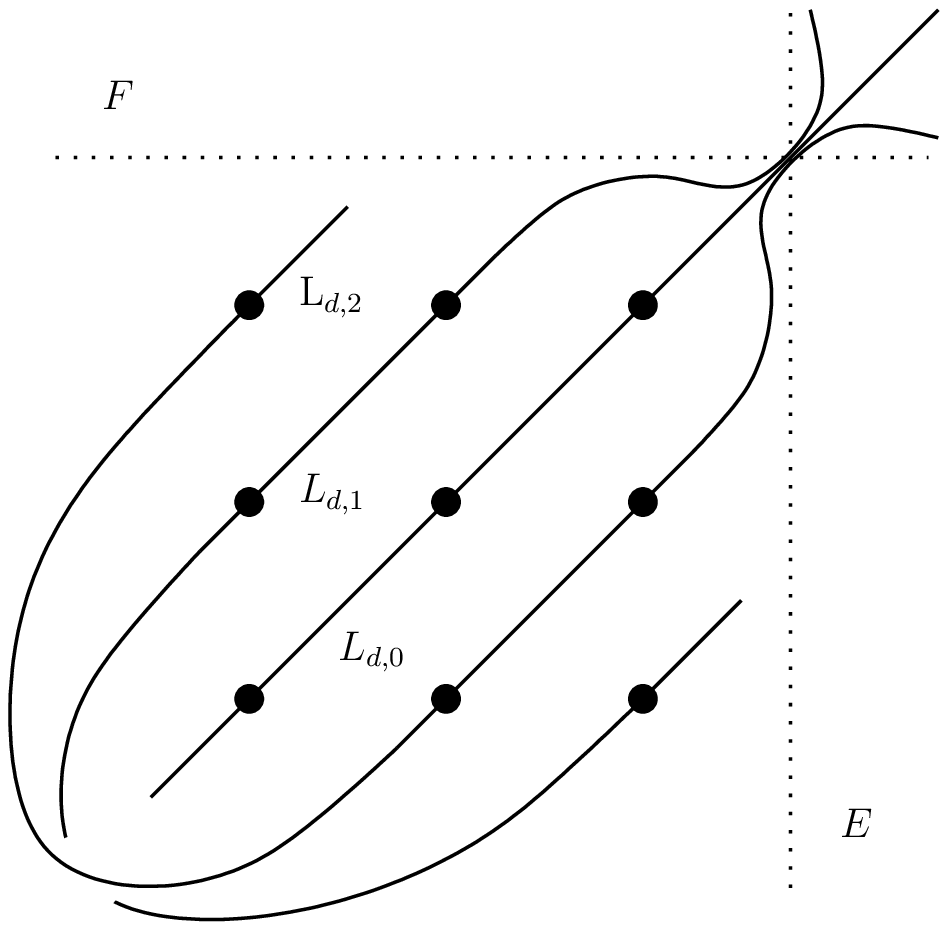}
\end{center}

\begin{rem}
Par le plongement de Segré, $\P^1\times \P^1$ s'identifie à une quadrique hyperbolique de $\P^3$. Les droites $L_{d, \alpha}$ sont les tirés en arrière des éléments d'un pinceau de coniques rationnelles obtenues par des sections de cette quadrique par des plans rationnels contenant tous une même droite tangente à la quadrique en un point.  
\end{rem}

Dans ce qui suit, les diviseurs $D_a$ et $D_b$ sont ceux qui ont été définis sur $\P^1 \times \P^1$ page \pageref{P1prodexdelta}. 

\begin{prop}
Si les entiers $m$ et $n$ sont tous deux compris entre $0$ et $q-1$,
on a alors les trois relations suivantes.
$$
\begin{array}{llcl}
(i) & C_{\Omega}(\Delta, D_a,D_d,G_{m,n}) & = & {\F}_q^q \otimes RS_q(q-2-n).  \\
(ii) & C_{\Omega}(\Delta,D_b,D_d,G_{m,n}) & = & RS_q(q-2-m) \otimes {\F}_q^q .\\
(iii) & {C_L(\Delta,G_{m,n})}^{\bot} & = & C_{\Omega}(\Delta,D_a,D_d ,G_{m,n})+ C_{\Omega}(\Delta,D_b, D_d,G_{m,n}).
\end{array}
$$
Par conséquent, le code $C_L (\Delta, G_{m,n})^{\bot}$ est une somme de deux codes différentiels.
\end{prop}

\begin{proof}
  D'après la relation (\ref{prout}) page \pageref{prout}, si $(i)$ et $(ii)$ sont vérifiées, alors $(iii)$ l'est également. De plus, par symétrie, $(i)$ et $(ii)$ sont équivalentes.
Reste donc à prouver que $(i)$ est vérifiée. Posons
$$
\nu:=\frac{ dx}{\prod_{\alpha \in \F_q}(x-u-\alpha)} \w
\frac{dy}{ \prod_{\beta \in \F_q}(u-\beta)}.
$$
Cette $2$-forme vérifie les trois propriétés du corollaire \ref{omega0}. Calculons le diviseur de $\nu$. Sur $U$, on a 
$$
(\nu_{|U})=-{D_a}_{|U}-{D_d}_{|U}.
$$
De plus, pour tout élément $\alpha$ de $\F_q$, la droite $L_{d,\alpha}$ est linéairement équivalente à $E+F$, donc
$$
D_d \sim q(E+F).
$$

\noindent Par un calcul analogue à celui qui a été effectué dans l'exemple précédent, on montre que
$$
(\nu)= (2q-2)E+(q-2)F -D_a-D_d.
$$ 
De fait,
$$
\begin{array}{rcl}
 C_{\Omega}(\Delta,D_a,D_d ,mE+nF) & = & C_L(\Delta,(2q-2-m)E+(q-2-n)F)\\
 & = & RS_q(2q-2-m)\otimes RS_q(q-2-n).
\end{array}
$$
Pour finir, il suffit de constater que si $m$ est compris entre $0$ et $q-1$, alors $2q-2-m$ est supérieur à $q-1$ et le code $RS_q(2q-2-m)$ est égal à $\F_q^q$.
\end{proof}

Cette stratégie de réalisation de l'orthogonal est celle qui va motiver le chapitre \ref{chapreal}, à savoir, si l'on ne peut réaliser l'orthogonal d'un code fonctionnel avec l'aide d'un seul code différentiel, il est peut-être possible de le réaliser comme somme de codes différentiels.

\begin{rem}
Noter que le contre-exemple ci-dessus s'étend aisément en tout dimension supérieure ou égale à $2$. En général, l'orthogonal d'un code fonctionnel sur un produit de droites projectives ne se réalise pas comme code fonctionnel sur cette variété.
\end{rem}

\section{Heuristique, est-ce un problème de super abondance?}\label{heuris}

Prenons le temps de commenter les phénomènes étudiés dans les exemples précédents. Pour ce faire, commençons par revenir  provisoirement au cas des codes géométriques construits sur des courbes.
Soient $X$ une courbe algébrique projective lisse géométriquement intègre sur $\F_q$ munie d'un diviseur $G$ et $D=P_1+\cdots +P_n$ une somme de points rationnels de $X$.
On sait que dans ce cas, on a systématiquement
$${C_L(D,G)}^{\bot}=C_{\Omega}(D,G).$$
On démontre l'inclusion $\supseteq$ avec la formule des résidus comme l'inclusion réciproque par un argument d'égalité de dimension. 
Voyons comment s'obtient cette égalité de dimension.
Considérons la suite exacte de faisceaux
$$
0\ \rightarrow \mathcal{L}(G-D)\ \rightarrow \mathcal{L}(G)\ 
\rightarrow \ \mathcal{L}(G)/\mathcal{L}(G-D)\ \rightarrow 0.
$$ 
Le terme le plus à droite de cette suite exacte est un faisceau gratte-ciel. Il est donc flasque et son $H^1$ est trivial (cf \cite{H} ch 3.5 théorème 5.1). On en déduit la suite exacte longue en cohomologie,
$$
0\ \rightarrow \ L(G-D)\ \rightarrow L(G)\ \rightarrow
\F_q^n\ \rightarrow H^1(X,\mathcal{L}(G-D))\ \rightarrow
H^1(X,\mathcal{L}(G))\ \rightarrow 0.
$$
Si l'on note $E^{\vee}$ le dual d'un espace vectoriel $E$, alors par dualité de Serre, on a la suite exacte
$$
0\ \rightarrow \ L(G-D)\ \rightarrow L(G)\ \rightarrow
\F_q^n\ \rightarrow {\Gamma (X,\Omega^1(G-D))}^{\vee}\ \rightarrow
{\Gamma (X,\Omega^1(G))}^{\vee}\ \rightarrow 0.
$$
La somme alternée des dimensions permet de conclure,
$$
\underbrace{\dim L(G-D)-\dim L(G)}_{=\dim C_L(D,G)} +n-
\underbrace{(\dim \Gamma (X,\Omega^1(G-D))-\dim \Gamma (X,\Omega^1(G))}_{=\dim C_{\Omega}(D,G)}=0.
$$

Revenons à présent aux surfaces. L'étude des deux exemples triviaux que sont le plan projectif et le produit de deux droites projectives peut encourager le raisonnement heuristique suivant.

\medbreak

\textit{Si l'on arrive a avoir l'égalité de dimension dans le cas o\`u $S=\mathbf{P}^2$ c'est parce que sur cette surface la super abondance\footnote{La dimension du $H^1$} d'un faisceau inversible est nulle (cf \cite{H} théorème III.5.1). Il est donc tentant de penser que l'écart de dimension entre le code différentiel et l'orthogonal du code fonctionnel est lié à la super abondance.}

\medbreak

En réalité ce raisonnement est trop rapide. 
Pour s'en convaincre nous allons essayer de reproduire le raisonnement effectué ci-dessus, dans le cas des surfaces.
Soit $\mathcal{I}_{\Delta}$ le faisceau d'idéaux associé à la sous-variété finie $\supp (\Delta)$.
Concernant la construction fonctionnelle il faut considérer la suite exacte courte de faisceaux suivante,
$$
0\rightarrow\ \mathcal{L}(G)\otimes \mathcal{I}_{\Delta}
\rightarrow\ \mathcal{L}(G)\ \rightarrow\
\mathcal{L}(G)/ \left( \mathcal{L}(G)\otimes \mathcal{I}_{\Delta} \right) \rightarrow 0.
$$

\noindent En remarquant qu'ici encore le dernier faisceau est un faisceau gratte-ciel on en déduit la suite exacte longue
$$
0\rightarrow\ L(G)_{\Delta}\ \rightarrow\ 
L(G)\ \rightarrow\ \F_q^n\
\rightarrow\ H^1(S, \mathcal{L}(G)\otimes \mathcal{I}_{\Delta})\
\rightarrow\ H^1(S, \mathcal{L}(G))\ \rightarrow 0,
$$
où $L(G)_{\Delta}$ désigne l'ensemble des fonctions de $L(G)$ qui s'annulent en tous les points du support de $\Delta$. 
Ici la dualité de Serre ne permet pas de traduire tous les $H^1$ sous formes d'espaces de $2$-formes différentielles. Il faut donc considérer une seconde suite exacte de faisceaux, à savoir:
$$
0\rightarrow\ \Omega^2(G-D)\otimes \mathcal{I}_{\Delta}
\rightarrow\ \Omega^2(G-D)\ \rightarrow\
\Omega^2(G-D)/\left( \Omega^2(G-D)\otimes \mathcal{I}_{\Delta} \right) \rightarrow 0.
$$

\noindent On en déduit la suite exacte longue en cohomologie
$$
\begin{array}{rcl}
0\rightarrow\ \Gamma(S,\Omega^2(G-D))_{\Delta}\  \rightarrow\ 
\Gamma(S,\Omega^2(G-D))\ \rightarrow & \F_q^n & \\
\rightarrow\ H^1(S, \Omega^2(G-D)\otimes \mathcal{I}_{\Delta}) 
& \rightarrow & H^1(S, \Omega^2(G-D))\ \rightarrow\  0,
\end{array}
$$

\noindent où $\Gamma(S,\Omega^2(G-D))_{\Delta}$ décrit l'ensemble des éléments de $\Gamma(S,\Omega^2(G-D))$ qui s'annulent en tous les éléments du support de $\Delta$.
Les faisceaux $\mathcal{L}(G)$ et $\Omega^2(G-D)$ sont inversibles. Donc, si la surface est $\mathbf{P}^2$, leurs $H^1$ sont nuls et les suites exactes longues donnent les égalités
$$
\begin{array}{rcl}
\dim H^1(\mathbf{P}^2, \mathcal{L}(G)\otimes \mathcal{I}_{\Delta})
 & = & \textrm{codim}\ C_L(\Delta,G) \\
\dim H^1(\mathbf{P}^2, \Omega^2(G-D)\otimes \mathcal{I}_{\Delta})
 & = & \textrm{codim}\ C_{\Omega}(D,G).
\end{array}
$$
Il faut ensuite réussir à prouver que la somme des dimensions de ces deux $H^1$ est égale à $n$.
Dans tous les cas, le fait que les super abondances des faisceaux inversibles soient nulles ne suffit pas pour démontrer l'égalité de dimension espérée. 

\medbreak

\textit{Ce qui semble réellement faire défaut à cette construction différentielle est moins la super abondance que l'asymétrie des constructions. Plus précisément, le fait qu'en dimension supérieure ou égale à $2$, les points et les diviseurs sont des objets de dimension différente.
Du fait de cette asymétrie, on doit introduire une paire de diviseur $\Delta$-convenable pour construire des codes différentiels, cette dernière étant
complètement absente dans la construction fonctionnelle.}

\section*{Conclusion}
Nous avons étendu la construction différentielle de codes correcteurs aux surfaces.
Nous avons également montré que, tout comme dans la cas des courbes, les codes fonctionnels et différentiels sur les surfaces appartiennent à la même classe de codes.
en d'autres termes, tout code fonctionnel sur une surface se réalise comme code différentiel sur cette même surface et réciproquement. 

La différence majeure avec la théorie des courbes est que, à $S$ et $\Delta$ fixés \textbf{l'orthogonal d'un code fonctionnel n'est pas fonctionnel (donc différentiel) en général}. Ces codes appartiennent à une ``classe différente'' de codes construits à partir de $S$.
Voloch et Zarzar avaient d'ailleurs déjà  constaté ce phénomène dans \cite{agctvoloch}. Dans cet article, les auteurs remarquent en effet que les codes sur les surfaces qu'ils étudient sont LDPC\footnote{\textit{Low Density Parity-Check}, c'est-à-dire admettant une matrice de parité creuse (voir chapitre \ref{chapldpc}).}.
De ce fait, ces codes sont très différents de leur orthogonal, ce qui n'est pas le cas des codes géométriques construits à partir de courbes algébriques.

Aussi, l'étude des codes différentiels sur les surfaces et des exemples que nous avons traités offrent des perspectives intéressantes. À ce titre, nous conclurons ce chapitre par deux questions.

\begin{ques}\label{QQk}
Peut-on estimer les paramètres des codes qui sont l'orthogonal de codes fonctionnels?  
\end{ques}

\begin{ques}\label{Qsom} Si l'orthogonal d'un code fonctionnel ne peut se réaliser comme un code différentiel associé à une paire de diviseurs $\Delta$-convenables, peut-on le réaliser comme somme de tels codes?
\end{ques}

La question \ref{QQk} donne lieu aux travaux présentés dans le chapitre \ref{chaporth}.
Concernant la question \ref{Qsom}, une réponse partielle sera donnée dans le chapitre \ref{chapreal}.


\newpage
\thispagestyle{empty}

\chapter{Théorème de réalisation}\label{chapreal}

\begin{flushright}
\begin{tabular}{p{6cm}}
\begin{flushright}
{\small \textit{``Entre les désirs et leurs réalisations
s'écoule toute une vie humaine.''}}
\medbreak
{\small Schopenhauer}
\end{flushright}
\end{tabular}
\end{flushright}

Dans ce chapitre, nous allons nous intéresser la question \ref{Qsom} posée à la fin du chapitre \ref{chapdiff}. À savoir: \textit{l'orthogonal $C_L(\Delta,G)^{\bot}$ d'un code fonctionnel  sur une surface $S$ se réalise-t-il comme somme de codes différentiels?} Une réponse positive à cette question sera donnée sous certaines conditions sur la surface $S$ et le diviseur $G$. Ces conditions sont décrites dans la section \ref{contextchap3} ci-dessous.

\section{Contexte}\label{contextchap3}
Dans ce chapitre, $S$ désigne une surface projective lisse géométriquement intègre au-dessus de $\F_q$.
On se donne également un diviseur $\F_q$-rationnel  $G$ sur $S$ et une famille $P_1, \ldots, P_n$ de points rationnels de $S$. On appelle $\Delta$, le $0$-cycle
$$
\Delta:=P_1+\cdots + P_n.
$$

\begin{nota}\label{L_X}
Soit $H$ un hyperplan de $\P^r$, pour toute sous-variété $X$ de $\P^r$ non contenue dans $H$, on note $L_X$ la classe d'équivalence linéaire du diviseur $\varphi^* H$ sur $X$, où $\varphi$ désigne l'injection canonique $\varphi: X \hookrightarrow \P^r$. De même, la classe canonique sur $X$ sera notée $K_X$.
\end{nota}

À partir de la section \ref{secreal}, on supposera (hypothèse \ref{secreal}) que $S$ est intersection complète dans un espace projectif $\P^r$ et que $G$ est linéairement équivalent à $mL_S$ pour un certain entier naturel $m$.

\section{Sous-$\Delta$-convenance}\label{secsous}

Dans cette section, nous allons définir une notion très proche de celle de $\Delta$-convenance et qui continue à vérifier le résultat du théorème d'orthogonalité (théorème \ref{orthocode}).

Commençons par justifier ce besoin d'introduire une nouvelle notion.
La question \ref{Qsom} posée à la fin du chapitre \ref{chapdiff} était en partie motivée par l'étude de la surface $\P^1 \times \P^1$.
En effet, on a vu en section \ref{P1P1} que l'orthogonal d'un code fonctionnel sur cette surface peut se réaliser comme somme de deux codes différentiels.
Tout cela encourage à essayer de construire l'orthogonal d'un code fonctionnel en ``plusieurs morceaux''.
Pour ce faire, on peut par exemple chercher des mots de code ou des sous-codes de $C_L(\Delta, G)^{\bot}$ dont le support est strictement contenu dans $\{1, \ldots, n\}$.
De plus l'exemple des quadriques elliptiques du chapitre précédent (section \ref{qq}) montre que la construction d'un diviseur $\Delta$-convenable devient vite ardue, lorsque la surface $S$ est plus compliquée que le plan projectif ou le produit de deux droites projectives.
Ces deux arguments motivent la définition de paires de diviseurs sous-$\Delta$-convenables qui suit.

\begin{defn}[Diviseurs sous-$\Delta$-convenables]\label{sousDconv}
Une paire $(D_a, D_b)$ est dite sous-$\Delta$-convenable si elle est $\Lambda$-convenable pour un $0$-cycle $\Lambda$ vérifiant
$$
0 \leq \Lambda \leq \Delta.
$$ 
\end{defn}

\begin{rem}\label{remsousD}
La sous-$\Delta$-convenance peut également s'énoncer de la façon suivante.
Soient $D_a$ et $D_b$ deux diviseurs dont l'intersection ensembliste des supports est finie et $D$ le diviseur $D:=D_a+D_b$. La paire $(D_a,D_b)$ est sous-$\Delta$-convenable si et seulement si elle vérifie les propriétés suivantes.
\begin{enumerate}
\item[$(i)$] Pour tout point géométrique $P\in \overline{S}$, l'application
$
\res^2_{D_a,P}: \overline{\Omega^2(-D)}_P\rightarrow \F_q
$
est $\mathcal{O}_{\overline{S},P}$-linéaire.
\item[$(ii)$] L'application ci-dessus est nulle pour tout point géométrique $P$ non contenu dans le support de $\Delta$.
\end{enumerate}
\end{rem}

Rappelons que les diviseurs $\Delta$-convenables ont été introduits pour obtenir une relation d'orthogonalité entre codes fonctionnels et codes différentiels (voir théorème \ref{orthocode}).
Le lemme qui suit et dont la démonstration est immédiate montre que les paires de diviseurs sous-$\Delta$-convenables sont en ce sens presque aussi intéressantes que les paires $\Delta$-convenables.

\begin{lem}\label{sousorth}
Soit $S$ une surface lisse géométriquement intègre au-dessus de $\F_q$ et munie d'un diviseur $G$ et d'un $0$-cycle $\Delta$ qui est la somme formelle de points rationnels de $S$.
Soit enfin $(D_a,D_b)$ une paire sous-$\Delta$-convenable de diviseurs, alors
$$
C_{\Omega}(\Delta,D_a,D_b,G) \subseteq C_L(\Delta,G)^{\bot}.
$$  
\end{lem}

\section{Sur les notions de réalisation}

La question de la réalisation de l'orthogonal d'un code fonctionnel en utilisant des $2$-formes rationnelles peut se poser de deux façons différentes. Il y a d'abord la question \ref{Qsom} posée à la fin du chapitre \ref{chapdiff} que nous rappelons ici.

\medbreak

\noindent \textbf{Question \ref{Qsom}.} \textit{Si l'orthogonal d'un code fonctionnel ne peut se réaliser comme un code différentiel associé à une paire de diviseurs (sous-)$\Delta$-convenable, peut-on le réaliser comme somme de tels codes?}

\medbreak

\noindent On peut également se poser une question sensiblement différente, à savoir la question \ref{Qsom}bis qui suit. Le théorème de réalisation énoncé en section \ref{secreal} y répondra sous certaines conditions.

\medbreak

\noindent \textbf{Question \ref{Qsom}bis.}
\textit{  Étant donné un mot de code $c$ appartenant à $C_{L,S}(\Delta,G)^{\bot}$, existe-t-il une paire de diviseurs (sous-) $\Delta$-convenable $(D_a,D_b)$ et une $2$-forme $\omega$ appartenant à  $\Gamma (S,\Omega^2(G-D_a-D_b))$ et telle que 
$$
c=\res^2_{D_a,\Delta}(\omega)\textrm{?}
$$ }

Le fait qu'un code est un espace vectoriel de dimension finie permet de montrer aisément qu'une réponse positive à la question \ref{Qsom}bis entraîne une réponse positive à la question \ref{Qsom}.
La réciproque de cette dernière assertion est également vraie, c'est ce que montre la proposition qui suit. Les problèmes posés par les questions \ref{Qsom} et \ref{Qsom}bis sont donc équivalents.

\begin{prop}\label{realsom} Soient $c_D$ et $c_E$ deux mots du code $C_L (\Delta,G)^{\bot}$. Supposons qu'il existe deux paires de diviseurs sous-$\Delta$-convenables $(D_a,D_b)$ et $(E_a,E_b)$ et deux $2$-formes rationnelles $\omega_D$ et $\omega_E$ appartenant respectivement aux espaces $\Gamma (S, \Omega^2 (G-D_a-D_b))$ et $\Gamma (S, \Omega^2 (G-E_a-E_b))$ et telles que
$$
c_D=\res^2_{D_a,\Delta}(\omega_D) \quad \textrm{et}\quad c_E=\res^2_{E_a,\Delta}(\omega_E).
$$
Alors, il existe une paire sous-$\Delta$-convenable de diviseurs $(F_a,F_b)$ et une $2$-forme rationnelle $\omega_F$ appartenant à $\Gamma (S, \Omega^2(G-F_a-F_b))$ telle que 
$$
c_D+c_E=\res^2_{F_a,\Delta}(\omega_F).
$$
\end{prop} 

\noindent La preuve de proposition \ref{realsom} est une conséquence des lemmes \ref{suppdis} et \ref{movingDel} qui suivent. 

 \begin{lem}\label{suppdis}
   Soient $(D_a, D_b)$ et $(E_a,E_b)$ deux paires sous-$\Delta$-convenables de diviseurs sur $S$ telles que les supports des diviseurs $D:=D_a+D_b$, $E:=E_a+E_b$ et $G$  n'ont pas de composante irréductible commune.
Soient également $\omega_D$ et $\omega_E$ deux $2$-formes rationnelles sur $S$ appartenant respectivement à $\Gamma(S, \Omega^2 (G-D))$ et $\Gamma (S, \Omega^2 (G-E))$. Alors, il existe une paire sous-$\Delta$-convenable de diviseurs $(F_a,F_b)$ telle que la $2$-forme $\omega_D+\omega_E$ appartienne à $\Gamma (S, \Omega^2 (G-F))$ où $F$ désigne le diviseur $F:=F_a+F_b$. De plus,
$$
\res^2_{F_a,\Delta}(\omega_D+\omega_E)=\res^2_{D_a, \Delta}(\omega_D)+\res^2_{E_a, \Delta}(\omega_E).
$$
\end{lem}

\begin{lem}\label{movingDel}
  Soit $(D_a, D_b)$ une paire de diviseurs sous-$\Delta$-convenable et $\omega$ un élément de $\Gamma (S, \Omega^2 (G-D))$ où $D$ désigne le diviseur $D:=D_a+D_b$.
Soient également $C_1, \ldots, C_s$ une famille de courbes irréductibles sur $S$ deux à deux distinctes.
Alors, il existe une paire sous-$\Delta$-convenable $(D_a',D_b')$ vérifiant les propriétés suivantes.
\begin{enumerate}
\item Les diviseurs $D_a$ et $D_a'$ (resp. $D_b$ et $D_b'$) sont linéairement équivalents.
\item Le support de $D':=D_a'+D_b'$ ne contient aucune des courbes $C_1, \ldots, C_s$.
\item\label{encule} Pour toute $2$-forme $\omega$ appartenant à $\Gamma (S, \Omega^2 (G-D))$, il existe une $2$-forme $\omega'$ appartenant à $\Gamma (S, \Omega^2 (G-D'))$ telle que 
$$
\res^2_{D_a,\Delta}(\omega)=\res^2_{D_a',\Delta}(\omega').
$$
\end{enumerate}
\end{lem}

\begin{proof}[\textsc{Preuve de la proposition \ref{realsom}}]
  Si les supports des diviseurs $D:=D_a+D_b$ et $E:=E_a+E_b$ sont sans composante commune, on applique le lemme \ref{suppdis}. Sinon on se ramène à cette situation grâce au lemme \ref{movingDel}.
\end{proof}

\begin{proof}[\textsc{Preuve du lemme \ref{suppdis}}]

\noindent \textbf{Étape 1. Construction de $\mathbf{(F_a,F_b)}$.}
Les diviseurs $D_a^+$, $E_a^+$, $D_b^+$, $E_b^+$ sont respectivement de la forme
$$
\begin{array}{ccccccc}
D_a^+ & := & m_1 V_1 + \cdots +m_k V_k & \quad & E_a^+ & := & n_1 W_1 + \cdots
+n_l W_l \\
D_b^+ & := & r_1 X_1 + \cdots +r_p X_p & \quad & E_b^+ & := & s_1 Y_1 + \cdots
+s_q Y_q, \\
\end{array}
$$

\noindent où les $V_i$, $W_i$, $X_i$, $Y_i$ sont des courbes $\F_q$-irréductibles. Par hypothèse, ces courbes sont deux à deux disjointes.
Nous allons construire une paire de diviseurs effectifs $(F_a^+,F_b^+)$. Le diviseur $F_a^+$ fera apparaître tous les $V_i$ (resp. $W_j$), pôles de $\omega_D$ (resp. $\omega_E$) avec pour coefficient l'ordre de ce pôle. Le diviseur $F_b^+$ est construit exactement de la même manière en remplaçant les $V_i$ par des $X_i$ et les $W_j$ par des $Y_j$. 
C'est-à-dire que l'on pose
\begin{equation}\label{Fa+}
F_a^+  :=  \sum_{i=1}^k \max (\ \!0\ \!, -\val_{V_i}(\omega_D))\ \!V_i  +   \sum_{j=1}^l \max (\ \!0\ \!, -\val_{W_j}(\omega_E))\ \!W_j
\end{equation}

\noindent et
\begin{equation}\label{Fb+}
F_b^+  :=  \sum_{i=1}^p \max (\ \!0\ \!, -\val_{X_i}(\omega_D))\ \!X_i   +   \sum_{j=1}^q \max (\ \!0\ \!, -\val_{Y_j}(\omega_E))\ \!Y_j.
\end{equation}

\noindent Soit $\omega_F$, la $2$-forme définie par $\omega_F:=\omega_D+\omega_E$.
On rappelle que les composantes irréductibles des supports des diviseurs $D_a,D_b,E_a$ et $E_b$ sont par hypothèse deux à deux disjointes.
Par conséquent, il existe un diviseur \textbf{effectif} $R$ dont le support n'a aucune composante irréductible commune avec les supports de $D_a^+,D_b^+, E_a^+$ et $E_b^+$ et tel que
$$
(\omega_F)=G+R-F_a^+-F_b^+.
$$

\noindent On pose enfin
$$
F_b^- := R, \quad\textrm{et}\quad F_a^-:=0
$$

\noindent et on a donc
$$
F_a:=F_a^+ \quad \textrm{et} \quad F_b:=F_b^+-F_b^-.
$$

Le diviseur $F:=F_a+F_b$ est donc construit de telle sorte que l'on ait exactement
$$
(\omega_F)=G-F.
$$
Le fait que $\omega_F$ est un élément de $\Gamma (S, \Omega^2 (G-F))$ est immédiat.
Il reste à montrer que la paire $(F_a,F_b)$ est sous-$\Delta$-convenable.

\medbreak

\noindent \textbf{Étape 2. Sous-$\mathbf{\Delta}$-convenance de $\mathbf{(F_a,F_b)}$.}
Soit $P$ un point de $\overline{S}$, il existe un germe de fonction $f_P$ appartenant à $\overline{\L (G)}_P$ tel que le diviseur de la $2$-forme $f_P\omega_F$ au voisinage de $P$ soit égal à $-F$.
Alors, le germe de $2$-forme $f_P\omega_F$ au voisinage de $P$ engendre la fibre $\overline{\Omega^2 (-F)}_P$ comme $\mathcal{O}_{\overline{S},P}$-module. 
En effet, la $2$-forme $f_P \omega_F$ est construite de telle sorte que pour tout germe de courbe $C$ au voisinage de $P$ on ait
$$
\val_C (f_P \omega_F)=\min_{\mu_P \in \overline{\Omega^2 (-F)}_P} \val_C(\mu_P).
$$

\noindent Un germe de $2$-forme $\mu_P \in \overline{\Omega^2 (-F)}_P$ s'obtient donc par multiplication de $f_P \omega_F$ par une fonction régulière au voisinage de $P$.

Montrons que l'application $\res^2_{F_a,P}$ restreinte à $\overline{\Omega^2(-F)}_P$ est $\mathcal{O}_{\overline{S},P}$-linéaire. Soit $\varphi \in \mathcal{O}_{\overline{S},P}$. On a
\begin{equation}\label{IDIE}
\res^2_{F_a,P}(\varphi f_P\ \! \omega_F)=\underbrace{\res^2_{F_a,P}(\varphi f_P\ \! \omega_D)}_{I_D}+\underbrace{\res^2_{F_a,P}(\varphi f_P\ \! \omega_E)}_{I_E}.
\end{equation}

Nous allons montrer que $I_D=\varphi(P)\res^2_{D_a,P}(f_P \omega_D)$. Commençons par faire deux remarques.
\begin{enumerate}
\item D'après la construction de $F_a^+$ en (\ref{Fa+}), certaines des courbes $V_i$ peuvent ne pas apparaître dans l'expression de ce diviseur.
C'est ce qui arrive pour une courbe $V_i$ donnée si la valuation de $\omega_D$ le long de $V_i$ est positive.
Dans ce cas, le $2$-résidu de $\omega_D$ en $P$ le long de $V_i$ est nul. 
De plus, on rappelle que $\varphi$ est régulière au voisinage de $P$ et que, par hypothèse, le support de $G$ n'a pas de composante commune avec ceux de $D$ et $E$.
Par conséquent, si $\omega_D$ n'a pas de pôle le long de $V_i$, alors le $2$-résidu en $P$ le long de $V_i$ de $\varphi f_P \ \!\omega_D$ est nul.
\item L'hypothèse ``les supports des diviseurs $D_a,D_b,E_a,E_b$ n'ont pas de composante irréductible commune'' implique que pour tout $i$, les $2$-formes $\omega_D$ et $\varphi f_P \ \!\omega_D$ n'ont pas de pôle le long de $W_i$, donc ont un $2$-résidu nul en $P$ le long de cette courbe.
\end{enumerate}

\noindent On déduit de ces deux remarques que
\begin{equation}\label{ID}
I_D = \sum_{i=1}^k \res^2_{V_i,P}(\varphi f_P \ \!\omega_D) .
\end{equation}

\noindent Enfin, on rappelle que la définition de $2$-résidu en un point le long d'un diviseur ne dépend que du support de ce dernier (voir définition \ref{resdiv} et la mise en garde qui suit).
Par conséquent,
\begin{equation}\label{ID2}
I_D= \res^2_{D_a,P}(\varphi f_P\ \! \omega_D)=\varphi (P) \res^2_{D_a,P}(f_P\ \! \omega_D),
\end{equation}

\noindent la seconde égalité étant une conséquence de la sous-$\Delta$-convenance de $(D_a,D_b)$ et du fait que $f_P \ \! \omega_D$ appartient à $\overline{\Omega^2(-D)}_P$ comme $\mathcal{O}_{\overline{S},P}$-module.
Le cas de la quantité $I_E$ de l'expression (\ref{IDIE}) se traite de façon rigoureusement identique. On obtient
\begin{equation}\label{IE2}
I_E=\varphi (P) \res^2_{E_a,P}(f_P\ \! \omega_E).
\end{equation}

\noindent En combinant les relations (\ref{IDIE}), (\ref{ID2}) et (\ref{IE2}), on aboutit à
$$
\res^2_{F_a,P}(\varphi f_P \ \! \omega_F)=\varphi (P)\res^2_{F_a,P} (f_P \ \! \omega_F),
$$
ce qui permet de conclure quant à l'$\mathcal{O}_{\overline{S},P}$-linéarité de l'application $\res^2_{F_a,P}$ restreinte à $\overline{\Omega^2(-F)}_P$.

Enfin, comme les paires $(D_a,D_b)$ et $(E_a,E_b)$ sont sous-$\Delta$-convenables, pour tout point $P$ de $\overline{S}$ n'appartenant pas au support de $\Delta$, les applications $\res^2_{D_a,P}$ et $\res^2_{E_a,P}$ sont identiquement nulles respectivement sur $\overline{\Omega^2 (-D)}_P$ et $\overline{\Omega^2 (-E)}_P$.
D'après (\ref{IDIE}), (\ref{ID2}) et (\ref{IE2}), on en déduit que l'application $\res^2_{F_a,P}$ est identiquement nulle sur $\overline{\Omega^2 (-F)}_P$. D'après la remarque \ref{remsousD}, la  paire $(F_a,F_b)$ est sous-$\Delta$-convenable. 
\end{proof}

Pour finir, il nous reste à démontrer le lemme \ref{movingDel}.

\begin{proof}[\textsc{Preuve du lemme \ref{movingDel}}]

\noindent \textbf{Étape 0. Mise en place.}

Quitte à réorganiser l'ordre des courbes $C_1, \ldots, C_s$, on peut supposer que $C_1, \ldots , C_l$ sont contenues dans le support de $D_a$, que $C_{l+1},\ldots, C_m$ sont dans le support\footnote{On rappelle que par définition de la sous-$\Delta$-convenance, les supports de $D_a$ et $D_b$ n'ont pas de composante irréductible commune. La courbe $C_i$ ne peut donc pas être contenue dans l'intersection ensembliste des supports de $D_a$ et $D_b$.} de $D_b$  et $C_{m+1}, \ldots, C_s$ ne sont contenues dans aucun des deux supports.
Nous allons montrer comment \textit{bouger} $D_a$ afin d'éviter ces courbes. On pourra alors conclure d'après la remarque \ref{DasymDb} en appliquant un raisonnement identique à $D_b$.

\medbreak

\noindent \textbf{Étape 1. Déplacement de $\mathbf{D_a}$.}

Pour \textit{déplacer} le support de $D_a$, nous allons utiliser la proposition \ref{movingperso} énoncée en annexe \ref{annexeapprox} qui est un analogue du \textit{moving lemma}.
Commençons par établir une liste de points à éviter.

Soit $\mathcal{C}$, l'ensemble des courbes formé de la réunion des composantes irréductibles des supports de $D_a$, $D_b$ et $G$ et des courbes $C_1, \ldots, C_s$.
L'ensemble $\mathcal{P}$ est un ensemble de points fermés de $S$ formé de tous les points d'intersection (ensembliste) de deux éléments de $\mathcal{C}$. Si l'un des éléments $C$ de $\mathcal{C}$ n'en intersecte aucun autre, on choisit arbitrairement un point de $C$ que l'on ajoute dans $\mathcal{P}$, afin que ce dernier contienne au moins un point de chaque courbe appartenant à $\mathcal{C}$.

Soit $i$ un entier appartenant à $\{1, \ldots, l\}$, 
on partitionne $\mathcal{P}$ en deux ensembles $\mathcal{P}_1^i$ et $\mathcal{P}^i_2$.
L'ensemble $\mathcal{P}_1^i$ désigne l'ensemble des points $P$ qui appartiennent à $C_i$ et $\mathcal{P}^i_2$ désigne son complémentaire dans $\mathcal{P}$.
D'après la proposition \ref{movingperso}, il existe une fonction $f_i$, vérifiant les propriétés suivantes.
\begin{enumerate}
\item[$(i)$] La fonction $f_i$ est une équation locale de $C_i$ au voisinage de tout point $P \in \mathcal{P}^i_1$.
\item[$(ii)$] Le support du diviseur de $f_i$ évite tout point $P \in \mathcal{P}^i_2$.
\end{enumerate}

\noindent On pose
$$
m_i:=\val_{C_i}(D_a),
$$

\noindent et on définit la fonction rationnelle $\varphi$ sur $S$ par
$$
\varphi:=f_1^{m_1} \ldots f_l^{m_l}.
$$

\noindent De la même manière, on pose
$$
\widetilde{D}:=D-(m_1C_1 + \cdots + m_l C_l)+C_{m+1}+ \cdots +C_s,
$$

\noindent et on partitionne $\mathcal{P}$ en $\mathcal{P}_1$ et $\mathcal{P}_2$, formés respectivement des points de $\mathcal{P}$ contenus et non contenus dans le support de $\widetilde{D}$. D'après la proposition \ref{movingperso}, il existe une fonction rationnelle $g$ qui est une équation locale de $\widetilde{D}$ au voisinage de tout élément de $\mathcal{P}_1$ et dont le support du diviseur évite tous les éléments de $\mathcal{P}_2$.
Posons alors, 
$$
h:=\frac{\varphi}{g+\varphi}, \quad  \textrm{et} \quad D_a':=D_a-(h).
$$

\noindent Montrons que le diviseur de la fonction $h$ est de la forme
$$
(h)=m_1C_1+ \cdots + m_l C_l+R,
$$

\noindent et que le support de $R$ ne contient aucun élément de $\mathcal{C}$.
Le diviseur $(h)$ est la différence des diviseurs $(\varphi)$ et $(\varphi+g)$. Par construction, le diviseur $(\varphi)$ est de la forme
$$
(\varphi)=m_1C_1+ \cdots + m_l C_l+R_1,
$$
où le support de $R_1$ évite tout élément de $\mathcal{P}$, donc ne contient aucun élément de $\mathcal{C}$.
Quant à la fonction $\varphi+g$, elle est par construction régulière au voisinage de tout élément de $\mathcal{P}$.
Elle l'est donc sur un ouvert admettant une intersection non vide avec tout élément de $\mathcal{C}$ et n'admet donc aucune de ces courbes comme pôle.
De plus, pour toute courbe $C$ appartenant à $\mathcal{C}$, l'une des fonctions $\varphi$ ou $g$ s'annule le long de $C$ et l'autre ne s'annule pas. De fait, $C$ ne peut être un zéro de $\varphi+g$.
Par conséquent, le support du diviseur $D_a'$ ne contient aucune des courbes $C_1, \ldots, C_s$ ni aucune composante des supports de $D_b$ et $G$.

\medbreak

\noindent \textbf{Étape 2. Sous-$\mathbf{Delta}$-convenance de $\mathbf{(D_a',D_b)}$.}

Soit $\omega$ une section globale de $\Omega^2 (G-D)$, on vérifie aisément que $h\omega$ est une section globale de $\Omega^2(G-D_a'-D_b)$. Nous allons montrer que pour toute $2$-forme $\omega$ appartenant à $\Gamma(S, \Omega^2 (G-D))$, on a 
\begin{equation}\label{graal}
\res^2_{D_a,\Delta}(\omega)=\res^2_{D_a', \Delta}(h \omega),
\end{equation}
ce qui nous permettra de démontrer à la fois la propriété \ref{encule} de l'énoncé et le fait que la paire $(D_a',D_b)$ est sous-$\Delta$-convenable.
Commençons par noter que, d'après la remarque \ref{sym} du chapitre \ref{chapdiff}, la relation (\ref{graal}) et équivalente à 
\begin{equation}\label{graal2}
\res^2_{D_b,\Delta}(\omega)=\res^2_{D_b, \Delta}(h \omega).
\end{equation}

Soient donc $\overline{C}$ une composante géométrique irréductible du support de $D_b$ et $u$ un élément de $\mathcal{O}_{\overline{S}, \overline{C}}$ dont la restriction à $\overline{C}$ est un élément séparant de $\overline{\F}_q(\overline{C})/\overline{\F}_q$. Nous allons montrer que 
\begin{equation}\label{graal3}
(u)\res^1_{\overline{C}}(\omega)=(u)\res^1_{\overline{C}}(h\omega).
\end{equation}
D'après la proposition \ref{1resprop2} du chapitre \ref{chapres}, si l'on montre que la relation (\ref{graal3}) ci-dessus est vérifiée par toute composante géométrique du support de $D_b$, alors la relation (\ref{graal2}) sera vérifiée.

Soit $v$ une uniformisante de l'anneau $\mathcal{O}_{\overline{S}, \overline{C}}$.
On rappelle que le complété $\mathfrak{m}_{\overline{S}, \overline{C}}$-adique de $\overline{\F}_q (\overline{S})$ s'identifie au corps $\mathcal{K}_u ((v))$, où $\mathcal{K}_u$ est une copie de $\overline{\F}_q (\overline{C})$ contenue dans
$\widehat{\mathcal{O}}_{\overline{S}, \overline{C}}$ (voir chapitre \ref{chapres} section \ref{seclolo2}).
Posons
$$
m:=\val_{\overline{C}}(D_b).
$$

\noindent Comme $\omega$ est un section globale de $\Omega^2 (G-D)$, on sait que sa valuation le long de $\overline{C}$ est supérieure à $-m$. La construction des fonctions $\varphi$ et $g$ nous assure que ces dernières sont de valuations respectives $0$ et $m$ le long de $\overline{C}$.
On en déduit donc
$$
h\omega=\frac{\varphi}{\varphi+g}\omega=\frac{1}{1+g\varphi^{-1}}\omega
=\left(1-g\varphi^{-1}+\cdots\right)\omega.
$$
La $2$-forme $\omega$ étant de valuation supérieure à $-m$ le long de $\overline{C}$ et la fonction $g\varphi^{-1}$ de valuation $m$, on en déduit que les termes de la forme $(-1)^n (g \varphi^{-1})^n \omega$ de la série ci-dessus sont de valuation positive le long de $\overline{C}$. Leur contribution dans le calcul du $(u)$-$1$-résidu de $\omega$ le long de $\overline{C}$ est donc nulle. La relation (\ref{graal3}) est bien vérifiée. 

On conclut la preuve en appliquant un raisonnement identique à $D_b$.
\end{proof}

\section{Construction de l'orthogonal d'un code fonctionnel}\label{secreal}

Le but de cette section est de démontrer le théorème suivant.
On se place dans le contexte donné en section \ref{contextchap3} et on suppose de plus vérifiée l'hypothèse suivante.

\medbreak

\noindent \textbf{Hypothèse \ref{secreal}.}
\textit{La surface $S$ est plongée dans un espace projectif $\P^r_{\F_q}$ dans lequel elle est \textbf{intersection complète}. De plus; le diviseur $G$ est linéairement équivalent à l'intersection de $S$ avec une hypersurface de $\P^r$. En d'autres termes et en utilisant la notation \ref{L_X}, il existe un entier naturel $m$ tel que $G\sim mL_S$.}

\begin{thm}[Théorème de réalisation]\label{thmreal}
Sous l'hypothèse \ref{secreal}, soit $c$ un mot du code $C_{L,S} (\Delta, G)^{\bot}$. Alors, il existe une paire de diviseurs $(D_a, D_b)$ et une $2$-forme $\omega$ appartenant à l'espace des sections globales $\Gamma(S,\Omega^2(G-D_a-D_b))$, tels que
$$
c=\res^2_{D_a, \Delta}(\omega).
$$
De plus, on peut choisir le couple $(D_a,D_b)$ de telle sorte que
\begin{enumerate}
\item la paire $(D_a,D_b)$ vérifie le critère de la proposition \ref{crit};
\item $D_a$ soit égal à une courbe lisse irréductible plongée dans $S$ et $D_a \sim n_a L_S$ pour un certain entier strictement positif $n_a$;
\item $D_b \sim n_b L_S$ pour un certain entier $n_b$.
\end{enumerate}
\end{thm}

\noindent Avant de démontrer ce théorème, énonçons un corollaire immédiat de ce dernier.

\begin{cor}\label{correal}
Sous l'hypothèse \ref{secreal},
il existe une famille finie $(D_a^{(1)}, D_b^{(1)}), \ldots ,$ $(D_a^{(s)}, D_b^{(s)})$ de paires de diviseurs sous-$\Delta$-convenables telles que
$$
C_{L,S}(\Delta, G )^{\bot} = \sum_{i=1}^s C_{\Omega,S} (\Delta, D_a^{(i)}, D_b^{(i)}, G).
$$
\end{cor}

\begin{proof}[\textsc{Preuve du corollaire \ref{correal}}]
L'inclusion vers la gauche est une conséquence immédiate du lemme \ref{sousorth}.
Pour ce qui est de l'inclusion réciproque, le théorème \ref{thmreal} implique que $C_{L,S}(\Delta, G)^{\bot}$ est égal à la somme de tous les codes de la forme $C_{\Omega,S}(\Delta, D_a, D_b, G)$ tels que $(D_a,D_b)$ est sous-$\Delta$-convenable.
Comme un code est un espace de dimension finie, on peut extraire de cette somme une somme finie.
\end{proof}

\noindent Le lemme qui suit est la première étape de la preuve du théorème \ref{thmreal}.

\begin{lem}\label{surj}
Sous l'hypothèse \ref{secreal}, soit $C$ une courbe lisse absolument irréductible plongée dans $S$ obtenue par l'intersection de $S$ avec une hypersurface de $\P^r$. On suppose également que $C$ n'est pas contenue dans le support\footnote{Notons que si $C$ est contenue dans le support de $G$, on peut remplacer ce dernier par un autre élément du système linéaire $|G|$, cette condition n'est donc pas vraiment problématique.
D'une façon générale, on peut éviter ce type de restriction en adoptant un langage plus \textit{faisceautique}.
En effet, le tiré en arrière de $G$ sur $C$ n'a pas de sens quand $C$ est contenue dans le support de $G$, le tiré en arrière de $\L (G)$ lui, est toujours bien défini (voir \cite{H} II.5).
Le défaut de ce point de vue est que dans ce cas, les sections de $i^*\L (G)$ ne peuvent plus être vues comme des restrictions à $C$ de fonctions sur $S$.
Nous avons donc préféré conserver une approche plus \textit{fonctionnelle}.} de $G$.
Soit $G^*$ le tiré en arrière de $G$ par l'inclusion canonique $C \hookrightarrow S$. Alors l'application de restriction à $C$
$$
r: \Gamma (S, \L (G)) \rightarrow \Gamma (C, \L (G^*))
$$
est surjective.
\end{lem}

\begin{proof}
  La courbe $C$ est lisse donc normale. C'est de plus une intersection complète dans $\P^r$. Ainsi, d'après \cite{H} II.8 ex 4, la courbe $C$ est projectivement normale.
Cela signifie par définition, que l'algèbre graduée des coordonnées homogènes de $C$ pour le plongement $i: C \hookrightarrow \P^r$, est intégralement close. 
D'après \cite{H} II.5 ex 14, cette algèbre graduée s'identifie à
$$
\bigoplus_{m \in \N} \Gamma (C, i^* \mathcal{O}_{\P^r} (m))
$$
et sa clôture intégrale à 
$$
\bigoplus_{m \in \N} \Gamma (C, \mathcal{O}_C (m)).
$$ 
La normalité projective de $C$ entraîne que pour tout entier naturel $m$, l'application de restriction
$$
\psi_m : \Gamma (\P^r, \mathcal{O}_{\P^r} (m)) \rightarrow \Gamma (C, \mathcal{O}_C(m))
$$
est surjective (cf \cite{H} II.5 ex 14 (d)).
Par ailleurs, le diviseur $G^*$ est linéairement équivalent à $mL_C$, donc le faisceau $\L (G^*)$ est isomorphe à $\mathcal{O}_C (m)$. Considérons le diagramme commutatif
$$
\xymatrix{\relax & \Gamma(S, \mathcal{O}_S(m)) \ar[dd]^{r_m} \\
\Gamma (\P^r, \mathcal{O}_{\P^r}(m)) \ar[ru]^{\phi_m} \ar[rd]_{\psi_m} & \\
 & \Gamma (C, \mathcal{O}_C (m)).
}
$$
La surjectivité de l'application $\psi_m$
entraîne celle de l'application $r_m$. 
\end{proof}

Le second ingrédient de la preuve du théorème \ref{thmreal} est le théorème 3.3 de l'article \cite{poon} de Poonen. Il s'agit d'un théorème ``à la Bertini'' pour des variétés sur des corps finis. Énonçons ce résultat.

\begin{thm}[Poonen 2004]\label{interpoon}
  Soit $X$ un sous-schéma quasi-projectif lisse de $\P^r$ de dimension $m\geq 1$ au-dessus de $\F_q$ et soit $F \subset X$ un ensemble fini de points fermés. Alors, il existe une hypersurface lisse et géométriquement intègre $H\subset \P^r$ telle que l'intersection $H\cap X$ est lisse, de dimension $m-1$ et contient F.
\end{thm}

\begin{rem}\label{connex}
À la suite de ce théorème, l'auteur remarque que, si $X$ est projective et géométriquement connexe, alors $H\cap X$ l'est également d'après \cite{H} corollaire III.7.9.   
\end{rem}

\begin{proof}[\textsc{Démonstration du théorème \ref{thmreal}}]
\textbf{Étape 1. Construction de $\mathbf{\omega}, \mathbf{D_a}$ et $\mathbf{D_b}$.}
Soient $i_1, \ldots, i_s$ les indices du support du mot de code $c$. D'après le théorème \ref{interpoon} et la remarque \ref{connex}, il existe une hypersurface $H$ contenue dans $\P^r$ telle que l'intersection $C:=H\cap S$ est une courbe projective lisse connexe contenant les points $P_{i_1}, \ldots, P_{i_s}$.
On note $G^*$ le tiré en arrière de $G$ sur $C$ par l'inclusion canonique $C \hookrightarrow S$ et $\Lambda_C$ le diviseur
$$
\Lambda_C:=P_{i_1}+ \cdots + P_{i_s}\ \in \textrm{Div}_{\F_q}(C).
$$
D'après le lemme \ref{surj}, l'application $\Gamma (S, \L (G)) \rightarrow \Gamma (C, \L (G^*))$ est surjective et induit donc une application surjective de codes
$$
r: C_{L,S} (\Delta, G) \rightarrow C_{L,C} (\Lambda_C,G^*).
$$

\noindent Soit à présent $c^*:=(c_{i_1}, \ldots , c_{i_s})$ le mot de code poinçonné obtenu en ne conservant que les coordonnées du mot $c$ d'indices $i_1, \ldots, i_s$.
La surjectivité de l'application $r$ entraîne que le mot $c^*$ est un élément de $C_{L,C}(\Lambda_C, G^*)^{\bot}$. On sait également que
$$
C_{L,C}(\Lambda_C, G^*)^{\bot}=C_{\Omega,C}(\Lambda_C, G^*).
$$

\noindent Par conséquent, il existe une $1$-forme $\mu$ sur $C$ appartenant à $\Gamma (C, \Omega^1(G^*-\Lambda_C))$ et telle que $c^*=\res_{\Lambda_C}(\mu)$, où $\res_{\Lambda_C}$ désigne l'application
$$
\res_{\Lambda_C}:\left\{\begin{array}{ccc}
  \Gamma (C, \Omega^1 (G^*-\Lambda_C)) & \rightarrow & \F_q^s \\
  \omega & \mapsto & (\res_{P_{i_1}}(\omega), \ldots, \res_{P_{i_s}}(\omega)).
\end{array}\right.
$$

\noindent Notons que, par hypothèse, les coordonnées du mot de code poinçonné $c^*$ sont toutes non nulles (il a été construit en éliminant les coordonnées nulles du mot $c$). De ce fait, 
\begin{equation}\label{valmu}
\forall k \in \{1, \ldots, s\},\quad \textrm{val}_{P_{i_k}}(\mu)=-1.  
\end{equation}

\noindent Soit à présent $\mu_*$ un relevé arbitraire de $\mu$ sur $S$, c'est-à-dire une $1$-forme rationnelle sur $S$ vérifiant 
$
{\mu_*}_{|C}=\mu.
$
Soit également $v$ une uniformisante de l'anneau $\OC$. Posons alors,
$$
\omega:= \mu_* \w \frac{dv}{v}.
$$
Le diviseur de $\omega$ est de la forme
$
(\omega)=-C+R,
$
où $R$ est un diviseur sur $S$ dont le support ne contient pas la courbe $C$. Pour finir, on pose
$$
D_a:=C \quad \textrm{et} \quad D_b:=G-R.
$$

\noindent Ainsi $\omega$ est bien un élément de $\Gamma (S, \Omega^2 (G-D_a-D_b))$. De plus, comme le $1$-résidu de $\omega$ le long de $C$ est $\mu$, on en déduit
$$
\res^2_{D_a, \Delta}(\omega)=\res^2_{C,\Delta}(\omega)=c.
$$

\medbreak

\noindent \textbf{Étape 2. Sous-$\Delta$-convenance de ($\mathbf{D_a,D_b}$).} 
Soit $\Lambda$, le $0$-cycle sur $S$ défini par
$$
\Lambda:=P_{i_1}+\cdots +P_{i_s}.
$$
Nous allons montrer que la paire $(D_a,D_b)$ est $\Lambda$-convenable. Pour ce faire, nous allons utiliser le critère de la proposition \ref{crit}.
Si l'on note $i$ l'inclusion canonique $i: C \hookrightarrow S$, d'après le lemme \ref{valres} on a l'égalité de diviseurs sur $C$:
$$
(\mu)= i^*(G-D_b).
$$
Comme $\mu$ est un élément de $\Gamma (C, \Omega^1 (G-\Lambda_C))$, on en déduit que $i^*D_b \leq \Lambda$, ce qui implique l'inégalité de $0$-cycles sur $S$ suivante:
\begin{equation}\label{zerocycle}
D_a\cap D_b \leq \Lambda.
\end{equation}
Soit alors
$P$, un point de $\overline{S}$ non contenu dans le support de $\Lambda$ et  n'appartenant pas à la courbe $\overline{C}$ (c'est-à-dire au support de $D_a$). Dans cette situation, $D_a$ peut jouer le rôle\footnote{Voir proposition \ref{crit} pour une description de $D^*$.} de $D_*$ en $P$ et comme $D_a$ est nul au voisinage de ce point, le critère y est trivialement vérifié.
Soit à présent un point $P$ de $\overline{C}$ non contenu dans le support de $\Lambda$. Comme $\overline{C}$ est une courbe irréductible lisse, le diviseur $D_a=C$ peut encore jouer le rôle de $D_*$.  
D'après la relation (\ref{zerocycle}), la multiplicité d'intersection de $D_a$ avec $D_b$ en $P$ est négative, le critère est donc vérifié en ce point.
En un point $P$ du support de $\Lambda$, le diviseur $D_a$ joue encore le rôle de $D_*$ et l'inégalité (\ref{zerocycle}) implique que la multiplicité d'intersection de $D_a$ et $D_b$ en $P$ est inférieure ou égale à $1$. 
D'après la relation (\ref{valmu}) et le lemme \ref{valres} l'inégalité est en fait un égalité. Le couple $(D_a,D_b)$ vérifie donc bien le critère de $\Lambda$-convenance.

\medbreak

\noindent \textbf{Étape 3. Classes d'équivalence linéaire de $\mathbf{D_a}$ et $\mathbf{D_b}$.}
D'après la construction de la courbe $C=D_a$ dans l'étape 1, on sait qu'il existe un entier naturel non nul $n_a$ tel que
$$
D_a \sim n_aL_S.
$$

\noindent D'après \cite{H} II.8 ex 4(e), la classe canonique d'une sous-variété $X$ intersection complète de $\P^r$ est de la forme $K_X \sim kL_X$ où $k$ dépend du degré des hypersurfaces dont l'intersection est égale à $X$.
Soit donc $k$ l'entier tel que $K_S \sim kL_S$.
\noindent Comme le diviseur de $\omega$ vérifie
$$
(\omega)=G-D_a-D_b,
$$
et que $G\sim mL_S$, on en déduit que 
$$
D_b \sim (m-k-n_a)L_S.
$$
\end{proof}

\begin{rem}\label{ducon}
  La courbe $C$ qui définit le diviseur $D_a$ dans la preuve du théorème \ref{thmreal} est construite de manière à interpoler les points correspondant au support du mot de code que l'on peut réaliser. Noter que l'on aurait pu tout aussi bien choisir une bonne fois pour toute une courbe interpolant tous les points du support de $\Delta$ et ne travailler que sur cette dernière.

Notons au passage qu'une telle approche permet de démontrer qu'un code fonctionnel construit sur $S$ à partir d'un diviseur $G\sim mL_S$ se réalise toujours comme code sur une courbe $C$ contenue dans $S$.
Ce fait n'a rien de nouveau, Pellikaan, Shen et  Van Wee montrent dans \cite{PelBZS} que tout code correcteur se réalise comme code sur une courbe.
L'exploitation potentielle de ce fait en vue d'une étude du code fonctionnel sera discutée en section \ref{bebert}.
\end{rem}

\section{Discussion autour du théorème de réalisation}\label{discuss}

Quelques commentaires s'imposent au sujet du théorème \ref{thmreal} et de sa démonstration.
D'abord, il est important de noter que la preuve de ce théorème de réalisation n'est malheureusement pas constructive.
En effet, cette dernière repose sur le théorème \ref{interpoon} de Poonen qui n'est lui-même qu'un résultat d'existence.
Ce dernier ne donne par exemple aucune information sur le degré minimal de l'hypersurface qui permet de construire le diviseur $D_a$.

Ensuite, on rappelle que le résultat n'est démontré que sous certaines conditions, à savoir que la surface $S$ est intersection complète et que le diviseur $G$ est linéairement équivalent à $mL_S$. En fait, ces conditions sont principalement là pour assurer la surjectivité de l'application
$$
\Gamma (S, \L (G)) \rightarrow \Gamma (C, \L(G^*)).
$$
Il s'avère que cette application est fréquemment surjective mais ce n'est pas systématique (un contre-exemple est donné en \ref{P1bu}). Les hypothèses du théorème assurent la surjectivité de l'application pour toute courbe lisse obtenue par intersection de $S$ avec une hypersurface.  
En conclusion, il s'agit de conditions suffisantes, mais absolument pas nécessaires. Il est fort possible que le résultat reste vrai en omettant ces hypothèses, nous n'avons cependant pas été à même de le démontrer dans un cas plus général.
L'exemple élémentaire présenté dans la section \ref{P1bu} va confirmer l'aspect non nécessaire de ces hypothèses.

Cela nous amène à poser la question ouverte suivante.

\begin{ques}\label{hypless}
Le résultat du théorème de réalisation reste-t-il vrai si l'on élimine les hypothèses que doivent vérifier $S$ et $G$?  
\end{ques}

\noindent Notons également que le théorème de réalisation (plus exactement le corollaire \ref{correal}) répond à la question \ref{Qsom} posée page \pageref{Qsom} sous certaines hypothèses sur $S$ et $G$. Cependant, si l'on sait que sous ces hypothèses l'orthogonal d'un code fonctionnel se réalise comme une somme de codes différentiels, la question suivante reste ouverte.

\begin{ques}\label{bornsom}
  Sous les conditions du corollaire \ref{correal}, peut-on estimer le nombre de minimal de codes différentiels dont la somme est égale à l'orthogonal d'un code fonctionnel en fonction d'invariants géométriques de la surface?
\end{ques}

\noindent L'exemple qui suit a été suggéré par Antoine Ducros.

\subsection{Un exemple de réalisation sans que les conditions du théorème de \ref{thmreal}  soient vérifiées}\label{P1bu}

Soient $S$ la surface obtenue par l'éclatement de $\P^2$ en un point $O$ et
$$
\pi: S \rightarrow \P^2,
$$
l'éclatement de $\P^2$ en $O$. 
Le diviseur $G$ est le diviseur exceptionnel de $S$ et le $0$-cycle $\Delta$, la somme des points rationnels de $S$ non contenus dans le support de $G$.
La surface $S$ peut être plongée dans $\P^5$ via le plongement de Segré (\cite{sch1} I.5.1). Pour ce plongement, $S$ est un intersection complète.
Cependant, le diviseur $G$ ne peut s'identifier à une section hyperplane de $S$ pour aucun plongement de cette surface. En effet, il est d'auto-intersection $-1$, donc ne vérifie pas le critère de Nakai-Moishezon (\cite{H} théorème V.1.10). 
L'espace $\Gamma(S, \L (G))$ est de dimension $1$ et ne contient que les fonctions constantes.
En effet, comme $G$ est d'auto-intersection négative, il est le seul élément du système linéaire $|G|$ qui est donc de dimension nulle. Par conséquent, la dimension de $\Gamma(S, \L (G))$ est égale à $1$. On vérifie ensuite que les constantes sont bien des éléments de cet espace.

Soit à présent $L$ la transformée stricte d'une droite de $\P^2$ passant par $O$. La courbe $L$ intersecte $G$ transversalement en un unique point $Q$. Le tiré en arrière $G^*$ de $G$ par l'inclusion canonique de $L$ dans $S$ est égal à $Q$. De fait, l'espace  $\Gamma (L, \L (G^*))$ est de dimension $2$ et donc l'application
$$
\Gamma (S, \L (G)) \rightarrow \Gamma (L, \L (G^*))
$$
n'est pas surjective.
Montrons maintenant que l'on peut tout de même réaliser tous les mots de $C_L (\Delta, G)^{\bot}$ comme résidus de $2$-formes sur $S$.

\medbreak

\subsubsection{Approche non constructive.}
Soit $c$ un mot de $C_L (\Delta, G)^{\bot}$ et soit $\Lambda$ le $0$-cycle de $S$ correspondant au support de $c$. Il existe une courbe irréductible lisse $C$ de $S$ qui contient tous les points du support de $\Lambda$ et dont l'intersection avec $G$ est vide. En effet, cela revient à construire une courbe lisse de $\P^2$ qui interpole une famille finie de points et évite le point $O$.
Le tiré en arrière $G^*$ de $G$ sur $C$ est nul et donc $\Gamma (C, \L (G ^*))$ est l'ensemble des fonctions constantes sur $C$. Par conséquent, l'application
$$
\Gamma (S, \L (G)) \rightarrow \Gamma (C, \L (G^*))
$$
est surjective et on peut effectuer la construction effectuée dans la démonstration du théorème \ref{thmreal}.

\medbreak

\subsubsection{Approche constructive.}
Le code $C_{L,S}(\Delta, G)$ est le code de répétition pure et de longueur $n=q^2+q$. Son orthogonal est donc un code de dimension $n-1$. 
On note $c_2, \ldots ,c_n$ les mots de la forme
$$
c_i:=(1, 0, \ldots, 0,-1,0, \ldots, 0),
$$
le ``$-1$'' apparaissant en $i$-ème position. La famille $(c_2, \ldots, c_n)$ est une base de $C_L (\Delta, G)^{\bot}$.
D'après la proposition \ref{realsom}, il suffit de réaliser ces $n-1$ mots pour montrer que tout mot de $C$ est réalisable.

\medbreak

\noindent \textbf{Étape 1.} Soit donc $i$ un entier compris entre $2$ et $n$ et supposons que les points $\pi(P_1)$ et $\pi( P_i)$ ne sont pas alignés avec $O$ dans $\P^2$. On appelle $C$, la droite de $\P^2$ reliant $\pi (P_1)$ et $\pi (P_i)$. On choisit deux droites $C_1$ et $C_i$ dans $\P^2$ distinctes de $C$ et contenant respectivement $\pi( P_1)$ et $\pi (P_i)$ et évitant le point $O$.

\begin{center}
\includegraphics[height=4.63cm, width=6cm]{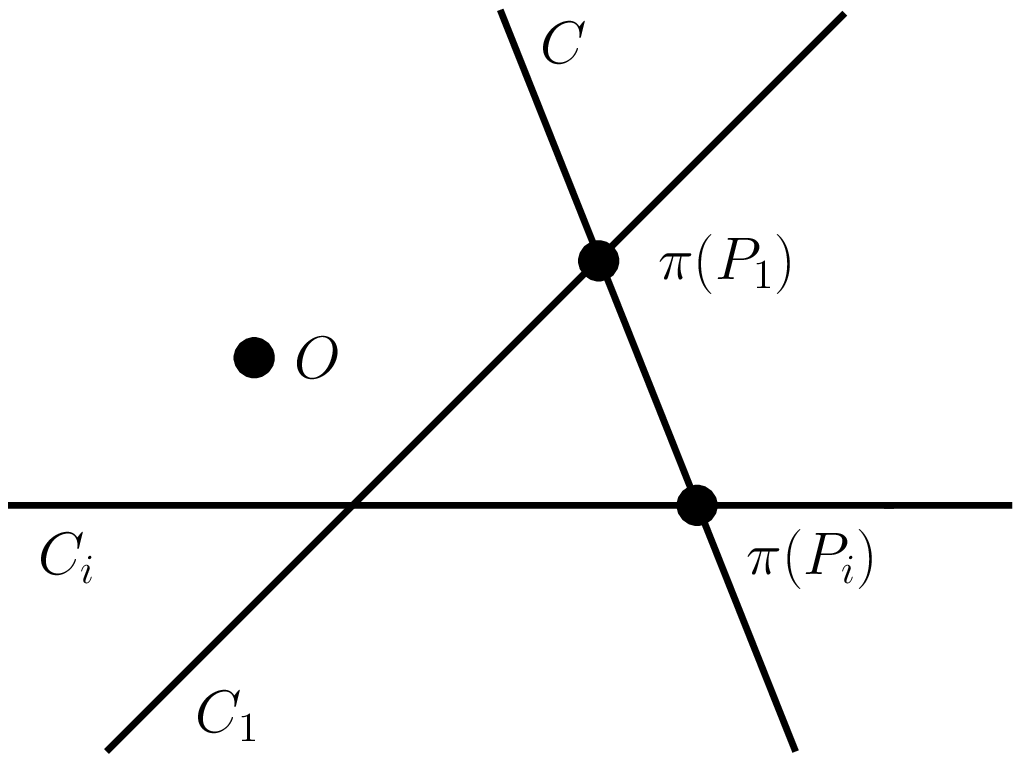}
\end{center}

\noindent On rappelle que la classe canonique dans $\P^2$ est égale à $-3L$, où $L$ désigne la classe d'équivalence linéaire des droites du plan projectif. De fait, le diviseur $-C-C_1-C_i$ est canonique, il existe donc une $2$-forme $\omega$ sur $\P^2$ telle que
$$
(\omega):=-C-C_1-C_2.
$$
D'après le lemme \ref{valres}, la $1$-forme $\res^1_C (\omega)$ sur $C$ n'a de pôles qu'en $\pi (P_1)$ et $\pi (P_i)$ et ces pôles sont simples. Elle a donc des résidus non nuls en ces points et d'après la formule des résidus, ils sont opposés. De ce fait, quitte à multiplier $\omega$ par un scalaire non nul, on a
$$
\res^2_{C,P_1}(\omega)=1 \quad \textrm{et} \quad \res^2_{C,P_i}(\omega)=-1.
$$
De plus, la $2$-forme $\omega$ n'a ni zéro ni pôle au voisinage de $O$. Donc, d'après le lemme \ref{valex} la $2$-forme $\pi^* \omega$ sur $S$ est de valuation $1$ le long du diviseur exceptionnel, on a donc
$$
(\pi^* \omega)=G-\widetilde{C}-\widetilde{C}_1-\widetilde{C}_2.
$$
On pose 
$$
\Lambda_i:=P_1+P_i, \quad D_a:=\widetilde{C} \quad \textrm{et} \quad D_b:=\widetilde{C}_1+\widetilde{C}_2
$$
et on vérifie aisément que $(D_a, D_b)$ est $\Lambda_i$-convenable (on peut par exemple voir qu'il satisfait le critère de la proposition \ref{crit}). De fait, le mot $c_i$ est réalisé par la $2$-forme $\pi^* \omega$.

\medbreak

\noindent \textbf{Étape 2.} Si maintenant les points $\pi (P_1)$ et $\pi (P_k)$ sont alignés avec $O$.
On choisit deux autres points rationnels $\pi (P_j)$ et $\pi (P_k)$ de $\P^2$ tels que les points $\pi (P_1),\ \pi (P_i),\ \pi (P_j)$ et $\pi (P_k)$ soient en position générale (trois d'entre eux ne sont pas alignés).
Il existe au moins deux coniques rationnelles $C$ et $C'$ distinctes interpolant ces quatre points et évitant le point $O$. En effet, le système linéaire des coniques interpolant ces points est de dimension $1$, donc même si le corps de base est $\F_2$, il y a au moins $3$ éléments dans ce système et un seul d'entre eux interpole $0$.
On appelle $C''$ la droite reliant $\pi (P_j)$ et $\pi (P_k)$.

\begin{center}
\includegraphics[height=6.4cm , width=10cm]{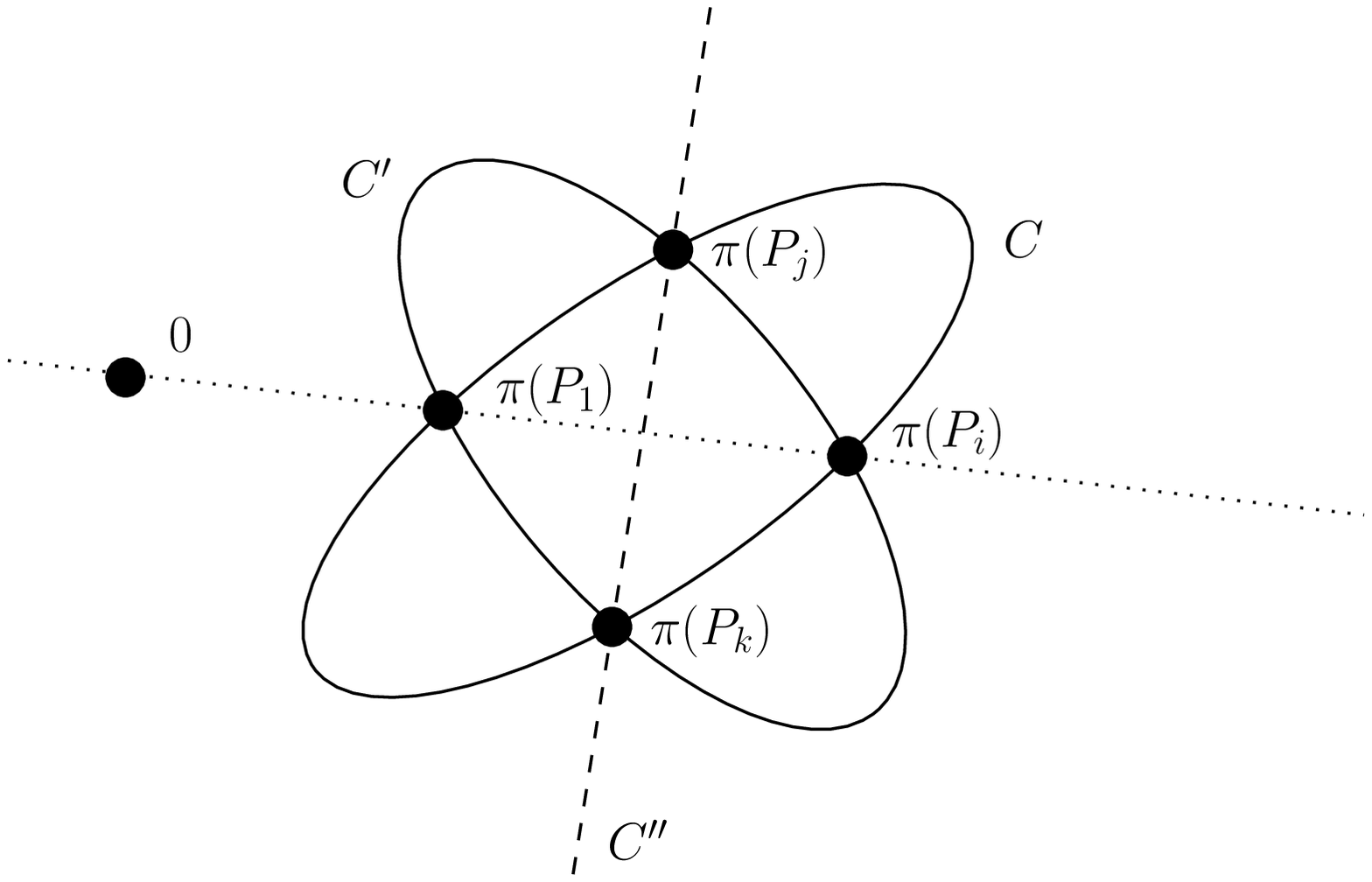}
\end{center}

\noindent Le diviseur $-C-C'+C''$ est linéairement équivalent à $-3L$, c'est donc un diviseur canonique et il existe une $2$-forme $\omega$ sur $\P^2$ vérifiant 
$$
(\omega)=-C-C'+C''.
$$
Avec le lemme \ref{valres} on montre que $\pi (P_1)$ et $\pi (P_i)$ sont les seuls pôles de la $1$-forme $\res^1_{C}(\omega)$ sur $C$. On en déduit que, quitte à multiplier $\omega$ par un scalaire inversible, on a  
$$
\res^2_{C,P_1}(\omega)=1 \quad \textrm{et} \quad \res^2_{C,P_i}(\omega)=-1.
$$
Par ailleurs, $\omega$ n'a ni zéro ni pôle au voisinage de $O$, donc $\pi^* \omega$ vérifie
$$
(\pi^* \omega)=G-\widetilde{C}-\widetilde{C}'+\widetilde{C}''.
$$
On finit en posant
$$
\Lambda_i:=P_1+P_i, \quad D_a:=\widetilde{C} \quad \textrm{et} \quad D_b:=\widetilde{C}'+\widetilde{C}''
$$
et en vérifiant (grâce au critère de la proposition \ref{crit}) que la paire $(D_a,D_b)$ ainsi construite est bien $\Lambda_i$-convenable.

\section{Une autre application possible des théorèmes ``à la Bertini''}\label{bebert}

Cette section, qui conclut le chapitre \ref{chapreal}, a pour but de montrer qu'une réponse à une certaine question ouverte pourrait permettre dans certaines situations de minorer la distance minimale du code fonctionnel $C_L(\Delta,G)$.
L'objectif est d'utiliser la constatation de la remarque \ref{ducon}.
Avant d'y arriver ouvrons une parenthèse historique sur la théorie des codes géométriques.

\subsection{Les travaux de Pellikaan, Shen et Wee}

Dans \cite{PelBZS}, les auteurs classent les codes correcteurs en WAG (\textit{Weakly Algebraic-Geometric}), AG (\textit{Algebraic-Geometric}) et SAG (\textit{Strongly  Algebraic-Geometric}).
Les codes WAG sont les codes $\Gamma$ admettant une représentation géométrique, c'est-à-dire les codes pour lesquels il existe une courbe $C$, un diviseur $G$ sur $C$ et une famille $P_1, \ldots, P_n$ de points rationnels de $X$ tels que
$$\Gamma=C_{L,C}(D,G)\quad  \textrm{avec} \quad D:=P_1+ \cdots +P_n.$$

\noindent Les codes AG sont les codes WAG qui admettent une représentation vérifiant $n > \deg (G)$. Quant aux SAG ce sont les WAG admettant une représentation telle que $n>\deg (G)>2g_C-2$.

L'un des résultats majeurs de l'article \cite{PelBZS} est le théorème 2 qui affirme que tout code est WAG.
Dans la suite de l'article, les auteurs donnent des exemples de codes qui ne le sont pas.
Ils signalent par exemple que les codes de Golay binaires ne sont pas AG (\cite{PelBZS} corollaire 9).

Le problème des codes est que toutes leurs réalisations comme codes sur des courbes donne une distance construite de Goppa nulle. 
Par conséquent, une représentation géométrique du code ne fournit aucune information sur sa distance minimale. 
En ce qui concerne les codes fonctionnels sur une surface algébrique, notre objectif va être de savoir si l'on peut disposer d'une représentation AG.

\subsection{Le cas des codes fonctionnels sur une surface}

Jusqu'à la fin du chapitre, nous nous plaçons sous l'hypothèse \ref{secreal} énoncée page \pageref{secreal}.
Reprenons la remarque \ref{ducon}.
Le théorème \ref{interpoon} de Poonen nous assure l'existence d'une courbe lisse géométriquement intègre $C\subset S$ obtenue par intersection de $S$ avec une hypersurface de $\P^r$ et qui interpole tous les points $P_1, \ldots, P_n$ du support de $\Delta$.
Puis, d'après le lemme \ref{surj}, l'application de restriction à $C$
$$
\Gamma (S, \L (G)) \rightarrow \Gamma (C, \L(G^*))
$$
est surjective.
Si l'on note $D$ le diviseur sur $C$ défini par $D:=P_1+\cdots + P_n$, alors la surjectivité de l'application ci-dessus entraîne que les codes $C_{L,S}(\Delta, G)$ et $C_{L,C}(D, G^*)$ sont identiques.
Le code $C_{L,S}(\Delta,G)$ s'identifie donc à un code sur une courbe algébrique.
Si $G$ est linéairement équivalent à $mL_S$, alors $G^*$ es linéairement équivalent à $mL_C$. Le tout est de savoir quel est le degré de $L_C$. Ce degré est le nombre de points géométriques de l'intersection de $C$ avec un hyperplan géométrique générique, c'est donc le degré de la courbe $C$ pour son plongement dans $\P^r$. 
Enfin, au vu de la construction de $C$, son degré n'est autre que le degré de $S$ multiplié par le degré de l'hypersurface $H$ telle que $C=H \cap S$.
Par conséquent, une estimation suffisamment fine du degré de l'hypersurface $H$ permettrait de minorer la distance minimale du code $C_{L,S}(\Delta, G)$.

\begin{ques}[Arithmétique]\label{qpoon}
Soient $X$ une variété projective lisse géométriquement intègre sur un corps fini $\F_q$ et $P_1, \ldots, P_n$, une famille de points fermés de $X$.
Peut-on évaluer explicitement ou majorer de façon précise le plus petit entier $d$ tel qu'il existe au moins une hypersurface définie sur $\F_q$ de degré inférieur ou égal à $d$ qui interpole tous les $P_i$ et dont l'intersection schématique avec $X$ soit une sous-variété lisse géométriquement intègre de codimension $1$?
\end{ques}

Dans \cite{poon}, Poonen montre que de telles hypersurfaces existent et forment même un sous-ensemble de densité positive dans l'ensemble des hypersurfaces de $\P^r$, mais il ne donne aucune information sur le degré minimal d'une telle variété.
On est de ce fait assurés de l'existence de l'entier $d$ mais ne dispose d'aucun procédé d'estimation explicite.

Notons que la question \ref{qpoon}A (Arithmétique) ne porte que sur les hypersurfaces définies sur $\F_q$. La remarque qui suit assure que l'on peut se poser la question pour les hypersurfaces définies sur $\overline{\F}_q$.

\begin{rem}\label{clotalg}
Soit $\F_{q^m}$ une extension de $\F_q$ et $S'$ la surface $S':=S \times_{\F_q} \F_{q^m}$. Notons $\Delta'$ et $G'$ pour les tirés en arrière respectifs de $\Delta$ et $G$ sur $S'$.
On dispose alors de l'égalité de codes
$$
C_{L,S}(\Delta, G) \otimes_{\F_q} \F_{q^m} = C_{L,S'}(\Delta', G').
$$

\noindent En particulier, ces codes ont la même distance minimale.
\end{rem}

Par conséquent, on peut chercher une réalisation du code $C_{L,S}(\Delta, G)\otimes \F_{q^m}$ pour une extension quelconque de $\F_q$, ce qui ramène notre problème à la question suivante.

\medbreak

\noindent \textbf{Question \ref{qpoon}} (Géométrique)\textbf{.}\label{qbert}
\textit{  Soit $X$ une variété projective irréductible lisse définie sur $\overline{\F}_q$ et $P_1, \ldots, P_n$ une famille de points de $X$.
Peut-on évaluer explicitement ou majorer de façon précise le plus petit entier $d$ tel qu'il existe au moins une hypersurface $H$ de degré inférieur ou égal à $d$, qui contienne tous les $P_i$ et telle que $H\cap X$ soit une sous-variété lisse de codimension $1$ de $X$?}

Cette dernière question ressemble fortement à un théorème ``à la Bertini''. En effet, si l'on note $\mathfrak{d}_d$ le système linéaire des sections hyperplanes de $X$ de degré $d$ et $\mathfrak{d}_d'$ le sous-système linéaire de $\mathfrak{d}_d$ des sections hyperplanes  de $X$ interpolant les $P_i$, alors la question \ref{qpoon}G (Géométrique) se traduit par: \textit{le système linéaire $\mathfrak{d}_d'$ possède-t-il un élément irréductible lisse?}

Les questions \ref{qpoon}A et \ref{qpoon}G restent ouvertes. Notons que l'article \cite{AK} d'Altman et Kleiman donne une piste pour tenter d'y répondre.
Les commentaires ci-dessous ont été suggérés par Steven L. Kleiman.

\begin{thm}[Altman, Kleiman 1979]
Soient $k$ un corps infini, $X$ un sous-schéma de l'espace projectif $\P^r_k$ et $Z$ un sous-schéma de $X$. Soit $\mathcal{I}_{\overline{Z}}$ le faisceau d'idéaux de $\mathcal{O}_{\P^r}$ associé à l'adhérence de $Z$ dans $X$.
Soit $d$ un entier tel que $\mathcal{I}_{\overline{Z}} (d)$ est engendré par ses sections globales.
Supposons que $X\smallsetminus \overline{Z}$ est lisse, alors l'intersection de $X$ par des hypersurfaces génériques indépendantes de degré $d+1$ est lisse hors de $\overline{Z}$.
\end{thm}

Dans notre situation, soit $Z$ la réunion des points $P_1, \ldots, P_n$.
Supposons que l'on connaisse un entier $d$ tel que $\mathcal{I}_{Z}(d)$ soit engendré par ses sections globales et que pour tout $i$, le faisceau $(\mathcal{I}_Z/\mathfrak{m}_{P_i}\mathcal{I}_Z) (d)$ soit engendré par les sections globales de $\mathcal{I}_Z(d)$.
Alors, une section globale générique de $\mathcal{I}_Z(d)$ sera envoyée sur un élément non nul de $\mathfrak{m}_{P_i}/\mathfrak{m}_{P_i}^2$ et sera donc non singulière en ce point, elle sera également lisse hors de $Z$ d'après le théorème ci-dessus.
Le problème reste en tous les cas de trouver un tel entier $d$ ou une borne supérieure pour celui-ci.


\newpage
\thispagestyle{empty}

\chapter{Orthogonal d'un code fonctionnel}\label{chaporth}

Dans ce chapitre, nous allons travailler sur le problème de la minoration de la distance minimale de l'orthogonal d'un code fonctionnel.
Nous allons présenter deux approches.
La première s'applique aux codes fonctionnels construits à partir de variétés projectives de dimension quelconque et pas seulement aux surfaces.
Elle est de plus indépendante de tous les résultats précédemment énoncés et ne fait pas intervenir la notion de formes différentielles.
La seconde approche, elle, utilise les résultats obtenus dans le chapitre \ref{chapreal}.

Pour le reste, ce chapitre ne peut être considéré comme totalement abouti.
Il ouvre cependant un certain nombre de problèmes de géométrie algébrique sur les systèmes linéaires dont la résolution permettrait d'obtenir des minorations de la distance minimale de l'orthogonal d'un code fonctionnel.

\section*{Notations}

Nous allons reprendre dans ce chapitre un certain nombre de notations utilisées dans le chapitre \ref{chapreal}.
En particulier, on rappelle que si $X$ est une sous-variété fermée d'un espace projectif, on note $L_X$ la classe d'équivalence linéaire d'une section hyperplane de $X$ et $K_X$ la classe canonique sur $X$. 
\section{Première approche}\label{approch1}

Cette approche consiste en fait à n'utiliser que des méthodes d'algèbre linéaire relativement élémentaires.
Dans cette section, $N$ désigne un entier naturel non nul et $k$ un corps quelconque.
On se donne une variété projective lisse géométriquement intègre $X$ intersection complète dans un espace projectif $\P^N$ et  munie d'un diviseur $G$ et d'un $0$-cycle $\Delta$.
On suppose également qu'il existe un entier naturel $m$ tel que $G\sim mL_X$ et que $\Delta$ est une somme de points rationnels de $X$ qui évitent le support de $G$.

\subsection{Notion de $m$-généralité}

Pour alléger les notations, on désignera
l'espace $\Gamma (\P_k^N, \mathcal{O}_{\P_k^N}(m))$ des formes homogènes de degré $m$ sur $\P^N_k$ par $\mathcal{F}_m^N$.
Enfin, pour tout point rationnel $P$ de $\P^N$, on note $\ev_P$ l'application d'\textit{évaluation} décrite dans la définition \ref{evproj} (voir annexe \ref{annexelachaud}).

\begin{defn}
On dit que les points $P_1, \ldots , P_r \in \P^N (k)$ sont liés en degré $m$ ou $m$-liés si les formes linéaires $\ev_{P_1}, \ldots, \ev_{P_r}$ sont liées dans le dual ${(\mathcal{F}^N_m)}^{\vee}$ de $\mathcal{F}^N_m$. 
Dans le cas contraire, si ces formes linéaires forment une famille libre de ${(\mathcal{F}^N_m)}^{\vee}$, on dit que ces points sont en position $m$-générale.
\end{defn}

La définition ci-dessus peut s'interpréter de façon géométrique, comme le montre le lemme qui suit.

\begin{lem}
  Une famille de $r$ points $P_1, \ldots, P_r$ de $\P^N$ est en position $m$-générale, si et seulement si pour tout entier $i \in\{1, \ldots, r\}$, il existe une hypersurface de degré $m$ qui contient les points $P_1, \ldots, \widehat{P_i}, \ldots, P_r$ et évite $P_i$.
\end{lem}

\begin{proof}
C'est un exercice élémentaire de dualité en algèbre linéaire.
\end{proof}

\begin{rem}
  La notion de $1$-généralité correspond à la définition classique de position générale. Le plus souvent, dans la littérature, une famille de points de $\P^N$ est dite en position générale si et seulement ces points forment un repère projectif de la variété linéaire projective qu'ils engendrent.
Si l'on se donne un entier $m$, alors une famille de points de $\P^N$ sont en position $m$-générale si et seulement si leurs images par le $m$-ème plongement de Veronèse\footnote{Voir \cite{sch1} I.4.4.} sont en position générale au sens classique (décrit ci-dessus) dans $\P^M$ avec $M=\left(
\begin{array}{c}
N+d \\
d
\end{array}
\right)-1$

\end{rem}

\begin{rem}
On aurait pu donner une définition plus générale de ces notions en ne considérant plus seulement des points rationnels de $\P^N$, mais des points fermés et même des points infiniment près\footnote{C'est-à-dire des points appartenant à une variété obtenue par une séquence d'éclatements de sous-variétés de $X$. Voir \cite{H} V.3.} de $\P^N$.
Un tel point de vue étant totalement inutile dans ce qui suit, nous avons choisi de nous restreindre au cas des points rationnels de $\P^N$.
\end{rem}

\begin{exmp}\label{gendim1}
 Supposons que $N=1$, on travaille donc sur la droite projective.
Dans ce cas, $r$ points deux à deux distincts $P_1, \ldots, P_r$ sont en position $m$-générale si et seulement si $r\leq m+1$. En effet, on peut construire une forme homogène de degré inférieur ou égal à $m$ ayant exactement $r-1$ racines données. 
\end{exmp}

\subsection{Systèmes linéaires de $\P^N$}
La notion de $m$-généralité peut se reformuler en termes de systèmes linéaires.
Pour ce faire, on adoptera les notations de \cite{H} V.4.
Soient $m$ un entier naturel, $\mathfrak{d}$ le système linéaire sur $\P^N$ des hypersurfaces de degré $m$ et
$P_1, \ldots, P_r$ une famille de points rationnels de $\P^N$. Pour tout entier naturel $1\leq i \leq r$, on note
$\mathfrak{d}_i$
le sous-système linéaire de $\mathfrak{d}$ des hypersurfaces de degré $d$ contenant les points $P_1, \ldots, \widehat{P}_i,$ $\ldots, P_r$.
Selon les notations de \cite{H} V.4,
$$\mathfrak{d}_i:=\mathfrak{d}-P_1-\cdots -\widehat{P}_i-\cdots -P_r.$$
En termes de systèmes linéaires, la $m$-généralité de $P_1, \ldots, P_r$ se formule de la façon suivante.
Les points $P_1, \ldots, P_r$ sont en position $d$-générale si et seulement si pour tout $i$, le point $P_i$ n'est pas un point base du système linéaire $\mathfrak{d}_i$.

\subsection{Lien avec les notions de distance minimale}

Munis de ces définitions, la question de la minoration de la distance minimale du code $C_{L,X} (\Delta, G)^{\bot}$ peut se traduire sous forme d'un problème géométrique.

\begin{prop}\label{dmindu}
Soit $m$ un entier tel que $G \sim mL_X$.
La distance minimale $d^{\bot}$ du code $C_{L,X} (\Delta, G)^{\bot}$ est  égale au nombre minimal $s$ de points $P_1, \ldots, P_s$ du support de $\Delta$ qui sont $m$-liés. 
\end{prop}

\begin{proof}
  Soit $c=(c_1, \ldots, c_n)$ un mot de $C_L (\Delta, G)^{\bot}$.
Cela signifie que l'application $\varphi_c:=c_1\ev_{P_1}+\cdots +c_n \ev_{P_n}$ est identiquement nulle sur $\Gamma(X, \L (G))$.
Comme $X$ est intersection complète, d'après \cite{H} II.8 ex 14, elle est projectivement normale et il y a donc une application surjective
$$
f:\Gamma (\P^N, \mathcal{O}_{\P^N}(m))\rightarrow \Gamma(X, \L (G)).
$$

\noindent L'application $\varphi_c\circ f$ définie sur $\Gamma (\P^N, \mathcal{O}_{\P^N}(m))$ est donc nulle, or ce morphisme n'est autre que $c_1\ev_{P_1}+\cdots +c_n \ev_{P_n}$ vu comme une forme linéaire sur $\Gamma (\P^N, \mathcal{O}_{\P^N}(m))=\mathcal{F}_m^N$.
\end{proof}

Cette reformulation du problème de la distance minimale de $C_L (\Delta, G)^{\bot}$, nous ramène à un problème qui est loin d'être aussi élémentaire qu'il en a l'air.
Autant il est aisé de savoir si une famille de points sont indépendants en dimension $1$ (voir exemple \ref{gendim1}), autant le problème se complique lourdement en dimension supérieure.
En d'autres termes, il est très difficile en dimension supérieure à $2$ de décider si une famille de points est en position $d$-générale ou, ce qui revient au même, de montrer qu'un système linéaire n'a pas d'autres points bases que ceux qu'on lui a assignés. Pour s'en convaincre on peut regarder les démonstrations du chapitre V.4 de \cite{H}.

\subsection{Minorations de la distance minimale de l'orthogonal d'un code fonctionnel}\label{ssecminor}

Nous allons utiliser la proposition \ref{dmindu} pour obtenir deux résultats de minoration de la distance minimale  du code $C_L(\Delta, G)^{\bot}$.

\begin{thm}\label{minor}
 On suppose $N$ supérieur ou égal à $2$.
  Soit $m$ un entier tel que $G\sim mL_X$, alors
  \begin{enumerate}
  \item\label{minor1}  la distance minimale $d^{\bot}$ du code $C_{L,X} (\Delta, G)^{\bot}$ vérifie $$d^{\bot} \geq m+2$$ et il y a égalité si et seulement si le support de $\Delta$ contient $m+2$ points alignés;
  \item\label{minor2} sinon, si le support de $\Delta$ ne contient pas $m+2$ points alignés, alors $$d^{\bot}\geq 2m+2$$ et il y a égalité si et seulement si le support de $\Delta$ contient $2m+2$ points sur une même conique plane.  
  \end{enumerate}
\end{thm}

La démonstration du (\ref{minor1}) de ce théorème fera appel aux lemmes \ref{trucalacon2} et \ref{kipu2} qui suivent et qui seront démontrés en annexe \ref{annexegen}.

\begin{lem}\label{trucalacon2}
Soient $r$ et $m$ deux entiers naturels avec $r\geq 1$, alors toute famille de $rm+2$ points rationnels distincts de $\P^N$ appartenant à une même courbe de degré $r$ est $m$-liée.  
\end{lem}

\begin{lem}\label{kipu2}
  Soit $m$ un entier naturel.
\begin{enumerate}
\item\label{my2} Si $m+2$ points rationnels distincts de $\P^N$ sont $m$-liés, alors ils sont alignés.
\item\label{cock2} Tout $r$-uplet de points rationnels deux à deux distincts de $\P^N$ avec $r\leq m+1$ est en position $m$-générale.
\end{enumerate}
\end{lem}

\begin{proof}[Démonstration du (\ref{minor1}) du théorème \ref{minor}]
La proposition \ref{dmindu} et la propriété (\ref{cock2}) du lemme \ref{kipu2} entraînent que la distance minimale du code $C_L(\Delta, G)^{\bot}$ est supérieure à $m+2$ et qu'une condition nécessaire pour qu'elle soit atteinte est que le support de $\Delta$ contienne $m+2$ points alignés.
D'après le lemme \ref{trucalacon2}, cette dernière condition est suffisante.
\end{proof}

\begin{exmp}
  Si $G\sim L_X$, alors la distance minimale de $C_{L,X}(\Delta, G)^{\bot}$ est minorée par $3$. Cette borne est atteinte dès que le support de $\Delta$ contient trois points alignés.
Remarquons que la borne est par exemple atteinte dès que $X$ contient une droite rationnelle.
\end{exmp}

\medbreak

Le point (\ref{minor2}) du théorème \ref{minor} se démontre de la même façon que le point (\ref{minor1})  en utilisant la proposition \ref{dmindu}, le lemme \ref{trucalacon2} et le lemme \ref{chiasse2} énoncé ci-dessous.
Nous renvoyons le lecteur à l'annexe \ref{annexegen} pour une démonstration de ce dernier.

\begin{lem}\label{chiasse2}
Soient $m$ et $r$ deux entiers naturels tels que $r\leq 2m+1$.
\begin{enumerate}
\item\label{butt} Une famille de $r$ points distincts de $\P^N$  telle que $m+2$ d'entre eux sont non alignés est en position $m$-générale.
\item\label{your} Soit $P_1, \ldots, P_{2m+2}$ un $(2m+2)$-uplet de points rationnels distincts de $\P^N$ tels que $m+2$ d'entre eux ne sont pas alignés.
Alors, ces points sont $m$-liés, si et seulement s'ils appartiennent à une même conique plane.
\end{enumerate}
\end{lem}

\subsubsection{Peut-on aller plus loin?}

Une généralisation naturelle (mais fausse) du théorème \ref{minor} serait:
\textit{``soient $m,s$ deux entiers naturels, supposons que pour tout $r<s$, un $(rm+2)$-uplet de points du support de $\Delta$ n'est jamais contenu dans une courbe de degré $r$, alors $d^{\bot}\geq sm+2$.''}

Malheureusement, ce résultat est faux.
En effet, d'après la proposition \ref{dmindu}, un tel résultat impliquerait que $sm+1$ points de $\P^N$ tels que pour tout $r<s$, un $(rm+2)$-uplet d'entre eux n'est jamais contenu dans une courbe de degré $r$, sont en position $m$-générale.
Or, si $s=3$ et $m=3$, cela signifierait que $10$ points de $\P^2$ tels que $5$ d'entre eux sont non alignés et $8$ d'entre eux ne sont pas sur une même conique sont toujours en position $3$-générale. 
Or d'après \cite{H} corollaire V.4.5, on peut construire un $9$-uplet de points $3$-liés vérifiant ces propriétés.
Un tel $9$-uplet de points est construit en prenant les points d'intersection de deux cubiques réduites sans composante irréductible commune.
Ces configurations de points provenant d'intersections de $N$ hypersurfaces dans $\P^N$ sont difficiles à repérer et compliquent les démonstrations de $m$-généralité lorsque l'on veut améliorer les lemmes \ref{kipu2} et \ref{chiasse2}.

\medbreak

En conclusion, on sait que les deux premières configurations minimales de points rationnels $m$-liés dans $\P^N$ sont
\begin{enumerate}
\item[$(i)$] $m+2$ points alignés;
\item[$(ii)$] $2m+2$ points sur une même conique plane.
\end{enumerate}

\noindent Nous laissons une question ouverte.

\begin{ques}\label{configmin}
Quelles sont les configurations minimales suivantes?  
\end{ques}

\subsection{Applications}

Le théorème \ref{minor} permet d'obtenir des minorations assez fines de la distance minimale de l'orthogonal du code $C_{L,X}(\Delta,G)^{\bot}$ dans le cas où l'entier $m$ tel que $G\sim mL_X$ est petit.
Commençons par étudier le cas bien connu où $X$ est une courbe.

\subsubsection{Courbes algébriques planes, comparaison avec la distance minimale construite de Goppa}

Soit $X$ une courbe algébrique projective plane lisse de degré $d\geq 2$ et définie sur $\F_q$.
Soient $m$ un entier naturel et $G$ un diviseur sur $X$ linéairement équivalent à $mL_X$.
Soient enfin $P_1, \ldots, P_n$ une famille de points rationnels de $X$ qui évitent le support de $G$ et $D$ le diviseur $D:=P_1+\cdots +P_n$.
On rappelle que le genre de $X$ s'obtient par la formule
$$g_X=\frac{(d-1)(d-2)}{2}$$
et que l'orthogonal du code fonctionnel $C_L (D,G)$ est le code différentiel $C_{\Omega}(D,G)$. Par ailleurs, on rappelle également que la distance minimale $d^{\bot}$ du code $C_{\Omega}(D,G)$ (qui est égal à $C_L (D,G)^{\bot}$) vérifie
$$
d^{\bot} \geq \deg (G)-(2g_X-2).
$$ 

\noindent La quantité $\deg (G)-(2g_X-2)$ est appelée \textit{distance construite de Goppa} et notée $\delta^{\bot}$. 

La courbe $X$ est supposée plane et lisse, elle est donc irréductible. Ainsi, comme $d\geq 2$, alors $X$ ne contient pas plus de $d$ points géométriques alignés.
Par conséquent, nous allons distinguer les cas $0 \leq m \leq d-2$ et $m\geq d-1$.
\medbreak

\begin{itemize}
\item[\textbullet] \textbf{Si $\mathbf{0 \leq m \leq d-2}$}, alors le théorème \ref{minor} (\ref{minor1}) nous fournit la minoration
$$
d^{\bot} \geq m+2.
$$

\noindent Quant à la distance construite de Goppa, on peut l'exprimer en fonction de $m$ et $d$.
En effet, comme $G \sim mL_X$, on en déduit que le degré de $G$ est $md$ et
$$
\delta^{\bot}=md-(d-1)(d-2)+2.
$$

\noindent Faisons la différence de ces deux minorants de $d^{\bot}$.
$$
m+2 -(md-(d-1)(d-2)+2) = m-md+(d-1)(d-2)=(d-1)(d-2-m).
$$

\noindent En conclusion, la minoration fournie par le théorème \ref{minor} (\ref{minor1}) est meilleure que la distance construite de Goppa si $m<d-2$. Elle est en particulier nettement meilleure lorsque $m$ est petit.

\item[\textbullet] \textbf{Si $\mathbf{m \geq d-1}$}, alors, d'après le théorème de Bezout, $m+2$ points de $X$ ne sont jamais alignés. Le théorème \ref{minor} (\ref{minor2}) fournit la minoration
$$
d^{\bot} \geq 2m+2.
$$ 

\noindent La différence entre ce minorant de $d^{\bot}$ et la distance construite de Goppa est 
$$
2m+2 - \delta^{\bot}=(d-2)(d-1-m).
$$

\noindent Comme $m$ est supposée supérieure à $d-1$, la différence ci-dessus est toujours négative et donc la distance construite de Goppa fournit une meilleure minoration de $d^{\bot}$.
\end{itemize}

\paragraph{Conclusion.} Dans ce contexte des courbes planes,
les techniques développées en section \ref{ssecminor}
fournissent une meilleure minoration de
la distance minimale de $C_L(D,G)^{\bot}=C_{\Omega}(D,G)$
que celle fournie par la distance construite
de Goppa si et seulement si $$m\leq d-2.$$

\subsubsection{Surfaces de $\P^3$}

Soit $S$ une surface de $\P^3$ de degré $d$ définie sur $\F_q$. Soient également
 $m$ un entier naturel, $G$ un diviseur sur $S$ tel que $G\sim mL_S$ et $\Delta$ un $0$-cycle de la forme $\Delta=P_1+ \cdots +P_n$ où les $P_i$ sont des points rationnels de $S$ qui évitent le support de $G$. 
On note de nouveau $d^{\bot}$, la distance minimale du code $C_{L,S}(\Delta, G)^{\bot}$.
Tout comme dans le paragraphe précédent, le théorème \ref{minor} fournit les minorations suivantes.
\begin{enumerate}
\item[$(i)$] Pour tout $m$, on a $d^{\bot}\geq m+2$.
\item[$(ii)$] Si de plus $m \geq d-1$ et que $S$ ne contient pas de droite rationnelle, alors $d^{\bot}\geq 2m+2$.
\end{enumerate}

\begin{exmp}\label{excub}
  Soit $S$, une surface cubique lisse de $\P^3$. Soit $L$ un diviseur sur $S$ donné par une section hyperplane de $S$ et $G:=mL$ avec $m\in \N$.
On choisit enfin comme $0$-cycle $\Delta$, la somme des points rationnels de $S$ qui évitent le support de $G$.
On note $d^{\bot}(m)$ la distance minimale du code $C_{L,S}(\Delta, mL)^{\bot}$. Les résultats de la section \ref{ssecminor} nous donnent
$$
\begin{array}{ccll}
  d^{\bot}(1) & \geq & \ \ \ \ \!3; & \\
  d^{\bot}(2) & \geq & \left\{
    \begin{array}{cl}
     4 & \textrm{si}\ S\ \textrm{contient une droite rationnelle};\\
     6 & \textrm{sinon.} 
    \end{array}
\right. 
\end{array}
$$

\noindent Dans le premier cas la borne est atteinte seulement si le support de $\Delta$ contient trois points alignés.
Nous verrons au chapitre \ref{chapldpc} que ce phénomène est très fréquent et que ce code possède en général de nombreux mots de poids $3$.
Dans le dernier cas, la borne inférieure n'est atteinte que si le support de $\Delta$ contient $6$ points appartenant à une même conique plane (éventuellement réductible).
\end{exmp}

\begin{rem}
  Remarquons que la classification des surfaces cubiques lisses réalisée par Swinnerton-Dyer dans \cite{swd} assure l'existence de cubiques ne contenant pas de droites rationnelles. 
Cela fait d'ailleurs partie des exemples introduits par Zarzar et Voloch dans \cite{agctvoloch}.
\end{rem}

\begin{exmp}
On reprend les mêmes notations que dans l'exemple \ref{excub}, mais cette fois, $S$ est une surface lisse de degré $4$. On obtient alors les minorations
$$
\begin{array}{rccll}
 (i) & d^{\bot}(1) & \geq & \ \ \ \ \!3; & \\
 (ii) &  d^{\bot}(2) & \geq & \ \ \ \ \!4; & \\
 (iii) & d^{\bot}(3) & \geq & \left\{
    \begin{array}{cl}
     5 & \textrm{si}\ S\ \textrm{contient une droite rationnelle et}\ \sharp\ \F_q\geq 5;\\
     8 & \textrm{sinon.} 
    \end{array}
\right. 
\end{array}
$$ 

\noindent Dans le cas $(i)$ (resp. $(ii)$), la borne est atteinte seulement si $S$ contient $3$ (resp. $4$) points alignés.
Dans le dernier cas, la borne n'est atteinte que si le support de $\Delta$ contient $8$ points appartenant à une même conique plane.
\end{exmp}

\begin{rem}
En ce qui concerne le cas $(iii)$ de l'exemple précédent, dans \cite{sch1} théorème I.6.9, on montre qu'une surface cubique contient toujours au moins une droite géométrique (elle en contient même $27$ quand elle est lisse) et qu'une surface générique de degré supérieur à $4$ ne contient pas de droite géométrique. Ainsi, en général, si $S$ est une surface de degré $4$, la distance minimale de $C_{L,S}(\Delta, 3L)$ est supérieure ou égale à $8$. 
\end{rem}

\section{Seconde approche, un problème ouvert}

Dans cette section nous revenons au contexte classique, à savoir celui des surfaces algébriques.
Cette \textit{seconde approche} ne fournira pas à proprement parler de minoration de la distance minimale de l'orthogonal d'un code fonctionnel.
Il ne s'agit donc pas d'une section réellement aboutie, mais d'une ouverture vers des problèmes de géométrie algébrique qui auraient d'intéressantes applications à la théorie des codes correcteurs d'erreurs.

L'objectif étant d'utiliser les résultats du chapitre \ref{chapreal}, nous allons nous replacer dans le contexte de ce dernier.
À savoir, $S$ désigne une surface projective lisse géométriquement intègre, qui est intersection complète dans un espace projectif $\P^r$.
On se donne également un entier naturel $m$ et un diviseur $G$ vérifiant $G\sim mL_S$, où $L_S$ désigne la classe d'équivalence linéaire d'une section de $S$ par un hyperplan de $\P^r$.
Enfin, $P_1, \ldots, P_n$ désignent des points rationnels de $S$ et $\Delta$ leur somme.

\subsubsection{Exploitation du théorème de réalisation}

D'après le théorème de réalisation (théorème \ref{thmreal}), pour tout mot de code $c\in C_L(\Delta, G)^{\bot}$, il existe un couple sous-$\Delta$-convenable de diviseurs $(D_a,D_b)$ et une $2$-forme $\omega$ vérifiant

\begin{equation}\label{omeho}
(\omega)=G-D_a-D_b\quad \textrm{et telle que}
\quad
c=\res^2_{D_a, \Delta}(\omega).
\end{equation}

\noindent De plus, le diviseur $D_a$ est une courbe lisse géométriquement intègre provenant de l'intersection de $S$ avec une hypersurface de $\P^r$. Il existe donc un entier naturel $n_a$ tel que $D_a \sim n_a L_S$.
Par conséquent, dans ce qui suit nous nous autoriserons l'abus de langage consistant à désigner par $D_a$ à la fois le diviseur et la courbe irréductible sous-jacente.
Enfin, le théorème de réalisation affirme qu'il existe un entier relatif $n_b$ tel que $D_b \sim n_b L_S$.

La $2$-forme $\omega$ est de valuation supérieure ou égale à $-1$ le long de la courbe $D_a$, le $1$-résidu de $\omega$ le long de $D_a$ est donc bien défini. On pose
$$\mu:= \res^1_{D_a}(\omega) \in \Omega^1_{\F_q (D_a)/ \F_q}$$

\noindent et d'après le lemme \ref{valres}, pour tout point géométrique $P$ de $D_a$, on a
$$
\val_P (\mu)=m_P (D_a, G-D_b).
$$

\noindent Par conséquent, si l'on note $D^*$ le tiré en arrière sur $D_a$ d'un diviseur $D$ sur $S$ dont le support ne contient pas $D_a$, alors le diviseur de $\mu$ s'écrit
\begin{equation}\label{divimu}
(\mu)=G^*-D_b^*.
\end{equation}

\noindent Soit $\Lambda_c$ le diviseur sur la courbe $D_a$ correspondant au support du mot de code $c$, on a 
$$
(\mu)\geq G^*-\Lambda_C
$$
et on déduit de cette inégalité et de (\ref{divimu}) que
$
\Lambda_C \geq D_b^*
$.
On a de plus,
$$
w(c)=\deg (\Lambda_C)\quad \textrm{et} \quad \deg (D_b^*) =D_a.D_b,
$$

\noindent où $w(.)$ désigne le poids de Hamming d'un mot de code et $D.D'$ le produit d'intersection de deux diviseurs sur $S$.
On en déduit la relation
\begin{equation}\label{relafond}
w(c) \geq D_a.D_b =n_an_b L_S^2.
\end{equation}

\noindent Soient $k$ et $m$ les entiers tels que
$$
K_S \sim kL_S \quad \textrm{et} \quad G \sim mL_S,
$$

\noindent où $K_S$ désigne la classe canonique sur $S$. D'après la relation (\ref{omeho}) page \pageref{omeho}, on a
$$
n_b=m-n_a-k.
$$ 

\noindent Si l'on injecte cette relation dans (\ref{relafond}), on obtient
\begin{equation}\label{salope}
  w(c)\geq n_a (m-k-n_a)L_S^2.
\end{equation}

\noindent Le souci est que la quantité avec laquelle on minore le poids de Hamming de $c$ est négative dès que $n_a\geq m-k$.
Partant de la discussion ci-dessus, 
le théorème \ref{pipeau} qui suit n'est pas réellement exploitable en l'état. 
Il offre cependant des perspectives de minoration de la distance minimale de $C_{L,S}(\Delta,G)^{\bot}$ sous réserve d'obtenir des réponses à une question ouverte qui sera posée plus loin (question \ref{Qpipeau}).
La preuve de ce théorème est suivie d'une discussion sur l'énoncé.

\begin{thm}\label{pipeau}
Soit $S$ une surface projective lisse intersection complète et géométriquement intègre. Soit $m$ un entier et $G$ un diviseur sur $S$ vérifiant $G \sim mL_S$. Soit enfin $\Delta$ une somme formelle de points rationnels de $S$ évitant le support de $G$.
Supposons qu'il existe un entier naturel $s$ vérifiant les conditions suivantes.
  \begin{enumerate}
  \item[$(i)$] $s$ est supérieur à la distance minimale $d^{\bot}$.
  \item[$(ii)$] Pour toute configuration $P_{i_1}, \ldots, P_{i_s}$ de points du support de $\Delta$, il existe une hypersurface $H$ de $\P^r$ \textbf{définie sur $\overline{\F}_q$} de degré inférieur à $m-k-1$ qui contient $P_{i_1}, \ldots, P_{i_s}$ et telle que $H \cap S$ est une courbe lisse. On rappelle que l'entier $k$ est celui qui vérifie $K_S \sim kL_S$ où $K_S$ désigne la classe canonique.
\end{enumerate}

\noindent Alors, la distance minimale $d^{\bot}$ du code $C_L(\Delta, G)^{\bot}$ vérifie
$$
d^{\bot} \geq (m-k-1)L_S^2.
$$
\end{thm}

\begin{proof}
D'après la remarque \ref{clotalg} page \pageref{clotalg}, la distance minimale d'un code géométrique est invariante par extension des scalaires.
Il suffit donc que l'on soit à même de réaliser géométriquement les mots de code de $C_{L,S}(\Delta, G)^{\bot}$ vus comme des mots de $C_{L, S'}(\Delta',G')$, où $S'$ désigne $S\times_{\F_q} \F_{q^l}$ pour un certain entier naturel $l$ et $\Delta'$ et $G'$ les tirés en arrière respectifs de $\Delta$ et $G$ sur $S'$.
C'est ce qui justifie dans le $(ii)$ le ``\textit{définie sur $\overline{\F}_q$}''.

Soient $s$ un entier vérifiant les condition $(i)$ et $(ii)$ de l'énoncé et $E$, l'ensemble des mots non nuls de $C_{L,S}(\Delta, G)^{\bot}$ de poids de Hamming inférieur ou égal à $s$. D'après $(ii)$ et (\ref{salope}), on a 
$$
\forall c\in E,\ w(c)\geq \min_{n_a=1}^{m-k-1} n_a(m-k-n_a)L_S^2. 
$$

\noindent La définition de $E$ entraîne que l'inégalité ci-dessus est en fait vérifiée par tous les mots non nuls de $C_{L,S}(\Delta,G)^{\bot}$. Par ailleurs, l'étude de la fonction $x\mapsto x(m-k-x)$ sur l'intervalle $[1, m-k-1]$ permet de voir que le minimum de l'expression $n_a(m-k-n_a)L_S^2$ est atteint pour $n_a=1$. On en déduit
$$
d^{\bot} \geq (m-k-1)L_S^2.
$$
\end{proof}

\paragraph{Discussion au sujet du théorème \ref{pipeau}.}
L'énoncé peut sembler troublant en ce sens où l'entier $s$ doit être supérieur à la quantité $d^{\bot}$ sur laquelle on cherche à s'informer.
Pour exploiter ce théorème il faut disposer d'une majoration \textit{à priori} de $d^{\bot}$. Voici un certain nombre de pistes, pour obtenir une telle majoration.
\begin{enumerate}
\item L'approche la plus naïve serait de majorer $d^{\bot}$ par la longueur du code.
\item Si l'on connaît la dimension du code on peut utiliser la borne de singleton  à savoir $d^{\bot}\leq n-\dim (C_{L,S}(\Delta,G)^{\bot})+1$.
\item Sous réserve de disposer d'une évaluation de la distance minimale d'un code fonctionnel, le théorème d'orthogonalité (théorème \ref{orthocode}) fournit une majoration de la distance minimale de $C_L(\Delta,G)^{\bot}$. 
En effet, comme pour toute paire de diviseurs (sous-) $\Delta$-convenable $(D_a,D_b)$ on a 
$$
C_{\Omega,S}(\Delta,D_a,D_b,G) \subseteq C_{L,S}(\Delta,G)^{\bot}.
$$
Par conséquent, la distance minimale du code $C_{L,S}(\Delta,G)^{\bot}$ est inférieure à celle du code $C_{\Omega,S}(\Delta,D_a,D_b,G)$.
Ce dernier code est fonctionnel d'après le théorème \ref{diff=fonc}. Si l'on et capable d'estimer la distance minimale d'un code fonctionnel sur $S$ on peut en déduire une bonne majoration à priori de $d^{\bot}$.
\end{enumerate}

Il s'agit là d'un piste à explorer dans le futur.
Notons tout de même que dans certaines situations,
 un tel entier $s$ n'existe pas, il suffit par exemple que $m-k-1$ soit négatif.
Avant de conclure, donnons deux exemples élémentaires. L'un assurant que l'entier $s$ existe (au moins dans certaines situations élémentaires) et le second présentant un cas où l'entier $s$ n'existe pas, bien que $m-k-1$ soit strictement positif.

\begin{exmp}[Un exemple où $s$ existe.]
On reprend l'exemple présenté dans la section \ref{P2exdelta} du chapitre \ref{chapdiff}.
  Soit $S$ le plan projectif sur un $\F_q$ de caractéristique différente de deux. Supposons par exemple que $m=1$. On a $k=-3$ et le code $C_L(\Delta,G)$ avec $G\sim L_{S}$ est de dimension $3$ et donc de longueur $q^2$.
On en déduit que le code $C_L(\Delta,G)^{\bot}$ est de dimension $q^2-3$ et la borne de singleton nous assure que sa distance minimale vérifie
$$
d^{\bot}\leq 4
$$
On a $m-k-1=3$, l'objectif est donc de montrer que pour tout quadruplet de points rationnels de $\P^2$, il existe une courbe lisse de degré inférieur ou égal à trois qui contient ces points.

Soient donc $P_1, P_2, P_3, P_4$ quatre points de $\P^2$. S'ils sont alignés c'est terminé, une droite étant une courbe lisse de degré $3$. Si $P_1, P_2, P_3$ appartiennent à une même droite $L$ qui ne contient pas $P_4$, il existe un système de coordonnées homogènes $(X,Y,Z)$ sur $\P^2$ et un élément $a\in \F_q \smallsetminus \{0,1\}$ tels que,
$$
\begin{array}{rclcrcl}
  P_1 & = & (0:0:1) & \quad & P_1 & = & (1:0:1)\\
  P_3 & = & (a:0:1) & \quad & P_4 & = & (0:1:0).
\end{array}
$$
La courbe elliptique d'équation 
$
Y^2Z=X(X-Z)(X-aZ)
$
est lisse et interpole ces quatre points.

Enfin, supposons que trois de ces points ne soient pas alignés, alors il existe au moins une conique lisse qui les contient.
L'entier $s=4$ vérifie donc les conditions $(i)$ et $(ii)$ du théorème \ref{pipeau}, ce dernier nous fournit donc une minoration de la distance minimale du code $C_{L,S}(\Delta,G)^{\bot}$, à savoir
$$
d^{\bot}\geq (m-k-1)L_S^2=3.
$$ 
D'après le théorème \ref{minor} (\ref{minor1}), le minorant obtenu est en fait exactement la distance minimale du code étudié.
\end{exmp}

\begin{exmp}[Un exemple où $m-k-1 \geq 0$ mais $s$ n'existe pas.]
Soit $S$ une surface cubique lisse de $\P^3_{\F_q}$ avec $q\geq 3$ et contenant au moins une droite rationnelle. Soit $\Delta$ la somme de tous les points rationnels de $S$ et $G$ un diviseur tel que $G\sim 2L_S$ et dont le support évite celui de $\Delta$ (un tel $G$ existe, voir annexe \ref{idem}).
On rappelle que les surfaces cubiques sont des surfaces de Del Pezzo. Donc $k=-1$ et $m-k-1=2$. Notons que, d'après le théorème \ref{minor} (\ref{minor2}), on sait que $d^{\bot}=4$, les mots de poids minimal étant ceux donc le support correspond à des points d'une droite rationnelle contenue dans $S$.

À présent, raisonnons par l'absurde,
 en supposant l'existence d'un entier $s$ vérifiant les conditions du théorème \ref{pipeau}. Alors, d'après ce théorème, la distance minimale de $C_{L,S}(\Delta, G)^{\bot}$ serait supérieure ou égale à $6$, ce qui est faux d'après les remarques ci-dessus.
\end{exmp}

\paragraph{Conclusion.} Le théorème \ref{pipeau} motive la question ouverte suivante.

\begin{ques}\label{Qpipeau}
  Soient $X$ une sous-variété irréductible lisse géométriquement intègre de $\P^r_{\overline{\F}_q}$
et $d$ un entier naturel.
Soient $P_1, \ldots, P_n$ une famille de points de $X$.
Sous quelles conditions sur $X$ et $P_1, \ldots , P_n$ a-t-on l'existence d'un entier $s$ tel que pour tout $s$-uplet de points parmi $P_1, \ldots, P_n$, il existe une hypersurface $H$ de degré $d$ contenant ce $s$-uplet de points et telle que $H\cap X$ soit une sous-variété lisse de codimension $1$ de $X$?
\end{ques}

Notons qu'une réponse à la question \ref{qpoon}G posée page \pageref{qbert}, fournirait sans doute des éléments de réponse, voire même une réponse complète à la question ci-dessus.
L'obtention d'un tel résultat ``à la Bertini'' nous donnerait de nombreuses informations, à la fois sur les codes fonctionnels et sur leurs orthogonaux. La question \ref{qpoon}G est donc un problème ouvert ouvrant de nombreuses perspectives d'application.




\newpage
\thispagestyle{empty}

\chapter{Constructions de mots de faible poids et codes LDPC}\label{chapldpc}

\begin{flushright}
\begin{tabular}{p{8cm}}
{\small Une idée reçue atteste que la creuse est le département le moins peuplé de France, ce qui est totalement faux.}
\medbreak
\hfill {\small \textsc{Wikipedia}} 
\end{tabular}
\end{flushright}

On signale dans la section \ref{discuss} du chapitre \ref{chapreal} que le théorème de réalisation n'est pas constructif. Aussi, ce chapitre est-il en partie consacré à la présentation de méthodes constructives de réalisation différentielle de mots de code appartenant à l'orthogonal d'un code fonctionnel.
Les mots de l'orthogonal d'un code fonctionnel qui vont nous intéresser et qui s'avèreront être les plus simples à calculer seront ceux dont le poids de Hamming est \textit{petit}.
Si ces mots engendrent le code qui les contient, on dit que ce code est LDPC (\textit{Low Density Parity Check}).
La première section de ce chapitre est une introduction à la théorie de ces codes.

\section{Introduction aux codes LDPC}\label{introldpc}

Un code LDPC est un code admettant une matrice de parité \textit{creuse}. En d'autres termes, c'est un code admettant une base duale composée de mots de petit poids de Hamming.

\subsection{Graphe de Tanner}  

\begin{defn}[Graphe biparti]
Un graphe biparti est la donnée de deux ensembles de sommets $V_1$ et $V_2$ et d'un ensemble d'arêtes $E$ tels que toute arête $a\in E$ relie un unique élément de $V_1$ avec un unique élément de $V_2$.  
\end{defn}

La définition qui suit a été introduite par R. Michael Tanner dans \cite{tanner}.

\begin{defn}[Tanner 1981]
Soient $C$ un code binaire de longueur $n$ et $H\in \mathfrak{M}_{r,n}(\F_2)$ une matrice de parité de $C$.
On appelle graphe de Tanner de $C$, le graphe biparti dont la première famille de sommets $V_1$ est indexée par les colonnes de $H$ et la seconde famille $V_2$ par les lignes. Une arête relie le $i$-ème sommet de la famille $V_1$ au $j$-ème de $V_2$ si et seulement si le coefficient $h_{i,j}$ de la matrice $H$ est non nul.
\end{defn}

\begin{rem}
  Remarquer qu'un code n'admet pas un unique graphe de Tanner. Aussi, on devrait parler du graphe de Tanner de $C$ associé à $H$ et non du graphe de Tanner de $C$. Dans la pratique, cet abus de langage est toléré et même fréquemment pratiqué.
\end{rem}

Dans ce qui suit, on représentera les sommets correspondant aux colonnes de la matrice par des \includegraphics[width=.4cm, height=.4cm]{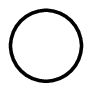} et on appellera ces sommets les \textit{n\oe uds de données} ou tout simplement les \textit{bits}. 
Les sommets correspondant aux lignes seront représentés par des
\includegraphics[width=.4cm, height=.4cm]{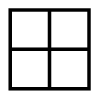} et on les appellera les \textit{n\oe uds de parité} ou les \textit{relations}\footnote{Dans la littérature anglophone, on parle de \textit{check nodes}, c'est-à-dire n\oe ud de contrôle. Nous avons préféré donner ce nom de \textit{relation}, car ces n\oe uds symbolisent une équation, donc une relation entre les bits qui lui sont voisins dans le graphe.}.

\begin{exmp}\label{exmptanner}
Soit $C$, le code de matrice de parité
$$
H=
\left(
\begin{array}{ccccccccc}
1 & 1 & 1 & 0 & 0 & 0 & 0 & 1 & 0 \\
0 & 0 & 1 & 1 & 1 & 0 & 0 & 0 & 0 \\
0 & 0 & 0 & 0 & 1 & 1 & 1 & 0 & 0 \\
0 & 0 & 0 & 0 & 1 & 0 & 0 & 1 & 1 
\end{array}
\right).
$$  

\noindent Le graphe de Tanner de $C$ correspondant à la matrice $H$ est de la forme suivante.

\begin{center}
\includegraphics[width=14cm, height=3.55cm]{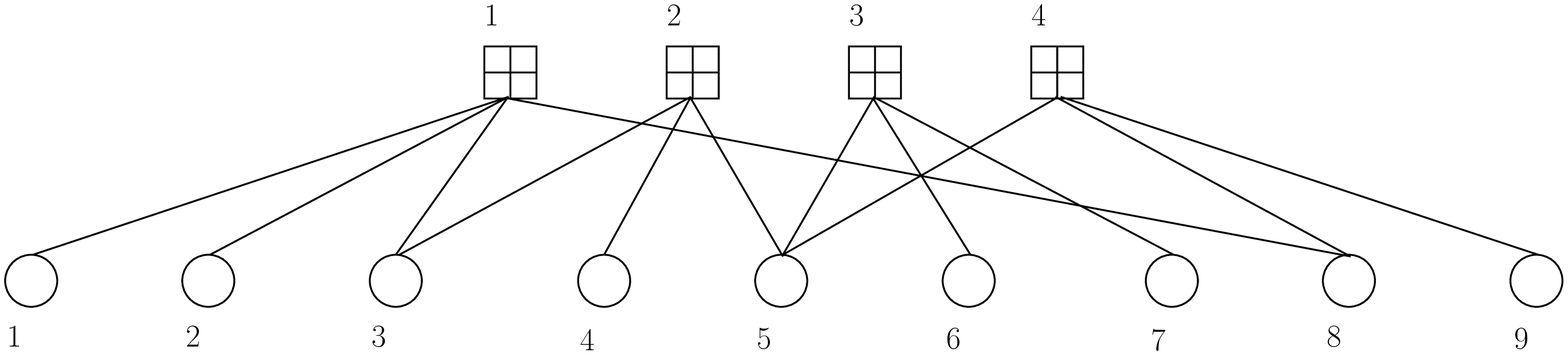}
\end{center}

\noindent On peut également essayer de l'\textit{étaler} afin d'y voir plus clair, sous réserve bien sûr que le graphe admette une représentation planaire.

\begin{center}
\includegraphics[width=10cm, height=6.1cm]{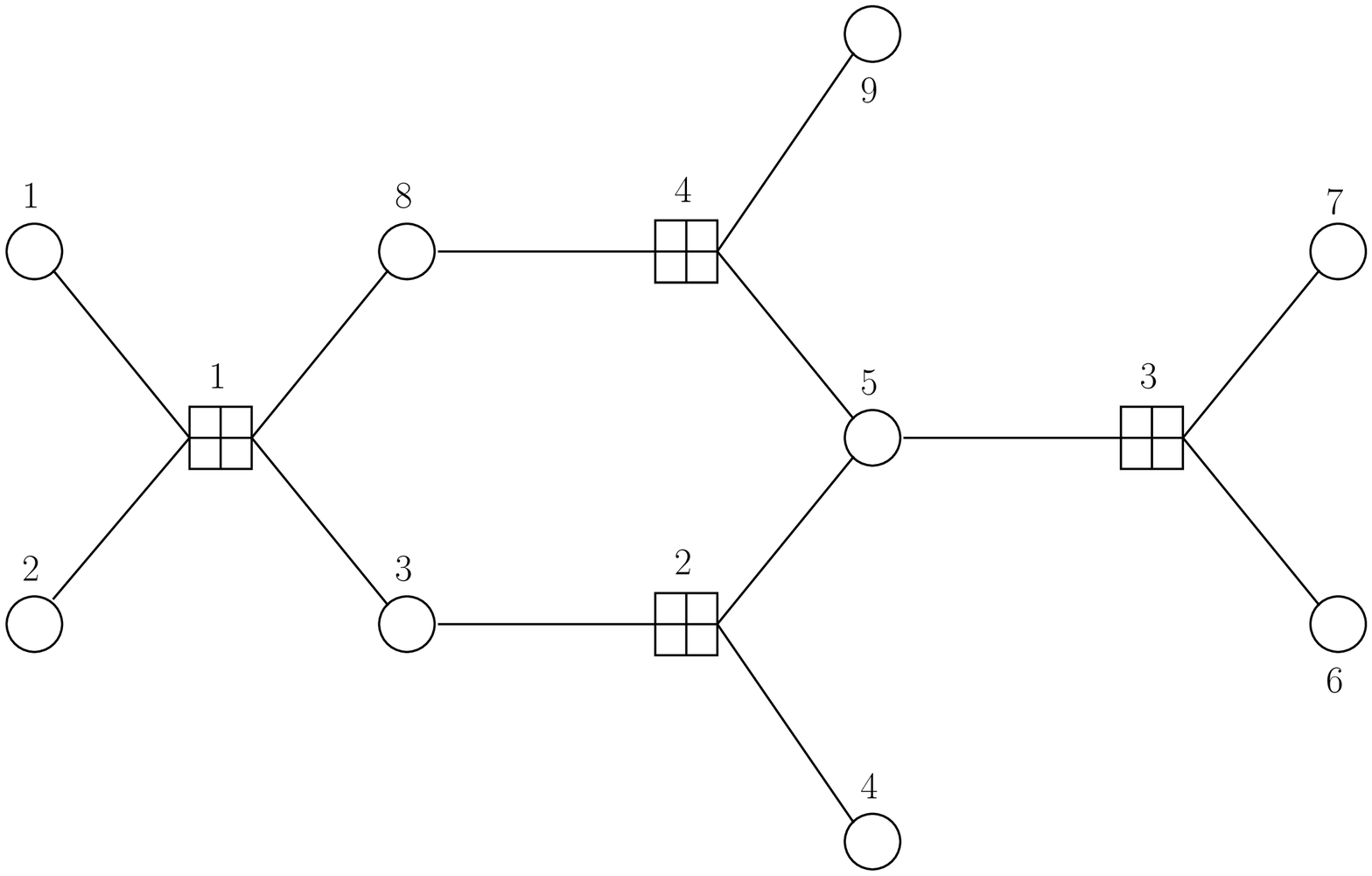}
\end{center}
\end{exmp}

Si le code n'est pas binaire on peut réaliser une construction semblable mais avec des arêtes \textit{pondérées}. Entre le bit $i$ et la relation $j$ on trace une arête pondérée par le coefficient $h_{i,j}$ de la matrice de parité si ce dernier est non nul et pas d'arête sinon. 

\begin{exmp}
Supposons que le corps de base soit $\F_5$ et considérons le code $C$ de matrice de parité 
$$
H=\left(
\begin{array}{ccccccccc}
2 & 0 & 3 & 1 & 0 & 0 & 0 & 0 & 0 \\
0 & 1 & 0 & 4 & 3 & 0 & 0 & 0 & 0 \\
0 & 0 & 3 & 0 & 0 & 1 & 2 & 0 & 0 \\
0 & 0 & 0 & 4 & 0 & 0 & 4 & 1 & 0 \\
0 & 0 & 0 & 0 & 1 & 0 & 0 & 2 & 1
\end{array}
\right).
$$

\noindent Le graphe de Tanner de $C$ associé à $H$ se représente de la façon suivante.

\begin{center}
\includegraphics[width=11cm, height=7.16cm]{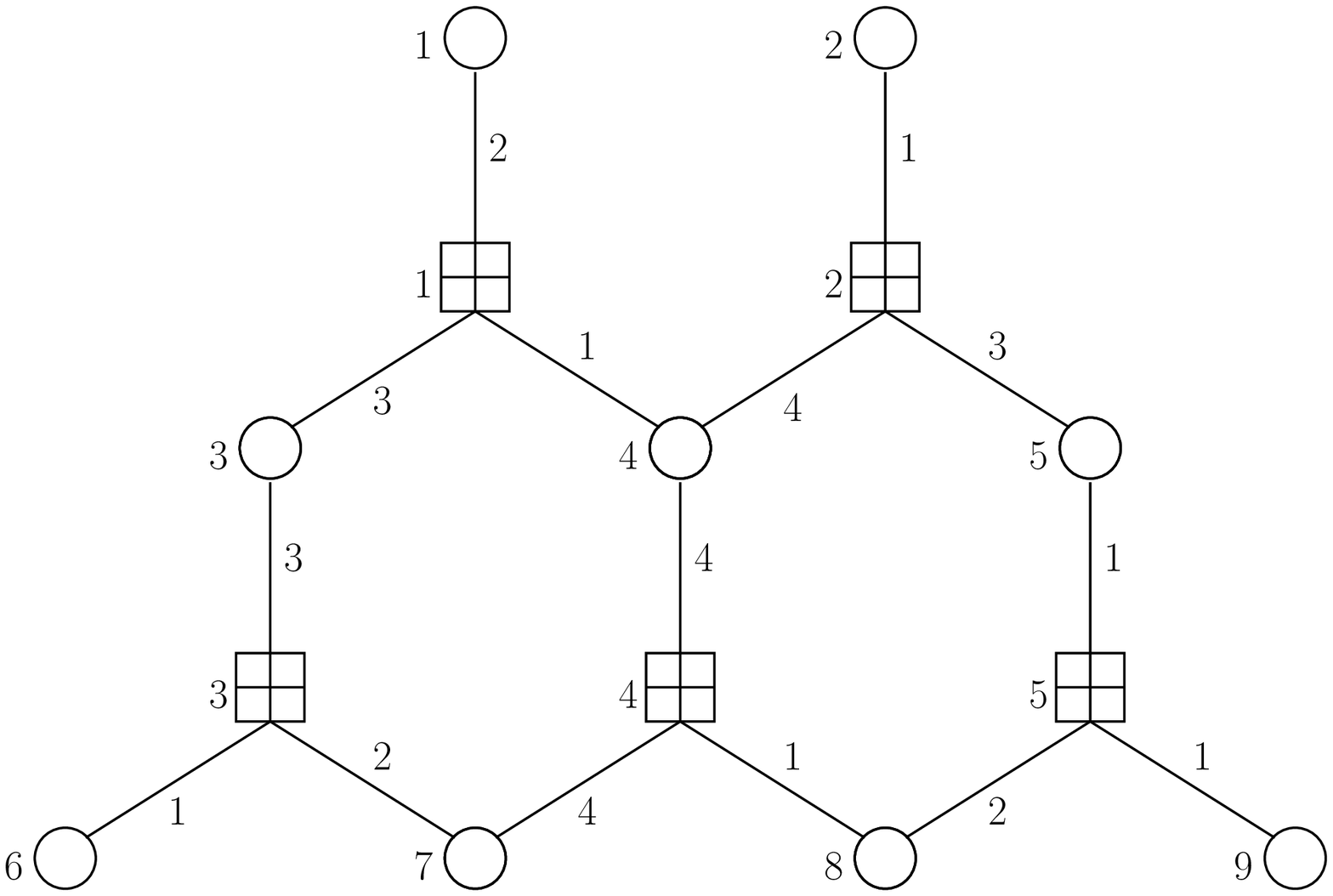}
\end{center} 
\end{exmp}

\subsection{Décodage itératif}

L'intérêt majeur de la représentation d'un code par un graphe de Tanner est le décodage itératif.
Le principe général consiste, étant donné un mot de code reçu $y$ à évaluer les coûts locaux d'assignation de chaque bit à une valeur prescrite.
Dans un second temps, par un principe de passage de messages dans le graphe, on actualise ces coûts en fonction du nombre de modifications qu'une assignation  d'un bit à une valeur prescrite entraînerait sur les bits voisins dans le graphe. 
De façon schématique, la répétition de ce procédé permet (à de nombreux détails près) de passer de coûts locaux à des coûts globaux. On choisit alors comme sortie de l'algorithme, le mot de code correspondant aux coûts globaux minimaux.

Il existe dans la littérature de nombreux algorithmes de décodage itératif. Celui que nous allons présenter porte en général le nom de \textit{algorithme min-somme}. Notons que l'on peut trouver une excellente présentation de cet algorithme dans la thèse de Niclas Wiberg \cite{wiberg}.

\paragraph{Description à partir d'un exemple.}

Le mécanisme d'un algorithme de décodage itératif, sans être très complexe, est relativement technique. Nous allons commencer par le décrire à l'aide d'un exemple élémentaire.

\begin{exmp}\label{exminsom}
Considérons le code binaire $C$ de matrice de parité
$$
H=\left(
\begin{array}{cccccc}
1 & 1 & 0 & 1 & 0 & 0 \\
0 & 1 & 1 & 0 & 1 & 0 \\
0 & 0 & 0 & 1 & 1 & 1
\end{array}
\right).
$$

\noindent Son graphe de Tanner associé à $H$ se présente sous la forme suivante.

\begin{center}
\includegraphics[height=6.5cm, width=8cm]{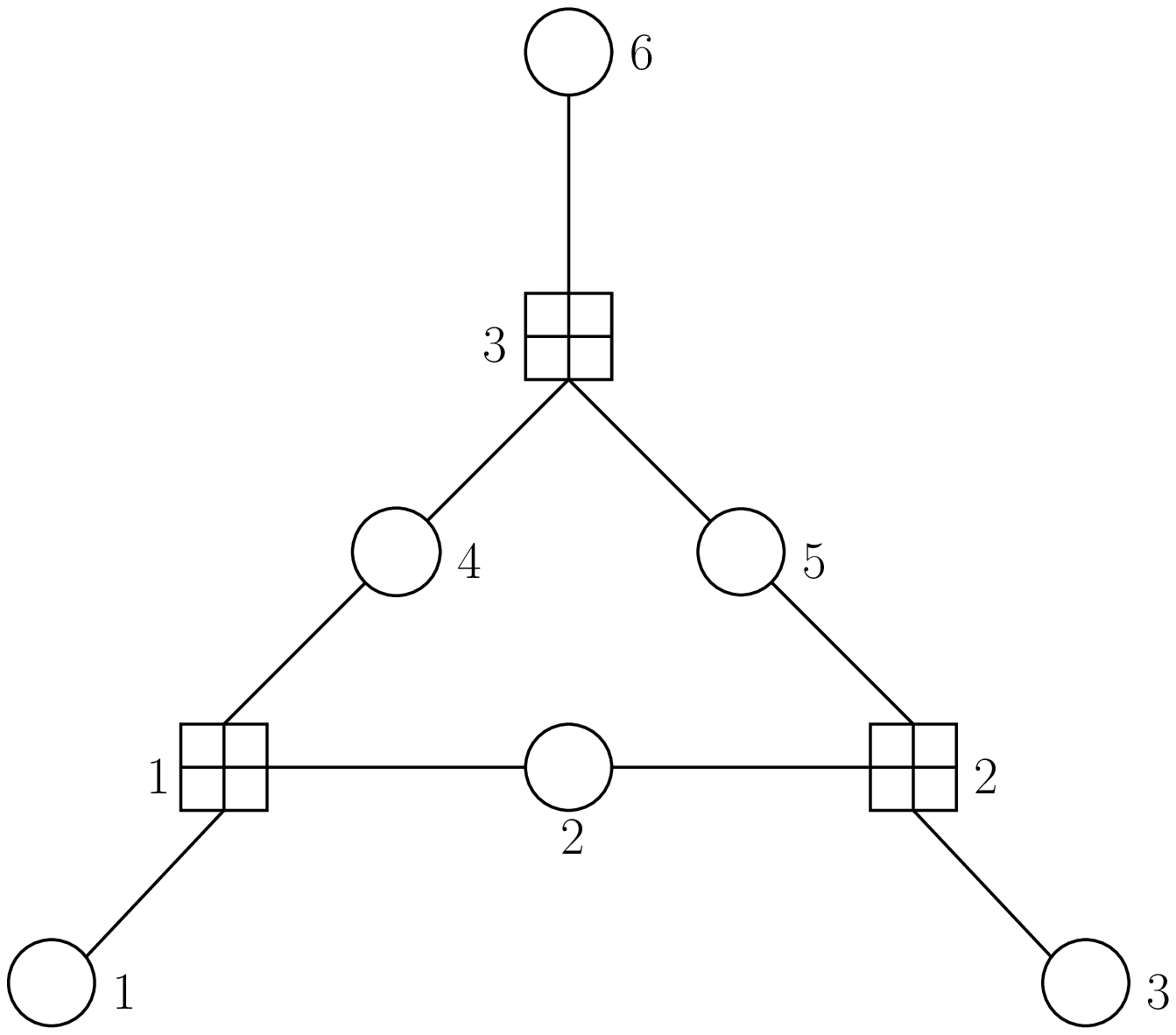}
\end{center}

\noindent Ce code est de distance minimale $3$, on peut donc corriger une erreur.
\end{exmp}

\noindent Considérons ce code $C$ et supposons que l'on ait reçu le mot
$$ 
y=(1,1,0,1,0,1).
$$
Le second bit est erroné.

\smallbreak

\noindent \textbf{Étape 1.} L'objectif de l'algorithme est d'évaluer le coût qu'aurait l'assignation d'un bit à une valeur prescrite dans $\F_2$. Pour ce faire, on commence par définir pour chaque bit une fonction de \textit{coût local}. C'est une fonction $C^{i}_{\textrm{loc}}: \F_2 \rightarrow \N$ telle que $C^i_{\textrm{loc}}(\alpha)$ est nul si $y_i=\alpha$ et égal à $1$ sinon. Par exemple $C^{1}_{\textrm{loc}}(0)=1$ et 
$C^{1}_{\textrm{loc}}(1)=0$.

\noindent Pour abréger, on note
$$
C^{1}_{\textrm{loc}}=[1,0].
$$ 

\noindent Cette fonction quantifie le nombre de changements qu'impliquerait l'assignation du bit $i$ à la valeur $\alpha$ sans chercher à vérifier les relations de parité. On voit facilement que 
$$
\begin{array}{rclcrclcrcl}
C^{1}_{\textrm{loc}} & = & [1,0] & \ & C^{2}_{\textrm{loc}} & = & [1,0] &
 & C^{3}_{\textrm{loc}} & = & [0,1] \\
C^{4}_{\textrm{loc}} & = & [1,0] & \ & C^{5}_{\textrm{loc}} & = & [0,1] &
 & C^{6}_{\textrm{loc}} & = & [1,0].
\end{array}
$$

\smallbreak

\noindent \textbf{Étape 2.} Dans un second temps, on va évaluer le nombre de changements qu'impliquerait l'assignation d'un bit à une valeur prescrite avec la contrainte de respecter les relations de parité voisines de ce bit. 

\medbreak

\noindent \textbf{Étude locale.} Focalisons nous sur le second bit. Il est voisin de deux relations de parité: la première et la seconde.
Supposons que l'on lui assigne la valeur $0$.
Alors, pour respecter la première équation de parité, on a deux possibilités.
\begin{enumerate}
\item Les bits $1$ et $4$ prennent tous deux la valeur $1$.
\item Les bits $1$ et $4$ prennent tous deux la valeur $0$.
\end{enumerate}
La première configuration est la moins coûteuse, elle n'implique aucun changement, c'est celle que l'on retient.
On en déduit que l'assignation du second bit à la valeur $0$ aura une répercussion de coût nul sur les autres bits voisins du premier n\oe ud de relation.
Si maintenant on assigne la valeur $1$ à ce bit on a également deux possibilités.
\begin{enumerate}
\item Le $1^{\textrm{er}}$ bit prend la valeur $1$ et le $4^{\textrm{e}}$ la valeur $0$.
\item Le $1^{\textrm{er}}$ bit prend la valeur $0$ et le $4^{\textrm{e}}$ la valeur $1$.
\end{enumerate}
Les deux configurations coûtent un changement.
On en déduit que l'assignation du second bit à la valeur $0$ a une répercussion de coût $1$ sur le premier n\oe ud de relation.

De la même manière, on montre que l'assignation du second bit à la valeur $0$ (resp. $1$) a une répercussion de coût $0$ (resp. $1$) sur les bits $3$ et $5$ voisins du second n\oe ud de relation.
\smallbreak

\noindent \textbf{Étape 3.} Au final l'assignation du second bit à la valeur $1$ coûte deux changements (un sur le bit $1$ ou $4$ et un autre sur le $2$ ou $5$) comme le sixième. Par contre l'assignation de ce bit à $0$ ne coûte qu'un seul changement, celui qui consiste à remplacer ce bit initialement à la valeur $1$ par un $0$. On est donc tentés d'assigner ce bit à $0$ ce qui corrige l'erreur.

\begin{rem}
Dans cet exemple, nous nous sommes focalisés sur un seul bit pour tenter de comprendre le mécanisme. En réalité, l'algorithme réalise en parallèle la même démarche pour chaque bit.
\end{rem}

\begin{rem}
Il est important de remarquer que nous n'avons pas vérifié si l'assignation du second bit à une valeur donnée avait des répercussions sur les bits plus éloignés.
On s'est limités au premier voisinage du second bit pour prendre notre décision. 
En général, on réitère le processus décrit dans l'étape 2 de façon à obtenir des informations sur les répercussions d'une assignation sur les bits \textit{éloignés}.
\end{rem}

\medbreak

\noindent L'idée de l'algorithme min-somme peut se résumer de la façon suivante.
\begin{enumerate}
\item On commence par compter le nombre de changements qu'impliquerait l'assignation du $i$-ème bit à une valeur donnée, sans tenir compte des relations de parité.
\item On compte ensuite le nombre de changements que cela impliquerait pour les autres bits reliés à $i$ par une relation de parité. C'est-à-dire le coût d'une telle assignation pour les bits étant dans le premier voisinage du $i$-ème bit.
\item En réitérant ce procédé on peut compter le nombre de changements qu'implique une telle assignation pour les bits appartenant au second voisinage du $i$-ème bit.
\item On réitère le processus...
\item Lorsque l'on dispose du coût d'assignation du $i$-ème bit à une valeur donnée pour un voisinage \textit{suffisamment grand} de ce dernier, on prend une décision sur la valeur à laquelle on l'assigne en choisissant bien sûr celle qui est la moins coûteuse.
\end{enumerate}

\noindent Encore une fois, les opérations sont réalisées en parallèle pour tous les bits.
 
\subsection{L'algorithme min-somme}

Nous allons à présent donner une description générale et rigoureuse de l'algorithme min-somme.

\medbreak

\noindent \textbf{Étape 1. Initialisation.}
À chaque bit, on associe une fonction de coût local $C_{\textrm{loc}}^i:\F_q \rightarrow \N$.
À l'état initial, la fonction de coût local du $i$-ème bit est extrêmement simple. Si la $i$-ème coordonnée du mot reçu $y$ est égale à $\alpha \in \F_q$, alors la fonction $f_i$ prend la valeur $0$ en $\alpha$ et la valeur $1$ en tous les autres éléments de $\F_q$. La valeur $C_{\textrm{loc}}^i(\beta)$ quantifie le coût d'assignation du $i$-ème bit à la valeur $\beta$ sans tenir compte des bits voisins.

Pour toute arête $(i,j)$ du graphe de Tanner, on définit les fonctions \textit{messages} $\mu_{i \rightarrow j}: \F_q \rightarrow \N$ et $\nu_{i \leftarrow j}: \F_q \rightarrow \N$. 
Ces fonctions peuvent être vues respectivement comme un message allant du bit $i$ vers la relation $j$ et réciproquement.
Ces fonctions sont des \textit{variables locales} de l'algorithme, c'est-à-dire qu'elles sont actualisées à chaque itération de l'algorithme.
Pour toute arrête $(i,j)$, ces fonctions sont initialement assignées à la fonction nulle
$$
\mu_{i \rightarrow j}:=0 \quad \textrm{et} \quad \nu_{i \leftarrow j}:=0.
$$

\medbreak

\noindent \textbf{Étape 2. Échanges de messages.}
Cette étape est itérée ``autant de fois que nécessaire''. Le nombre d'itérations sera discuté en section \ref{niter}.

\smallbreak

\noindent \textit{Étape 2a. Messages données $\rightarrow$ relations.}
Dans cette étape, on actualise les messages $\mu_{i \rightarrow j}$ des données vers les relations en tenant compte des nouvelles informations fournies par les messages  $\nu_{k \leftarrow l}$. Étant donnés un bit $i$ et une relation $j$, on note $j_1, \ldots, j_s$ les relations voisines de $i$ autres que $j$. Le n\oe ud de données $i$ centralise les informations transmises par les relations $j_1, \ldots, j_k$ et les envoie vers la relation $j$.
Le message $\mu_{i \rightarrow j}$ devient alors
$$
\mu_{i \rightarrow j}:=C^i_{\textrm{loc}}+\sum_{k=1}^s \nu_{i \leftarrow j_k}.
$$

\begin{center}
\includegraphics[width=7cm, height=6cm]{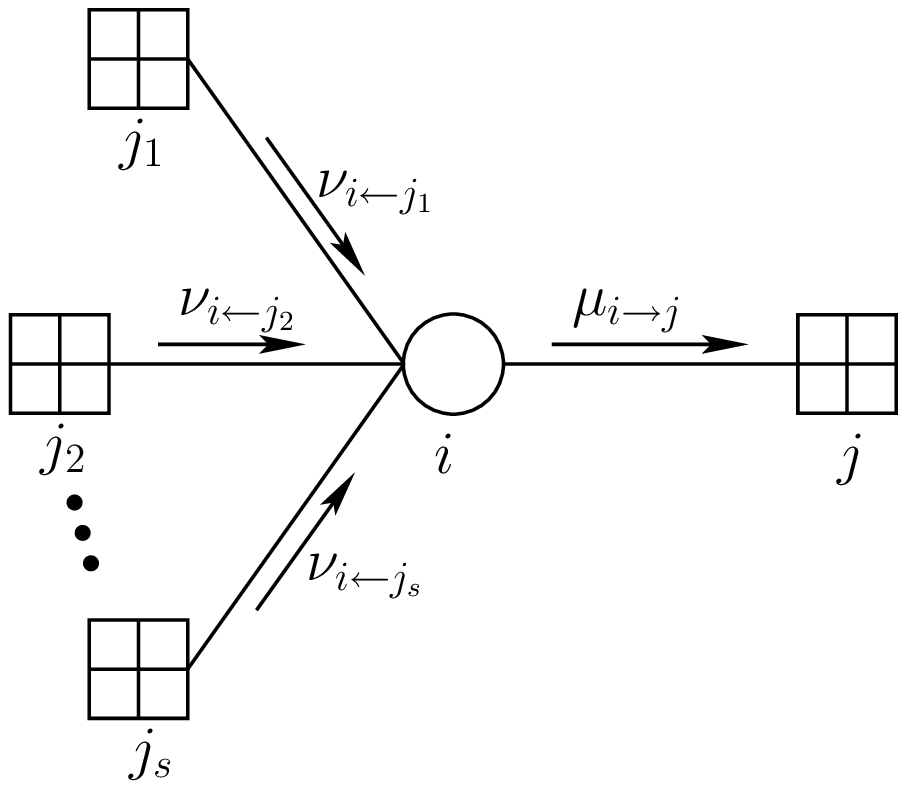}
\end{center}

\begin{rem}
Lors de la première itération de l'algorithme, la fonction $\mu_{i \rightarrow j}$ initialement assignée à la fonction nulle devient égale à la fonction $C^i_{\textrm{loc}}$.   
\end{rem}

\smallbreak

\noindent \textit{Étape 2b. Messages relations $\rightarrow$ données.}
Dans cette étape, on actualise les messages $\nu_{i \leftarrow j}$ en tenant compte des informations fournies par les messages $\mu_{i \rightarrow j}$.
Un noeud de relation $j$ centralise les informations fournies par les fonctions message $\mu_{i_k \rightarrow j}$ provenant des bits voisins autres que $i$ et les redirige vers ce dernier.
Soient donc $i_1,\ldots, i_r$ les n\oe uds de donnée voisins du n\oe ud de relation $j$ autres que $i$.
La fonction $\nu_{i \leftarrow j}$ est définie par
$$
\forall \alpha \in \F_q, \quad \nu_{i \leftarrow j}(\alpha):= \min
\left\{
\sum _{k=1}^r \mu_{i_k \rightarrow j}(\alpha_k)\ \left| \ \begin{array}{c}
(\alpha_1, \ldots , \alpha_r)   \in  \F_q^r, \\
h_{i,j}\alpha  +  h_{i_1,j}\alpha_1 + \cdots +h_{i_r,j}\alpha_r = 0 
\end{array}
\right.
\right\}.
$$ 

\noindent On calcule le coût minimal d'une configuration vérifiant la relation $j$ et telle que le bit $i$ vaille $\alpha$.
On rappelle que les coefficients $h_{i,j}$ sont les coefficients de la matrice de parité $H$ qui pondèrent les arêtes du graphe de Tanner. Dans la figure qui suit, ils n'ont pas été indiqués de façon à alléger la représentation.

\begin{center}
\includegraphics[width=7cm, height=6cm]{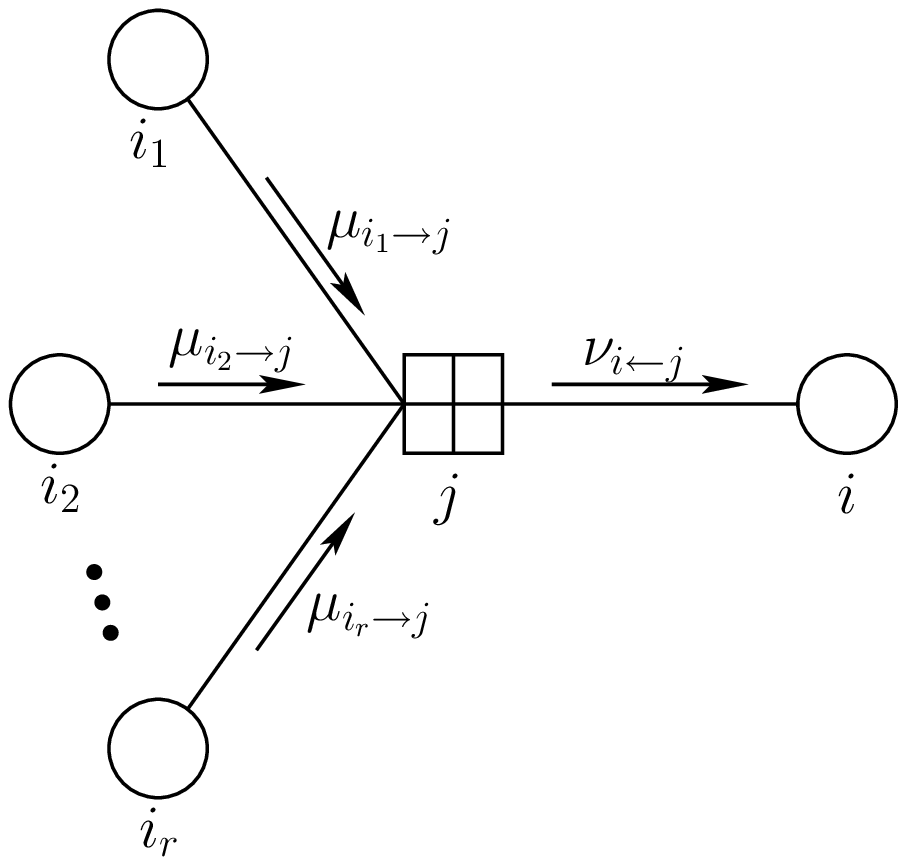}
\end{center}

\begin{rem}\label{iter1} Lors de la première itération de l'algorithme, l'entier $\nu_{i \leftarrow j}(\alpha)$ quantifie le nombre minimal de changements qu'entraînerait l'assignation du $i$-ème bit à la valeur $\alpha$ pour les autres bits intervenant dans la relation $j$.
\end{rem}

\begin{exmp}\label{exb2c}
Dans l'exemple \ref{exminsom} que nous avons étudié précédemment les fonctions $\nu_{2 \leftarrow 1}$ et $\nu_{2 \leftarrow 1}$ ont été calculées, on avait obtenu
$$
\nu_{2 \leftarrow 1}=[0,1]\quad \textrm{et}\quad \nu_{2 \leftarrow 2}=[0,1].
$$
Par le calcul on obtient également (la vérification est laissée au lecteur),
$$
\begin{array}{rclcrclcrclcrcl}
\nu_{1 \leftarrow 1} & = & [0,1] & & \nu_{3 \leftarrow 2} & = & [1,0] & &
\nu_{4 \leftarrow 1} & = & [0,1] & & \nu_{4 \leftarrow 3} & = & [1,0] \\
\nu_{5 \leftarrow 2} & = & [1,0] & & \nu_{5 \leftarrow 3} & = & [0,1] & &
\nu_{6 \leftarrow 3} & = & [1,0]. & & & & 
\end{array}
$$
\end{exmp}

\medbreak

\noindent \textbf{Étape finale. Décision.}
Chaque n\oe ud de donnée évalue ses coûts \textit{globaux} d'assignation avec l'aide des fonctions $\nu_{i \leftarrow j}$. Soient $i$ un noeud de données et $j_1, \ldots, j_t$ l'ensemble des n\oe uds de relation voisins de $i$. La fonction de coût \textit{global} $C^i_{\textrm{glob}}$ est définie par
$$
\forall \alpha \in \F_q, \quad C^i_{\textrm{glob}}(\alpha):=C^i_{\textrm{loc}}(\alpha)+ \sum_{k=1}^t \nu_{i \leftarrow j_k}(\alpha).
$$ 

\noindent On regarde ensuite s'il existe un élément $\alpha$ qui minimise la fonction $C^i_{\textrm{glob}}$. Si oui, on assigne la valeur $\alpha$ au bit $i$.

\begin{exmp}
Dans l'exemple \ref{exminsom}, si l'on prend une décision après une itération du processus d'échanges de messages, à partir des résultats de l'exemple \ref{exb2c}, on obtient
$$
\begin{array}{rclcrclcrcl}
C^1_{\textrm{glob}} & = & [1,1] & & C^2_{\textrm{glob}} & = & [1,2] & & C^3_{\textrm{glob}} & = & [1,1] \\
C^4_{\textrm{glob}} & = & [2,1] & & C^5_{\textrm{glob}} & = & [1,2] & & C^6_{\textrm{glob}} & = & [2,0].
\end{array}
$$
On ne peut donc pas prendre de décision quant à l'assignation finale des bits $1$ et $3$. Il ne fallait pas évaluer les fonctions de coûts globaux à cette étape mais réitérer le processus.
À la seconde itération, l'actualisation des fonctions $\mu_{i \rightarrow j}$ donne
$$
\begin{array}{rclcrclcrcl}
\mu_{1 \rightarrow 1} & = & [1,0] & & \mu_{2 \rightarrow 1} & = & [1,1] & & 
\mu_{2 \rightarrow 2} & = & [1,1] \\
\mu_{3 \rightarrow 2} & = & [0,1] & & \mu_{4 \rightarrow 1} & = & [2,0] & & 
\mu_{4 \rightarrow 3} & = & [1,1] \\
\mu_{5 \rightarrow 2} & = & [0,2] & & \mu_{5 \rightarrow 3} & = & [1,1] & & 
\mu_{6 \rightarrow 3} & = & [1,0]. \\
\end{array}
$$
Après quoi, l'actualisation des fonctions $\nu_{i \leftarrow j}$ donne
$$
\begin{array}{rclcrclcrcl}
\nu_{1 \leftarrow 1} & = & [1,1] & & \nu_{2 \leftarrow 1} & = & [0,1] & & 
\nu_{2 \leftarrow 2} & = & [0,1] \\
\nu_{3 \leftarrow 2} & = & [1,1] & & \nu_{4 \leftarrow 1} & = & [1,1] & & 
\nu_{4 \leftarrow 3} & = & [1,1] \\
\nu_{5 \leftarrow 2} & = & [1,1] & & \nu_{5 \leftarrow 3} & = & [1,1] & & 
\nu_{6 \leftarrow 3} & = & [2,2]. \\
\end{array}
$$
Si l'on évalue les coûts globaux à la fin de cette seconde itération, on obtient
$$
\begin{array}{rclcrclcrcl}
C^1_{\textrm{glob}} & = & [2,1] & & C^2_{\textrm{glob}} & = & [1,2] & & C^3_{\textrm{glob}} & = & [1,2] \\
C^4_{\textrm{glob}} & = & [3,2] & & C^5_{\textrm{glob}} & = & [2,3] & & C^6_{\textrm{glob}} & = & [3,2].
\end{array}
$$
On peut donc prendre une décision, on choisit comme mot décodé le mot
$$
c=(1,0,0,1,0,1),
$$
qui est bien le mot le plus proche du mot reçu $y$ pour la distance de Hamming.
\end{exmp}

\subsection{Discussion sur l'algorithme}

Avant de rentrer dans des considérations plus techniques, commençons par quelques remarques concernant cet algorithme.
\begin{itemize}
\item Le nom de l'algorithme provient bien sûr de l'étape 2b et d'une façon plus générale, du fait que les seules opérations effectuées sont des sommes et des calculs de minima.
\item Il existe également un algorithme appelé somme-produit dont le fonctionnement est assez comparable. Moralement le min-somme calcule des coûts, alors que le somme-produit évalue des probabilités.
\end{itemize}

\subsubsection{Nombre d'itérations}\label{niter}

Concrètement, au départ, chaque bit ne possède que l'information qui le concerne,
à savoir son coût d'assignation à une valeur donnée, sans tenir compte de ses voisins.
C'est l'information qu'il va transmettre à tous les n\oe uds de relation voisins sous la forme des fonctions $\mu_{i\rightarrow j}$ durant la première itération (voir remarque \ref{iter1}).
À la fin de la première itération on peut savoir le coût qu'aurait l'assignation d'un bit à une valeur prescrite pour ce bit et ses voisins (c'est-à-dire à une distance de deux arêtes).
Si l'on réitère le processus $p$ fois, on dispose de toutes les informations provenant des bits à une distance inférieure à $p$ du $i$-ème bit. 
Aussi il semble raisonnable de choisir comme nombre d'itérations la distance maximale entre deux bits dans le graphe de Tanner.

\subsubsection{Le problème des cycles}

Le problème majeur de ces algorithmes est qu'ils agissent localement, sans tenir compte de la géométrie du graphe.
Pour le comprendre reprenons l'exemple \ref{exminsom} et supposons que l'on a effectué deux itérations de la phase d'échanges de messages (étape 2).
Si l'on évalue la fonction de coût global $C^2_{\textrm{glob}}$ après ces deux itérations, le coût évalué prend en compte la contribution de tous les bits qui sont à distance inférieure à $2$ du $2^{\textrm{e}}$. 
Le souci est que, en partant du $2^{\textrm{e}}$ bit, le $4^{\textrm{e}}$ (ainsi que le $5^{\textrm{e}}$) peut être atteint par un chemin de longueur\footnote{Il s'agit d'un graphe biparti, aussi on appelle \textit{chemin de longueur $n$} un chemin de $2n$ arêtes.} deux de deux façons différentes, comme le montre la figure ci-dessous.
De fait, dans le calcul du coût global $C^2_{\textrm{glob}}$, la contribution du quatrième et du cinquième bit est comptée deux fois, ce qui peut biaiser la décision finale.
\begin{center}
\includegraphics[height=6.5cm, width=8cm]{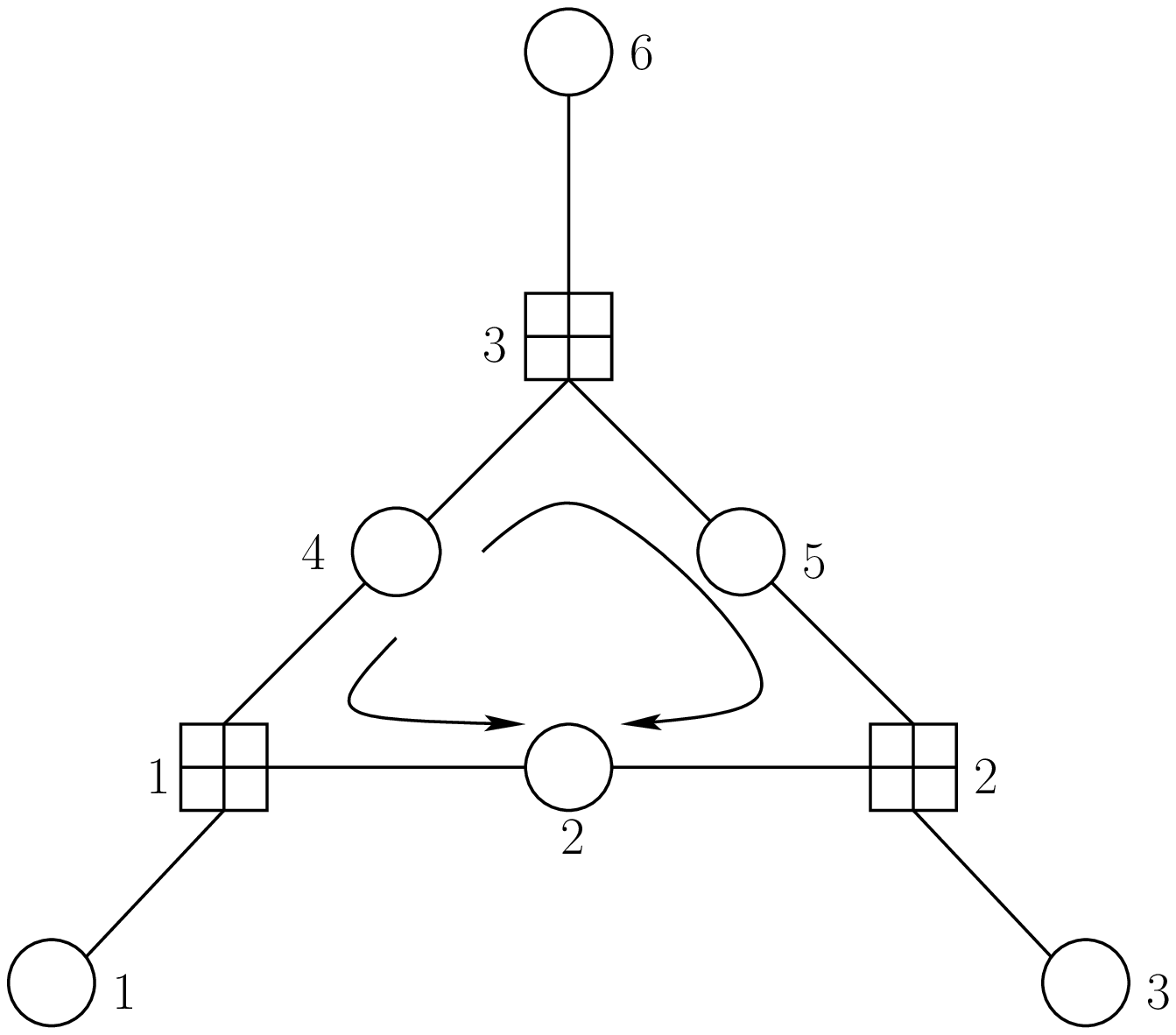}
\end{center}

Les résultats connus sur l'efficacité de l'algorithme sont que le coût global réel d'assignation d'un bit ne peut être calculé exactement par cet algorithme que si le graphe de Tanner est sans cycles. Cependant, les codes dont le graphe de Tanner est sans cycles sont peu intéressants (mauvais paramètres).

De fait, on ne dispose pas réellement de résultats sur l'efficacité d'un tel algorithme. On dispose cependant d'une constatation empirique, à savoir que si le graphe de Tanner n'a pas trop de ``petits cycles'', alors les algorithmes de décodage itératif sont extrêmement efficaces. Ils permettent en particulier de corriger un grand nombre d'erreurs en un temps relativement limité, \textbf{à condition que le code soit LDPC}.

\paragraph{Conclusion.}
L'étape réellement coûteuse est la seconde qui est exponentielle en la valence des n\oe uds de relation. C'est la raison pour laquelle, si l'on travaille sur un code pour lequel cette valence est bornée par une petite valeur, alors l'algorithme tournera rapidement.

\section{Codes LDPC et surfaces de petit degré}

Soient $N$ un entier supérieur ou égal à $3$ et $X$ une hypersurface projective lisse géométriquement intègre de degré $d$ de $\P^N$.
Tout comme dans les chapitres précédents, on note $L_X$, la classe d'équivalence linéaire d'une section hyperplane de $X$.
On se donne également $G$, un diviseur sur $X$ tel que $G\sim mL_X$ pour un certain entier $m$ et $\Delta$, la somme formelle de tous les points rationnels de $X$ qui évitent le support de $G$.
Notons que si $m$ est strictement négatif, le code $C_L(\Delta, G)$ est nul. On peut donc supposer $m$ positif ou nul.
D'après le théorème \ref{minor} (\ref{minor1}) (chapitre \ref{chaporth}), on sait l'orthogonal d'un code fonctionnel sur $X$ a une distance minimale supérieure ou égale à $m+2$ et que cette borne est atteinte dès que le support de $\Delta$ contient $m+2$ points alignés. Cette borne inférieure ne peut donc être atteinte que dans deux situations. 
\begin{enumerate}
\item\label{PicBien} Le degré de $X$ est supérieur ou égal à $m+2$ et il existe une droite de $P^N$ dont l'intersection avec $X$ contient au moins $m+2$ point rationnels. 
\item\label{ContDroite} L'hypersurface $X$ contient une droite rationnelle et toute droite contient au moins $m+2$ point rationnels, ce qui revient à dire que $\sharp \F_q \geq m+1$.
\end{enumerate}

\subsection{Objectifs}
Dans ce qui suit, notre but est de chercher des codes construits sur des surfaces et dont l'orthogonal est engendré par des mots de petit poids, voire de poids minimal. 
Dans cette optique, la situation \ref{PicBien} ci-dessus est en fait la plus intéressante. En effet, la situation \ref{ContDroite} est en général assez rare.
Par exemple, dans le cas où $X$ est une surface ($N=3$), d'après \cite{sch1} théorème I.6.9, une surface générique de degré supérieur ou égal à $4$ ne contient pas de droites  et une surface générique de degré $3$ n'en contient qu'un nombre fini ($27$ si elle est lisse).
De plus, ce dernier résultat est géométrique, ce qui signifie que les droites sur une surface cubique peuvent ne pas être rationnelles.
Les mots de code provenant de la situation \ref{ContDroite}, seront donc en général peu nombreux et engendreront un code dont le support sera souvent strictement contenu dans $\{1, \ldots, n\}$ où $n$ désigne la longueur du code $C_L (\Delta, G)$.
En effet, si $P$ est un point de $\supp (\Delta)$ qui n'est contenu dans aucune droite rationnelle contenue dans $X$, alors l'indice correspondant n'est dans le support d'aucun mot de code provenant de la situation \ref{ContDroite}.

Dans ce qui suit, nous allons nous intéresser aux surfaces fournissant un grand nombre de mots de codes provenant de la situation \ref{PicBien}.
N'ayant pas obtenu de résultat théorique permettant d'orienter cette recherche, nous avons fait appel à l'outil informatique (plus précisément le logiciel \textsc{Magma}).
Dans la section \ref{testsdeouf}, nous allons présenter des résultats expérimentaux effectués sur des surfaces cubiques de $\P^3$.
Auparavant, nous allons nous intéresser au calcul explicite de ces mots en utilisant des résidus.

\section{Calcul explicite de mots de codes de petit poids}\label{secexpl}

Dans ce qui suit, $S$ désigne une surface projective lisse plongée dans $\P^3$ (elle est donc absolument irréductible). On note $d$ le degré de la surface et on suppose que $d \geq 3$.
L'espace projectif $\P^3$ est muni de coordonnée homogènes $(X,Y,Z,T)$.
Le plan d'équation $T=0$ est appelé \textit{plan à l'infini} et noté $\Pi$.
La carte affine $\{T\neq 0\}$ de $\P^3$ est notée $\mathcal{U}_T$ et son intersection avec $S$ est appelée $U_t$. On note $x$, $y$ et $z$ les fonction rationnelles sur $\P^3$ suivantes
$$
x:=\frac{X}{T}, \quad y:=\frac{Y}{T}, \quad \textrm{et} \quad z:=\frac{Z}{T}.
$$

\noindent Ces trois fonctions forment un système de coordonnées affines dans la carte affine $\mathcal{U_T}$ de $\P^3$.
Par ailleurs, on suppose que la surface $S$ n'est contenue dans aucun plan de $\P^3$ et on note $L_{\infty}$ le tiré en arrière du plan à l'infini $\Pi$ sur $S$ via l'injection canonique $S\hookrightarrow \P^3$.
Pour finir, on se donne un entier naturel $m$, on pose 
$$
G:=mL_{\infty}
$$
et on appelle $\Delta$ la somme formelle des points rationnels de $S$ qui évitent le support de $G$. En d'autres termes, $\Delta$ est la somme de tous les points rationnels de la carte affine $U_t$ de $S$.
Nous allons présenter une méthode de calcul explicite des mots de poids minimal du code $C_L(\Delta,G)^{\bot}$ provenant des situations \ref{PicBien} et \ref{ContDroite} signalées page \pageref{PicBien}.

\subsection[Mots provenant de droites non contenues dans $S$]{Mots provenant de droites non contenues dans $\mathbf{S}$}

Nous commençons par considérer un cas simple, à savoir $d=m+2$. Pour des raisons que nous énoncerons plus loin c'est cette situation que nous traiterons plus en détail dans la section \ref{testsdeouf}. Plus précisément nous utiliserons l'outil informatique pour étudier spécifiquement les cas des surfaces cubiques avec $G=L_{\infty}$.

\subsubsection{Le cas $\mathbf{d=m+2}$}
Soit $F$ une droite non contenue dans $S$ et contenant exactement $d$ points $P_1, \ldots , P_d$ appartenant au support de $\Delta$. D'après le théorème de Bezout, le schéma $F\cap S$ est réduit et son ensemble sous-jacent est égal à la réunion des points $P_1, \ldots, P_d$.

\begin{fait}
Étant donné que les points du support de $\Delta$ évitent le support de $G$, ils sont tous contenus dans la carte affine $\mathcal{U}_T$ de $\P^3$. La droite $F$ n'est donc pas contenue dans l'hyperplan à l'infini.  
\end{fait}

Quitte à faire un changement de coordonnées, on peut supposer que la droite $F$ est définie dans la carte affine $\mathcal{U}_T$ par
$$
F_{|\mathcal{U}_T} = \{x=0, \ y=0\}.
$$

\begin{fait}
Soit $G(x,y,z)$ l'équation de $S$ dans $\mathcal{U}_T$. D'après \cite{sch1} III.6.4, la $2$-forme sur $S$
$$
\omega:= \frac{1}{\left(\frac{\displaystyle \p G}{\displaystyle \p z} \right)}. dx \w dy
$$
est régulière sur $U_t$ et ne s'annule en aucun point de cet ouvert. Plus précisément, son diviseur est de la forme $(\omega)=(d-4)L_{\infty}$.
\end{fait}

Soient $D_a$ et $D_b$ les diviseurs très amples définis respectivement par
$$
D_a:=i^*\{X=0\}\quad \textrm{et} \quad D_b:=i^* \{Y=0\},
$$

\noindent où $i$ désigne l'injection canonique $i:S\hookrightarrow \P^3$. On a
$$
(x)=D_a-L_{\infty} \quad \textrm{et} \quad (y)=D_b-L_{\infty}.
$$

\begin{fait}
La $2$-forme sur $S$
\begin{equation}\label{omegap}
\omega':= \frac{1}{\left(\frac{\displaystyle \p G}{\displaystyle \p z} \right)}. \frac{\displaystyle dx }{\displaystyle x} \w \frac{\displaystyle dy}{\displaystyle y}
\end{equation}

\noindent vérifie $(\omega')=(d-2)L_{\infty}-D_a-D_b$. Or, on rappelle que l'on a supposé $d=m+2$ (avec $G=mL_{\infty}$), donc 
$$
(\omega')=G-D_a-D_b.
$$
\end{fait}

\noindent Il reste à vérifier que la paire $(D_a,D_b)$ est sous-$\Delta$-convenable.
Les supports de ces diviseurs s'intersectent seulement en les points $P_1, \ldots, P_d$ et, du fait que le schéma $F \cap S$ est réduit, on en déduit que $D_a$ et $D_b$ sont lisses et s'intersectent transversalement en chacun de ces points.
Si l'on pose $\Lambda:=P_1+\cdots+ P_d$, on a $0 \leq \Lambda \leq \Delta$ et $(D_a,D_b)$ est $\Lambda$-convenable, donc sous-$\Delta$-convenable.

\paragraph{Calcul explicite du mot de code correspondant.}

L'objectif est de calculer de façon explicite les $2$-résidus
$$
\res^2_{D_a,P_i}(\omega'), \quad \textrm{pour}\quad i\in \{1, \ldots, d\},
$$

\noindent où $\omega'$ est la $2$-forme sur $S$ définie dans l'expression (\ref{omegap}).

\noindent On a vu qu'en tout point $P_i$ pour $i\in \{1, \ldots, d\}$, les diviseurs $D_a$ et $D_b$ se croisent transversalement. De plus, au voisinage de ces points, ces deux diviseurs sont respectivement définis par les équations locales $x=0$ et $y=0$. D'après le lemme \ref{valres}, on a pour tout $i\in \{1, \ldots, d\}$
\begin{equation}\label{res_Db}
\res^2_{D_b,P_i}(\omega')=\frac{\displaystyle 1}{\left(\frac{\displaystyle \p G}{\displaystyle \p z} \right)(P_i)}=-\res^2_{D_a,P_i}(\omega').
\end{equation}

 \noindent En effet, soit $C$, la composante de $D_b$ qui passe par $P_i$. La $2$-forme $\omega'$ a un pôle simple le long de $C$, donc
$$
 \res^1_C(\omega')=
\frac{1}{\left(\frac{\displaystyle \p G}{\displaystyle \p z} \right)_{|C}}
. \frac{d\bar{x}}{\bar{x}}
$$
 et le calcul du résidu de cette $1$-forme en $P_i$ donne $\res^2_{C,P_i}(\omega')$ qui est en fait égal à $\res^2_{D_b,P_i}(\omega')$ (voir définition \ref{resdiv}). On en déduit donc la relation (\ref{res_Db}).

\begin{rem}
 Notons que sur la carte affine $U_t$ de $S$, le lieu d'annulation de la fonction $\left(\frac{\p G}{\p z} \right)$ est le lieu des points de branchement du morphisme de projection de $S$ sur le plan d'équation $z=0$.
C'est également le lieu d'annulation de la $2$-forme $dx\w dy$.
Par conséquent, en un point $P_i$ en lequel $D_a$ et $D_b$ s'intersectent transversalement, le couple $(x,y)$ est un système de paramètres locaux et $dx\w dy$ ne s'annule pas.
L'expression (\ref{res_Db}) est donc bien définie. 
\end{rem}

Pour finir, remarquons que l'expression (\ref{res_Db}) s'obtient à condition d'avoir bien effectué un changement de coordonnées pour lequel la droite $F$ est définie par les équations $x=0$ et $y=0$.
Or, si l'on veut réaliser un programme calculant tous les mots de code de $C_L(\Delta, G)^{\bot}$ provenant de droites intersectant $S$ en exactement $d$ points, il sera malcommode de réaliser le changement de variables pour chaque droite.
Une alternative à ce changement de variables, consiste à choisir des équations de $F$ de la forme $L_{|\mathcal{U}_T}\{f(x,y,z)=0,\ g(x,y,z)=0\}$ et un vecteur directeur $\mathbf{v}$ de $F$. On considère alors la $2$-forme sur $S$,

\begin{equation}\label{omegas}
\omega'':= \frac{\displaystyle 1}{\displaystyle \langle\mathbf{grad}(G),\mathbf{v}\rangle}. \frac{\displaystyle df}{f} \w \frac{dg}{g}.
\end{equation}

\noindent
Pour un changement de coordonnées affines de $\mathcal{U}_T$ adapté,
la $2$-forme $\omega''$ ci-dessus coïncide avec la $2$-forme $\omega'$ de l'expression (\ref{omegap}).
L'intérêt de l'expression (\ref{omegas}) est qu'elle fournit une méthode de calcul explicite des mots de $C_L(\Delta,G)^{\bot}$ associés à des droites qui intersectent $S$ en exactement $d$ points rationnels distincts, sans avoir à effectuer de changement de coordonnées. On en déduit le lemme suivant.

\begin{lem}\label{explicit}
  Soit $G$ une équation de $S$ dans la carte affine $\mathcal{U}_T$ et soit $F$ une droite de $\P^3$ qui intersecte $S$ en exactement $d$ points $P_{i_1},\ldots, P_{i_d}$. Alors le code $C_L(\Delta, G)^{\bot}$ contient le mot
$$
c:=\res^2_{D_b,\Delta}(\omega'')\quad \textrm{tel que}\quad
c_i=\left\{
  \begin{array}{ll}
    0 & \textrm{si}\ i \notin \{i_1, \ldots, i_d\},\\
    {\langle \mathbf{grad}_{P_i}(G), \mathbf{v} \rangle}^{-1} & \textrm{sinon} .
  \end{array}
\right.
$$ 
\end{lem}

Nous allons à présent considérer le cas $d<m+2$.

\subsubsection{Le cas $\mathbf{d<m+2}$}

Soient $P_1, \ldots, P_{m+2}$ une famille de points alignés de $\supp (\Delta)$ et soit $F$ la droite les contenant.
Tout comme dans le cas précédent, on va supposer que la droite $F$ est définie par les équations $x=0$ et $y=0$ (ce qui sera toujours vrai après avoir effectué un changement de coordonnées adapté).
On note $\Lambda$, le $0$-cycle défini par $\Lambda:=P_1+\cdots +P_{m+2}$.
Comme $d<m+2$, le $0$-cycle d'intersection $F \cap S$ vérifie
$$\Lambda \leq F \cap S.$$

\noindent Il peut y avoir dans le support de $F \cap S$ d'autres points que les $P_i$ éventuellement de degré supérieur à $1$, certains points peuvent également apparaître avec multiplicité supérieure ou égale à $1$.
Notre objectif est de construire une paire de diviseurs $\Lambda$-convenable $(D_a,D_b)$ pour laquelle l'espace $\Gamma (S, \Omega^2 (G-D_a-D_b))$ est non nul.

Soient $\Pi_1$ et $\Pi_2$ les plans d'équations respectives $x=0$ et $y=0$. On pose
$$
D_a^+:=i^* \Pi_1 \quad \textrm{et} \quad D_b^+:=i^* \Pi_2,
$$ 
où $i$ désigne l'injection canonique $i:S \hookrightarrow \P^3$. Pour tout point géométrique $P$ de $S$ appartenant au support de $L \cap S$, on note $z_P \in \overline{\F}_q$, la coordonnée suivant $z$ de ce point et on pose
\begin{equation}\label{rP}
r_P:=m_P(F,S).
\end{equation}

\begin{rem}
  Dans ce qui suit, afin d'éviter d'alourdir on notera indifféremment par ``$m_P(\ \!.\ \!,.\ \!)$'', la multiplicité d'intersection d'un diviseur et d'un $0$-cycle de $\P^3$ et la multiplicité de deux diviseurs de $S$. L'expression (\ref{rP}) correspond à une multiplicité d'intersection dans $\P^3$. Quant à l'expression (\ref{rP2}) qui suit elle correspond à une multiplicité d'intersection dans $S$.
\end{rem}

D'après les définitions de $D_a^+$ et $D_b^+$, on montre aisément que 
\begin{equation}\label{rP2}
  m_P(D_a^+,D_b^+)=r_P.
\end{equation}

\noindent On définit ensuite pour tout point géométrique $P$ appartenant à $\supp (F\cap S- \Lambda)$ le coefficient
\begin{equation}\label{sP}
s_p:=\left\{
\begin{array}{ccc}
r_P & \textrm{si} & P \notin \supp (\Lambda) \\ 
r_P-1 & \textrm{sinon} & 
\end{array}
\right.
\end{equation}

\noindent Notons que $s_P$ n'est autre que le coefficient de $P$ dans le $0$-cycle 
$\overline{F}\cap \overline{S}$ sur $\overline{S}:=S\times_{\F_q} \overline{\F}_q$.
Comme $F\cap S - \Lambda$ est un $0$-cycle $\F_q$-rationnel, on en déduit que l'ensemble $\{z_P\ |\ P\in \supp (F\cap S - \Lambda)\}$ est invariant sous l'action de $\textrm{Gal}(\overline{\F}_q/\F_q)$. Par conséquent, la fonction 
$$
h:=\prod_{P \in \supp (F\cap S -\Lambda)}(z-z_P)^{s_P}
$$ 
est définie sur $\F_q$ et son degré est égal à celui du $0$-cycle $F\cap S- \Lambda$. On a donc
\begin{equation}\label{degh}
(h) \sim (d-m-2)L_{\infty}.
\end{equation}
 
\noindent Pour finir, posons
\begin{equation}\label{machintruc}
D_a:=D_a^+, \quad D_b:=D_b^+ - (h)^+ \quad \textrm{et}\quad D:=D_a+D_b.
\end{equation}

\begin{lem}
La paire $(D_a,D_b)$ décrite ci-dessus est $\Lambda$-convenable.
\end{lem}

\begin{proof}
  Nous allons utiliser le critère de la proposition \ref{crit}.

\medbreak

\noindent \textbf{Étape 1.} Soit $P$ un point géométrique de $S$ non contenu dans le support de $\Lambda$. Si l'un des diviseurs $D_a$ ou $D_b$ ne contient pas $P$ dans son support, alors le critère est trivialement vérifié en ce point (voir remarque \ref{Det}). Sinon, si le point $P$ fait partie des points géométriques de $\supp (F \cap S)$ autres que $P_1, \ldots , P_{m+2}$.
Comme les diviseurs $D_a^+$ et $D_b^+$ sont obtenus à partir de sections planes de $S$ et que $S$ est lisse, elle ne peut donc pas avoir deux plan tangents distincts en $P$. Ainsi, $\supp(D_a^+)$ ou $\supp (D_b^+)$ est lisse en $P$. 
Supposons que ce soit $\supp (D_a^+)$, il faut alors étudier la multiplicité d'intersection $m_P(D_a^+,D-D_a^+)$ qui n'est autre que $m_P(D_a,D_b)$. Or, d'après~(\ref{machintruc}) et (\ref{rP2}), on vérifie aisément que
$$
\begin{array}{rcl}
m_P(D_a,D_b) & = & m_P(D_a,D_b^+)-m_P(D_a,(h)^+)\\
 & \leq & r_P - s_P.
\end{array}
$$

\noindent Or, comme $P$ n'est pas dans le support de $\Lambda$, d'après (\ref{sP}), on a $s_P=r_P$ et 
$$
m_P(D_a,D_b) \leq 0.
$$

\medbreak

Si maintenant c'est $\supp(D_b^+)$ qui est lisse en $P$, il faut étudier la multiplicité d'intersection $m_P(D_b^+, D_a^+-(h)^+)$. Par un raisonnement identique on montre que cette multiplicité est négative ou nulle.

\medbreak

\noindent \textbf{Étape 2.}
Supposons maintenant que $P$ soit contenu dans le support de $\Lambda$. Par un raisonnement analogue à celui qui a été effectué dans l'étape précédente, on sait que $\supp (D_a^+)$ ou  $\supp (D_b^+)$ est lisse en $P$.
Supposons que $\supp (D_a^+)$ soit lisse en $P$. Il faut calculer 
\begin{equation}\label{acalcu}
m_P(D_a^+, D-D_a^+)=m_P(D_a^+, D_b)=m_P(D_a^+,D_b^+)-m_P(D_a^+,(h)^+).
\end{equation}

\noindent D'après (\ref{rP2}), le terme $m_P(D_a^+,D_b^+)$ est égal à $r_P$. Nous allons montrer que $$m_P(D_a^+, (h)^+)=r_P-1.$$

\smallbreak

\noindent \textit{Étape 2a.} Si $r_P=1$, alors d'après la définition de $h$, le diviseur $(h)^+$ est nul au voisinage de $P$ et 
$$
m_P(D_a^+,(h)^+)=0=r_P-1.
$$ 

\smallbreak

\noindent \textit{Étape 2b.} Si $r_P \geq 2$, alors sur un voisinage $V$ de $P$, on a 
$$
(h)^+_{|V}=((z-z_P)^{s_P})_{|V}=((z-z_P)^{r_P-1})_{|V}.
$$
Comme le plan $\Pi_0$ d'équation $z=z_P$ ne contient pas la droite $F$ et que cette dernière est par hypothèse tangente à $S$ en $P$, on en déduit que le plan $\Pi_0$ est non tangent à $S$ en $P$. Par conséquent, soit $C$ le tiré en arrière de $\Pi_0$ sur $S$. Sur un voisinage de $P$, on a $C=\supp ((h)^+)$ et cette courbe est lisse au voisinage de $P$, de plus il intersecte $\supp (D_a^+)$ transversalement en ce point. En conclusion, sur un voisinage $V$ de $P$, on a 
$$
(h)^{+}_{|V}=s_P C_{|V}.
$$
et 
$$
m_P(D_a^+, (h)^+)=s_P m_P(D_a^+,C)=s_P=r_P-1.
$$

\smallbreak

\noindent Ainsi, quelle que soit la valeur de $r_P$, la relation (\ref{acalcu}) donne
$$
m_P(D_a^+,D_b)=1
$$

\medbreak

Si maintenant, c'est $\supp (D_b^+)$ qui est lisse en $P$, on effectue le même raisonnement en partant de la relation
$$
m_P(D_b^+,D-D_b^+)=m_P(D_b^+,D_a^+)-m_P(D_b^+,(h)^+)
$$
et en montrant que cette multiplicité d'intersection est égale à $1$.

\medbreak

\noindent \textbf{Conclusion.} Le couple $(D_a, D_b)$ vérifie le critère de la proposition \ref{crit}, il est donc $\Lambda$-convenable.
\end{proof}

Soit donc $\omega$ la $2$-forme sur $S$ définie par 
$$
\omega:= \frac{\displaystyle h}{\displaystyle \left( \frac{\displaystyle \p G}{\p z} \right)}. \frac{\displaystyle dx}{\displaystyle x} \w \frac{\displaystyle dy}{\displaystyle y}.
$$

\noindent Un bref calcul permet de montrer que cette $2$-forme vérifie
$$
(\omega)= G-D_a-D_b
$$
et le mot $c:=\res^2_{D_a,\Delta}(\omega)$ appartient à $C_L (\Delta, G)^{\bot}$ et a pour support $\{ 1, \ldots, d \}$.

\begin{rem}
Il est plus délicat de donner une formule simple pour calculer le mot de code $c$ comme dans le lemme \ref{explicit}.
La difficulté vient de ce que les fonctions $h$ et $\langle\mathbf{grad}(G),\mathbf{v}\rangle$ s'annulent toutes deux en $P$. Cependant, si $\supp (D_a^+)$ (resp.$\supp (D_b^+)$) est lisse en $P$ les fonctions $h_{|\supp (D_a^+)}$ et $\langle\mathbf{grad}(G),\mathbf{v}\rangle_{|\supp (D_a)^+}$ ont même valuation en $P$, on peut donc donner un sens à l'évaluation en $P$ de leur rapport.
\end{rem}

\subsection[Mots provenant de droites contenues dans $S$]{Mots provenant de droites contenues dans $\mathbf{S}$}

On suppose dans cette sous-section que le cardinal du corps de base $\F_q$ est supérieur ou égal\footnote{Dans l'introduction de cette section page \pageref{secexpl}, on demande que le cardinal du corps de base soit supérieur ou égal à $m+1$.
En effet, dans cette introduction, la seule contrainte à laquelle on est soumis pour que le support de $\Delta$ puisse contenir $m+2$ points alignés est que le nombre de points rationnels d'une droite projective soit supérieur ou égal à $m+2$.
Maintenant que l'on a précisé le contexte, il faut faire plus attention, car même si la droite projective sur $\F_q$ a $q+1$ points, les points de $\supp (\Delta)$ sont tous par hypothèse contenus dans une carte affine de $S$. Par conséquent, le nombre maximal de points alignés de $\supp (\Delta)$ est au plus égal à $q$. C'est ce qui explique cette hypothèse ``$q$ supérieur ou égal à $m+2$''.} à $m+2$.
Soit $F$ une droite contenue dans $S$ et contenant une famille de points $P_1, \ldots, P_l$ appartenant au support de $\Delta$.
Une fois de plus, quitte à faire un changement de variables, on peut supposer que la droite $F$ est définie sur l'ouvert affine\footnote{Voir page \pageref{secexpl} pour une définition de l'ouvert affine $\mathcal{U}_T$ ainsi que des fonctions $x$, $y$ et $z$.} $\mathcal{U}_T$ de $\P^3$ par les équations $x=0$ et $y=0$.

Soit $\Pi_0$ le plan d'équation $x=0$.
La courbe $C$ définie par l'intersection schématique $C:=\Pi_0 \cap S$ est une courbe plane (car contenue dans $\Pi_0$). Elle est de plus réunion de $F$ et d'une courbe $C'$ de degré $d-1$.
Soit $E(\bar{y}, \bar{z})$ une équation de la courbe $C'$ dans le plan $\Pi_0$. On relève cette fonction en une fonction rationnelle $E(x,y,z)$ sur $\P^3$ qui ne dépend pas de $x$.

Soit enfin $P_1, \ldots, P_{m+2}$ une famille de points de $\supp (\Delta)$ contenus dans $F$. On note $z_1, \ldots, z_{m+2}$ les coordonnées respectives de ces points suivant $z$. On pose
$$
h:=\prod _{i=1}^{m+2} (z-z_i)
$$
et
$$
D_a:=F,\quad D_b:=(h)_0. 
$$

\begin{lem}
Soit $\Lambda$ le $0$-cycle défini par $\Lambda:=P_1+\cdots +P_{m+2}$.
Alors, la paire $(D_a,D_b)$ est $\Lambda$-convenable.
\end{lem}

\begin{proof}
Le diviseur $D_a$ est une droite, c'est donc une courbe lisse. Et le $0$-cycle d'intersection de ces diviseurs est exactement $\Lambda$. De cette dernière assertion, on déduit aisément que cette paire vérifie le critère de la proposition \ref{crit}. Elle est donc $\Delta$-convenable. 
\end{proof}

\noindent Soit alors
$$
\omega:= \frac{\displaystyle E}{\displaystyle h}. \frac{\displaystyle 1}{\displaystyle \left( \frac{\displaystyle \p G}{\p y} \right)}. \frac{\displaystyle dx}{x}\w dz.
$$

\noindent Calculons le diviseur de cette $2$-forme sur $S$.

\smallbreak

\noindent \textbf{Étape 1.} D'après \cite{sch1} III.6.4, on a 
$$
\Bigg(\frac{\displaystyle 1}{\displaystyle \left( \frac{\displaystyle \p G}{\p y} \right)}. du\w dz \Bigg)=(d-4)L_{\infty},
$$

\noindent où l'on rappelle que $L_{\infty}$ désigne la section plane à l'infini.

\smallbreak

\noindent \textbf{Étape 2.} On rappelle que le plan $\Pi_0$ d'équation $x=0$ intersecte $S$ suivant une courbe $F\cup C'$ où $C'$ est une courbe plane de degré $d-1$. De fait,
$$
(u)=F+C'-L_{\infty}.
$$

\smallbreak

\noindent \textbf{Étape 3.} Par construction, la fonction $E$ s'annule suivant la courbe $C'$. On rappelle également que $E$ est un polynôme de degré $d-1$ en $y$ et $z$. Il existe donc un diviseur effectif $C''$ sur $S$ vérifiant
$$
(E)=C'+C''-(d-1)L_{\infty}.
$$

\smallbreak

\noindent \textbf{Étape 4.} La fonction $h$ est un polynôme de degré $m+2$. On a donc
$$
(h)=(h)_0 - (m+2)L_{\infty}.
$$ 

\smallbreak

\noindent \textbf{Étape finale.}
On en déduit donc le calcul du diviseur de $\omega$,
$$
\begin{array}{rcl}
(\omega) & = & (d-4)L_{\infty}-F-C+L_{\infty}+C'+C''-(d-1)L_{\infty}-(h)_0+(m+2)L_{\infty} \\
 & = & mL_{\infty}+C''-F-(h)_0 \\
 & = & C''+G-D_a-D_b \\
 & \geq & G-D_a-D_b.  
\end{array}
$$

\medbreak

\noindent \textbf{Conclusion.} Le mot de code
$$
c:=\res^2_{D_a, \Delta}(\omega)
$$
appartient au code $C_L (\Delta,G)^{\bot}$. De plus, son support correspond exactement aux points $P_1, \ldots, P_{m+2}$.

\section{Expérimentations avec \textsc{Magma}}\label{testsdeouf}

\subsection{Codes sur des surfaces cubiques}
Dans cette section $S$ est une surface cubique lisse de $\P^3$. On munit $\P^3$ d'un système de coordonnées homogènes $(X,Y,Z,T)$. On note $L_{\infty}$ le diviseur défini par l'intersection schématique de $S$ avec le plan \textit{à l'infini} d'équation $T =0$. Soient $P_1, \ldots, P_n$ les points de $S$ qui évitent le support de $L_{\infty}$, on pose alors
$$
\Delta:=P_1+ \cdots +P_n \quad \textrm{et} \quad G:=L_{\infty}.
$$
On rappelle que $U_t$ désigne la carte affine $\{T \neq 0\} \cap S$ de $S$.

\subsubsection{Codes fonctionnels}

Dans le contexte ci-dessus, l'espace $\Gamma (S, \L (G))$ s'identifie à l'espace des polynômes en $x$, $y$ et $z$ de degré inférieur ou égal à $1$. Il est donc de dimension $4$. Le code fonctionnel $C_L (\Delta,G)$ est donc un code de longueur $n$ et de dimension $4$.
La distance minimale de ce type de code dépend du fait que $S$ contienne ou non des droites rationnelles. 
En effet, la distance minimale de ce code est minorée par $n-l_S$, où $l_S$ désigne le nombre maximal de points rationnels d'une section hyperplane de $U_t$.
Pour minorer cette distance minimale, il faut majorer $l_S$. Or, $l_S$ est le nombre maximal de points rationnels d'une courbe affine plane de degré $3$.
Les courbes de degré $3$ ayant le plus grand nombre de points rationnels sont les réunions de trois droites rationnelles concourantes (Voir la lettre de Serre à M.Tsfasman \cite{lettre}).

Ainsi, si la surface $S$ ne contient pas de droite rationnelle (ce qui est possible, voir \cite{swd}), alors la distance minimale du code fonctionnel est plus grande et le code fonctionnel est meilleur. C'est l'approche adoptée par Voloch et Zarzar dans \cite{zarzar} et \cite{agctvoloch}\footnote{Ces articles sont cités dans leur ordre d'écriture.}.

\subsubsection{Orthogonaux des codes fonctionnels}

Dans \cite{agctvoloch}, les auteurs proposent une méthode pour construire de nombreux mots appartenant à l'orthogonal d'un code fonctionnel.
Pour ce faire, ils se donnent un famille de courbes lisses $C_1, \ldots, C_r$ tracées sur $S$.
Ensuite, ils considèrent les codes $C_{L,C_i}(D_i, G_i)$ où $D_i$ est la somme des points de $\supp (\Delta)$ qui appartiennent à $C_i$ et $G_i$ est le tiré en arrière de $G$ sur $C_i$. 
Enfin, pour chacun de ces codes fonctionnels sur des courbes, ils en construisent l'orthogonal et en déduisent des mots dans le code $C_L (\Delta, G)^{\bot}$.
Ces mots ont un support contenu dans l'ensemble des indices correspondant aux points de $\supp (D_i)$. Ils ont donc un poids petit par rapport à la longueur du code.
De cette manière, ils peuvent décoder le code $C_L(\Delta,G)$ par le biais d'un algorithme de décodage itératif du type de celui qui nous avons présenté en section \ref{introldpc}. L'algorithme qu'ils utilisent est présenté dans l'article \cite{luby} de Luby et Mitzenmacher.

Dans ce qui suit, nous allons utiliser les résultats théoriques présentés précédemment pour calculer un grand nombre de mots de poids minimal du code $C_L (\Delta, G)^{\bot}$. Nous chercherons ensuite à savoir si la famille de mots ainsi construite engendre le code $C_L (\Delta, G)^{\bot}$.
Le cas échéant, l'algorithme min-somme pourra être utilisé pour décoder le code $C_L (\Delta, G)$.

\subsubsection{Construction d'un graphe de Tanner}

Notre objectif est de construire un graphe de Tanner en vue d'un décodage itératif. D'après ce qui a été vu en section en section \ref{introldpc}, le cahier  des charges pour un bon graphe de Tanner repose sur deux contraintes.
\begin{enumerate}
\item Éviter les noeuds de relation dont la valence est trop importante car ils augmentent lourdement la complexité de l'algorithme de décodage.
\item  Éviter les \textit{petits cycles}.
\end{enumerate}

Concernant la première contrainte, d'après le théorème \ref{minor} (\ref{minor1}), la distance minimale de $C_L(\Delta, G)^{\bot}$ est supérieure ou égale à $3$. De plus les mots de poids $3$ de ce code ont pour support les indices de trois points alignés du support de $\Delta$.

On dispose par ailleurs d'une formule explicite (lemme \ref{explicit}) pour calculer les mots correspondant à des triplets de points appartenant à des droites non contenues dans $S$. Pour ce qui est des autres mots, il est préférable de les éviter afin de répondre à la seconde contrainte du cahier des charges.
En effet, soit $F$ une droite rationnelle contenue dans $S$. Alors, à tout triplet de points de $F$ contenus dans $\supp (\Delta)$, on associe un mot de poids $3$ dans $C_L (\Delta, G)^{\bot}$. Si l'on note $P_1, \ldots, P_s$ les points de $F(\F_q) \cap \supp (\Delta)$ et que l'on suppose de $s > 3$, alors les mots correspondant par exemple aux triplets $(P_1, P_2, P_4)$ et $(P_1, P_2, P_4)$ donneront un cycle de longueur\footnote{On rappelle que la longueur d'un chemin dans un graphe biparti est la moitié du nombre d'arêtes composant ce chemin.} $2$ dans le graphe de Tanner.

\begin{center}
\includegraphics[width=8cm, height=5cm]{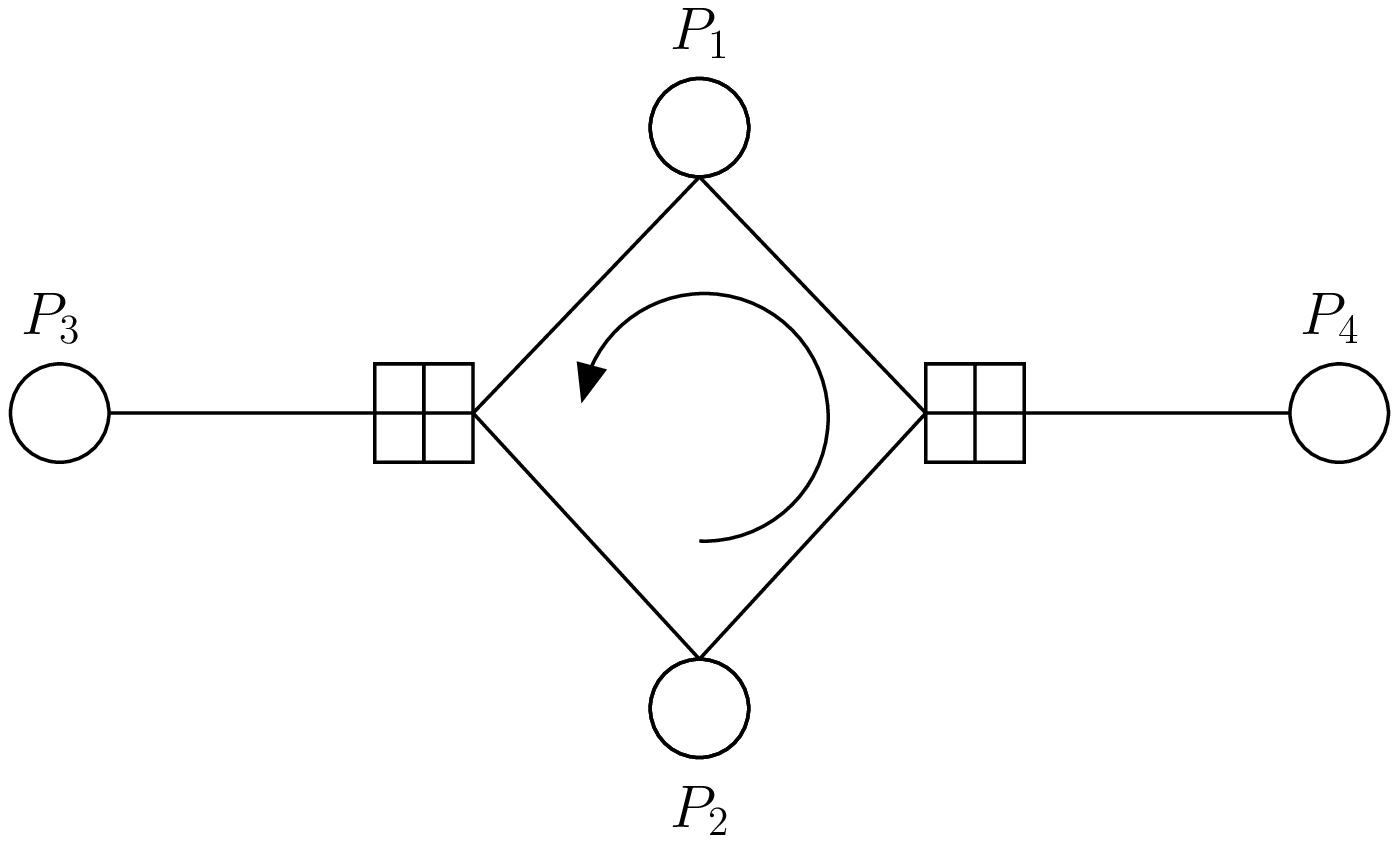}
\end{center}

Il est donc préférable de ne pas choisir tous les mots construits à partir de cette droite. Nous avons fait le choix de n'en retenir aucun.
Nous allons voir que dans la pratique, à condition que le corps $\F_q$ soit assez grand, les mots issus de triplets de points d'une droite non contenue dans $S$ suffisent à engendrer le code $C_L(\Delta, G)^{\bot}$.
Si de plus on ne considère que ces mots là, on évite les cycles de longueur $2$.
Les cycles minimaux seront alors de longueur $3$ et correspondront à la donnée de trois droites coplanaires non contenues dans $S$, non concourantes et telles que les points d'intersection de deux d'entre elles sont dans $\supp (\Delta)$.

\subsection{Implémentation}

Pour construire un graphe de Tanner, nous allons utiliser l'algorithme suivant.

\medbreak

\noindent \textbf{Algorithme de construction d'équations de parité.}

\noindent \textit{Entrées:} Une surface cubique lisse de $\P^3$.
\textit{Sorties:} Une matrice.

\begin{enumerate}
\item On se donne une liste de mots de codes $M$, initialement vide ($M:=[\ ]$).
\item On crée l'ensemble $Points$ des points rationnels de la carte affine $U_t$ de $S$.
\item On crée un second ensemble PointsBis qui est initialisé à $PointsBis:=Points$.
\item Pour $P \in Points$

\begin{itemize}
\item[\textbullet] Enlever $P$ de $PointsBis$.

\item[\textbullet] Pour $Q\in PointsBis$,

\begin{itemize}
\item[$\circ$] Si la droite $(PQ)$ n'est pas contenue dans $S$ et contient exactement trois éléments de
$Points$, alors on construit le mot de code $c$ correspondant par la formule  du lemme \ref{explicit} et on ajoute $c$ dans la liste $M$.
\item[$\circ$] Sinon, on ne fait rien.
\end{itemize}
\end{itemize}

\item On construit une matrice dont les lignes sont les éléments de $M$.
\end{enumerate}

\begin{rem}
  Un programme Magma de cet algorithme est donné en annexe \ref{prgmldpc}.
\end{rem}

Une question se pose ensuite, \textit{La matrice ainsi construite est-elle une matrice génératrice  de $C_L (\Delta, G)^{\bot}$?}
Remarquons que le code $C_L(\Delta, G)$ est de dimension $4$. Aussi, pour vérifier qu'une matrice obtenue par l'algorithme ci-dessus est bien génératrice de $C_L(\Delta, G)^{\bot}$, il suffit de vérifier qu'elle est de rang $n-4$.
Les expérimentations présentées dans le tableau ci-dessous montrent qu'en général la matrice obtenue est bien une matrice génératrice de l'orthogonal. Nous ne sommes toutefois pas parvenus à fournir une preuve théorique de ce fait.

\medbreak

Dans le tableau qui suit, nous présentons une série d'expériences.
Le test de base est le suivant.
Pour une surface cubique lisse sur un corps $\F_q$ avec $q>2$ on calcule une matrice en utilisant l'algorithme ci-dessus.
Ensuite, on calcule le rang de cette matrice et le compare avec $n-4$. S'il y a égalité le test est positif sinon il est négatif.

Voici les résultats de tests sur des surfaces choisies de façon aléatoire pour différents corps de base.

$$
\begin{array}{|c|c|c|c|c|c|}
\hline
\textrm{\textbf{Corps}} & \textrm{\textbf{Nombre de tests}} & 
\textrm{\textbf{Nombre de tests}} & \mathbf{\%} &
\textrm{\textbf{Écart}} & \textrm{\textbf{Longueur}} \\
\textrm{\textbf{de base}} & \textrm{\textbf{effectués}} & 
\textrm{\textbf{positifs}} & & \textrm{\textbf{moyen}} & \textrm{\textbf{moyenne}}\\
\hline
\F_3 & 10000 & 1139 & 11,39 \% & 2,15 & 8,80 \\
\hline
\F_4 & 10000 & 6274 & 62,64 \% & 1,96 & 15,88 \\
\hline 
\F_5 & 10000 & 9763 & 97,63 \%  & 1,82 & 24,90 \\
\hline
\F_7 & 1000 & 1000 & 100 \% & - & 48,66 \\
\hline
\F_8 & 500 & 500 & 100 \% & - & 64,05 \\
\hline
\F_9 & 500 & 500 & 100 \% & - & 81,008 \\
\hline

\end{array}
$$

\medbreak

\noindent Les deux dernières colonnes fournissent les quantités suivantes.

\begin{itemize}

\item[\textbullet] \textbf{Écart moyen.} C'est la moyenne sur l'ensemble des tests négatifs de la quantité $n-4-\textrm{Rang}(M)$, où $M$ désigne la matrice construite grâce à l'algorithme ci-dessus.

$$
\textrm{\textbf{\'Ecart moyen}}= \frac{n-4-\textrm{Rg}(M)}{ \textrm{\textbf{Nombre de tests négatifs}}}.
$$

\item[\textbullet] \textbf{Longueur moyenne.} C'est la longueur moyenne des codes construits, c'est-à-dire le nombre moyen de points rationnels des cartes affines $U_t$ des surfaces cubiques testées.

$$
\textrm{\textbf{Longueur Moyenne}}= \frac{\textrm{\textbf{Longueur du code}}}{ \textrm{\textbf{Nombre de tests effectués}}}.
$$

On remarque dans le tableau ci-dessus que cette longueur moyenne est proche de $q^2$, ce qui est naturel puisque le nombre de points rationnels d'une surface cubique affine lisse est lui-même proche de $q^2$.
\end{itemize}

\begin{rem}
  Lorsque la taille du corps grandit, le nombre moyen de points rationnels d'une surface cubique augmente de façon quadratique en la taille du corps, ce qui contribue à augmenter lourdement la complexité de l'algorithme. 
C'est la raison pour laquelle le nombre de tests est moins important lorsque le corps de base est plus grand.
\end{rem}

\paragraph{Conclusion.}
Il semble très probable que pour $q\geq 7$, le code $C_L(\Delta, G)^{\bot}$ ainsi construit soit engendré par ses mots de poids $3$.
Il serait d'ailleurs intéressant d'obtenir une démonstration mathématique de ce résultat (si du moins il est vrai).

\subsection{Codes sur des surfaces quartiques}

Nous avons réalisé le même type d'expérience dans le cas où $S$ est une surface de degré $4$ et $G\sim 2L_S$.
Voici les résultats de l'expérience.

$$
\begin{array}{|c|c|c|c|c|c|}
\hline
\textrm{\textbf{Corps}} & \textrm{\textbf{Nombre de tests}} & 
\textrm{\textbf{Nombre de tests}} & \mathbf{\%} &
\textrm{\textbf{Écart}} & \textrm{\textbf{Longueur}} \\
\textrm{\textbf{de base}} & \textrm{\textbf{effectués}} & 
\textrm{\textbf{positifs}} & & \textrm{\textbf{moyen}} & \textrm{\textbf{moyenne}}\\
\hline
\F_4 & 1000 & 40 & 4 \% &   5.11 &
15.9\\
\hline 
\F_5 & 1000 & 0 & 0 \% & 
 9.66 &  
24.57 \\
\hline
\F_7 & 1000 &  204 & 20,4 \% &  8.87 &
48.71 \\
\hline
\F_8 & 1000 & 633 & 63,3 \% & 5,37 & 64,14 \\
\hline
\F_9 & 1000 & 894 & 89,4 \% & 2,7 & 80,98 \\
\hline
\F_{11} & 1000 & 999 & 99,9 \% & 1 & 121,233 \\
\hline
\F_{13} & 1000 & 1000 & 100 \% & - & 168,711 \\
\hline

\end{array}
$$

La conclusion est sensiblement la même que pour l'expérience sur les cubiques. Il semble que l'orthogonal du code fonctionnel soit engendré par ses mots de poids $4$ à condition que le cardinal du corps de base soit suffisamment grand.

\subsection{Utilisation de l'algorithme min-somme pour le décodage de ces codes.}

Le graphe de Tanner construit de cette manière offre de bonne perspectives de décodage.
Reprenons par exemple la surface donnée par Voloch et Zarzar dans \cite{agctvoloch}, c'est-à-dire la surface $S$ sur $\F_3$ d'équation
$$
X^3+Y^3+Z^3-ZX^2-XY^2-YZ^2+XZ^2+T^3.
$$
Les auteurs montrent dans l'article cité ci-dessus que le code $C_{L,S}(\Delta, G)$, où $G$ est la section plane à l'infini, est un code de longueur $13$, de dimension $4$ et de distance minimale $7$.

Les points rationnels de cette surface sont
$$
\begin{array}{ccccccccc}
P_1 & = &  (2 : 0 : 0 : 1) & P_2 & = &  (1 : 0 : 1 : 1) & P_3 & = &  (0 : 0 : 2 : 1)\\
P_4 & = & (1 : 0 : 2 : 1) & P_5 & = &  (2 : 1 :1 : 1) & P_6 & = &  (0 : 1 : 2 : 1)\\
P_7 & = &  (2 : 1 : 2 : 1) & P_8 & = & (0 : 2 : 0 : 1) & P_9 & = & (1 : 2 : 0 : 1) \\
P_{10} & = & (2 :2 : 0 : 1) & P_{11} & = & (2 : 2 : 1 : 1) & P_{12} & = &  (0 : 2 : 2 : 1)\\
P_{13} & = &  (2 : 2 : 2 : 1) & & & & & & 
\end{array}
$$
et ils évitent tous le support de $G$. Si maintenant, on applique l'algorithme, les droites qui coupent $S$ en exactement trois points rationnels sont les treize droites ci-dessous.
$$
\begin{array}{cclccl}
L_1 & = & \{ x + t=0,\ y + z=0\} & L_2 & = & \{x+t=0, y+2z=0\}\\
L_3 & = & \{ x+z+t=0,\ y =0\}    & L_4 & = & \{x+z+t=0, y+2z+t=0\}\\
L_5 & = & \{ x + 2z=0,\ y + 2z+t=0\} & L_6 & = & \{x=0, z+t=0\}\\
L_7 & = & \{ x+2z + 2t=0,\ y + z+t=0\} & L_8 & = & \{x+2y+2t=0, z+t=0\}\\
L_9 & = & \{ x +y+2t=0,\  z+t=0\} & L_{10} & = & \{x+z=0, y+z+t=0\}\\
L_{11} & = & \{ y + t=0,\ z=0\} & L_{12} & = & \{x+2z+2t=0, y+t=0\}\\
L_{13} & = & \{ x + t=0,\ y + t=0\} & &  & \\
\end{array}
$$

\noindent Enfin, par la formule du lemme \ref{explicit}, on obtient la matrice

$$
M=\left(\begin{array}{ccccccccccccc}
2 & 0& 0& 0& 0& 0& 2& 0& 0& 0& 2& 0& 0 \\
2& 0& 0& 0& 2& 0& 0& 0& 0& 0& 0& 0& 2 \\
2& 2& 2& 0& 0& 0& 0& 0& 0& 0& 0& 0& 0 \\
0& 2& 0& 0& 0& 0& 2& 2& 0& 0 &0& 0& 0 \\
0& 1& 0& 0& 0& 1& 0& 0& 0& 1& 0& 0& 0 \\
0& 0& 2& 0& 2& 0& 0& 0& 2& 0& 0& 0& 0 \\
0& 0& 1& 0& 0& 1& 0& 0& 0& 0& 0& 1& 0 \\
0& 0& 0& 1& 0& 0& 1& 0& 0& 0& 0& 1& 0 \\
0& 0& 0& 1& 0& 1& 0& 0& 0& 0& 0& 0& 1 \\
0& 0& 0& 1& 1& 0& 0& 1& 0& 0& 0& 0& 0 \\
0& 0& 0& 0& 0& 0& 0& 2& 2& 2& 0& 0& 0 \\
0& 0& 0& 0& 0& 0& 0& 0& 1& 0& 1& 1& 0 \\
0& 0& 0& 0& 0& 0& 0& 0& 0& 2& 2& 0& 2 
\end{array}\right)
$$

On vérifie que la matrice $M$ est de rang: $9=13-4$. Elle vérifie bien la relation $\textrm{Rg}(M)=n-\dim C_L$, la matrice $M$ est donc bien génératrice du code $C_L(\Delta, G)^{\bot}$. 
On peut l'utiliser pour implémenter un algorithme de décodage min-somme pour le code $C_L (\Delta,G)$.


\newpage
\thispagestyle{empty}

\chapter*{Conclusion}\label{chapconclu}

Pour construire des codes différentiels à partir de surfaces algébriques, nous avons développé le matériel théorique nécessaire à l'obtention d'une formule de sommation de résidus en dimension $2$.
Ce résultat était déjà connu et dans un contexte plus général que celui des surfaces (voir \cite{harRD}, \cite{parshin} et \cite{lip}). Cependant, l'approche adoptée dans le premier chapitre fournit des constructions explicites et une démonstration plus accessible de cette formule.

Dans un second temps, on montre qu'un certain nombre de propriétés vérifiées par les codes différentiels construits sur les courbes s'étendent aux codes différentiels construits sur les surfaces.
En fait, seule la relation d'orthogonalité ne s'étend pas parfaitement. 
Le théorème de réalisation est, en un certain sens, une manière de remédier à ce défaut d'inclusion réciproque dans la relation 
``$
C_{\Omega} \subset C_L^{\bot}
$''.

Pour le reste, cette absence d'inclusion réciproque rend, d'une certaine manière, l'étude des codes géométriques construits sur des surfaces plus riche que celle des codes construits sur des courbes.
En effet, dans le contexte des surfaces algébriques, les codes fonctionnels n'appartiennent plus en général à la même classe de codes que leurs orthogonaux.

Pour finir, rappelons que l'étude de ces deux classes de codes ouvre d'intéressants problèmes de théorie des codes et de géométrie algébrique que nous rappelons une dernière fois afin de conclure cette thèse. Commençons par énoncer les différentes questions posées tout au long de ce texte.

\medbreak

\noindent \textbf{Question \ref{QQk}.}
\textit{Peut-on estimer les paramètres des codes qui sont l'orthogonal de codes fonctionnels?  }

\medbreak

\noindent \textbf{Question \ref{Qsom}.} \textit{Si l'orthogonal d'un code fonctionnel ne peut se réaliser comme un code différentiel associé à une paire de diviseurs (sous-)$\Delta$-convenables, peut-on le réaliser comme somme de tels codes?}

\medbreak

\noindent \textbf{Question \ref{Qsom}bis.}
\textit{  Étant donné un mot de code $c$ appartenant à $C_{L,S}(\Delta,G)^{\bot}$, existe-t-il une paire de diviseurs (sous-) $\Delta$-convenable $(D_a,D_b)$ et une $2$-forme $\omega$ appartenant à  $\Gamma (S,\Omega^2(G-D_a-D_b))$ et telle que 
$$
c=\res^2_{D_a,\Delta}(\omega)\textrm{?}
$$ }

\medbreak

\noindent \textbf{Question \ref{hypless}.} \textit{Le résultat du théorème de réalisation(théorème \ref{thmreal}) reste-t-il vrai si l'on élimine l'hypothèse \ref{secreal} sur $S$ et $G$ ($S$ est intersection complète et $G$ est linéairement équivalent à l'intersection de $G$ avec une hypersurface)?}

\medbreak

\noindent \textbf{Question \ref{bornsom}.} \textit{Sous les conditions du corollaire \ref{correal}, peut-on estimer le nombre de minimal de codes différentiels dont la somme est égale à l'orthogonal d'un code fonctionnel en fonction d'invariants géométriques de la surface?}

\medbreak

\noindent \textbf{Question \ref{qpoon}} (Arithmétique)\textbf{.}
\textit{Soient $X$ une variété projective lisse géométriquement intègre sur un corps fini $\F_q$ et $P_1, \ldots, P_n$, une famille de points fermés de $X$.
Peut-on évaluer explicitement ou majorer de façon précise le plus petit entier $d$ tel qu'il existe au moins une hypersurface définie sur $\F_q$ de degré inférieur ou égal à $d$ qui interpole tous les $P_i$ et dont l'intersection schématique avec $X$ soit une sous-variété lisse géométriquement intègre de codimension $1$?}

\medbreak

\noindent \textbf{Question \ref{qpoon}} (Géométrique)\textbf{.}
\textit{Soit $X$ une variété projective irréductible lisse définie sur $\overline{\F}_q$ et $P_1, \ldots, P_n$ une famille de points de $X$.
Peut-on évaluer explicitement ou majorer de façon précise le plus petit entier $d$ tel qu'il existe au moins une hypersurface $H$ de degré inférieur ou égal à $d$, qui contienne tous les $P_i$ et telle que $H\cap X$ soit une sous-variété lisse de codimension $1$ de $X$?}
 
\medbreak

\noindent \textbf{Question \ref{configmin}.}
\textit{Sachant que les deux premières configurations minimales de points rationnels $m$-liés dans $\P^N$ sont la donnée de $m+2$ points alignés et $2m+2$ points sur une même conique plane, quelles sont les configurations minimales suivantes.}

\medbreak

\noindent \textbf{Question \ref{Qpipeau}.}
\textit{Soient $X$ une sous-variété irréductible lisse géométriquement intègre de $\P^r_{\overline{\F}_q}$
et $d$ un entier naturel.
Soient $P_1, \ldots, P_n$ une famille de points de $X$.
Sous quelles conditions sur $X$ et $P_1, \ldots , P_n$ a-t-on l'existence d'un entier $s$ tel que pour tout $s$-uplet de points parmi $P_1, \ldots, P_n$, il existe une hypersurface $H$ de degré $d$ contenant ce $s$-uplet de points et telle que $H\cap X$ soit une sous-variété lisse de codimension $1$ de $X$?}

\thispagestyle{myheadings}
\markboth{}{} 

\medbreak

Le diagramme suivant représente les relations entre ces questions ainsi que certaines parties ou résultat de cette thèse.
Une flèche $A\rightarrow B$ doit se lire ``$A$ motive $B$'' lorsque $B$ est une question et ``$B$ répond à $A$'' lorsque $A$ est une question. les flèches pointillées signifient que les questions/résultats/parties ne sont qu'indirectement liés.

\bigbreak

 \xymatrix{& & & & \includegraphics[width=1.42cm, height=1cm]{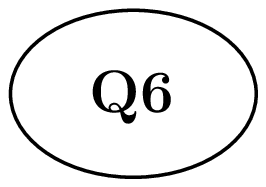}\\
 & & \includegraphics[width=1.42cm, height=1cm]{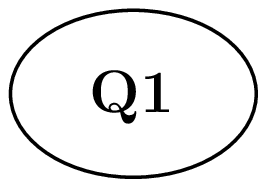} \ar[r] \ar@{.>}[rrd] & \includegraphics[width=1.5cm, height=0.75cm]{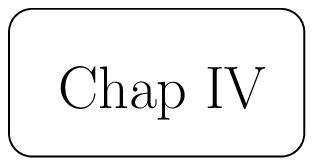} \ar[ru] \ar[r] \ar@{.>}[rd] & \includegraphics[width=1.42cm, height=1cm]{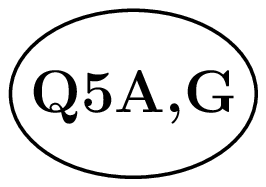} \ar[l] \\
 &\includegraphics[width=2cm, height=0.66cm]{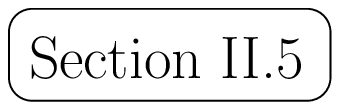}\ar[ru] \ar[rd] & & & \includegraphics[width=1.42cm, height=1cm]{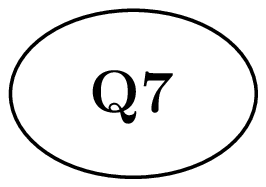} \\
 & &\includegraphics[width=1.42cm, height=1cm]{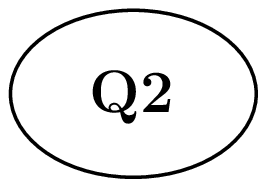} \ar[r] \ar@{=>}[d] & \includegraphics[width=2cm, height=0.75cm]{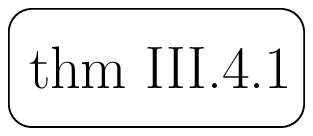} \ar@{.>}[ru] \ar[r] \ar[rd] & \includegraphics[width=1.42cm, height=1cm]{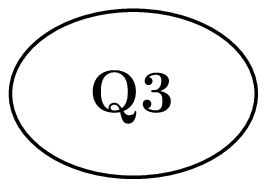} \\
 & & \includegraphics[width=1.42cm, height=1cm]{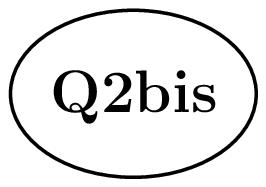} \ar@{=>}[u] & & \includegraphics[width=1.42cm, height=1cm]{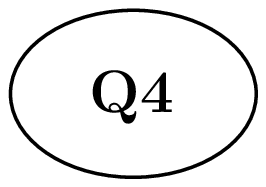}}

\medbreak Nous conclurons en rappelant les perspectives qu'ouvriraient certaines questions ou problèmes ouverts.
En particulier, une solution aux problèmes portant sur des systèmes linéaires de type Bertini posés par les questions \ref{qpoon} A et G fournirait une méthode d'estimation de la distance minimale de codes fonctionnels sur des surfaces ou même des variétés de dimension supérieures.
Quant à la question \ref{Qpipeau} qui peut être vue comme une variante des \ref{qpoon} A et G, une réponse à cette dernière fournirait une élégante méthode d'estimation de la distance minimale de l'orthogonal d'un code fonctionnel sur une surface algébrique.


\part*{Annexes}

\appendix

\chapter{Séries de Laurent}\label{annexe_lolo}

\begin{flushright}
\begin{tabular}{p{6.5cm}}
\begin{flushright}
{\small \textit{Les maths, ça s'écrit comme du français...  ...les formules en plus}.}
\end{flushright}
\medbreak
\hfill Marc Perret
\end{tabular}
\end{flushright}

Cette première annexe contient toutes les démonstrations techniques relatives aux séries de Laurent et aux formes différentielles formelles. 

\section{Sur les modules de différentielles relatives}\label{ann_dif_lolo}

Le but est de démontrer le lemme \ref{reldif} sur les modules de différentielles relatives. Le résultat peut sembler assez élémentaire, cependant, les opérations de complétions et de passage au module des différentielles relatives ne commutent pas en général (voir \cite{eis} ex 16.14). Nous avons donc choisi d'en donner une preuve détaillée faute de référence. 

On se place dans le cadre décrit en \ref{cadre}. On dispose de plus d'une $(P,C)$-paire faible  sur $S$ (voir définition \ref{PCfaible}).

\medbreak

\noindent \textbf{Étape 1.} Commençons par montrer que
$$
\Omega^1_{k((u))/k} \cong \Omega^1_{k(C)/k}\otimes_{k(C)} k((u)).
$$
On sait que $\Omega^1_{k(C)/k}$ est un $k(C)$-espace vectoriel de dimension $1$. De plus, $\bar{u}$ étant une uniformisante de $\mathcal{O}_{C,P}$, c'est un élément séparant de $k(C)/k$, donc la forme $d\bar{u}$ est non nulle sur $C$ et engendre donc $\Omega^1_{k(C)/k}$ sur $k(C)$.
Il suffit donc de montrer que $du$ engendre $\Omega^1_{k((u))/k}$ sur $k((u))$.
D'après \cite{matsumura} ex 25.3, l'application $d: k((u)) \rightarrow \Omega^1_{k((u))/k}$ est continue pour la topologie $(u)$-adique. Cela implique que pour tout $f\in k((u))$, on a $df=f'(u)du$ où $f'$ est la dérivée formelle de $f$ par rapport à $u$. Par conséquent, tout élément de $\Omega^1_{k((u))/k}$ étant une somme finie d'éléments de la forme $fdg=fg'du$, on en déduit que $\Omega^1_{k((u))/k}$ est engendré par $du$ sur $k((u))$.  

\medbreak

\noindent \textbf{Étape 2.} Montrons maintenant que $\Omega^1_{k((u))((v))/k}\cong \Omega^1_{k(S)/k}\otimes_{k(S)} k((u))((u))$.
Comme $\Omega^1_{k(S)/k}$ est librement engendré par $du$ et $dv$ sur $k(S)$, le module $\Omega^1_{k(S)/k}\otimes_{k(S)} k((u))((u))$ est librement engendré par $du$ et $dv$ sur $k((u))((v))$.
Considérons l'application
$$
\delta: \left\{
\begin{array}{ccc}
k((u))((v)) & \rightarrow & \Omega^1_{k(S)/k}\otimes_{k(S)} k((u))((u)) \\
f & \mapsto & \frac{\p f}{\p x}dx + \frac{\p f}{\p y}dy 
\end{array}
\right. .
$$
Cette application est une dérivation. Donc, d'après la propriété universelle du module des différentielles, il existe un unique morphisme $\bar{\delta}$ qui fasse commuter le diagramme suivant
$$
\xymatrix{ 
k((u))((v)) \ar[r]^-{\delta} \ar[d]^d & \Omega^1_{k(S)/k}\otimes_{k(S)} k((u))((u)) \\
\Omega^1_{k((u))((v))/k} \ar[ru]^-{\bar{\delta}}
}.
$$
De fait, les différentielles $du$ et $dv$ appartenant à $\Omega^1_{k((u))((v))/k}$ sont respectivement envoyées par $\bar{\delta}$ sur $du$ et $dv$ appartenant à $\Omega^1_{k(S)/k}\otimes_{k(S)} k((u))((v))$.
Par conséquent, ces deux formes différentielles sont linéairement indépendantes sur $k((u))((v))$ dans $\Omega^1_{k((u))((v))/k}$. Ce dernier est donc de dimension au moins $2$ sur $k((u))((v))$, mais, d'après \cite{matsumura} théorème 25.1, les injections successives
$$
k \rightarrow k((u)) \rightarrow k((u))((v))
$$
donnent une suite exacte
$$
\xymatrix{
\Omega^1_{k((u))/k} \otimes_{k((u))} k((u))((v))
\ar[r] &  \Omega^1_{k((u))((v))/k} \ar[r] & \Omega^1_{k((u))((v))/k((u))} \ar[r] & 0
}.
$$

\begin{rem}
  Cette suite exacte est en général appelée première suite exacte fondamentale (\cite{H} proposition II.8.3.A).
\end{rem}

\noindent D'après l'étape $1$, l'espace $\Omega^1_{k((u))((v))/k((u))}$ est de dimension $1$ sur $k((u))((v))$ et engendré par $dv$. De même, $\Omega^1_{k((u))/k} \otimes_{k((u))} k((u))((v))$ est de dimension $1$ et engendré sur $k((u))((v))$ par $du$. L'espace $\Omega^1_{k((u))((v))/k}$ est donc de dimension au plus $2$. On conclut qu'il est de dimension $2$ et librement engendré par $du$ et $dv$.

\medbreak

\noindent \textbf{Étape 3.} Pour conclure quant à l'isomorphisme $\Omega^2_{k((u))((v))/k}\cong \Omega^2_{k(S)/k}\otimes_{k(S)}k((u))((v))$, il suffit de remarquer que
$$\Omega^{2}_{k((u))((v))/k}= \bigwedge^2\ \Omega^1_{k((u))((v))/k}.$$

\section{Démonstration du lemme \ref{jacob}}\label{demojacob}

Pour commencer nous allons introduire deux applications.
La première est notée $J$ (comme Jacobien), et est définie par
$$
J :\left\{\begin{array}{ccc}
k((u))((v))^2 & \rightarrow & k((u))((v)) \\
(A,B) & \mapsto & \frac{\partial A}{\partial u} \frac{\partial
  B}{\partial v} - \frac{\partial A}{\partial v} \frac{\partial B}{\partial u}
\end{array}\right.
.
$$ 
La seconde est notée $\rho$ est définie par
$$\rho :\left\{\begin{array}{ccc}
k((u))((v)) & \rightarrow & k((u)) \\
\sum_{i\geq -n} h_i(u)v^i & \mapsto & h_{-1}(u)
\end{array}\right. .$$

\begin{lem}
Soient $A,B$ deux éléments de $k((u))((v))$, alors il existe une série de Laurent $\phi \in k((u))$ telle que
$$\rho \circ J(A,B)=\phi'(u),$$  

\noindent où $\phi'(u)$ désigne la dérivée formelle d'une série de Laurent $\phi \in k((u))$.
\end{lem}

\begin{proof}
Comme les applications $\rho$ et $J$ sont respectivement $k$-linéaire et $k$-bilinéaire antisymétrique,
on peut démontrer le lemme en  ne considérant que les trois situations ci-dessous. Le résultat s'en déduira en utilisant ces propriétés de linéarité et d'antisymétrie.

\begin{enumerate}
\item $A$ et $B$ sont éléments de $k((u))[[v]]$.
\item $A\in k((u))[[v]]$ et $B$ est de la forme $B=\frac{b(u)}{v^{n}}$
  avec $n\in \mathbf{N}^*$ et $b(u)\in k((u))$.
\item $A=\frac{a(u)}{v^{m}}$ et $B=\frac{b(u)}{v^{n}}$ avec $m,n \in
  \mathbf{N}^*$ et $a(u),b(u)\in k((u))$.
\end{enumerate}

\noindent Traitons séparément ces trois situations.

\medbreak

\noindent\textbf{Cas 1.} Les séries $A$ et $B$ n'ont pas de pôle suivant la variable $v$, leurs dérivées partielles non plus, donc
$\rho \circ J(A,B)=0$.

\medbreak

\noindent\textbf{Cas 2.} La série $A$ est de la forme, $A=\sum_{i\geq 0}a_i(u)v^i$. Le calcul de $J(A,B)$ donne
$$
J (A,B) =  \sum_{i\geq 0} a_{i}'(u)b(u)(-n)v^{i-n-1} - \sum_{i\geq 0} a_i(u)b'(u)iv^{i-n-1}.
$$
On a donc
$$
\rho (J (A,B))  =  -n(a_n'(u)b(u)+a_n(u)b'(u))=(-na_n(u)b(u))'.
$$

\medbreak

\noindent\textbf{Cas 3.} On a,
$$J (A,B)=\left(-n \frac{a' (u)b (u)}{v^{m+n+1}} - (-m)\frac{a (u)b'(u)}{v^{m+n+1}} \right).$$

\noindent Comme $m$ et $n$ sont supposés strictement positifs, il n'y a pas de terme en $v^{-1}$ et $\rho (J (A,B))=0$.

\end{proof}

Par conséquent, soit $(A,B)$ une paire d'éléments de $k((u))((v))$. Un calcul simple montre que
$
dA\w dB=J(A,B)
$.
De ce fait, 
$$
(u,v)\res^1(dA\w dB)=\rho (J(A,B))du. 
$$
D'après le lemme précédent, il existe $\phi \in k((u))$ tel que ce $1$-résidu est égal à $\phi'(u)du$.  Cette $1$-forme n'a donc pas de terme en $du/u$ et on a 
$$
(u,v)\res^2 (dA \w dB)=0.
$$

\noindent Pour finir, noter que si une $2$-forme formelle $\omega\in \Omega^{2}_{k((u))((v))/k}$ est de la forme $\omega=dA \w dB$ pour $A,B \in k((u))((v))$ et que $(x,y)$ est un couple de séries lié à $(u,v)$ par un changement de variables de le forme (\ref{cv})\footnote{Voir lemme \ref{lemcv}.}, alors $\omega$ est de la forme $dA' \w dB'$ pour $A',B' \in k((x))((y))$. Les séries $A'$ et $B'$ ne sont autres que $A(f(x,y), g(x,y))$ et $B(f(x,y), g(x,y))$. De fait le travail effectué ci-dessus permet de déduire que
$$
(x,y)\res^2 (\omega)=0
$$

\noindent et ce, pour tout couple $(x,y)$ lié à $(u,v)$par un changement de variables de la forme (\ref{cv}).
\section{Topologie de $k((u))[[v]]$}\label{ann_topo}

Le but de cette section est de prouver le lemme \ref{defcv}. Pour ce faire, nous allons introduire quelques notions de topologie sur $k((u))[[v]]$. La première question à se poser est: \emph{de quelle topologie doit-on munir  $k((u))[[v]]$?}
Il pourrait sembler logique de le munir de la topologie associée à la valuation $(v)$-adique, c'est-à-dire, la topologie rendant l'addition continue  et telle que les idéaux $\{(v^n), n\in \N \}$ forment une base de voisinage de $0$. Le défaut d'un tel choix est que pour cette topologie, la suite de terme général ${(u^n)}_{n\in \N}$ diverge. On souhaiterait donc munir $k((u))[[v]]$ d'une topologie qui tiendrait compte à la fois de la valuation $(v)$-adique mais également de la valuation $(u)$-adique sur $k((u))$. 
Pour ce faire, on rappelle que $k((u))[[v]]$ est une limite projective
$$
k((u))[[v]] \cong \lim_{\longleftarrow} k((u))[v]/(v^n).
$$
De fait, si l'on munit $k((u))$ de sa topologie $(u)$-adique, on définit une topologie de limite projective sur $k((u))[[v]]$.
Les ensembles suivants fournissent une base de voisinages de $0$ pour cette topologie.
$$
V_{i_0,\ldots, i_r}:=\left\{
s=\sum_{j\geq 0} s_j(u)v^j \in k((u))[[v]],\ \  \textrm{val}_{(u)} (s_{k}) \geq i_k,\ \  \forall k\in \{0,\ldots, r\}
\right\}.
$$
On rappelle que $\textrm{val}_{(u)}$ désigne la valuation $(u)$-adique sur $k((u))$.
Pour cette topologie, une suite ${(s^{(n)})}_{n\in \N}$ de séries converge vers $0$ si et seulement si elle converge vers $0$ \textit{coordonnée par coordonnée}.
C'est-à-dire:
$$
\lim_{n\rightarrow +\infty} s^{(n)}=0 \quad \Longleftrightarrow
\quad \forall j \in \N,\ \lim_{n\rightarrow +\infty}  s_j^{(n)}(u)=0.
$$

\begin{lem}\label{serie_proj}
Pour cette topologie, une série est convergente si et seulement si son terme général tend vers $0$. 
\end{lem}

\begin{proof}
Soit $(s^{(n)})_{n\in \N}$ une suite qui tend vers $0$ pour la topologie de la limite projective.
Cela signifie que pour tout entier naturel $j$, la suite $(s_j^{(n)})$ d'éléments de $k((u))$ converge vers $0$ pour la topologie $(u)$-adique. 
La topologie $(u)$-adique provenant d'une norme ultramétrique, on en déduit que pour tout entier naturel $j$, la série de terme général $s_j^{(n)}$ converge. Donc la suite des sommes partielles $(\sum_{k=0}^n s^{(k)})_{n\in \N}$ converge coordonnée par coordonnée, elle converge donc pour la topologie de la limite projective.
\end{proof}

\begin{rem}\label{topofine}
La topologie de limite projective est moins fine que la topologie $(v)$-adique. On note par exemple que la topologie $(v)$-adique induit sur $k((u))$ une topologie discrète. De fait, une suite $(s^{(n)})_n$ qui converge $(v)$-adiquement vers une certaine limite $s$, converge vers cette même limite pour la topologie de la limite projective.
\end{rem}

\noindent Nous avons à présent les cartes en main pour démontrer le lemme \ref{defcv}.

\begin{proof}[\textsc{Démonstration du lemme \ref{defcv}}]
Pour commencer, remarquons qu'il suffit de prouver que changement de variable est bien défini sur $k((u))[[v]]$ et induit un isomorphisme local $k((u))[[v]] \rightarrow k((x))[[y]]$. La propriété universelle des corps de fractions permettra ensuite de conclure.

\medbreak  

\noindent \textbf{Étape 1.} Nous allons montrer que la suite $(f^n)_{n\in \N}$ converge vers $0$ pour la topologie de limite projective.
Rappelons que $f$ est un élément de $k((x))[[y]]$ de la forme
$$f=f_0(x)+f_1(x)v+\cdots$$
et que la valuation $(x)$-adique de $f_0$ est égale à $1$.
Soient $n\in \N$ et $k\leq n$, on a 
$$
f^n=f_0^n +f_0^{n-1} f_1 y+f_0^{n-2} (f_1^2+f_2)y^2+\cdots + f_0^{n-k} P_k (f_1,\ldots, f_k)y^k + \cdots,
$$
où $P_k$ désigne un polynôme en les séries de Laurent $f_0, \ldots, f_k$. Ce polynôme ne dépend pas de $n$. Nous donnons ci-dessous, les premiers termes de cette suite de polynômes.
$$
\begin{array}{rclcrcl}
  P_0 & = & 1 & \quad & P_3 & = & f_1^3+ 2f_1 f_2+ f_3 \\
  P_1 & = & f_1 & & P_4 & = & f_1^4+ 2f_1 f_3+ f_2^2+f_4 \\
  P_2 & = & f_1^2 +f_2 & &  P_5 & = & f_1^5 + 2f_1 f_4 + 2f_2 f_3 +f_5 .
\end{array}
$$
Ainsi, étant donné $i\in \N$, pour $n$ assez grand, le coefficient de $y^i$ dans $f^n$ sera égal à $f_0^{n-i}P_i(f_0,\ldots ,f_i)$. Comme $f_0$ est de valuation $(x)$-adique $1$, le coefficient de $y^i$ tend vers $0$ quand $n$ tend vers l'infini. La suite $(f^n)_n$ converge donc vers $0$ pour la topologie de limite projective.
Par conséquent, pour toute série $\varphi(u)\in k((u))$, la série $\varphi(f(x,y))$ converge dans $k((x))[[y]]$.
Remarquons enfin que le coefficient en $y^0$ de la série $\varphi (f(x,y))$ est $\sum_i \varphi_{i}f_0^i$. Cette dernière série est non nulle, car $f_0$ est de valuation $(x)$-adique $1$. De fait, la valuation $(y)$-adique de $\varphi (f(x,y))$ est nulle.

\medbreak

\noindent \textbf{Étape 2.} 
Soit $\psi (u,v) \in k((u))[[v]]$ de la forme $\psi=\sum_{j \geq 0} \psi_j(u)v^j$.
D'après l'étape $1$, pour tout $j\in \N$, la série $\psi_j(f(x,y))$ est bien définie.
Ensuite, comme $g$ est de valuation $(y)$-adique $1$, la série $\psi_j(f)g^j$ est de valuation $(y)$-adique supérieure ou égale à $j$.
Donc, d'après le lemme \ref{serie_proj} et la remarque \ref{topofine}, la série de terme général $\psi_j(f)g^j$ converge dans $k((x))[[y]]$.
Le changement de variables est donc bien défini.
Par ailleurs, on a vu à la fin de l'étape $1$ que la valuation $(y)$-adique d'un coefficient $\psi_j(f(x,y))$ est nulle.
Donc, comme la série $g$ est de valuation $(y)$-adique $1$, on déduit que la valuation $(y)$-adique de $\psi (f(x,y),g(x,y))$ est égale à la valuation $(v)$-adique de $\psi (u,v)$.

\medbreak

\noindent \textbf{Étape 3.} Il nous reste à traiter le cas des $2$-formes différentielles formelles. Soit donc une forme différentielle formelle $\omega=h(u,v)du\w dv$ de valuation $(v)$-adique $n$, montrons que la valuation $(y)$-adique de $h(f,g)df \w dg$ est également $n$. En utilisant les étapes précédentes, cela revient à montrer que $df \w dg$ est de valuation $(y)$-adique nulle. On a,
$$
df\w dg= \left( \frac{\p f}{\p x} \frac{\p g}{\p y} -
\frac{\p f}{\p y} \frac{\p g}{\p x} \right) dx\w dy.
$$

\noindent Calculons les premiers termes de ces produits de dérivées partielles
$$
\begin{array}{rcl}
\frac{\displaystyle \p f}{\displaystyle \p x} \frac{\displaystyle \p g}{\displaystyle \p y} & = & f_0'(x)g_1(x)+ (2f_0'(x)g_2(x)+g_1(x)f_1'(x))y+\cdots \\
 & & \\
\frac{\displaystyle \p f}{\displaystyle \p y} \frac{\displaystyle \p g}{\displaystyle\p x} & = & f_1(x)g_1'(x)y+ (2f_2(x)g_1'(x)+f_1(x)g_2'(y))y^2+\cdots
\end{array}
$$ 
Par définition du changement de variables (\ref{cv}), les séries $f_0'$ et $g_1$ sont non nulles, on en déduit que le jacobien $ \frac{\p f}{\p x} \frac{\p g}{\p y} -
\frac{\p f}{\p y} \frac{\p g}{\p x}$ est de valuation $(y)$-adique nulle, ce qui achève cette démonstration.
\end{proof}

\section{Démonstration du théorème \ref{inv2res} en caractéristique positive}\label{demochiante}

Cette section concerne les formes différentielles formelles. Les notations et définitions utilisées proviennent des sections \ref{ratformel} et \ref{resformel}.

\begin{rem}
Dans cette section nous utilisons le système d'indexage de la notation \ref{indices}.
\end{rem}

Commençons par étudier quelles parties de la preuve du théorème \ref{inv2res},  nécessitent vraiment le fait que le corps $k$ est de caractéristique nulle.
Les lemmes \ref{lem_invCV1} et \ref{valCV} ainsi que la remarque \ref{invCV1} sont valables en caractéristique quelconque. Seule la preuve à proprement parler du théorème, qui commence page \pageref{preuveinv} et se termine page \pageref{finpreuve} fait intervenir des primitives formelles qui n'existent pas toujours en caractéristique positive.
Nous allons donc reprendre l'étude du comportement sous l'action de (CV2) de différentielles de la forme
$$
\omega=\phi  (u)du \w \frac{dy}{y^{n+1}}\quad \textrm{où} \phi \in k((u)) \quad \textrm{et} \quad n\geq 1.
$$

\noindent Soit $N\in \N$, considérons un changement de variables de la forme (CV2):
$$
u=f(x,y)\quad \textrm{avec}\quad f=\sum_{j\geq 0}f_j(x)y^j
$$

\noindent où $f_0$ est de valuation $(x)$-adique $1$.
On suppose de plus que
\begin{equation}\label{minval}
\min_{k=1\ldots n} \{\textrm{val}_{(x)}(f_k)\} =-N,
\end{equation}

\noindent où $\textrm{val}_{(x)}$ désigne la valuation $(x)$-adique sur $k((x))$.

\medbreak

\noindent \textbf{Étape 1.} Si $\omega$ est de la forme $\omega = u^m du\w \frac{dy}{y^{n+1}}$ avec $m\in \N$. On a alors,
$$
\omega = \underbrace{(f'_0(x)+f'_1(x)y+\cdots)}_{\frac{\p f}{\p x}}
\underbrace{(f_0(x)+f_1(x)y+\cdots)^m}_{f^m} dx \w \frac{dy}{y^{n+1}}.
$$
Le $(x,y)$-$1$-résidu de $\omega$ est le coefficient en $y^{n-1}$ de la série $f^m \p f/\p x$. Ce résidu est de la forme
\begin{equation}\label{diese}
(x,y)\res^1_{C,P}(\omega)=P_{m,n}(f_0,\ldots,f_n,f'_0,\ldots,f'_n)dx,
\end{equation}
où $P_{m,n} \in \Z[X_0, \ldots X_n, Y_0,\ldots , Y_n]$ est un polynôme qui ne dépend pas du corps de base, il ne dépend en fait que de $m$ et $n$. Par un raisonnement analogue, le coefficient $p_{m,n}$ en $x^{-1}$ de $P_{m,n}$ est une expression polynomiale en les $f_{i,j}$ avec
$$
-N \leq i \leq N+1 \quad \textrm{et} \quad 0 \leq j \leq n.
$$

\noindent En effet, $P_{m,n}$ est un polynôme en les $f_j$ et $f_j'$ pour $j\in \{0,\ldots, n\}$ ce qui explique l'encadrement de $j$. Pour ce qui en de l'encadrement de $i$, on rappelle que, d'après (\ref{minval}), les séries de Laurent $f_0, \ldots, f_n$ sont de valuation supérieure ou égale à $-N$.
Leurs dérivées sont donc de valuation $(x)$-adique supérieure ou égale à $-N-1$ et donc les termes de degré en $x$ maximal intervenant dans cette expression sont ceux de degré en $x$ égal à $N$. Ces termes peuvent être des $f_{N,j}x^N$ provenant de $f_j(x)$ ou des $(N+1)f_{N+1,j}x^N$ provenant de $f_j'(x)$.
Ainsi, l'indice $i$ est donc toujours inférieur à $N+1$.
Pour finir, d'après la preuve du théorème \ref{inv2res} en caractéristique nulle, on sait que  l'expression polynomiale $p_{m,n}$ s'annule sur l'ensemble $\{f_{1,0} \neq 0\}$. Ainsi, d'après le théorème de prolongement des identités algébriques (\cite{bou} IV.2.3 théorème 2), ce polynôme est nul. Donc, le $(x,y)$-$2$-résidu de $\omega$ est nul.

\medbreak

\noindent \textbf{Étape 2.} Supposons à présent que $\omega$ soit de la forme $\omega=\phi (u)du \w \frac{dy}{y^{n+1}}$ où $\phi=\sum_{m\geq 0}\phi_i u^i$ est une série de Taylor (un élément de $k[[u]]$).
D'après le travail effectué dans l'étape 1, on a 
\begin{equation}\label{seriegolo}
(x,y)\res^1(\omega)=\sum_{m\geq 0} \phi_m P_{m,n}(f_0,\ldots, f_n,f'_0,\ldots,f'_n)dx,
\end{equation}
où les $P_{m,n}$ sont les polynômes définis dans la relation (\ref{diese}) de l'étape 1.
Le $(x,y)$-$1$-résidu de $\omega$ est bien défini, donc la série apparaissant en (\ref{seriegolo}) converge $(x)$-adiquement dans $k((x))$. Par conséquent, la valuation $(x)$-adique de ses termes tend vers l'infini quand $m$ tend vers l'infini, elle est donc positive à partir d'un certain rang $M$.
On a donc
$$
(x,y)\res^1_{C,P}(\omega)=\sum_{m=0}^M \phi_m P_m dx + \underbrace{\sum_{m=M+1}^{+\infty} \phi_m P_m dx}_{\val_{(x)} \geq 0}.
$$ 
Le reste de la série étant de valuation $(x)$-adique positive, son résidu est nul. Quant à la somme de $0$ à $M$, son résidu est nul d'après le résultat obtenu dans l'étape 1 et étendu par linéarité. Ici encore, le $(x,y)$-$2$-résidu de $\omega$ est nul.

\medbreak

\noindent \textbf{Étape 3.}
Supposons maintenant que $\omega$ est de la forme $\omega=\frac{du}{u^m} \w \frac{dy}{y^{n+1}}$, avec $m\in \N^*$. Après changement de variables, on obtient
$$
\omega  = \frac{1}{f^m} \frac{\p f}{\p x}  dx \w \frac{dy}{y^{n+1}}.
$$
Le $(x,y)$-$1$-résidu de $\omega$ est égal au coefficient en $y^n$ de $\frac{1}{f^m}\frac{\p f}{\p x}$ multiplié par $dx$. Nous devons donc étudier la série $\frac{1}{f^m}\frac{\p f}{\p x}$. Commençons par travailler sur $\frac{1}{f^m}$. On a
$$
\begin{array}{rcl}
\frac{\displaystyle 1}{\displaystyle f^m} & = &\frac{\displaystyle 1}{\displaystyle f_0^m {\left( 1+ \frac{f_1}{f_0}y+ \frac{f_2}{f_0}y^2+\cdots
\right)}^m}\\
 & & \\
 & = &  \frac{\displaystyle 1}{\displaystyle f_0^m} {\left(
1+ \frac{\displaystyle U_{1}(f_0,f_1)}{\displaystyle f_0}y+ \cdots + \frac{\displaystyle U_{p}(f_0,\ldots,f_p)}{\displaystyle f_0^2} y^2+\cdots
\right)}^m,
\end{array}
$$
où $U_{p} \in \Z[X_0, \ldots, X_p]$ est un polynôme homogène de degré $p$ qui ne dépend que de $p$. Voici les premiers termes de cette suite de polynômes
$$
\begin{array}{rcl}
U_{1}(X_0,X_1) & = & -X_1 \\
U_{2}(X_0,X_1,X_2) & = & -X_0X_2+X_1^2 \\
U_{3}(X_0,X_1,X_2,X_3) & = &- X_0^2 X_3 +2X_0X_1X_2 -X_1^3 \\
U_{4}(X_0,X_1,X_2,X_3,X_4) & = & -X_0^3X_4+2X_0^2X_1X_3+ X_0^2X_2^2 -3X_0X_1^2X_2+
 X_1^4 .
\end{array}
$$
On développe ensuite le terme à la puissance $m$,
$$
\frac{1}{f^m}= \frac{1}{f_0^m}\left(
1+ \frac{mU_1(f_0,f_1)}{f_0}+ \frac{\left(
    \begin{array}{c}
      m \\ 2
    \end{array}
\right) U_1(f_0,f_1)^2+mU_2(f_0,f_1,f_2)}{f_0^2}y^2+\cdots
\right).
$$

\noindent On introduit alors les notations suivantes 
\begin{equation}\label{Vij}
 \frac{1}{f^m}  = \frac{1}{f_0^m}
\left(1+ \frac{V_{m,1}(f_1)}{f_0}y+\cdots+\frac{V_{m,p}(f_1,\ldots,f_p)}{f_0^p}+\cdots \right),
\end{equation}
où les polynômes $V_{m,p} \in \Z[X_1,\ldots, X_p]$ sont des polynômes homogènes de degré $p$ qui ne dépendent que de $m$ et $p$. Le coefficient de $y^n$ de $\frac{1}{f^m}\frac{\p f }{\p x}$ est donc
$$
C_{m,n}(x):= \frac{1}{f_0^{m}}\left(f'_n+ f'_{n-1}\frac{V_{m,1}(f_0,f_1)}{f_0}+\cdots+f'_0 \frac{V_{m,n}(f_0,\ldots,f_n)}{f_0^n} \right).
$$

\noindent Pour tout entier $k$ appartenant à $\{1,\ldots, n\}$, on pose 
$$S_{m,n,k}(f_0,\ldots , f_k):=f_0^{n-k}V_{m,k}(f_0, \ldots, f_k)$$
et $S_{m,n,0}(f_0):=f_0^n$. Les polynômes, $S_{m,n,k}$ sont homogènes de degré $n$ et
\begin{equation}\label{Cmn}
C_{m,n}(x):=\underbrace{\frac{1}{f_0^{m+n}}}_{A_{m,n}(x)}
\underbrace{\sum_{k=0}^n f_{n-k}' S_{m,n,k}(f_0,\ldots, f_k)}_{B_{m,n}(x)}.
\end{equation}

\noindent Nous allons étudier $A_{m,n}$ séparément. On rappelle que $f_0$ était de valuation $(x)$-adique $1$, c'est-à-dire que $f_0(x)=f_{1,0}x+f_{2,0}x^2+\cdots$. En reprenant le calcul effectué précédemment pour $\frac{1}{f^m}$, on obtient
$$
\begin{array}{rcl}
A_{m,n}(x)  & = & \frac{\displaystyle 1}{\displaystyle  f_{1,0}^{m+n}} \Bigg(
1+ \frac{\displaystyle V_{m+n,1}(f_{1,0},f_{2,0})}{\displaystyle f_{1,0}}x
+ \cdots  \\
 & & \\
 & & \qquad  \cdots +
\frac{\displaystyle V_{m+n,p-1}(f_{1,0},\ldots , f_{p,0})}{\displaystyle f_{1,0}^{p-1}}x^{p-1}+\cdots
\Bigg),
\end{array}
$$
où les polynômes $V_{i,j}$ sont ceux introduits dans l'expression (\ref{Vij}).
Rappelons que l'objectif initial est de montrer le $(x,y)$-$2$-résidu, qui est le coefficient $c_{m,n,-1}$ de $x^{-1}$ dans $C_{m,n}$, s'obtient comme une expression polynomiale en un nombre fini des coefficients $f_{i,j}$ de $f$.
Nous allons voir quels coefficients interviennent.

\medbreak

\noindent \textbf{Dans $A_{m,n}$.} Comme les polynômes $S_{m,n,k}$ sont de degré $n$ pour tout entier $k\in \{1,\ldots, n\}$, la valuation $(x)$-adique de $B_{m,n}$ vérifie
$$
\textrm{val}_{(x)} (B_{m,n}) \geq -nN -(N+1) = -(n+1)N-1.
$$ 
Le $-nN$ est la contribution de $S_{m,n,k}$ et le $-(N+1)$ celle de $f_{n-k}'$.
Par conséquent, les termes de $A_{m,n}$ intervenant dans le calcul de $c_{m,n,-1}$ sont de degré au plus $N(n+1)$. En reprenant l'expression de $A_{m,n}$ ci-dessus, on voit que les termes de degré inférieur à $N(n+1)$ font intervenir les $f_{i,0}$ avec $-N \leq i\leq N(n+1)+1$.

\medbreak

\noindent \textbf{Dans $B_{m,n}$.}
Comme la série de Laurent $A_{m,n}$ est de valuation $-m-n$, les termes de $B_{m,n}$ intervenant dans le calcul de ce $2$-résidu $c_{m,n,-1}$ sont de degré au plus $m+n-1$.
D'après l'expression de $B_{m,n}$, ces termes font intervenir les $f_{i,j}$ avec $i\leq m+n$ (on ajoute $1$ du fait de la présence des dérivées $f_i'$ dans l'expression). 

\medbreak

\noindent \textbf{Conclusion.} Les coefficients $f_{i,j}$ intervenant dans le calcul de $c_{m,n,-1}$ sont ceux dont les indices vérifient
$$
\left\{ \begin{array}{ccccl}
-N & \leq & i & \leq & \max\{m+n, (n+1)N +1 \}   \\
0 & \leq & j & \leq & n .
\end{array} \right.
$$
On note $I$ cet ensemble de paires d'indices. Il existe donc un polynôme $Q_{m,n}(X_{i,j}) \in \Z[X_{i,j}|(i,j)\in I]$ et un entier $M$ qui ne dépendent pas du corps de base $k$ et tels que
$$
c_{m,n,-1}=(x,y)\res^2(\omega)=\frac{1}{f_{1,0}^M}Q(f_{i,j} | (i,j)\in I).
$$
Par un raisonnement analogue à celui effectué à la fin de l'étape 1 et faisant intervenir le théorème de prolongement des identités, on montre que ce résidu est nul.

Ainsi, nous avons montré l'invariance du $2$-résidu d'une $2$-forme formelle sous l'action d'un changement de variables de la forme $u=f(x,y)$ tel que les $n+1$ premiers coefficients $f_j(u)$ de $f$ sont de valuation supérieure à $-N$. Ce résultat est valable pour tout entier $N$, donc pour tout changement de variables de la forme (CV2), ce qui conclut la démonstration.

\chapter{Indépendance des valuations}\label{annexeapprox}

L'objectif de ce chapitre est de donner une preuve de la proposition \ref{approx}. Nous allons en fait démontrer un résultat en dimension quelconque, ce qui n'augmente pas le niveau de difficulté de la preuve.

\begin{prop}\label{movingperso}
  Soient $X$ une variété quasi-projective lisse irréductible de dimension supérieure ou égale à $2$ au-dessus d'un corps quelconque $k$ et $Y$ une sous-variété de $X$ de codimension $1$. Soient $P_1,\ldots, P_n$ une famille de points fermés de $X\smallsetminus Y$ et $Q_1, \ldots, Q_s$ une famille de points fermés de $Y$.
Alors, il existe une uniformisante $v\in \mathcal{O}_{X,Y}$ telle que le support du diviseur principal $(v)$ évite les points $P_1, \ldots, P_n$ et le support de $(v)-Y$  évite $Q_1, \ldots, Q_s$.
\end{prop}

\noindent La démonstration suivante a été suggérée par Gerhard Frey.

\begin{proof}
Commençons par remarquer que l'on peut se contenter de démontrer le résultat dans le cas où $Y$ est irréductible. Le cas général se déduira de ce cas particulier en raisonnant sur les composantes irréductibles de $Y$. 
On suppose donc désormais que $Y$ est irréductible.

Dans un premier temps nous allons nous ramener au cas d'une variété affine.
Il existe un ouvert affine de $X$ contenant tous les points $P_1,\ldots,\ \! P_r,\ \!  Q_1, \ldots ,\ \!  Q_s$. L'intersection de cet ouvert avec $Y$ est non vide car contient les points $Q_i$.
Aussi, quitte à se restreindre à cet ouvert, on peut supposer que $X$ est affine.
Dans ce qui suit, nous noterons $I_{Y}$ l'idéal de $k[X]$ associé à $Y$ et $\mathfrak{m}_{P_1},\ldots,\ \!  \mathfrak{m}_{P_r},\ \!  \mathfrak{m}_{Q_1}, \ldots,\ \!  \mathfrak{m}_{Q_s}$ les idéaux maximaux correspondant respectivement aux points $P_1,\ldots,\ \! P_r$ et $Q_1, \ldots,\ \!  Q_s$.
On notera également $m_1, \ldots, m_s$ les multiplicités respectives de $Y$ en $Q_1, \ldots, Q_s$. On rappelle que la multiplicité de $Y$ en un point $Q$ est le plus petit entier $m$ tel que
$$
I_Y\subset \mathfrak{m}_Q^m.
$$

\medbreak

\noindent \textbf{Reformulation du problème.}
On cherche une fonction $v\in k[X]$ vérifiant les propriétés suivantes.
\begin{enumerate}
\item\label{crut} La fonction $v$ appartient à $I_Y\smallsetminus I_Y^2$. 
\item\label{crat} Pour tout entier $i$ appartenant à $\{1, \ldots, r\}$, la fonction $v$ n'appartient pas à $\mathfrak{m}_{P_i}$.  
\item\label{crot} Pour tout entier $j$ appartenant à $\{1, \ldots, s\}$, la fonction $v$ n'appartient pas à $\mathfrak{m}_{Q_j}^{m_j+1}$.
\end{enumerate}

\noindent Le premier critère permet d'affirmer que $v$  est bien une uniformisante de $\mathcal{O}_{X,Y}$, le second signifie que $v$ ne s'annule en aucun des $P_i$, et le troisième équivaut au fait que $Y$ soit l'unique zéro de $v$ au voisinage de $Q_i$.
Dans ce qui suit, lorsque $P$ est un point fermé de $X$ et $f\in \mathcal{O}_{X,P}$ on notera $f(P)$ l'image de $f$ dans le corps résiduel $k(P)$.

\medbreak

\noindent \textbf{Étape 0.} Pour toute famille finie de points fermés deux à deux distincts $A_1,\ldots, A_n$ de $X$ et toute famille d'entiers naturels $p_2, \ldots, p_n$, il existe une fonction $g$ appartenant à l'idéal produit $\mathfrak{m}_{A_2}^{p_2}\cdots \mathfrak{m}_{A_{n}}^{p_n}$ telle que 
$$g(A_1)=1.$$
Pour le montrer, il suffit de constater que les idéaux $\mathfrak{m}_{A_1}$ et $\mathfrak{m}_{A_2}^{p_2}\cdots \mathfrak{m}_{A_n}^{p_n}$ sont étrangers, c'est-à-dire que
$$\mathfrak{m}_{A_1}+\mathfrak{m}_{A_2}^{p_2}\cdots
\mathfrak{m}_{A_n}^{p_n}=k[X].$$
De fait, il existe une fonction $f$ appartenant à $\mathfrak{m}_{A_1}$ et une fonction $g$ appartenant à $\mathfrak{m}_{A_2}^{p_2}\cdots
\mathfrak{m}_{A_n}^{p_n}$ telles que $f+g=1$. Ainsi, on a bien $g(A_1)=1$.

\medbreak

\noindent \textbf{Étape 1.}\ (\emph{Construction d'une fonction qui vérifie \ref{crot}}).
 Soit $j_0$ un entier appartenant à $\{1,\ldots, s\}$. Par définition de la multiplicité, il existe une fonction $f$ appartenant à l'idéal $I_Y$ et n'appartenant pas à $\mathfrak{m}_{Q_{j_0}}^{m_{j_0}+1}$.
Soit $f_{j_0}$ une telle fonction, d'après l'étape précédente il existe une fonction $h_{j_0}$ régulière sur $X$ telle que $h_{j_0}$ est un élément de l'idéal produit $\mathfrak{m}_{Q_1}^{m_1}\cdots \widehat{\mathfrak{m}_{Q_{j_0}}}\cdots \mathfrak{m}_{Q_s}^{m_s}$ et $h_{j_0}$ n'appartient pas à l'idéal $\mathfrak{m}_{Q_{j_0}}$.
Posons alors
$$v_{j_0}:=f_{j_0} h_{j_0} \quad \textrm{et} \quad
v:=\sum_{j=1}^s v_{j}. 
$$
La fonction $v$ appartient à l'idéal $I_Y$, mais n'appartient à aucun $\mathfrak{m}_{Q_j}^{m_j+1}$. En effet, soit $j_0$ un entier appartenant à $\{1,\ldots,s\}$, alors
$$
v_{j_0}\in \mathfrak{m}_{Q_{j_0}}^{m_{j_0}}\smallsetminus \mathfrak{m}_{Q_{j_0}}^{m_{j_0}+1}\quad \textrm{et}
\quad \forall j\neq j_0, \ v_j \in \mathfrak{m}_{Q_{j_0}}^{m_{j_0}+1}.
$$
Remarquons au passage que la fonction ainsi construite n'appartient pas à $I_Y^2$, car sinon elle appartiendrait à $\mathfrak{m}_{Q_j}^{2m_j}$ pour tout entier $j$ appartenant à $\{1, \ldots, s\}$.

\medbreak

\noindent \textbf{Étape 2.} Nous allons montrer par récurrence sur $r$ que, quitte à ajouter à $v$ un élément de $I_Y^2$, ce qui n'aura aucune conséquence sur les propriétés acquises dans l'étape précédente, on peut faire en sorte que $v$ ne s'annule en aucun des $P_i$. 

\smallbreak

\noindent \textbf{Pour $\mathbf{r=1}$},
si $v(P_1)$ est non nul, c'est terminé. Sinon, comme
 $P_1$ n'est pas contenu dans  $Y$, l'idéal $I_Y$ ne contient pas $\mathfrak{m}_{P_1}$.
De même, l'idéal $I_Y^2$ ne contient pas $\mathfrak{m}_{P_1}$, car $\mathfrak{m}_{P_1}$ est radical. De fait, les idéaux $I_Y^2$ et $\mathfrak{m}_{P_1}$ sont étrangers, il existe donc une fonction $f$ appartenant à $I_Y^2$ et une fonction $g$ appartenant à $\mathfrak{m}_{P_1}$ telles que $f+g=1$.
On remplace alors $v$ par $v+f$ et cette nouvelle fonction ne s'annule plus en $P_1$.

\smallbreak

\noindent \textbf{Soit $\mathbf{r\geq 1}$}, supposons la propriété vraie au rang $r$ et montrons qu'elle est vérifiée au rang $r+1$. Si $v(P_{r+1})$ est non nul, c'est terminé. Sinon, les idéaux $I_Y^2 \mathfrak{m}_{P_1}\cdots \mathfrak{m}_{P_r}$ et $\mathfrak{m}_{P_{r+1}}$ sont étrangers, on en déduit l'existence d'une fonction $f$ appartenant à $I_Y^2 \mathfrak{m}_{P_1}\cdots \mathfrak{m}_{P_r}$ qui ne s'annule pas en $P_{r+1}$ et on remplace $v$ par $v+f$.
\end{proof}

\chapter{Complément d'algèbre linéaire}\label{linalg}

Le but de cette annexe est de fournir quelques résultats d'algèbre linéaire utilisés à la fin du chapitre \ref{chapdiff}.
Ces résultats sont des exercices relativement élémentaires. Nous avons cependant choisi de les démontrer ici, faute de références.

Dans ce qui suit, $E$ et $F$ désignent des espaces vectoriels sur un corps quelconque $k$. On dispose de plus sur de formes bilinéaires $E$ et $F$ notées respectivement $\langle .\ ,\ .\rangle_E$ et $\langle  .\ ,\ .\rangle_F$.
On peut alors construire une forme bilinéaire sur $E\otimes_k F$ qui soit canoniquement associée à $\langle   .\ ,\ . \rangle_E$ et $\langle   .\ ,\ . \rangle_F$. Pour ce faire, on la définit sur les tenseurs élémentaires par
$$
\langle  .\ ,\ . \rangle_{E\otimes F}\ : \left\{
\begin{array}{ccc}
E \otimes F \times E \otimes F &  \rightarrow & k \\
(a\otimes b, a' \otimes b') & \mapsto & \langle  a,a'\rangle_E . \langle  b,b'\rangle_F 
\end{array}
\right.
$$
et on l'étend ensuite par linéarité.
On constate immédiatement que si $(e_i)_{i \in I}$ et $(f_j)_{j\in J}$ sont des bases orthogonales respectives de $E$ et $F$, alors $(e_i \otimes f_j)_{(i,j)\in I \times J}$ est une base orthogonale de $E\otimes F$.

\begin{lem}
Si les formes bilinéaires $\langle .\ ,\ .\rangle_E$ et $\langle .\ ,\ .\rangle_F$ sont non dégénérées, alors il en est de même pour $\langle .\ ,\ .\rangle_{E\otimes F}$. De même, si $\langle .\ ,\ .\rangle_E$ et $\langle .\ ,\ .\rangle_F$ sont symétriques, alors $\langle .\ ,\ .\rangle_{E\otimes F}$ l'est également.
\end{lem}

\begin{proof}
  Le fait que la symétrie de $\langle .\ ,\ .\rangle_E$ et $\langle .\ ,\ .\rangle_F$ entraîne celle de $\langle .\ ,\ .\rangle_{E\otimes F}$ est évident. Supposons que $\langle .\ ,\ .\rangle_E$ et $\langle .\ ,\ .\rangle_F$ soient non dégénérées. Soient $(e_i)_{i \in I}$ et $(f_j)_{j\in J}$ des bases orthogonales respectives de $E$ et $F$ et soit $u$ un élément de $E \otimes F$ tel que la forme linéaire $\langle u,\ .\rangle_{E\otimes F}$ soit identiquement nulle.
Le vecteur $u$ est de la forme
$$
u= \sum_{i,j}u_{i,j}e_i \otimes f_j.
$$
Par ailleurs, pour tout couple $(i,j)$, on a
$$
\langle u, e_i \otimes f_j\rangle= 0, \quad \textrm{or} \quad \langle u, e_i \otimes f_j\rangle=u_{i,j}.
$$
On en déduit que $u$ est nul et que la forme bilinéaire $\langle .\ ,\ .\rangle_{E\otimes F}$ est non dégénérée.
\end{proof}

On suppose désormais que les formes bilinéaires $\langle  .\ ,\ . \rangle_E$ et $\langle  .\ ,\ . \rangle_E$ sont non dégénérées.

\begin{lem}\label{orthotens}
Supposons que $E$ et $F$ soient de dimension finie.
Soient $A$ et $B$ deux sous-espaces vectoriels respectifs de $E$ et $F$, alors 
$$
(A\otimes B)^{\bot}= A^{\bot} \otimes F + E \otimes B ^{\bot}.
$$ 
\end{lem}

\begin{proof}
L'inclusion $\supseteq$ est élémentaire. En effet on prend $a,a',b,f$ appartenant respectivement à $A, A^{\bot},B$ et $F$. On a
$$
\langle a\otimes b, a' \otimes f\rangle_{E \otimes F} = \underbrace{\langle a,a'\rangle}_{=0} \langle b,f\rangle=0.
$$

\noindent On en déduit que $A^{\bot}\otimes F$ est inclus dans $(A\otimes B)^{\bot}$. Par un raisonnement identique, on montre ensuite que $E\otimes B^{\bot}$ est inclus dans $(A\otimes B)^{\bot}$.

Pour l'inclusion réciproque, commençons par montrer que
$$
A^{\bot}\otimes F \cap E \otimes B^{\bot} = A^{\bot}\otimes B^{\bot}.
$$

L'inclusion $\supseteq$ est immédiate. Montrons donc l'inclusion réciproque.
Soient $(e_i)_{i\in I_0}$ et $(f_j)_{j \in J_0}$ des bases respectives de $A^{\bot}$ et $B^{\bot}$ complétées en des bases $(e_i)_{i \in I}$ et $(f_j)_{j \in J}$ de $E$ et $F$.
Soit $s$ un vecteur de $A^{\bot}\otimes F \cap E \otimes B^{\bot}$. Décomposons le dans la base $(e_i\otimes f_j)_{(i,j)\in I \times J}$.
$$
s=\sum_{(i,j)\in I \times J}s_{ij}\ \! e_i \otimes f_j.
$$

\noindent Comme $s\in A^{\bot}\otimes F$ (resp. $s\in E \otimes B^{\bot}$), les $s_{ij}$ sont nuls pour $i \notin I_0$ (resp. pour $j \notin J_0$), ce qui prouve l'inclusion réciproque.

Pour finir, posons
$$
\begin{array}{rclcrcl}
m & := & \dim(E) &  & a & := & \dim (A) \\
n & := & \dim(F) &  & b & := & \dim (B),
\end{array}
$$

\noindent et calculons la dimension de $A^{\bot}\otimes F + E \otimes B^{\bot}$.
$$
\begin{array}{rcl}
  \dim (A^{\bot}\otimes F + E \otimes B^{\bot}) & = & \dim (A^{\bot}\otimes F)+
\dim (E \otimes B^{\bot})- \dim (A^{\bot}\otimes B^{\bot})\\
 & = & n(m-a)+m(n-b)-(m-a)(n-b) \\
 & = & mn-ab \\
 & = & \dim (A^{\bot}\otimes B^{\bot}).
\end{array}
$$
\end{proof}

Nous pouvons à présent énoncer un résultat fort utile dans le chapitre \ref{chapdiff}, puisque c'est celui qui nous permet de prouver que l'orthogonal d'un code fonctionnel sur une surface n'est pas toujours fonctionnel. En particulier, il ne l'est jamais lorsque la surface est un produit de deux courbes.

\begin{cor}\label{tens}
Si $A$ et $B$ sont des sous-espaces non triviaux\footnote{Par \emph{non triviaux}, on entend que $A$ (resp. $B$) est non nul et strictement inclus dans $E$ (resp. $F$).} respectifs de $E$ et $F$, alors le sous-espace $(A\otimes B)^{\bot}$ de $E \otimes F$ n'est pas un produit tensoriel élémentaire $U\otimes V$ de sous-espaces de $E$ et $F$. 
\end{cor}

\begin{proof}
C'est une conséquence immédiate du lemme \ref{orthotens} et du lemme \ref{nonelem} qui suit.
\end{proof}

\begin{lem}\label{nonelem}
Soient $E$ et $F$ deux espaces vectoriels sur un corps $k$ quelconque.
Soient $A$ et $B$ des sous-espaces respectifs non triviaux de $E$ et $F$. Alors, le sous-espace $A\otimes_k F+E \otimes_k B$ de $E \otimes_k F$ ne peut pas être écrit sous la forme d'un produit tensoriel élémentaire $U \otimes_k V$.  
\end{lem}

\begin{proof}
Raisonnons par l'absurde et supposons qu'il existe $U$ et $V$, sous-espaces 
respectifs de $E$ et $F$ tels que 
\begin{equation}\label{sodome}
A\otimes_k F+E \otimes_k B=U \otimes_k V.
\end{equation}

Supposons que $U$ ne soit pas inclus dans $A$.
Soient $(e_i)_{i \in I_0}$ et $(f_j)_{j\in J_0}$ des bases respectives de $A$ et $B$ complétées en des bases $(e_i)_{i \in I}$ et $(f_j)_{j\in J}$ de $E$ et $F$.
Soit également $u \in U \smallsetminus A$. Pour tout vecteur $v\in V$, on a
$$
u \otimes v= \sum_{(i,j)\in I \times J}u_i v_j\ \! e_i\otimes f_j.
$$
D'après (\ref{sodome}), le produit $u_iv_j$ est nul pour tout couple $(i,j)$ tel que $i\notin I_0$ et $j\notin J_0$. Comme par hypothèse $u$ n'est pas dans $A$, il existe au moins un $i_1\notin I_0$ tel que $u_{i_1}$ est non nul. Donc pour tout $j \notin J_0$, on a $u_{i_1}v_j=0$, donc $v_j=0$.
Par conséquent, $v$ appartient à $B$ et ce pour tout $v\in V$. Donc
\begin{equation}\label{raslecul}
V \subseteq B, \quad \textrm{donc}\quad U\otimes V \subseteq E\otimes B.
\end{equation}

\noindent Soient maintenant $f$ un vecteur de $F$ n'appartenant pas à $B$ et $a$ un vecteur non nul de $A$. Alors $a\otimes f$ appartient à $A\otimes F$ mais pas à $E\otimes B$
et d'après (\ref{raslecul}), il n'appartient pas à $U\otimes V$ ce qui contredit l'hypothèse de départ (\ref{sodome}).

Si maintenant $U$ est inclus dans $A$, alors $U \otimes V$ est inclus dans $A \otimes F$ et on aboutit à une contradiction en réalisant le raisonnement symétrique de celui qui vient d'être effectué.
\end{proof}

\newpage
\thispagestyle{empty}

\chapter{Construction de codes fonctionnels}\label{annexelachaud}

Dans \cite{lachaud2}, Lachaud présente une autre procédé de construction de codes fonctionnels sensiblement différent de celui qui est présenté dans le chapitre \ref{chapdiff}.
Les codes ainsi construits sont en général appelés codes de Reed-Müller projectifs.

\section{Construction}
On se place dans l'espace projectif $\P^r_{\F_q}$ muni d'un système de coordonnées homogènes $(X_0, \ldots, X_n)$. Pour tout entier naturel $d$
on note $\mathcal{F}_d$, l'espace
$$
\mathcal{F}_d:=\left\{f \in \F_q[X_0, \ldots , X_n],\ f\ \textrm{homogène de degré}\ d  \right\}\cup \{0\}.
$$
En d'autres termes, $\mathcal{F}_d$ est l'espace des sections globales du faisceau $\mathcal{O}_{\P^r}(d)$.
Rappelons que, s'il est possible de définir le lieu d'annulation dans $\P^r$ d'un élément non nul de $\mathcal{F}_d$, l'évaluation d'un élément de $\mathcal{F}_d$ en un point de $\P^r$ n'a pas de sens. D'une certaine façon, la définition qui suit consiste à lui en donner un.

\begin{defn}[Lachaud 1988]\label{evproj}
Soit $P$ un point rationnel de $\P^r$ de coordonnées homogènes $(x_0: \cdots : x_n)$. Soit $h$ le plus petit entier compris entre $0$ et $n$ tel que $x_h$ soit non nul. Pour tout entier naturel $d$, on définit l'application d'évaluation en $P$ par
$$
\textrm{ev}_P:\left\{
\begin{array}{ccc}
\mathcal{F}_d & \rightarrow & \F_q \\
 f & \rightarrow & \frac{\displaystyle f(x_0,\ldots ,x_n)}{\displaystyle x_h^d} .
\end{array}
\right.
$$ 
\end{defn}

\noindent L'application est bien définie, mais n'est pas canonique, elle dépend du choix d'un système de coordonnés homogènes.
À partir de cette application, on peut construire des codes correcteurs d'erreurs de la façon suivante.

\begin{defn}\label{codelachaud}
Soient $X$ une sous-variété fermée lisse absolument irréductible de $\P^r_{\F_q}$ et $P_1, \ldots, P_n$, l'ensemble des points rationnels de $X$. Pour tout entier naturel $d$, le code $C_d (X)$ est l'image de l'application
$$
c: \left\{
\begin{array}{ccc}
\mathcal{F}_d & \rightarrow & \F_q^n \\
 f & \mapsto & (ev_{P_1}(f), \ldots, ev_{P_n}(f)).     
\end{array}
\right.
$$
\end{defn}

Bien que moins canonique que la construction présentée dans le chapitre \ref{chapdiff}, cette approche présente de nombreux avantages. Tout d'abord, elle permet une évaluation en \textbf{tous les points de la variété}, sans avoir à se soucier des questions de définitions. Par rapport à la définition du chapitre \ref{chapdiff}, on évite ici les restrictions du type ``les points $P_1, \ldots, P_n$ évitent le support de $G$''.
Ensuite, l'évaluation de la distance minimale et de la distribution des poids du code $C_d(X)$ se traduit sous la forme d'un problème de comptage explicite de tous les points rationnels de variétés projectives\footnote{Éventuellement réduites ou réductibles.}.


Toutefois, pour le cadre qui nous intéresse, à savoir celui de pouvoir construire l'orthogonal du code fonctionnel à partir de différentielles, cette construction n'est pas optimale. En effet, le défaut de canonicité des applications d'évaluation empêche toute possibilité de reproduction du raisonnement de la preuve du théorème \ref{orthocode}. L'obstruction est liée au fait que les éléments de $\mathcal{F}_d$ ne sont pas des fonctions. On peut localement les identifier à des fonctions rationnelles sur $X$, c'est ce qui est fait dans la définition de l'application $\textrm{ev}_P$, mais cette identification dépend du point $P$. D'une certaine façon, c'est comme si l'on évaluait en des points différents des fonctions différentes. On ne peut donc plus reproduire le raisonnement consistant à appliquer la formule des résidus à la $2$-forme $f\omega$.

\section{Essentiellement, c'est la même chose}\label{idem}

Par la méthode de construction présentée dans le chapitre \ref{chapdiff}, on peut construire un code fonctionnel isomorphe au code $C_d(X)$.
Pour ce faire, on définit $H_{\infty}$, l'hyperplan de $\P^r$ d'équation $X_0=0$.
On appelle $i$ l'inclusion canonique $i: X \hookrightarrow \P^r$
et on pose
$$L:=i^*H_{\infty} \in \div (X).$$
Pour tout entier naturel $d$, les faisceaux $\L (dL)$ et $i^*\mathcal{O}_{\P^r} (d)$ sont isomorphes. Aussi, si les points $P_1, \ldots, P_n$ évitent le support de $L$, alors les codes $C_{L,X} (P_1+ \cdots +P_n, dL)$ et $C_d(X)$ sont isomorphes.
Pour se défaire de la condition \textit{évitent les point du support de $L$}, il suffit d'utiliser le \textit{moving lemma} (\cite{sch1}  III.1.3 théorème 1 et annexe \ref{annexeapprox}). En vertu de ce résultat, on sait qu'il existe $G\in \div (X)$ dont le support évite les points $P_1, \ldots, P_n$ et linéairement équivalent à $dL$.
Cette fois-ci, on est assuré de la bonne définition du code $C_{L,X} (P_1+ \cdots +P_n, G)$ et le faisceau $\L (G)$ est isomorphe à $i^* \mathcal{O}_{\P^r} (d)$.
On obtient donc l'isomorphisme de codes
$$
C_{L,X}(P_1+ \cdots, P_n,G) \cong C_d (X).
$$

\noindent On peut donc construire un code fonctionnel (provenant de la définition de \cite{manin} donnée dans le chapitre \ref{chapdiff}) isomorphe à $C_d(X)$. Malheureusement la construction d'un diviseur $G$ dont le support évite tous les points rationnels de $X$ n'est pas toujours évidente. La définition \ref{codelachaud} reste donc commode car elle fournit une construction explicite simple.
\begin{rem}
Notons que l'on a identifié les éléments de $\mathcal{F}_d$ restreints à $X$ à des sections globales du faisceau $i^* \mathcal{O}_{\P^n}(d)$. Comme on l'a vu dans le chapitre \ref{chapreal}, si $X$ est projectivement normale\footnote{Par exemple si c'est une intersection complète lisse de $\P^r$.} en respect au plongement $X \hookrightarrow \P^r$, alors le faisceau $i^* \mathcal{O}_{\P^r}(d)$ s'identifie au faisceau $\mathcal{O}_X (d)$ (voir \cite{H} II.5 ex 14).   
\end{rem}

\chapter{Points en position générale}\label{annexegen}

Cette annexe contient les démonstrations de lemmes techniques utilisés dans la section \ref{approch1} du chapitre \ref{chapreal}.

\medbreak

\noindent \textbf{Lemme \ref{trucalacon2}.}
\textit{Soient $r$ et $m$ deux entiers naturels avec $r\geq 1$, alors toute famille de $rm+2$ points rationnels distincts de $\P^N$ appartenant à une même courbe de degré $r$ est $m$-liée.}

\begin{proof}
  Soient $P_1, \ldots, P_{rm+2}$ une telle famille de points et $C$ la courbe de degré $r$ qui les contient.
D'après le théorème de Bezout, une hypersurface de degré $m$ qui ne contient pas $C$ l'intersecte en au plus $rm$ points géométriques.
Par conséquent, il n'existe pas d'hypersurface contenant tous ces points sauf un.
\end{proof}

\medbreak

\noindent \textbf{Lemme \ref{kipu2}.}
\textit{  Soit $m$ un entier naturel.
\begin{enumerate}
\item\label{my} Si $m+2$ points rationnels distincts de $\P^N$ sont $m$-liés, alors ils sont alignés. 
\item\label{cock} Tout $r$-uplet de points rationnels deux à deux distincts de $\P^N$ avec $r\leq m+1$ est en position $m$-générale.
\end{enumerate}}

\begin{proof}
\textbf{Preuve de (\ref{my}).}
Soient $P_1,\ldots, P_{m+2}$ des points distincts de $\P^N$ qui sont $m$-liés.
Supposons qu'ils ne soient pas alignés.
Il existe donc un hyperplan $H$ qui évite $P_{m+2}$ et contient au moins deux points parmi $P_1, \ldots, P_{m+1}$.
Quitte à réordonner les indices des points on peut supposer qu'il existe $l\geq 2$ tel que $P_1, \ldots, P_l$ appartiennent à $H$.
\begin{itemize}
\item[\textbullet] Si $l=m+1$, on dispose d'une hypersurface (en l'occurrence l'hyperplan $H$) qui contient tous les $P_i$ sauf $P_{2m+2}$.

\item[\textbullet] Sinon, pour tout $l+1\leq j \leq m+1$, il existe un hyperplan $H_j$ qui contient $P_j$ et évite $P_{m+2}$.
\end{itemize}

\smallbreak

\noindent L'hypersurface $H\cup H_{l+1} \cup \cdots \cup H_{m+1}$ est de degré inférieur ou égal à $m$ (car $l \geq 2$), contient tous les $P_i$ sauf $P_{m+2}$.
Par symétrie, on en déduit l'existence pour tout $i_0 \in \{1, \ldots, m+2\}$ d'une hypersurface de degré $m$ qui contient tous les $P_i$ sauf $P_{i_0}$.

\medbreak

\noindent \textbf{Preuve de (\ref{cock}).}
Soient $P_1, \ldots, P_r$ des points distincts de $\P^N$ avec $r\leq m+1$.
En réutilisant la technique présentée dans la preuve de (\ref{my}), on peut construire à l'aide de réunions d'hyperplans, des hypersurfaces de degré $m$ qui interpolent tous les points sauf un. On en déduit que ces points sont en position $m$-générale.
\end{proof}

\medbreak

\noindent \textbf{Lemme \ref{chiasse2}.}
\textit{Soient $m$ et $r$ deux entiers naturels.
\begin{enumerate}
\item\label{butt2} Si $r\leq 2m+1$, alors une famille de $r$ points distincts de $\P^N$  telle que $m+2$ d'entre eux sont non alignés est en position $m$-générale.
\item\label{your2} Soient $P_1, \ldots, P_{2m+2}$ un $(2m+2)$-uplet de points rationnels distincts de $\P^N$ tels que $m+2$ d'entre eux ne sont pas alignés. Alors ces points sont $m$-liés, si et seulement s'ils appartiennent à une même conique plane.
\end{enumerate}}

\begin{proof}
\textbf{Preuve de (\ref{butt2}).}
Nous allons démontrer le résultat par récurrence sur $m$. 

Pour $m=0$, c'est évident car un point de $\P^N$ est toujours en position $0$-générale.

Soit $m\geq 0$, supposons que la propriété soit vérifiée au rang $m$ et montrons qu'elle l'est au rang $m+1$.
Soit $P_1, \ldots, P_{2m+3}$ une famille de points rationnels de $\P^N$ telle que $m+3$ d'entre eux ne sont pas alignés.
Soit $\mathfrak{L}$ l'ensemble des droites contenant $P_{2m+3}$ et au moins un point parmi $P_1, \ldots, P_{2m+2}$.
On choisit un élément $L_1$ de $ \mathfrak{L}$ contenant un nombre maximal de points parmi les $P_i$.
On choisit ensuite une droite $L_2$ de $ \mathfrak{L}\smallsetminus L_1$ contenant un nombre maximal de points parmi les $P_i$. 
Les droites $L_1$ et $L_2$ sont distinctes et se croisent en $P_{2m+3}$.
Quitte à réordonner les points $P_i$, on peut supposer que $P_{1}\in L_1$ et $P_{2}\in L_2$.

Ensuite, montrons que, par hypothèse, les droites $L_1$ et $L_2$ contiennent chacune au plus $m+2$ points parmi les $P_i$ et que les autres droites de $\mathfrak{L}$ en contiennent au plus $m+1$.
On distingue deux situations.
\begin{enumerate}
\item Les droites $L_1$ et $L_2$ contiennent toutes deux $m+2$ points parmi les $P_i$. Comme par définition, ces droites contiennent toutes deux le point $P_{2m+3}$, l'ensemble des $P_i$ contenus dans $L_1 \cup L_2$ est égal à l'ensemble ${P_1, \ldots, P_n}$.
Aussi, dans cette situation l'ensemble $\mathfrak{L}$ est égal à $ \{L_1,L_2\}$ et le résultat attendu est trivialement vérifié.
\item La droite $L_2$ contient moins de $m+1$ points parmi les $P_i$ et par définition de $L_2$, les éléments de $\mathfrak{L}\smallsetminus \{L_1,L_2\}$ contiennent tous moins de $m+1$ points parmi les $P_i$.
\end{enumerate}

Par conséquent, d'après l'hypothèse de récurrence, les points $P_3, \ldots, P_{2m+3}$ sont en position $m$-générale, il existe donc une hypersurface $H'$ de degré $m$ qui interpole $P_3, \ldots, P_{2m+2}$ et évite le point $P_{2m+3}$.
On choisit de plus un hyperplan $H$ qui contient $P_{1}, P_{2}$ et évite $P_{2m+3}$.
L'hypersurface $H\cup H'$ est de degré $m+1$ et interpole tous les $P_i$ sauf $P_{2m+3}$.
On en déduit que pour tout $i_0\in \{1, \ldots, 2m+3\}$, il existe une hypersurface de degré $m+1$ qui interpole tous les $P_i$ sauf $P_{i_0}$.

\bigbreak

\noindent \textbf{Preuve de (\ref{your2}).}

\noindent\textbf{Étape 2a.}
Si $m=0$, alors deux points distincts sont toujours $0$-liés et également tou--jours contenus dans une courbe de degré $2$, la propriété est donc trivialement vérifiée dans ce cas.
On supposera donc désormais que $m\geq 1 $ et on se donne un $(2m+2)$-uplet de points $P_1, \ldots, P_{2m+2}$ $m$-liés et deux à deux distincts.

\medbreak

\noindent \textbf{Étape 2b.}
Montrons que si les $P_i$ ne sont pas coplanaires et que $m+1$ d'entre eux sont alignés, alors les $m+1$ restants ne le sont pas.
Pour ce faire, on doit supposer que $\P^N$ est de dimension $N\geq 3$.

Raisonnons par l'absurde et supposons qu'il existe deux droites disjointes $L_1$ et $L_2$ contenant respectivement les points $P_1, \ldots, P_{m+1}$ et $P_{m+2}, \ldots, P_{2m+2}$.
Il existe un unique plan contenant $P_{m+2}$ et la droite $L_{1}$. Ce plan évite les points $P_{m+3}, \ldots, P_{2m+2}$.
On en déduit l'existence d'un hyperplan $H$ vérifiant les mêmes propriétés.
D'après la propriété (\ref{cock}) du lemme \ref{kipu2}, il existe une hypersurface $H'$ de degré $m-1$ qui contient les points $P_{m+3}, \ldots,P_{2m+1}$ et l'hypersurface $H\cup H'$ est de degré $m$ et contient tous les $P_i$ sauf $P_{2m+2}$.
De même, pour tout $i_0$, on peut construire une hypersurface contenant tous les $P_i$ sauf $P_{i_0}$, ce qui contredit le fait que les $P_i$ sont liés.

\medbreak

\noindent \textbf{Étape 2c.} Montrons que les $P_i$ sont coplanaires.
Ici encore, nous allons raisonner par l'absurde en supposant qu'ils ne le sont pas.
Soit $\mathfrak{D}$, l'ensemble des droites de $\P^N$ contenant au moins deux points distincts parmi les $P_i$.
Soit $L_1$ une droite de $\mathfrak{D}$ contenant un nombre maximal de ces points. Quitte à réorganiser les indices, on peut supposer que $P_1 \in L_1$.
Par hypothèse, les $P_i$ sont non coplanaires, on peut donc supposer que les points $P_1,P_2,P_3$ et $P_{2m+2}$ sont non coplanaires.

De plus, comme $P_1$ est dans $L_1$, deux situations sont possibles.

\begin{itemize}

\item[\textbullet] Soit $L_1$ contient moins de $m$ points, moyennant quoi, par définition de $L_1$, il n'y a pas de $(m+1)$-uplet de $P_i$ alignés.

\item[\textbullet] Soit $L_1$  contient $m+1$ points et d'après l'étape 2b, les $m+1$ points restants sont non alignés.
\end{itemize}

\medbreak

Ainsi, les points $P_4, \ldots, P_{2m+2}$ forment un $(2m-1)$-uplet de points tel que $m+1$ d'entre eux ne sont pas alignés. 
D'après la preuve la propriété \ref{butt2}, ces points sont en position $(m-1)$-générale, il existe donc une hypersurface $H'$ de degré $m-1$ qui contient $P_4, \ldots, P_{2m+1}$ et évite $P_{2m+2}$.
Il existe de plus un hyperplan $H$ contenant $P_1, P_2, P_3$ et évitant $P_{2m+2}$.
L'hypersurface $H\cup H'$ est donc de degré $m$ et contient tous les $P_i$ sauf $P_{2m+2}$.

En appliquant aux points  $P_1, \ldots, P_{2m+1}$, le raisonnement que l'on vient d'appliquer à $P_{2m+2}$, on conclut que ces points sont en position $m$-générale.
Il y a contradiction, les points $P_1, \ldots, P_{2m+2}$ sont donc coplanaires.

\medbreak

\noindent \textit{Étape 2d.} Maintenant que l'on sait que les $P_i$ sont coplanaires, montrons qu'un $(2m+2)$-uplet de points de $\P^2$ est $m$-lié seulement si ces points sont sur une même conique. Nous allons traiter séparément les cas $m=1$ et $m=2$.

Si $m=1$, quatre points du plan sont toujours $1$-liés et sont également toujours contenus dans une même conique.
Si $m=2$ et que les six points ne sont pas tous sur une même conique, alors pour tout $1 \leq i_0 \leq 6$, il existe une conique\footnote{On rappelle que par cinq points de $\P^2$ passe toujours au moins une conique. Il y a d'ailleurs unicité si et seulement si quatre d'entre eux ne sont pas alignés, c'est-à-dire s'ils sont en position $2$-générale.} contenant tous les $P_i$ sauf $P_{i_0}$.
Il y a contradiction avec le fait que les $P_i$ sont $2$-liés. Ils sont donc bien sur une même conique.

Supposons maintenant que $m\geq 3$ et $P_1, \ldots, P_{2m+2}$ soit un $(2m+2)$-uplet de points $m$-liés dans $\P^2$ tel que $m+2$ d'entre eux ne sont pas alignés. Supposons que ces points ne soient pas tous sur une même conique.
On note de nouveau $\mathfrak{D}$ l'ensemble des droites de $\P^2$ contenant au moins deux points distincts parmi les $P_i$.
Soit $L_1$ une droite contenant un nombre maximal de ces points.
On peut supposer, quitte a changer l'ordre des indices, que $P_1\in L_1$.
Montrons qu'il existe quatre points $P_{i_1}, \ldots, P_{i_4}$ parmi $P_2, \ldots, P_{2m+1}$ tels qu'il existe une conique contenant $P_1,P_{i_1}, \ldots, P_{i_4} $ et évitant $P_{2m+2}$.

Dans le cas contraire, pour tout $i \in \{2, \ldots, 2m-2\}$ une conique $C_i$ contenant $P_1$,$P_i$,$P_{i+1}$,
$P_{i+2}$ et $P_{i+3}$ contiendrait $P_{2m+2}$.
Soit alors $i\in \{2, \ldots, 2m-3\}$, les coniques $C_i$ et $C_{i+1}$ ont cinq points d'intersection, d'après le théorème de Bezout, elles ont donc une composante irréductible commune $\Gamma$ qui est de degré un ou deux.
Par récurrence, on montre que tous les $P_i$ appartiennent à cette courbe $\Gamma$ ce qui contredit le fait que les $P_i$ ne sont pas tous sur une même conique.

Ainsi, quitte à réorganiser les indices, on peut supposer qu'il existe une conique $C$ contenant les points $P_1, \ldots, P_5$ et évite $P_{2m+2}$.
On rappelle que $P_1$ est un contenu dans une droite $L_1$ contenant un nombre maximal de points parmi les $P_i$.
Par un raisonnement identique à celui qui a été effectué dans l'étape $2c$, on montre que les points $P_6, \ldots , P_{2m+2}$ forment un $(2m-3)$-uplet de points tel que $m$ d'entre eux ne sont pas alignés.
D'après la propriété \ref{butt2} appliquée à $m-2$ et $r=2m-3= 2(m-2)+1$, les points $P_6, \ldots, P_{2m+2}$ sont en position $(m-2)$-générale.
Il existe donc une courbe $C'$ de degré $m-2$ qui contient les points $P_6, \ldots, P_{2m+1}$ et qui évite $P_{2m+2}$. La courbe $C\cup C'$ est de degré $m$ et contient tous les $P_i$ sauf $P_{2m+2}$.
On en déduit la $m$-généralité de ces points, il y a contradiction.
Les points $P_i$ sont donc bien tous sur une même conique plane.
\end{proof}

\newpage
\thispagestyle{empty}

\chapter{Programmes Magma}

\section{Diviseurs $\Delta$-convenables}\label{prgmDconv}

Dans cette section, nous allons présenter un programme magma qui permet de calculer une paire de diviseur $\Delta$ convenable sur une surface lisse $S$ munie d'un $0$-cycle $\Delta$.
La paire ainsi construite vérifie les propriétés suivantes:
\begin{itemize}
\item elle satisfait le critère de la proposition \ref{crit};
\item le diviseur $D_a$ est effectif;
\item les diviseurs $D_a$ et $D_b$ sont tous deux linéairement équivalents à des sections de la surface $S$ avec une hypersurface de son espace ambiant.
\end{itemize}

Pour ce faire on utilise essentiellement la méthode présentée dans la démonstration du lemme \ref{Deltacex}. À savoir: les diviseurs $D_a$ et $D_b^+$ sont construits par interpolation des points du support de $\Delta$. Le diviseur $D_b^-$ est construit en interpolant les points où les deux autres diviseurs ne doivent pas se croiser ou se croisent avec une multiplicité trop importante et en évitant les points où ces diviseurs se croisent bien.

Ce programme se décompose en trois grosses parties que sont le calcul de $D_a$, de $D_b^+$ et de $D_b^-$. Le principe est à chaque fois le même, on considère un système linéaire de courbes qui vérifient de bonnes propriétés et on cherche un candidat dans ce dernier. Suivant sa dimension (et donc son nombre d'éléments définis sur $\F_q$), la recherche se fera de façon exhaustive ou aléatoire.

\paragraph{Préliminaires.} On a besoin d'un certain nombre de scripts pour y parvenir.

\begin{lstlisting}[frame=single]
// Cette fonction est une variation de la fonction "IsReduced"
// de Magma qui a tendance a ramer dans le cas ou le schema
// considere n'est pas une hypersurface de son espace (affine
// ou projectif) ambiant.


function EstReduit(V)
d:=Dimension(V);
A:=SingularSubscheme(V);
dd:=Dimension(A);
delete A;
if dd eq d then return false;
   else return true;
end if;
end function;
\end{lstlisting}

\newpage

\begin{lstlisting}[frame=single]

// Ce programme permet d'evaluer l'entier n minimal pour lequel
// l'interpolation d'une famille finie de points par une
// hypersurface de degre n est possible.


function MinDim(V,Points)

// L'entier $n$ va designer un degre que l'on va
// incrementer au debut de chaque iteration. 

n:=0;
d:=-1;

Proj:=AmbientSpace(V);
Liste:=[Proj!p: p in Points];

while d lt 0 do
      n+:=1;
      L:=LinearSystem(Proj,n);
      L1:=LinearSystem(L,Liste);
      L2:=LinearSystemTrace(L1,V);      
      d:=Dimension(L2);
end while;

return [*n,d,L2,Proj,Liste*];

end function;


\end{lstlisting}

\medbreak

\begin{lstlisting}[frame=single]
// Cette fonction permet de calculer un vecteur tangent
// en un point d'une courbe. Il est exprime dans le
// systeme de coordonnees de la carte affine
// AffinePatch(Proj,P).

function TangentVector(C,P)

assert P in C;
assert Dimension(C) eq 1;
assert IsNonSingular(C,P);

Proj:=AmbientSpace(C);
CAff:=AffinePatch(C,P);
Aff:=AmbientSpace(CAff);
d:=Dimension(Aff);
Q:=Aff!P;

T:=TangentSpace(CAff,Q);
TP:=RationalPoints(T);

if TP[1] eq Q then 
   R:=TP[2];
   else R:=TP[1];
end if;

QQ:=Coordinates(Q);
RR:=Coordinates(R);
return [RR[i]-QQ[i]: i in {1..d}];

end function;
\end{lstlisting}

\medbreak

\begin{lstlisting}[frame=single]
// Cette fonction permet de calculer un systeme lineaire
// ayant un vecteur tangent impose en un point base rationnel.

// ATTENTION: Le vecteur $v$ est un objet "affine".
// Il est defini dans la carte affine naturelle
// Choisie pour P dans "AffinePatch(Proj,P)";

function LinearSystemTangent(Proj,L,P,v)

// On commence par verifier que la longueur de $v$ est la bonne.

d:=Dimension(Proj);
F:=BaseField(Proj);
assert #v eq d;

// On se donne une carte affine naturellement adaptee a $P$.

Aff:=AffinePatch(Proj,P);
p:=Aff!P;
h:=ProjectiveClosureMap(Aff);
LL:=LinearSystem(L,P);
LLL:=Pullback(h,LL);
S:=Sections(LL);
T:=Sections(LLL);
m:=#S;
assert #T eq m;

// On commence par travailler localement, donc sur le systeme
// lineaire affine. Pour chaque section $f$ de notre systeme
// lineaire affine on evalue $\langle grad_p (f), v \rangle$.

Liste:=[&+[Evaluate(Derivative(f,j),p)*v[j]: j in {1..d}]: f in T];
   
M:=Matrix(F,m,1,Liste);
B:=Basis(NullSpace(M));
delete Liste;
delete M;
Liste2:=[ &+[B[i][j]*S[j]: j in {1..m}]: i in {1..#B}];
Lpv:=LinearSystem(LL,Liste2);
return Lpv;

end function;
\end{lstlisting}

\newpage

\paragraph{Calcul de $\mathbf{Da}$.}

Ce programme fait appel à eux sous-programmes ``sectionreduite'' et ``randsectionreduite''. Dont nous donnerons le code source plus loin.

\medbreak

\begin{lstlisting}[frame=single]

// Pour des raisons de commodite, nous allons faire
// en sorte de choisir $D_a$ reduit.

function Da(V,Points,M,T,r)

R:=MinDim(V,Points);
n:=R[1];
d:=R[2];
L:=R[3];
Proj:=R[4];
Liste:=R[5];
Poly:=Random(L);
F:=BaseField(V);
q:=#F;

printf "Le degre minimal d'interpolation pour D_a est %o. ",R[1];

while q^d lt M do
      A:=SectionReduite(V,L);
      if A ne 0
      	 then printf "Le resultat obtenu pour
              D_a a ete trouve de maniere deterministe,
              avec une hypersurface de degre %o. ",n;
	 return A;
      end if;

// En cas d'echec, on augmente le degre d'interpolation
// et on reinitialise notre syst\`eme lineaire.

      n+:=1;      
      L:=LinearSystem(Proj,n);
      L:=LinearSystem(L,Liste);
      L:=LinearSystemTrace(L,V);
      d:=Dimension(L);
end while;


print "Passage en recherche aleatoire pour D_a. ";
delete d;


for i:=1 to r do
      printf "On demarre les recherches avec des
             hypersurfaces de degre %o pour D_a. ",n;
      A:=RandSectionReduite(V,L,T);
      if A ne 0     
      	 then printf "On obtient une hypersurface de
              degre %o fournissant D_a. ",n;
	 return A;
      end if;
      n+:=1;
      L:=LinearSystem(Proj,n);
      L:=LinearSystem(L,Liste);
      L:=LinearSystemTrace(L,V);
end for;

return "Echec";

end function;
\end{lstlisting}

\medbreak

\begin{lstlisting}[frame=single]

// Cette fonction fait une recherche deterministe
// d'un element reduit d'un systeme lineaire.

function SectionReduite(V,L)

assert IsNonSingular(V);
assert IsIrreducible(V);
F:=BaseField(V);
d:=Dimension(L);

// On doit parametrer le systeme lineaire par
// un espace projectif de meme dimension.
// Le probleme est que MAGMA, ne comprend
// la notion d'espace projectif de dimension nulle.
// Il faut donc commencer par traiter le cas ou 
// le systeme lineaire n'a qu'un seul element.

if d eq 0
   then Coeffts:=[[1]];
   else  Pd:=ProjectiveSpace(F,d);
   Coeffts:=RationalPoints(Pd);
end if;
Sec:=Sections(L);


// On teste tous les elements du systeme lineaire jusqu'a
// ce que l'on en trouve un reduit.

for i :=1 to #Coeffts do
    Poly:=&+[Coeffts[i][j]*Sec[j]: j in {1..d+1}];
    D:=Scheme(V,Poly);
    if EstReduit(D)
       then return Poly;
    end if;			
end for;

print "Aucun element reduit trouve de maniere deterministe";
return 0;

end function;
\end{lstlisting}

\newpage

\begin{lstlisting}[frame=single]
// Cette fonction effecute une recherche aleatoire 
// d'un element reduit d'un systeme lineaire.

function RandSectionReduite(V,L,m)

for i:=1 to m do
    
// Il faut eviter le diviseur nul dans notre systeme lineaire.
// Voila comment nous allons proceder.
    
    Poly:=Random(L);
    while Poly eq 0 do
    	  Poly:=Random(L);
    end while;
 
    D:=Scheme(V,Poly);  
    if EstReduit(D) then return Poly;
    end if; 
end for;

print "Aucun element reduit trouve de maniere aleatoire";
return 0;
end function;
\end{lstlisting}

\paragraph{Calcul de $\mathbf{D_b}$.}

Nous donnons les deux programmes permettant de calculer respectivement
$D_b^+$ et $D_b^-$. Tout comme le programme précédent, ces programmes appellent chacun deux sous-programmes effectuant des recherches exhaustives ou aléatoires.
Ces sous-programmes ressemblant fortement aux programmes ``sectionreduite'' et ``randsectionreduite'', nous n'avons pas jugé nécessaire de les ajouter dans cette annexe.

De même, nous ne présenterons pas les programmes dans leur intégralité. Car certaines parties sont un copié/collé d'un des programme déjà présenté

\medbreak

\begin{lstlisting}[frame=single]
// Cette fonction permet de calculer la partie 
// effective d'un diviseur $D_b$ lorsque
// l'on a calcule $D_a$.

// Pour des raisons de simplicite, nous allons imposer
// a $D_b^+$ d'eviter le lieu singulier de $D_a$ hors
// du support de $\Delta$.

function Dbplus(V,Points,M,T,r,Ha)

// On doit se donner un degre minimal d'interpolation
// On se donne celui de $H_a$ qui est tout indique.
	 
Proj:=AmbientSpace(V);
DA:=Scheme(Proj,Ha);
Da:=Scheme(V,Ha);
F:=BaseField(V);
q:=#F;
Liste:=[Proj!p: p in Points];


// Il faut ensuite construire le schema a eviter.
// C'est-a-dire tous les points geometriques singuliers
// de $D_a$ qui sont hors du support de $\Delta$.
// C'est pour cela que nous allons construire le schema
// "evite".

// Par ailleurs il faut memoriser les points du support
// de $\Delta$ qui sont singuliers pour $D_a$
// car en ces derniers, $D_b^+$ devra etre lisse!
// C'est pourquoi nous allons construire la liste DaSing,
// puis le schema "evite2".


UA:=SingularSubscheme(Da);
assert Dimension(UA) eq 0;
evite:=UA;

MultiA:=[];
DaSing:={@ @};
N:=#Points;

for i:=1 to N do
    mi:=Multi(UA,Liste[i]);
    Append(~MultiA,mi);

// Il faut enlever les points du support de $\Delta$
// que l'on ne doit bien sur pas eviter.
    
    if mi gt 0 then
        DaSing join:= {@Liste[i]@};   
    	evite:=DiffIter(evite,Liste[i],mi);
    end if;
end for;

// Il y a deux facons de programmer le schema vide.
// C'est le schema d'equation 1 et celui d'equations x,y,z,t.
// Le second se comporte mieux avec la fonction meet donc:

if IsEmpty(evite) then
   evite:=EmptySubscheme(Proj);
end if;


if IsEmpty(DaSing) then
   evite2:=EmptySubscheme(Proj);
   else evite2:=Cluster(DaSing);
end if;


// On construit ensuite notre systeme lineaire.

d:=-1;

// On initialise $n$ a degre de $H_a$ moins
// un car il sera incremente des le debut de la boucle.

deg:=TotalDegree(Ha);
n:=deg-1;


while d eq -1 do
      n+:=1;
      L:=LinearSystem(Proj,n);
      L:=LinearSystem(L,Liste);
      L:=Complement(L,DA);
      L:=LinearSystemTrace(L,V);

// Afin de gagner du temps,
// on verifie que le schema a eviter ne rencontre
// pas l'ensemble des points bases de $L$.
// On doit aussi verifier que las points de evite2
// ne sont pas generiquement de multiplicite 
// superieure ou egale a $2$ dans $L$.
      
      R:=BaseScheme(L);
      if IsEmpty(DaSing) then
      	 z:=0;
      	 else z:=&+[Multiplicity(L,P)-1: P in DaSing];      
      end if;

	if z eq 0 and IsEmpty(evite meet R) then
      	 d:=Dimension(L);      
	 else d:=-1;
      end if;

end while;

delete UA;
delete R;

printf "Le degre minimal d'interpolation pour D_b^+ est %o. ",n;

// Ensuite, on fait comme d'habitude...

\end{lstlisting}

\medbreak

\begin{lstlisting}[frame=single]
// Cette fonction fournit $D_b^-$.

function Dbmoins(V,Points,M,T,r,Ha,Hb)

Proj:=AmbientSpace(V);
DA:=Scheme(Proj,Ha);
Da:=Scheme(V,Ha);
DB:=Scheme(Proj,Hb);
Db:=Scheme(V,Hb);
F:=BaseField(V);
q:=#F;
Liste:=[Proj!p: p in Points];
N:=#Liste;


// Il faut construire le schema des points
// bases. 
// Nous allons par ailleurs lister trois types de
// Points:
// 1) les points que $D_b^-$ doit eviter;
// 2) les points que $D_b^-$ doit interpoler
// avec une multiplicite fixee;
// 3) ceux que $D_a^-$ doit interpoler
// avec un vecteur tangent fixe.

U:=Da meet Db;
SchemaBase:=U;
Base:=[];
NonBase:=[];
BaseSpe:=[];
mult:=[]; 


// Les elements de Base sont la donnee d'un point et de la
// multiplicite de $U$ en ce point. Les elements de baseSpe sont la
// donnee d'un point et d'un vecteur tangent en ce point.

for j:=1 to N do
    p:=Liste[j];
    mj:=Multi(U,p);
    Append(~mult,mj-1);
    
    if mj eq 0 then printf "ERREUR: D_a et D_b^+ ne se croisent
                    pas en %o. ",p;

    elif mj eq 1 then
       SchemaBase:=Difference(SchemaBase,Cluster(p));
       Append(~NonBase,p);

    elif mj eq 3 then
       SchemaBase:=DiffIter(SchemaBase,p,2); 
       if IsNonSingular(Da,p) 
       	  then Append(~BaseSpe,[*p,TangentVector(Da,p)*]);
	  else Append(~BaseSpe,[*p,TangentVector(Db,p)*]);
       end if;      

    else  SchemaBase:=DiffIter(SchemaBase,p,mj-1);
       Append(~Base,[*p,mj*]);

    end if;
end for;

print "D_b^-, premiere serie de calculs effecutee. ";

delete U;

// Il faut calculer le degre minimal du systeme lineaire.

n:=0;
d:=-1;

while d eq -1 do
      n+:=1;
      L:=LinearSystem(Proj,n);
      L:=LinearSystem(L,SchemaBase);
      for a in BaseSpe do
      	  L:=LinearSystemTangent(Proj,L,a[1],a[2]);
      end for;      
      
      if IsEmpty(NonBase) then z:=0;
      else z:=&+[Multiplicity(L,p): p in NonBase];
      end if;

      if IsEmpty(Base) then z:=z;
      else z:=z+ &+[Multiplicity(L,a[1])-a[2]+1:a in Base];
      end if;
      
      t:=0;
      
      for a in BaseSpe do
      	  if Multiplicity(L,a[1]) gt 2 then t+:=1;
	  end if;	  
      end for;   				     

      if t eq 0 and z eq 0 then
      	 d:=Dimension(L);
      else d:=-1;
      end if;
	 
end while;

printf "Le degre minimal d'interpolation pour D_b^- est %o. ",n;

// Ensuite c'est comme d'habitude...
\end{lstlisting}

\section{Calculs de matrices de parité de codes LDPC}\label{prgmldpc}

\begin{lstlisting}[frame=single]
// Cette fonction permet de repertorier separement
// toutes les droites rationnelles de P3 qui sont contenues
// dans une cubique donnee et toutes celles qui l'intersectent
// en exactement trois points distincts. Elle compte egalement
// le nombre de droites de la seconde categorie qui est issue
// de chaque point.

function PickInc(F,Form)

Surf:=Scheme(P3,Form);
AffSurf:=AffinePatch(Surf,1);
Points:=RationalPoints(AffSurf,F);

A1<t>:=AffineSpace(F,1);
LP:=[{@Parent(AffSurf)| @}: i in {1..#Points}];
LD3:={@Parent(AffSurf)| @};
LDinc:={@Parent(AffSurf)| @};
PointsBis:=IndexedSetToSet(Points);


for p in Points do
L1:=Coordinates(p);
Exclude(~PointsBis,p);
	for q in PointsBis do
	L2:=Coordinates(q);

// On construit notre droite.        

	h:=map<A1 -> A3 | [t*L1[k]+(1-t)*L2[k]:
                  k in {1..3}]>;
	D:=Image(h);
	Sec:=AffSurf meet D;
	PtSec:=RationalPoints(Sec,F);
		if D eq Sec then Include(~LDinc,D);
		elif #PtSec eq 3 then Include(~LD3,D);
			Include(~LP[Index(Points,p)],D);
			Include(~LP[Index(Points,q)],D);
		end if; 
	end for;
end for;
LPcard:=[#LP[i]: i in {1..#LP}];
return [* LPcard, LP ,LD3,LDinc *];
end function;
\end{lstlisting}

\medbreak

Nous allons ensuite présenter la fonction qui permet de calculer une matrice de parité creuse. Cette dernière fait appel à une fonction dont nous ne 
donnerons pas le code source et qui à une droite affine associe un vecteur directeur.

\medbreak

\begin{lstlisting}[frame=single]
// Cette fonction calcule une matrice de parite creuse.

function GenMatDiff(Form)

// On prend notre liste de droites interessantes.

L:=PickInc(F,Form)[3];
Surf:=Scheme(P3,Form);
AffSurf:=AffinePatch(Surf,1);
Points:=RationalPoints(AffSurf,F);
k:=#L;
n:=#Points;

// Il faut qu'on sache par quels points passe chaque droite:

IncidenceCtB:=[];
for i:=1 to k do
S:=[];


	for j:=1 to n do
		if Points[j] in Scheme(A3,Equations(L[i]))
                   then Append(~S,j);
		end if;
	end for;
Append(~IncidenceCtB,S);
end for;
UU:=[Equations(L[i]): i in {1..#L}];
U:=[];

// On calcule un vecteur directeur pour chaque droite.

for j :=1 to #UU do
Append(~U,DirVec(UU[j]));
end for;
delete UU;

// On calcule quelques derivees partielles...

AffForm:=UnHom(Form);
fx:=Derivative(AffForm,1,x);
fy:=Derivative(AffForm,1,y);
fz:=Derivative(AffForm,1,z);

// Et on va construire notre matrice.

MM:=[];
for j:=1 to k do
	ZZ:=[F| 0 : o  in {1..n}];
	a:=IncidenceCtB[j];
	for i:=1 to 3 do
		p:=Coordinates(Points[a][i]);
		d:=1/(Evaluate(fx,p)*U[j][1]+Evaluate(fy,p)*U[j][2]
		+ Evaluate(fz,p)*U[j][3] );
		ZZ[a[i]]:=d;
	end for;
	Append(~MM,ZZ);
end for;
Matmat:=Matrix(F,k,n,MM);
return Matmat;
end function;
\end{lstlisting}


\bibliographystyle{alpha}
\bibliography{biblio}

\end{document}